\documentclass[msom,nonblindrev]{informs4}

\OneAndAHalfSpacedXI
\usepackage{natbib}
\usepackage{diagbox}

 \bibpunct[, ]{(}{)}{,}{a}{}{,}%
 \def\bibfont{\small}%

\usepackage{amssymb}
 \usepackage{booktabs}
 \usepackage[table]{xcolor}
\usepackage{bigstrut}
\usepackage{float}
\usepackage{mathtools}
\usepackage{amsmath} 
\usepackage{algorithm}
\usepackage[noend]{algpseudocode}
\usepackage{optidef}
\usepackage{multirow}
\usepackage{comment} 
\usepackage{longtable}
\usepackage{ltablex}
\usepackage{makecell} 
\usepackage{mathrsfs} 
\usepackage{caption}
\usepackage{subcaption} 
\usepackage{graphicx}
\graphicspath{ {./figures/} }
\usepackage{epstopdf}
\AppendGraphicsExtensions{.pdf} 
\usepackage{enumerate}
\usepackage{bbm}
\usepackage{bbold}
\newcommand\numberthis[1][]{%
    \refstepcounter{equation}%
    \ifx#1\empty\else\label{eq:#1}\fi%
    \tag{\theequation}%
}

\makeatletter
\usepackage{tikz}
\usetikzlibrary{snakes} 
\usetikzlibrary{decorations.pathreplacing}
\definecolor{myLightGray}{RGB}{191,191,191}
\definecolor{myGray}{RGB}{160,160,160}
\definecolor{myDarkGray}{RGB}{144,144,144}
\definecolor{myDarkRed}{RGB}{167,114,115}
\definecolor{myRed}{RGB}{255,58,70}
\definecolor{myGreen}{RGB}{0,255,71}
\definecolor{myBlue}{rgb}{0.0, 0.0, 1.0} 

\usetikzlibrary{shapes.geometric, arrows} 
\tikzstyle{startstop} = [rectangle, rounded corners, minimum width=3cm, minimum height=1cm, text centered, text width=7cm, draw=black, fill=red!30]
\tikzstyle{io} = [trapezium, trapezium left angle=70, trapezium right angle=110, minimum width=3cm, minimum height=1cm, text centered, text width=7cm, draw=black, fill=blue!30]
\tikzstyle{decision} = [diamond, minimum width=3cm, minimum height=1cm, text centered, draw=black, fill=green!30] 
\tikzstyle{arrow} = [thick,->,>=stealth]
\tikzstyle{process} = [rectangle, minimum width=3cm, minimum height=1cm, text centered, text width=7cm, draw=black, fill=orange!30]

\usetikzlibrary{automata}

\usetikzlibrary{positioning}  %                 ...positioning nodes
\usetikzlibrary{arrows}       %                 ...customizing arrows
\tikzset{node distance=4.5cm, % Minimum distance between two nodes. Change if necessary.
         every state/.style={ % Sets the properties for each state
           semithick,
           fill=gray!10},
         initial text={},     % No label on start arrow
         double distance=4pt, % Adjust appearance of accept states
         every edge/.style={  % Sets the properties for each transition
         draw,
           ->,>=stealth',     % Makes edges directed with bold arrowheads
           auto,
           semithick}}
 \usetikzlibrary{arrows.meta, automata,
                positioning,
                quotes}
 \usetikzlibrary{calc}               
\makeatother 

\usepackage{xcolor}
\usepackage{adjustbox} %this
\usepackage{lipsum}% example text

\usepackage{float}

\TheoremsNumberedThrough     % Preferred (Theorem 1, Lemma 1, Theorem 2)
\ECRepeatTheorems

\EquationsNumberedThrough    % Default: (1), (2), ...
\MANUSCRIPTNO{MSOM-2025}

\begin{document}
%%%%%%%%%%%%%%%%

%\setcounter{page}{1}
%\DOI{10.1287/xxxx.0000.0000}%


%\RUNTITLE{Tis a Butter Place}

% \TITLE{Bundling Information Goods of Decreasing Value}
\TITLE{Dynamic Redeployment of Nurses Across Hospitals: A Sample Robust Optimization Approach}  

\ARTICLEAUTHORS{%
\AUTHOR{Wei Liu}
\AFF{Faculty of Business in Scitech, School of Management, University of Science and Technology of China, \EMAIL{liuweimn@ustc.edu.cn}} %, \URL{}}
\AUTHOR{Tianchun Li}
\AFF{Elmore Family School of Electrical and Computer Engineering, Purdue University, \EMAIL{li2657@purdue.edu}}
\AUTHOR{Mengshi Lu}
\AFF{Mitchell E. Daniels, Jr. School of Business, Purdue University, \EMAIL{mengshilu@purdue.edu}}
\AUTHOR{Pengyi Shi}
\AFF{Mitchell E. Daniels, Jr. School of Business, Purdue University, \EMAIL{shi178@purdue.edu}}
} % end of the block

\ABSTRACT{
{ \bf \textit{Problem definition:}} 
We study a workforce redeployment problem in hospital networks, where clinical staff, such as nurses, are temporarily reassigned from overstaffed to understaffed sites to address short-term imbalances. This practice of ``internal travel,'' which gained traction during the COVID-19 pandemic to tackle nurse shortages, presents new operational challenges that require tailored analytical support. Key requirements such as advance notice and short-term secondments must be incorporated. Moreover, in rapidly evolving environments, reliance on historical data leads to unreliable forecasts, limiting the effectiveness of traditional sample-based methods. 
% We consider a nurse transfer problem across a network of hospitals. 
% The core issue addressed is the redistribution of nurses from hospitals with surpluses to those experiencing shortages, ensuring a balance between demand and supply across the system. 
% Specifically, we consider a fully connected network, allowing transfers between any two hospitals. 
% Transferred nurses might need to stay at the destination hospital for a few days, known as secondment, to avoid long commutes. 
% Additionally, the nurses need to receive advance notice of their transfers. 
% However, due to the rapid spread of COVID-19, demand has shifted from historical data and demand predictions based on historical data may not be accurate. As a result, the traditional sample average approximation (SAA) method, which approximates the demand distribution using historical data or predictions,  may not work well. 
{ \bf \textit{Methodology:}} 
We formulate the problem as a stochastic dynamic program and incorporate demand uncertainty via a sample robust optimization (SRO) framework. Using linear decision rule approximation, we reformulate the problem as a tractable linear program. 
%We first model the nurse transfer problem as a stochastic dynamic programming problem. To address the potential deviation from historical data or predictions, we adopt a sample robust optimization (SRO) approach.  However, due to its complexity, this is challenging to solve.  To simplify, we reformulate the problem into a more manageable linear optimization problem using linear decision rules. 
{ \bf \textit{Results:}} 
We evaluate the impact of key network design components on system performance. Network connectivity has the largest effect in reducing the total cost, number of redeployments, and travel distance, but its benefits depend on aligning the secondment duration with the network structure. Full connectivity without proper secondments can be counterproductive. The SRO approach outperforms the traditional sample-average method in the presence of demand surges or under-forecasts by better anticipating emergency redeployments.  
% We show that the network design has a more significant impact on reducing total costs and other related metrics including the transfer volumes and transfer distances compared to the secondment scenarios and the methods applied (SAA and SRO).   
% Furthermore, our analysis reveals that longer secondments notably reduce these metrics as well. 
% Moreover, we show that the SRO performs better in an increasing demand pattern or when actual nurse demand is underestimated in a decreasing or stable demand pattern, as it can effectively reduce emergency transfers. 
{ \bf \textit{Managerial implications:}} 
Internal travel programs offer a promising strategy to alleviate workforce shortages in healthcare systems. Our results highlight the importance of network design, aligning secondment durations with the network structure, and adopting planning methods that are robust to demand surges or inaccurate predictions. 
%Together, these elements help reduce system costs, improve responsiveness, and minimize costly last-minute redeployments. 
%This underscores the importance of considering a fully connected network design, secondments, and the SRO approach to enhance cost efficiency, reduce emergency transfers, and improve the system's responsiveness to demand changes.  
}% 


%\KEYWORDS{butter, margarine, silliness}
% \HISTORY{This paper was first submitted on April 12, 1922 and has been with the authors for 83 years for 65 revisions.}

\maketitle
%%%%%%%%%%%%%%%%%%%%%%%%%%%%%%%%%%%%%%%%%%%%%%%%%%%%%%%%%%%%%%%%%%%%%%

% Samples of sectioning (and labeling) in MNSC
% NOTE: (1) \section and \subsection do NOT end with a period
%       (2) \subsubsection and lower need end punctuation
%       (3) capitalization is as shown (title style).
%
%\section{Introduction.}\label{intro} %%1.
%\subsection{Duality and the Classical EOQ Problem.}\label{class-EOQ} %% 1.1.
%\subsection{Outline.}\label{outline1} %% 1.2.
%\subsubsection{Cyclic Schedules for the General Deterministic SMDP.}
%  \label{cyclic-schedules} %% 1.2.1
%\section{Problem Description.}\label{problemdescription} %% 2.

% Text of your paper here 

% Brusie, C. These Hospitals Offer Internal, Direct-Hire, Travel Nurse Programs (List).
% h;ps://nurse.org/arHcles/hospitals-internal-travel-nurse-programs/ (2023).
% 7. Tarone, D. M. et al. Selected long abstracts from the st. luke’s university health network quality awards program 2021–2022). Int J Acad Med 9, 73–115 (2023). 

% cite more literature 

\section{Introduction}
\label{sec:Introduction}

The decades-long healthcare workforce shortage has escalated into a global crisis, with the United States projected to face significant deficits across multiple frontline roles---including nurses, nurse assistants, and clinical technicians---by 2030 (\citealt{Lyons2023nursing,Hoover2024data} and \citealt{PRSG2024navigating}).  
%\red{cite more here}. 
These shortages have serious implications, including increased staff burnout, higher risks of medical errors, and adverse patient outcomes (\citealt{Augustine_Health_Sciences2021,carayon2018nursing}). In response, a growing number of hospital systems have adopted ``internal travel'' programs, which temporarily redeploy medical staff across hospital sites to manage demand imbalances~(\citealt{brenner2025combating}, and \citealt{yuan2025managing}). This practice, which emerged during the COVID-19 pandemic, introduces new operational challenges that demand tailored analytical support. Our community partner, Indiana University Health (IU Health), is one such system that has implemented the workforce redeployment initiative, termed the Delta Coverage program (\citealt{helm2024delta}). 
Motivated by this context, our research in this paper focuses on nurse redeployment as a representative case, with the problem setup applicable to other roles such as nurse assistants and technicians. 
Next, we introduce the background and outline the key analytical challenges using IU Health as a motivating example, while highlighting the relevance in the broader healthcare context. We then present an overview of our approach and key contributions.  

\begin{figure}[h]
\centering
\includegraphics[width=.4\textwidth]{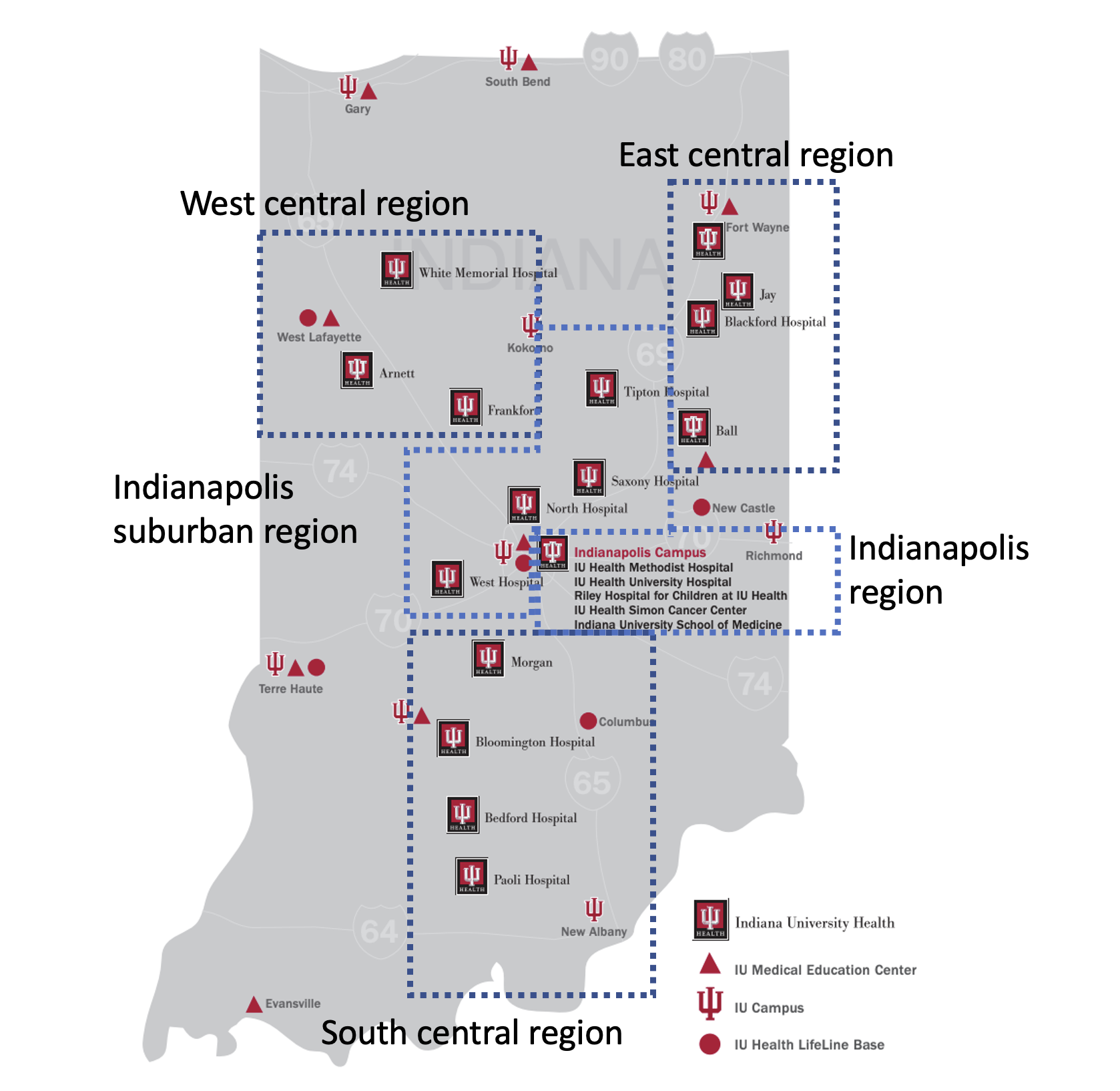}
\caption{Locations of each hospital in the IU Health network.} 
\label{fig_hospital_region}
\end{figure}

\subsection{Problem Background and Operations Challenges}

IU Health is the largest healthcare system in Indiana, operating 16 hospitals.
Figure \ref{fig_hospital_region} shows the distribution of these hospitals throughout the state, with four main regions segregated based on their proximity to Indianapolis (\citealt{IU_Health_fiveregions2023}). 
The central regions, including Indianapolis and its suburban areas, have higher population densities. During the early stages of the COVID-19 pandemic, they experienced more confirmed cases and hospitalizations (per 10,000 residents, \citealt{Indiana_Government_2023}), and hence, faced severe nurse shortages at the onset of the outbreak. 
In contrast, rural regions with lower population densities initially saw fewer cases but faced their unique challenges. With smaller staffing pools and limited surge capacity, these hospitals struggled to absorb day-to-day demand fluctuations. When COVID-19 eventually spread to these areas, they experienced even more severe staffing shortfalls than their urban counterparts. This stochastic and regionally asynchronous demand pattern was a key driver behind IU Health’s implementation of the nurse redeployment program.
%\red{Shall we call it ``nurse transfer'' or ``redeployment'' problem?}
%{\color{blue}: In the writing, I use ``nurse transfer''.} 
%That is, floating nurses across different hospitals to balance the stochastic demand and achieve resource pooling. 
While the idea behind workforce redeployment is intuitive (i.e., resource pooling), its implementation raises several key operations management (OM) challenges, including network design, advance-notice deployment decisions, and planning under demand uncertainty. We elaborate on each of these elements below, which is broadly relevant to hospital systems facing similar staffing shortages.

The \textbf{first element} is network design. Given the statewide distribution of hospitals, two structures are possible. The current model at IU Health follows a {hub-and-spoke} design (\citealt{vstevarova2018performance}), where nurses travel between central and rural hospitals and return home the same day. An alternative is a {fully connected} (or point-to-point) network (\citealt{cook2008airline}), allowing deployments between any hospitals, including rural-to-rural ones. IU Health has been actively exploring this design to increase flexibility and better respond to asynchronous demand surges across regions. However, longer distances between rural sites make daily commuting impractical, which leads to the second design element. %how to provide sufficient advance notice when nurses must travel over 100 miles and stay multiple days at the receiving site. 
%``secondment'' planning that allows a nurse to stay multiple days in the remote location.   

The \textbf{second element} involves accounting for long-distance redeployments through advance notice and secondments. \textit{Secondment} refers to a multi-day assignment at the receiving hospital, during which the nurse remains on-site to avoid prolonged and potentially hazardous daily commutes. It often comes with additional pay to compensate for remote travel.

IU Health initially implemented a real-time redeployment model, where decisions were made on the day of need, e.g., at 5 am each day before the morning shifts start. However, participation was minimal, as nurses were unwilling to engage without adequate notice. This led to the development of a two-stage framework: during the on-call planning phase (two weeks in advance), tentative assignments are made to provide nurses with early notification; during the deployment phase (each morning), these plans can be revised, either by canceling the planned deployment (with a cancellation fee) or initiating emergency deployment (at a premium cost).

This redesign highlights the \textit{human-centric} nature of the problem: nurses are not items, unlike in traditional inventory transshipment or physical resource repositioning models (e.g., \citealt{shu2013models}, and \citealt{he2020robust}). Operational plans must incorporate advance notice, consider the travel burden, and account for the cost of modifying the plans. Nurses are also not interchangeable physical resources—--each nurse has a home location. Nurses need to return to their home location after the secondment, making this fundamentally different from a network flow problem. Detailed discussion on the difference from the transshipment and repositioning literature is in Section~\ref{sec:LiteratureReview}.
Secondment further complicates the problem structure by introducing temporal dependencies within the deployment phase; that is, a decision to deploy a nurse today affects their availability for multiple future days. This necessitates a multi-stage decision problem with interdependent constraints across time. The resulting complexity exceeds what ad hoc or manual planning can effectively manage, calling for new analytical tools to support efficient nurse redeployment planning.

The resulting multi-stage problem naturally fits into a stochastic dynamic programming framework. This introduces the \textbf{third element}: how to make decisions under demand uncertainty, especially when assignments must be made up to two weeks in advance, before actual demand is observed. 
A common approach is the sample average approximation (SAA), which relies on historical data to estimate future demand. However, SAA becomes unreliable during disruptive events like a pandemic, where demand is highly non-stationary. This challenge is amplified in our human-centric setting: both understaffing and unnecessary travels can lead to negative consequences for nurse well-being and patient outcomes (\citealt{meredith2024nurse}). As such, the model must account for adverse scenarios, not just average-case performance.
%For example, understaffing can lead to higher nurse burnout, medication errors, and increased mortality rates among patients, while unnecessary transfers can cause disruptions in care continuity, stress for nurses, and increased operational costs.   
% \red{Add some discussions that we are considering human being here - so we need to account for the bad scenarios, as understaffing and unnecessary transfers both impose negative outcomes to patients and nurses.}  

\subsection{Overview and Contributions}  

We make contributions across modeling, analytical framework, and operational insights. 

First, we formulate a new multi-stage stochastic optimization model to support nurse redeployment with advance notice and secondments, particularly, capturing temporal correlations induced by the secondments; see Section~\ref{sec:formulation}. This bridges the gap between practical nurse redeployment and the literature on transshipment and repositioning by incorporating human-centric constraints such as advance notice and multi-stage plan adjustments. 

Second, to address dynamic decision-making under uncertain and non-stationary demand, we adopt the SRO framework that accounts for deviations from historical data or forecasts and improves the out-of-sample performance. To ensure tractability, we approximate the deployment policy using linear decision rules, which allow us to reformulate the problem as a linear optimization model (see Section~\ref{sec:DRO}). Additionally, we implement a rolling-horizon scheme to incorporate real-time demand updates across the planning horizon.
%Additionally, as this is a multi-stage optimization problem, we use a rolling-horizon approach to incorporate the updated demand observed at the start of each day. More details can be found in Section~\ref{sec:DRO}.   

Third, we generate actionable insights from an extensive case study in Section~\ref{sec:CaseStudy}. Among all design levers, network structure has the largest impact on system cost---fully connected networks outperform hub-and-spoke models by allowing direct transfers and reducing unnecessary routing. However, this benefit depends critically on secondment alignment: short secondments for long-distance deployments can lead to excessive travel costs, while overly long secondments reduce responsiveness to changing demand. We also show that the SRO approach improves performance primarily in environments with rapidly increasing demand or under-predicted demand by proactively allocating additional staff and reducing emergency deployments. 
%While its value is limited in stable environment, SRO offers the resilience under demand shocks, making it a practical tool for future disruptions such as pandemics. 
These findings highlight the need to jointly consider network flexibility, secondment planning, and robust decision-making in workforce redeployment strategies.

\section{Literature Review}
\label{sec:LiteratureReview}

Our research is related to two streams of research. The first stream of research, reviewed in Section~\ref{sec:LiteratureReview_sub1}, includes inventory transshipment and repositioning of physical resources such as vehicles. The second stream, reviewed in Section~\ref{sec:LiteratureReview_sub2},  relates to the healthcare staffing. 
%We introduce these two streams of research in Sections \ref{sec:LiteratureReview_sub1} and \ref{sec:LiteratureReview_sub2}, respectively.

\subsection{Transshipment and Repositioning}\label{sec:LiteratureReview_sub1}

Inventory transshipment involves moving inventory stock across locations within the same echelon to meet local demand, and has been extensively studied; see \cite{paterson2011inventory} for a comprehensive survey. A typical assumption is that inventory is consumed upon meeting demand—once used, it exits the system. In contrast, our nurse redeployment problem involves ``reusable'' human resources. That is, nurses who are temporarily deployed to another location return to their home hospital after completing their secondment and become available again (with the dynamics of availability dependent on prior assignments, not consumptions). This introduces different state dynamics and temporal dependencies not present in traditional inventory models. 

% between different locations of the same echelon to satisfy demand at each location. There is extensive research in this area, e.g., see a comprehensive survey in \cite{paterson2011inventory}. In this category of research, the inventory levels in the system depends on customer consumption, i.e., the items are not reusable. This differentiates from our nurse transfer problem, where nurses who are deployed to another location can be ``reused'' after the deployment period is done. 
%Inventory transshipment can be conducted either before observing the demand (\citealt{agrawal2004dynamic,abouee2015optimal,meissner2018approximate}, and \citealt{zhou2022replenishment}), or after observing the demand at each location (\citealt{ozdemir2006multi,archibald2009index,olsson2010inventory}, and \citealt{bhatnagar2019joint}).

Repositioning of physical resources involves moving items across locations to rebalance supply and demand. Examples include the repositioning of vehicles (e.g., bicycles, cars, emergency vehicles), empty containers, and other mobile assets~(\citealt{song2014empty,belanger2019recent}, and \citealt{xie2021analytics}). We focus on the vehicle repositioning literature, which studies how to relocate shared vehicles across stations to maintain service balance, as it is most relevant to our setting. This area has been extensively explored; see, for instance, \cite{nair2011fleet}, \cite{shu2013models}, \cite{lu2018optimizing}, and \cite{benjaafar2022dynamic}, which typically assume full knowledge of the demand distribution. More recent work considers uncertainty via robust optimization. Among them, the following two are the closest to our research.

\cite{hao2020robust} address idle vehicle pre-allocation under uncertain demand and covariate information, using a distributionally robust optimization (DRO) framework with moment and covariate-based ambiguity sets. \cite{he2020robust} also use a DRO approach to study fleet repositioning in free-float vehicle-sharing systems under temporal demand uncertainty. Our work differs from them in several important ways. (i) \cite{he2020robust} assume that both customer trips and repositioning trips can be completed within a single period, with vehicle availability determined by when they are next used or relocated. \cite{hao2020robust} focus on a single-period setting. In contrast, nurses in our setting are assigned to multi-day secondments, which introduces \textit{intertemporal dependencies} that significantly expand the state space and increase problem complexity.
(ii) Unlike vehicles, nurses are not identical units, in the sense that each nurse has a designated home hospital and must return after secondment. This makes their origin a key attribute in planning, and the traditional network flow problem solutions do not apply here due to the \textit{non-interchangeable} nature. 
(iii) As discussed earlier, nurse redeployment decisions must be communicated in advance to ensure preparedness and acceptance, unlike vehicle repositioning, which is typically decided and executed at the start of each period. This advance-notice requires a multi-stage planning structure with commitment and adjustment layers. (iv) We adopt a sample robust optimization (SRO) method that does not require knowledge of moment-based parameters, in contrast to the moment-based ambiguity sets used in \cite{hao2020robust} and \cite{he2020robust}. This allows us to deal with non-stationary demand without relying on parametric assumptions.

\subsection{Healthcare Staffing}\label{sec:LiteratureReview_sub2} 

Healthcare staffing in specific hospital units and emergency departments (EDs) has been extensively studied; see, for example, \cite{yom2014erlang,cho2019behavior},  \cite{ding2020parallel,prabhu2020team,razak2020modelling} and \cite{hu2025prediction}. Advanced planning for healthcare staff assignment across units has also been explored in 
\cite{huh2013multiresource}, \cite{rath2023multilocation,yao2024multi}.  %\cite{yuan2025managing}, and \cite{ryu2025nurse}. 
Our work is most closely related to \cite{yuan2025managing} and \cite{ryu2025nurse}.
\cite{yuan2025managing} propose a multistage staffing strategy for a centralized service system with multiple downstream units. They first set long-term base staffing levels, then determine staffing from pooled or dedicated agency staff at the start of each contract period. The pooled agency staff, deployable across locations, resemble our inter-facility float nurses. Their model is a stylized queueing network solved via stochastic fluid approximation, assuming full knowledge of the demand distribution. In contrast, we focus on daily deployment decisions with richer operational features: (i) advance notice requirements before redeployment, (ii) home locations with distance-dependent transfer costs, and (iii) uncertain, non-stationary demand with an unknown distribution, addressed via SRO. \cite{ryu2025nurse} develop a two-stage staffing model for float pools \emph{within} a hospital, where the first stage allocates staff across units and the float pool, and the second stage reassigns float nurses within a shift to meet realized demand.
Our setting instead addresses cross-hospital floating, which introduces fundamentally different operational challenges, particularly the advance notice and multi-day secondments, which are not relevant in within-hospital float settings. Moreover, while \cite{ryu2025nurse} use a two-stage DRO model with moment constraints, we employ an SRO framework that avoids parametric and moment-based assumptions, making it better suited for non-stationary demand in multi-hospital networks.

% For instance, \cite{vericourt2011nurse} proposed an asymptotic staffing rule based on a predetermined probability of excessive delay by considering the reentrant patients in a single-station closed queueing system. 
%  \cite{yom2014erlang}  further analyzed the staffing problem in a two-station open queueing network, and \cite{van2016restricted} extended the research of  \cite{yom2014erlang}  by considering staffing rule in a restricted queueing system. 
% Recently, \cite{hu2021prediction}  proposed a two-stage prediction-driven staffing framework. In the first stage, they determine the base staffing level weeks in advance, while in the second stage, they make nearly real-time surge staffing decisions in the ED. 
% However, for the research in this area,  they do not consider the assignment of healthcare staff. 
%  from the float pool to each individual unit or hospital. 
%
%Other works on nurse staffing estimation or prediction include studies by \cite{vericourt2011nurse}, \cite{yom2014erlang} and \cite{van2016restricted}.   
% Advanced planning of assignment of healthcare staff are studied in~\cite{campbell2011two}, \cite{rath2023multilocation},  and  \cite{ryu2025nurse}.   
% %\cite{easton2014service},   \cite{wang2014nurse},

Our study setting is based on the inter-facility float practice at IU Health, with two prior works studying this novel practice. \cite{shi2022operations} use nurse redeployment as an application to demonstrate their workload prediction framework, but their focus is on predictive modeling rather than planned decisions. \cite{helm2024delta} develop a stochastic optimization model for nurse transfer with one-day secondments, assuming accurate demand predictions. In contrast, our work addresses multi-period redeployment under limited and uncertain demand information, explicitly accounting for prediction errors using an SRO framework. This allows for more robust and adaptive staffing decisions over time.

\section{Model Formulation}  
\label{sec:formulation} 

In this section, we formulate the model for the nurse deployment problem. We consider a finite-horizon planning problem with the horizon length being $T$ days. Assume that there are $L$ hospitals (locations) in the network. Throughout the paper, we use $[\hspace{2pt}]$ to represent the set of running indices; for example, $[L]=\{1,2,\cdots, L\}$ represents the set of locations. 
%We use hospital and location interchangeably.  

The nurse demands are assumed to be random. Let $\xi_{t}^i\in \mathbb{R} $ be the nurse demand on day $t$ at location~$i$, where $ i\in [L]$ and  $t\in [T]$. 
Let ${\boldsymbol \xi}_{t} =[\xi_{t}^i, 1\le i \le L]$ be the demand vector at all locations on day $t$, and $ \boldsymbol{ \xi}_{[t]} = ( {\boldsymbol \xi}_{1} ,\cdots,{\boldsymbol \xi}_{t})$ be the historical demands up to day $t$ for $t \in [T]$. The uncertain demands observed over the entire planning horizon is a stochastic process, denoted as ${\boldsymbol \xi}_{[T]}=({\boldsymbol \xi}_1,\cdots,{\boldsymbol \xi}_{ T})$ with a joint probability distribution $\mathbb{P}$. 
Without loss of generality, we assume the support set of ${\boldsymbol \xi} $ is nonnegative, and the support set is given by 
 $$ \Xi= \{ \boldsymbol\xi \in \mathbb{R}^{L\times T} :  \xi_{t}^i \ge 0 ,  1\le i \le L, 1\le t\le T \} .  $$    
The central planner (e.g., the resource nurse manager) determines how many nurses to deploy from one location to another through a two-step decision process. First, an initial planned decision, commonly referred to as the on-call decision in hospital practice, is made at the start of the planning horizon (before any demand is realized) to provide nurses with advance notice. Then, daily deployment adjustments are made based on observed demand.

\begin{figure}[h]
\centering
\begin{tikzpicture} [scale=2] [%
    every node/.style={
        font=\scriptsize,
        % Better alignment, see https://tex.stackexchange.com/questions/315075
        text height=1ex,
        text depth=.25ex,
    },
]

%% draw horizontal line   
%\draw[->] (0,0) -- (8.5,0);
    % draw horizontal line   
    \draw (0,0) -- (2,0);
    \draw[snake] (2,0) -- (3,0);
    \draw (3,0) -- (4,0);
%    \draw [snake] (6,0) -- (7,0);
%      \draw  (7,0) -- (8,0); 
%      \draw (8,0) -- (8.5,0); 
%      \draw[ ->]  [snake]  (8.5,0) -- (9.3,0); 

% draw vertical lines
\foreach \x in {1,2,3}{
    \draw (\x cm,0.1) -- (\x cm,0pt);
}

\foreach \x in {0,4}{
    \draw (\x cm,0.2) -- (\x cm,0pt);
}

   \draw[ ->]   (0,0.6) -- (0,0);   
   \draw[ ->]   (0.3,-0.5) -- (0.3,0);   
   \draw[ ->]   (0.6,0.2) -- (0.6,0);   
      \draw[ ->]   (1.3,-0.5) -- (1.3,0);   
   \draw[ ->]   (1.6,0.2) -- (1.6,0);   
   \draw[ ->]   (3.3,-0.5) -- (3.3,0);   
   \draw[ ->]   (3.6,0.2) -- (3.6,0);   
   
%      \draw[ ->]   (0+4,0.6) -- (0+4,0);   
%   \draw[ ->]   (0.3+4,-0.5) -- (0.3+4,0);   
%   \draw[ ->]   (0.6+4,0.2) -- (0.6+4,0);   
%      \draw[ ->]   (1.3+4,-0.5) -- (1.3+4,0);   
%   \draw[ ->]   (1.6+4,0.2) -- (1.6+4,0);   
%   \draw[ ->]   (3.3+4,-0.5) -- (3.3+4,0);   
%   \draw[ ->]   (3.6+4,0.2) -- (3.6+4,0);   
%            \draw[ ->]   (4+4,0.6) -- (4+4,0);   

% place axis labels
\node[anchor=north] at (0,0) {0};
\node[anchor=north] at (1,0) {1};
%\node[anchor=north] at (4,0) {$t-1$};
%\node[anchor=north] at (5,0) {$t$};
\node[anchor=north] at (2,0) {2};
\node[anchor=north] at (3,0) {$T-1$};
\node[anchor=north] at (4,0) {$T$};
%\node[anchor=north] at (5,0) {D+2};
%%\node[anchor=north] at (4,0) {$t-1$};
%%\node[anchor=north] at (5,0) {$t$};
%\node[anchor=north] at (6,0) {D+3};
%\node[anchor=north] at (7,0) {2D}; 
%\node[anchor=north] at (8,0) {2D+1}; 

\node[anchor=north] at (0,1) {$\textcolor{red}{ {\boldsymbol a}_{[T]}}$}; 
\node[anchor=north] at (0.3,-0.4) {$ \textcolor{orange}{{\boldsymbol\xi}_{1}}$}; 
\node[anchor=north] at (0.65,0.6) {$ \textcolor{blue}{{\boldsymbol b}_{1}}$};  
\node[anchor=north] at (1.3,-0.4) {$ \textcolor{orange}{{\boldsymbol\xi}_{2}}$}; 
\node[anchor=north] at (1.6,0.6) {$ \textcolor{blue}{{\boldsymbol b}_{2}}$};  
\node[anchor=north] at (2.5,-0.4) {$ \textcolor{orange}{\cdots}$}; 
\node[anchor=north] at (2.5,0.6) {$ \textcolor{blue}{ \cdots}$}; 
\node[anchor=north] at (3.3,-0.4) {$ \textcolor{orange}{{\boldsymbol\xi}_{T}}$}; 
\node[anchor=north] at (3.6,0.6) {$ \textcolor{blue}{{\boldsymbol b}_{T}}$};  

%\node[anchor=north] at (0+4,1) {$\textcolor{red}{ a_{[(D+1):2D]}}$}; 
%\node[anchor=north] at (0.3+4,-0.4) {$ \textcolor{orange}{\xi_{D+1}}$}; 
%\node[anchor=north] at (0.6+4,0.6) {$ \textcolor{blue}{b_{D+1}}$};  
%\node[anchor=north] at (1.3+4,-0.4) {$ \textcolor{orange}{\xi_{D+2}}$}; 
%\node[anchor=north] at (1.6+4,0.6) {$ \textcolor{blue}{b_{D+2}}$};  
%\node[anchor=north] at (2.5+4,0.6) {$ \textcolor{blue}{\cdots}$}; 
%\node[anchor=north] at (2.5+4,-0.4) {$ \textcolor{orange}{\cdots}$}; 
%\node[anchor=north] at (3.3+4,-0.4) {$ \textcolor{orange}{\xi_{2D}}$}; 
%\node[anchor=north] at (3.6+4,0.6) {$ \textcolor{blue}{b_{2D}}$};  
%
%\node[anchor=north] at (4+4,1) {$\textcolor{red}{ a_{[(2D+1):3D]} }$}; 

\node[anchor=north] at (-1,1) { \textcolor{red}{Planned Decision}  }; 
\node[anchor=north] at (-1,0.6) { \textcolor{blue}{Deployment Decision} };  
\node[anchor=north] at (-1,-0.4) { \textcolor{orange}{Random Demand} };

\end{tikzpicture}
\caption{ Timeline of the nurse transfer problem}  
\label{timeline:original}
\end{figure}
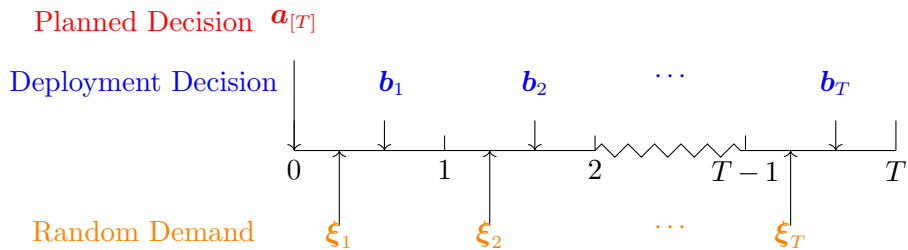

The overall decision timeline is illustrated in Figure~\ref{timeline:original}, and we elaborate on each step below. 
All notations are summarized in Table~\ref{tab:notation} in Appendix~\ref{app:table}.
Throughout the paper, we use the terms ``on-call decision'' and ``planned decision'' interchangeably, as well as ``nurse transfer'' and ``nurse (re)deployment.''

% The notations are summarized in Table~\ref{tab:notation} in Appendix \ref{app:table}. 
% Throughout the paper, we use the terms ``on-call decision'' and ``planned decision'' interchangeably. We also use the terms ``nurse transfer'' and ``nurse deployment/redeployment'' interchangeably. 

%. \red{We call planning decision as ``on-call'' in the text so we should be consistent.} 
%\red{similar as the winter sim paper, create a notation table in the appendix.}
%\blue{Done}

\textbf{Decision Variables. }
At the beginning of the planning horizon, the planner creates a tentative transfer plan--the \emph{planned decision}--for the entire time horizon. 
Let $a_{t}^{ij}$ denote the number of nurses tentatively planned to transfer from location $i$ to location $j$ on day $t$, where $i,j\in [L]$  and $t\in [T]$. By definition, we set $a_{t}^{ii} = 0$, as nurses are not transferred to their ``home'' location. For notational convenience, define ${\boldsymbol a}_t= [a_t^{ij}, 1\le i,j \le L]$ as the vector of on-call decisions for day $t$, and ${\boldsymbol a}_{[t]} = ({\boldsymbol a}_{1} ,\cdots,{\boldsymbol a}_{t})$ denote the sequence of on-call plans from day $1$ to day $t$.

%The purpose of the on-call decision is to determine the number of nurses to be transferred to other locations so that they can receive advance notice of the transfer. 

After observing the realized demand up to day $t$, the planner makes the \textit{deployment decision}, determining the actual number of nurses to be transferred between locations on day t. Let $b_{t}^{ij}$ denote the number of nurses deployed from location $i$ to location $j$ on day ${t}$, where $i,j\in [L]$. As before, we set $b_{t}^{ii} = 0$, and we define ${\boldsymbol b}_t= [b_t^{ij}, 1\le i,j \le L]$ as the vector of deployment decisions for day $t$, and let ${\boldsymbol b}_{[t]} = ({\boldsymbol b}_{1} ,\cdots,{\boldsymbol b}_{t})$ denote the sequence of deployment decisions up to day $t$. Importantly, each deployment decision $b_{t}^{ij}$ is adapted to the realized demand information ${\boldsymbol \xi}_{[t]}$ observed up to day $t$.

\textbf{Costs. }
When the nurses are planned to work at a location other than their home hospital, and the actual deployment occurs, they receive additional compensation. This includes a daily premium pay for working away from home and a relocation bonus upon transfer.
Note that the deployed nurses are required to stay at the destination location for a minimum number of days, referred to as the \emph{secondment}. %(After completing their secondment, nurses must return to their home location.) 
The length of the secondment affects the incurred cost and it may vary depending on the travel distance and the timing within the planning horizon. 
Let $\mu^{ij}(t)= \omega^{ij} \wedge ( T- t+1 )$ denote the secondment length for a nurse transferred from location $i$ to location $j$ on day $t$, where $\wedge$ denotes the minimum operator. The parameter $\omega^{ij}$ represents the minimum required stay when sufficient time remains in the horizon; if the remaining number of days $(T - t + 1)$ is less than $\omega^{ij}$, the nurse will stay at the destination hospital until the end of the planning horizon. Thus, the secondment length is given by the minimum of the two. 
%the required stay $\omega^{ij}$ and the number of days remaining in the horizon $(T - t + 1)$.

Let $p$ denote the daily premium for working at a non-home location, and let $\tau^{ij}$ represent the additional compensation for traveling from location $i$ to location $j$. The planned cost associated with the on-call decision ${\boldsymbol a}_t$ on day $t$ is given by:
\begin{equation}
\begin{aligned} 
c_t^p ({\boldsymbol a}_t )=   \sum_{i=1}^L\sum_{j=1}^L (p \mu^{ij}(t) +\tau^{ij} )a_{t}^{ij}, 
  \end{aligned}  
\end{equation}    
which accounts for both the total premium pay over the secondment period and the one-time bonus per nurse transferred.

Once the deployment decision is made, it incurs a deployment cost consisting of three components: (i) emergency transfer cost, incurred when a transfer is initiated without prior notice; (ii) cancellation cost, incurred when a previously planned transfer is canceled; and (iii) shortage cost, incurred when local nurse demand exceeds available staffing. We specify each of these cost components more precisely below. 

%We first introduce the emergency cost. 
%When the actual number of transfers required at a particular location exceeds the originally planned number for a given day, the planner needs to request an emergency transfer. 
%That is, w
\paragraph{Emergency Transfer Cost. } 
An emergency transfer occurs when the actual number of nurses deployed from location $i$ to $j$ on day $t$ exceeds the planned (on-call) amount, i.e., when $b_{t}^{ij}>a_{t}^{ij}$. In such cases, the excess $(b_{t}^{ij} - a_{t}^{ij})$ nurses are deployed without advance notice and incur an additional cost. Let $\theta_t$ denote the premium multiplier applied to the daily wage $p$ for emergency transfers on day $t$. The total emergency transfer cost for day $t$ is given by: 
\begin{equation} 
\begin{aligned} 
 \sum_{i=1}^{L} \sum_{j=1}^L (\theta_{t} p \mu^{ij}(t) + \tau^{ij} )  (b_{t}^{ij} - a_{t}^{ij}   )^+ , 
  \end{aligned}  
\end{equation}    
where $(x)^+ = \max\{x, 0\}$ denotes the positive part. This cost component accounts for the additional wage premiums (and travel bonuses) associated with emergency deployments.

\paragraph{Cancellation Cost. } 
A cancellation occurs when a previously planned transfer is not executed, i.e., when $a_{t}^{ij} > b_{t}^{ij}$. In this case, $(a_{t}^{ij} - b_{t}^{ij})$ planned transfers are canceled, and each incurs a cancellation fee. 
We denote the percentage fee for canceling a planned transfer by $\eta$. 
The total cancellation cost for day $t$ is given by:
\begin{equation}
\begin{aligned} 
 \sum_{i=1}^{L} \sum_{j=1}^L  (\eta-1) (p \mu^{ij}(t) +\tau^{ij} )  (a_{t}^{ij} - b_{t}^{ij}   )^+ , 
  \end{aligned}   
\end{equation}  
which captures the lost travel bonuses and wages associated with canceling planned transfers.

\paragraph{Shortage Cost. } We use $\delta_{t}^i$ to denote the imbalance between the demand and available nurses at location $i$ on day $t$: 
\begin{equation}
\label{eq:original_shortagecost2} 
\begin{aligned} 
  \delta_{t}^i    =  \xi_{t}^i - \bigg( K^i - \sum_{j=1}^L\sum_{k=({t}-\omega^{ij}+1) \vee  1 }^{t } b_{k}^{ij}   +  \sum_{j=1}^L\sum_{k=({t}-\omega^{ji} +1) \vee 1  }^{t} b_{{k}}^{ji}        \bigg),
  \end{aligned} 
\end{equation}   
where $K^i$ is the initial nurse capacity at location $i$, $\vee$ denotes the maximum operator,
$\sum_{j=1}^L\sum_{k=({t}-\omega^{ij}+1) \vee  1 }^{t } b_{k}^{ij}$ represents the number of nurses from location $i$ who are still on secondment elsewhere on day $t$,
and $\sum_{j=1}^L\sum_{k=({t}-\omega^{ji} +1) \vee 1}^{t} b_{{k}}^{ji}$ represents the number of nurses from other locations who are currently on secondment at location $i$. Then, $(\delta_{t}^i  )^+$ corresponds to the nurse shortage (understaffing) at location $i$ on day $t$, which occurs when available staffing falls short of local demand $\xi_t^i$. Let $s_t^i$ be the unit cost of shortage at location $i$ on day $t$. Then, the total shortage cost incurred on day $t$ is $\sum_{i=1}^{L} s_{t}^i ( \delta_{t}^i  )^+$. 
% \begin{equation}
% \begin{aligned}
%  \sum_{i=1}^{L} s_{t}^i   ( \delta_{t}^i  )^+  .
%   \end{aligned} 
% \end{equation}    

In summary, the total deployment cost incurred on day $t$ is given by  
{
\begin{equation}
\begin{aligned}
 c^d_t ( {\boldsymbol a}_t, {\boldsymbol b}_{t} , {\boldsymbol \xi}_{t}   )     = &  \sum_{i=1}^{L} \sum_{j=1}^L (\theta_{t} p \mu^{ij}(t) + \tau^{ij} )  (b_{t}^{ij} - a_{t}^{ij}   )^+ +  \sum_{i=1}^{L} \sum_{j=1}^L  (\eta-1) (p \mu^{ij}(t) +\tau^{ij} )  (a_{t}^{ij} - b_{t}^{ij}   )^+ +  \sum_{i=1}^{L} s_{t}^i   ( \delta_{t}^i  )^+  .  
\end{aligned} 
\end{equation}   
}
\textbf{State Transitions. }
Before presenting the optimization formulation, we need to describe the state transitions. These transitions are much more complicated than those in the transshipment or repositioning literature due to the secondment requirement and the fact that nurses are non-interchangeable—each has a designated home location and must return after deployment. As a result, we must explicitly track the number of deployed nurses by both their origin-destination pair and the remaining duration of their secondment. 
%This leads to a complex state space that evolves over time and across locations.

Let ${ z}_{t}^{ij}  (k) $ be the number of nurses transferred from location $i$ to location $j$, having the number of remaining secondment days as $k$ at the beginning of day $t$ (before the actual transfer), $ 1\le k\le \omega -1 $, where  $\omega= \max_{i,j}  \omega^{ij} $. 
Let $${\boldsymbol z}_{t} =[ { z}_{t}^{ij}  (k), 1\le i,j\le L , 1\le k\le \omega -1  ].$$ It is easy to see that $ {\boldsymbol z}_{1} = \boldsymbol {0}_{ L\times L \times (\omega -1)  } $.   
%Here, $t=0$ represents the time when we make the planning decision at the beginning of planning horizon. 
After observing the demand ${\boldsymbol \xi}_{t}  $ on day $t$, the planner makes the deployment decision ${\boldsymbol b}_{t}  $.  
Then on  day $t+1$ ($1\le t\le T-1$), we have 
\begin{equation} 
\begin{aligned}
{\boldsymbol z}_{t+1}^{ij}  (k)  &= {\boldsymbol z}_{t}^{ij}  (k+1)   +  b_t^{ij} \mathbbm{1}_{ \mu^{ij} (t)  = k+1 }  , \quad   1\le i , j \le L , 1\le k \le \omega-2,\\ 
{\boldsymbol z}_{t+1}^{ij}  ( \omega-1)  &=  b_t^{ij } \mathbbm{1}_{ \mu^{ij} (t)  = \omega } , \quad   1\le i,j  \le L . 
  \end{aligned}  
\end{equation}   
It should be noted that ${\boldsymbol z}_{T+1} =  \boldsymbol {0}_{ L\times L \times (\omega -1)  } $ since the nurses will return to their home location by the end of planning horizon.

Using ${\boldsymbol z}_{t}$'s, we can rewrite the imbalance quantity at location $i$ on day $t$ as 
\begin{equation}\label{eq:original_shortagecost3} 
\begin{aligned} 
 \delta_{t}^i &=  \xi_{t}^i - \bigg( K^i - \sum_{j=1}^L\sum_{k= 1 }^{ \omega -1 } z_{t}^{ij} (k)  -\sum_{j=1}^L b_t^{ij} + \sum_{j=1}^L \sum_{k= 1 }^{ \omega -1 } z_{t}^{ ji } (k)    + \sum_{j=1}^L b_t^{ji}     \bigg). 
  \end{aligned} 
\end{equation}   
%\red{$b^{ij}$ and $b^{ji}$ missing summation over $L$ locations.}

%Without loss of generality, we assume that $\boldsymbol{\xi}_{[0]}=\emptyset$, indicating the absence of prior demand information at the outset of the first day. 
%Then we have the following optimality equations. 
% It is easy to see that the nurse transshipment problem is a bilevel optimization problem. 
% In the upper-level  optimization problem,  we need to decide the planning decision  $\boldsymbol a_{[T]} $ at the beginning of the planning horizon. 
% In the lower-level  optimization problem,  we optimize the problem sequentially to determine the deployment decision $ {\boldsymbol b}_{t}$ for a given $\boldsymbol a_{[T]} $. 
% Let  $(\boldsymbol a_{[T]}, {\boldsymbol z}_{t}  )$ be the system state for the lower-level  optimization problem.  

\subsection{Optimization Formulation}

The planned decision $\boldsymbol a_{[T]}$ is made at the beginning of the planning horizon, before any demand realization. Subsequently, on each day $t$, the deployment decision $ {\boldsymbol b}_{t}$ is optimized sequentially based on realized demand and the pre-specified plan. Let $v_t  (\boldsymbol a_{[T]}, \boldsymbol z_{t})$ denote the optimal expected cost from day $t$ to the end of the horizon, given the advance plan $\boldsymbol a_{[T]}$ and the state ${\boldsymbol z}_{t}$. Without loss of generality, we set the terminal cost as $v_{T+1}  ( \boldsymbol a_{[T]}, {\boldsymbol z}_{T}  )  = 0$. 
This results in a nested optimization structure, where the outer problem determines the on-call plan: 
\begin{eqnarray} 
& &  \min_{ \boldsymbol a_{[T]} }  \sum_{t=1}^T c_t^p (\boldsymbol a_{t})  + v_1(\boldsymbol a_{[T]}, \boldsymbol z_{1}  )   \\ 
& &  \mathrm{s.t.} \quad \sum_{j=1}^{L} a_{t}^{ij}    \le K^i - \sum_{j=1}^L\sum_{k=(t-\omega^{ij} +1)\vee 1 }^{t-1  } a_{k}^{ij}  ,  {\quad} {1 \le i \le L, 1 \le t \le T,}  \label{orginal_model:eq2a} \\   
%& & \sum_{j=1}^{L} b_{1}^{ij}    \le K^i    {,\quad}{ 1 \le i \le L,} \label{orginal_model:eq3a} \\
 & &  \quad\quad\quad  a_{t}^{ij}  \ge 0 {,\quad}{ 1\le i ,j\le L, 1 \le t \le T . }  
 % & &  b_{1}^{ij}  \ge 0 {,\quad}{ 1\le i ,j\le L,}   
\end{eqnarray}    
Constraint~\eqref{orginal_model:eq2a} is the capacity constraint -- it ensures that the number of nurses planned for transfer from location $i$ on day $t$ does not exceed the number available at that location, accounting for those already on secondment. 

%where the expectation is taken over by the probability distribution of ${\boldsymbol \xi}_{ 1} $, and 
%  \begin{equation} 
% \begin{aligned} 
% J_1  (  {\boldsymbol b}_{1},  {\boldsymbol z}_{1}, {\boldsymbol \xi}_{ [1] }   )    =  c^s_1 ({\boldsymbol b}_{1}, {\boldsymbol \xi}_{ 1 }  )+ v_{2}   ( {\boldsymbol z}_{2}, {\boldsymbol \xi}_{ [1] }   )   . 
%   \end{aligned} 
% \end{equation}     

The inner problem is solved recursively through dynamic programming. Given $\boldsymbol a_{[T]}$ and state ${\boldsymbol z}_{t}$, the expected cost-to-go function is:
$$v_t ( \boldsymbol a_{[T]}, {\boldsymbol z}_{t}  )  = \mathbb{ E} v_t^\xi  ( \boldsymbol a_{[T]}, {\boldsymbol z}_{t}, \boldsymbol \xi_t  ) , $$ 
where $v_t^\xi$ is defined as the solution to the following problem (for realized $\xi_t)$:
\begin{eqnarray} 
v_t^\xi  ( \boldsymbol a_{[T]}, {\boldsymbol z}_{t}, \boldsymbol \xi_t  )   =  \min_{ \boldsymbol b_{t}  } & &  c^d_t ( {\boldsymbol a}_t, {\boldsymbol b}_{t} , {\boldsymbol \xi}_{t}   )    + v_{t+1}(\boldsymbol a_{[T]}, \boldsymbol z_{t+1} )   \\ 
    \mathrm{s.t.} & & \sum_{j=1}^{L} b_{t}^{ij}    \le K^i  -   \sum_{j=1}^L \sum_{k=1 }^{\omega -1}  z_{t}^{ij}(k)    {,\quad}{ 1 \le i \le L, }  \label{orginal_model:eq3a}  \\
 & & {\boldsymbol z}_{t+1}^{ij}  (k) = {\boldsymbol z}_{t}^{ij}  (k+1)   +  b_t^{ij} \mathbbm{1}_{ \mu^{ij} (t)  = k+1 }  , \quad   1\le i , j \le L , 1\le k \le \omega-2,\\ 
 & & {\boldsymbol z}_{t+1}^{ij}  ( \omega-1) =  b_t^{ij } \mathbbm{1}_{ \mu^{ij} (t)  = \omega } , \quad   1\le i,j  \le L , \\ 
 & & b_{t}^{ij}  \ge 0 {,\quad}{  1\le i ,j\le L.}
 \end{eqnarray}     
Constraint~\eqref{orginal_model:eq3a} is the capacity constraint for the inner problem. It ensures that deployment decisions do not exceed the number of nurses currently available at their home location, after accounting for those already on assignment.

Note that nurses are not transferred to their home location, and such transfers do not affect costs or feasibility. Thus, by construction, we have:
$a_t^{ii} = b_t^{ii} = z_t^{ii}(k) = 0, ~ \forall i \in [L],\, t \in [T],\, k \in [\omega - 1]$, and these terms are excluded from the optimization model.

\subsection{Stochastic Programming Formulation}

The nurse transfer problem can also be formulated as a multi-stage stochastic optimization model. 
% The dynamic programming problem in Section~\ref{sec:formulation} can be reformulated as a multistage stochastic optimization problem. 
% The aim of this stochastic optimization model is to obtain the optimal planned and deployment decisions that minimize the total expected cost.  
In this formulation, the deployment decision $\boldsymbol{b}_{t} ({\boldsymbol \xi}_{[t]} )$  specifies a decision rule that determines the actual nurse transfers to be made on day $t$ as a function of all the realized demand on or before day $t$, ${\boldsymbol \xi}_{[t]}$.  
With this decision rule definition, we formulate the following stochastic optimization programming. 

% Alternatively, one can formulate the multi-period nurse transshipment problem as a stochastic optimization problem.
% That is, 
\begin{equation} \label{eq:SO}
\begin{aligned}
  \min_{\substack{ \boldsymbol a_{[T]} , \\ \boldsymbol b_{[T]} (\boldsymbol \xi_{[T]})} }&   { \sum_{t=1}^{T}  \sum_{i=1}^L\sum_{j=1}^L (p \mu^{ij}(t) +\tau^{ij} ) a_{t}^{ij}} + \mathbb{E} \bigg(  \sum_{t=1}^{T}\sum_{i=1}^{L} \sum_{j=1}^L  \bigg( {  (\theta_{t} p \mu^{ij}(t) +  \tau^{ij} ) (b_{t}^{ij}  (\boldsymbol \xi_{[t]})- a_{t}^{ij})^+ } +  \\ 
 & (\eta -1) { (p \mu^{ij}(t) +\tau^{ij} )  (a_{t}^{ij} - b_{t}^{ij}  ({\boldsymbol \xi}_{[t]} )  )^+}  \bigg ) +\sum_{t=1}^{T}\sum_{i=1}^{L}   {  s_{t}^i(\delta_{t}^i  ({\boldsymbol \xi}_{[t]} )   )^+  }  \bigg)   \\  
   \mathrm{s.t.} \quad & \sum_{j=1}^L\sum_{k=(t-\omega^{ij} +1)\vee 1 }^{t  } a_{k}^{ij}  \le K^i {,\quad}{  1 \le i \le L, 1 \le t  \le T,  }  \\ 
& \sum_{j=1}^L\sum_{k=({t}-\omega^{ij}+1)\vee 1 }^{t}  b_{k}^{ij} ({\boldsymbol \xi}_{[k]} )      \le K^i {,\quad}{  1 \le i \le L, 1 \le t  \le T,  } \\
 & \delta_{t}^i   ({\boldsymbol \xi}_{[t]} )    =  \xi_{t}^i - \bigg( K^i - \sum_{j=1}^L\sum_{k=({t}-\omega^{ij}+1) \vee  1 }^{t } b_{k}^{ij} ({\boldsymbol \xi}_{[k]} )   +  \sum_{j=1}^L\sum_{k=({t}-\omega^{ji} +1) \vee 1  }^{t} b_{{k}}^{ji}   ({\boldsymbol \xi}_{[k]} )      \bigg),   
    { 1 \le t \le T, 1\le i \le L, }   \\    
 & a_{t}^{ij}, b_{t}^{ij} ({\boldsymbol \xi}_{[t]} )  \ge 0 {,\quad}{ 1\le i ,j\le L, 1 \le t \le T.}  
 \end{aligned}
\end{equation}  

%The aim of this stochastic optimization model is to obtain feasible planned and deployment decisions that minimize the total expected cost. 
%In other words, we seek a sequence of decision rules that determine the deployment decisions to be made on day $t$ as a function of all the realized demand until day $t$, including the demand on day $t$ at all locations, namely, ${\boldsymbol \xi}_{[t]}$, in addition to the planning decision made at the beginning of the planning horizon. 
%Consequently, the deployment decisions are non-anticipative.  

% The aim of this stochastic optimization model is to obtain the optimal planned and deployment decisions that minimize the total expected cost.  
% Specifically, we seek a sequence of decision rules that determine the deployment decisions to be made on day $t$ as a function of all the realized demand until day $t$, including the demand on day $t$ at all locations, namely, ${\boldsymbol \xi}_{[t]}$, in addition to the planning decision made at the beginning of the planning horizon. 
% We see that the decisions are non-anticipative. 

One major challenge in applying the stochastic optimization model is that the nurse demand distribution is seldom known in practice.  
One approach is to apply SAA to solve Problem~\eqref{eq:SO} by approximating the demand distribution using the empirical distribution of ${\boldsymbol \xi}_{[T]}$. However, SAA exhibits poor out-of-sample performance when there are deviations from historical data or predictions, as is often the case with nurse demand during the COVID-19 pandemic (\citealt{bertsimas2022two}). %\citealt{van2021data} and 
To address this issue, we adopt the sample robust optimization approach to obtain a robust solution to the nurse transfer problem in Section \ref{sec:DRO}.  

% Stochastic optimization requires full information on the nurse demand, which is seldom available in practice.  
% One approach is to estimate the demand distribution using the empirical distribution of ${\boldsymbol \xi}_{[T]}$, employing the sample average approximation (SAA) to solve Problem \ref{eq:SO}. However, SAA exhibits poor out-of-sample performance when the historical data is limited, as is often the case with nurse demand during the pandemic (\citealt{van2021data} and \citealt{bertsimas2022two}). To address this issue, we consider a sample robust optimization approach to tackle the nurse transshipment problem, which we will discuss in Section \ref{sec:DRO}.  

 \section{Sample Robust Optimization Approach} \label{sec:DRO}  
 In this section, we present the sample robust optimization approach to the multi-period nurse transfer problem. 
 %In the robust optimization approach, we solve the nurse transfer problem by considering potential deviations of demand sample paths, to yield better out-of-sample performance compared to traditional sample-based approaches. 
 Compared to conventional sample average approximation, the sample robust approach incorporates potential deviations from available data samples to improve out-of-sample performance. This would allow better staffing decisions in the presence of demand ambiguity.
 We introduce the uncertainty set for sample deviations in Section \ref{sec:Uncertainty_Set} and the corresponding sample robust optimization model in Section \ref{sec:Sample_Robust_Optimization}. 
To address the computational challenge, we employ the linear decision rule and show a tractable reformulation in Section~\ref{sec:Linear_decision_rules}. 
%However, in practice, we do not use the optimal linear decision rule as the policy for the implementation of a multi-stage optimization problem. 
For practical implementation, we further develop a rolling horizon approach in Section \ref{sec:Rolling_horizon}.
%This involves utilizing the observed demand at the beginning of each day to re-optimize the problem, leading to more informed and improved deployment decisions. 

%The timeline of the two-stage stochastic optimization problem is shown in \ref{timeline:approximation}. 
%In  the two-stage stochastic optimization problem, the planner first makes the planning decision $\boldsymbol a_{[T]}$. Then they observe all the demands in the entire planning horizon $\boldsymbol \xi_{[T]}$. After observing the demand, the planner optimizes the deployment decision $\boldsymbol b_{[T]}$. 

\subsection{Uncertainty Set for Sample Deviation}  \label{sec:Uncertainty_Set}

Assume the true distribution $\mathbb{P}$ of the uncertain nurse demand ${\boldsymbol \xi}_{[T]}$ is unknown.
There are $N$ available sample paths ${\boldsymbol \xi}_{[T]}^1, {\boldsymbol \xi}_{[T]}^2, \cdots, {\boldsymbol \xi}_{[T]}^N$.  
Each sample path consists of nurse demands across all locations over $T$ periods. 
%We assume that these sample paths are independent and identically distributed from distribution $\mathbb{P}$.
We refer to the available demand data as sample paths, though they can also be predictions generated using historical data. 
%It should be noted that ${\boldsymbol \xi}_{[T]}^n$ includes the demand realizations in all $L$ locations within $T$ days ($1\le n\le N$).

To incorporate potential perturbations of the realized demand from historical samples or predictions, we construct an uncertainty set $\mathcal{U}_N^{n}$ around each available sample ${\boldsymbol \xi}_{[T]}^n$ for $n=1,\ldots,N$,
\begin{equation}\label{uncertaintyset}
\begin{aligned}
\mathcal{U}_N^{n}  =\{  {\boldsymbol \zeta}_{[T]}  \in \Xi: \|  {\boldsymbol \zeta}_{[T]}   - {\boldsymbol \xi}_{[T]}^n \|_\infty \le \epsilon_N \}. 
\end{aligned}
\end{equation}
Note that both $\boldsymbol{\zeta}_{[T]}$ and $\boldsymbol{\xi}_{[T]}^n$ are $L \times T$ matrices, of which each column represents the demand vector on a particular day.
The matrix norm $\|\cdot\|_\infty$ is the entry-wise infinity norm, which is equal to the infinity norm of the vectorized form of the matrix.
%To apply the infinity norm, we flatten the matrices into vectors by stacking the demand vectors over time (i.e., column-wise vectorization). That is, $\|\cdot\|_\infty$ is taken over the vectorized form of the matrix, defined as the maximum absolute deviation across all elements in the $L \times T$ matrix.
In \eqref{uncertaintyset}, the potential realized demand ${\boldsymbol \zeta}_{[T]}$ lies within an infinity-norm ball centered at the available sample ${\boldsymbol \xi}_{[T]}^n$ with radius $\epsilon_N$.
The uncertainty set radius $\epsilon_N$, which we also refer to as the \emph{robust parameter}, determines the level of conservativeness of the robust model.
A larger $\epsilon_N$ incorporates larger deviations from the available samples and, therefore, generates more conservative solutions.
The choice of $\epsilon_N$ also depends on the sample size $N$.
A larger sample provides more accurate demand distribution information, and consequently, one can choose a smaller $\epsilon_N$.
Note that we adopt the infinity norm primarily for tractability, as it allows an exact reformulation of the robust model, which we detail later. While alternative norms yield similar managerial insights, they often lead to intractable or overly conservative formulations with significantly higher computational costs. Thus, the infinity norm offers a practical and efficient choice for our application context. It has also been advocated in the literature (e.g., \citealt{jiang2018risk}, \citealt{bertsimas2019}, and \citealt{wang2024wasserstein}).

For each $n = 1, \ldots, N$, the uncertainty set in \eqref{uncertaintyset} is equivalent to 
% $$\max_{1\le i\le L,1\le t\le T}  |  \zeta_{t}^i   - \xi_{t}^{i,n}   | \le \epsilon_N  .$$ 
%That is, 
% $$   \underline{\zeta}_{t}^{i,n}   = \xi_{t}^{i,n} - \epsilon_N  \le  \zeta_{t}^i   \le \xi_{t}^{i,n}  + \epsilon_N =  \bar \zeta_{t}^{i,n} $$ 
 $ \xi_{t}^{i,n} - \epsilon_N  \le  \zeta_{t}^i   \le \xi_{t}^{i,n}  + \epsilon_N $,
for all $t = 1, \ldots, T$ and $i = 1, \ldots, L$.   
For ease of notation, let $\underline{\zeta}_{t}^{i,n}   = \xi_{t}^{i,n} - \epsilon_N $, $\bar \zeta_{t}^{i,n} = \xi_{t}^{i,n}  + \epsilon_N $, $  \underline{\boldsymbol \zeta}_t^n =[\underline{ \zeta}_t^{i,n},  1\le i \le L ]  $, and $  \bar{\boldsymbol \zeta}_t^n =[\bar{ \zeta}_t^{i,n},  1\le i \le L ]  $. 
% Furthermore, we let $   \underline{\boldsymbol \zeta}^n = [  \underline{\boldsymbol \zeta}_1^n,\cdots, \underline{\boldsymbol \zeta}_T^n ] \in \mathbb{R}^{L\times T } $ and $   \bar{\boldsymbol \zeta}^n = [ \bar{\boldsymbol \zeta}_1^n,\cdots, \bar{\boldsymbol \zeta}_T^n ] \in \mathbb{R}^{L\times T } $. 
% These matrices represent the lower and upper bounds of the uncertainty set $ \mathcal{U}_N^{n} $. 

\subsection{Sample Robust Optimization Model}    \label{sec:Sample_Robust_Optimization}

In the sample robust optimization model, we 
%consider possible perturbations from the sample path in advance by 
minimize the average of the worst-case costs under each uncertainty set, which is formulated as 
%\begin{equation} 
\begingroup
\allowdisplaybreaks 
\begin{align*}\label{eq:DRO} 
  \min_{\substack{ \boldsymbol a_{[T]} , \\ \boldsymbol b_{[T]} (\boldsymbol \zeta_{[t]})} }&   { \sum_{t=1}^{T}  \sum_{i=1}^L\sum_{j=1}^L (p \mu^{ij}(t) +\tau^{ij} ) a_{t}^{ij}} + \frac{1}{N}    \sum_{n=1}^N  \sup_{ {\boldsymbol \zeta}_{[T]}  \in \mathcal{U}_{N}^{n } }  \bigg(  \sum_{t=1}^{T}\sum_{i=1}^{L} \sum_{j=1}^L  \bigg( {  (\theta_{t} p \mu^{ij}(t) +  \tau^{ij} ) (b_{t}^{ij} (\boldsymbol \zeta_{[t]})- a_{t}^{ij})^+ } +  \\ 
 & (\eta -1) { (p \mu^{ij}(t) +\tau^{ij} )  (a_{t}^{ij} - b_{t}^{ij} (\boldsymbol \zeta_{[t]})  )^+}  \bigg ) +\sum_{t=1}^{T}\sum_{i=1}^{L}   {  s_{t}^i(\delta_{t}^i  (\boldsymbol \zeta_{[t]})  )^+  }  \bigg)   \\  
\mathrm{s.t.} \quad & \sum_{j=1}^L\sum_{k=(t-\omega^{ij} +1)\vee 1 }^{t  } a_{k}^{ij}  \le K^i {,\quad}{  1 \le i \le L, 1 \le t  \le T,  }  \\ 
& \sum_{j=1}^L\sum_{k=({t}-\omega^{ij}+1)\vee 1 }^{t}  b_{k}^{ij}  (\boldsymbol \zeta_{[k]})    \le K^i {,\quad}{  1 \le i \le L, 1 \le t  \le T,  }  \numberthis \\
 & \delta_{t}^i  (\boldsymbol \zeta_{[t]})   =  \zeta_{t}^i - \bigg( K^i - \sum_{j=1}^L\sum_{k=({t}-\omega^{ij}+1) \vee  1 }^{t } b_{k}^{ij}(\boldsymbol \zeta_{[k]})   +  \sum_{j=1}^L\sum_{k=({t}-\omega^{ji} +1) \vee 1  }^{t} b_{{k}}^{ji}  (\boldsymbol \zeta_{[k]})      \bigg),   
    { 1 \le t \le T, 1\le i \le L, }   \\    
 & a_{t}^{ij}, b_{t}^{ij}  (\boldsymbol \zeta_{[t]}) \ge 0 {,\quad}{ 1\le i ,j\le L, 1 \le t \le T,}  \\  
 & \forall{\boldsymbol \zeta}_{[T]} \in \mathcal{U}_N^{n}, 1 \le n \le N.    
\end{align*}
\endgroup
The nonlinear terms in \eqref{eq:DRO} can be linearized, which we present in the following proposition. 
Specifically, we introduce two auxiliary variables 
$x_t^{ij}(\boldsymbol{\zeta}_{[t]})=(b_t^{ij}(\boldsymbol{\zeta}_{[t]}) - a_t^{ij}(\boldsymbol{\zeta}_{[t]}))^+$, 
and $y_t^i(\boldsymbol{\zeta}_{[t]})=(\delta_t^i(\boldsymbol{\zeta}_{[t]}))^+$. 
Note that exact linearization is due to the use of the infinity norm in the uncertainty set in \eqref{uncertaintyset}.
Under other norms, similar linearization will result in more conservative reformulations.
\begin{proposition}\label{prop:reformulation}
Problem  \eqref{eq:DRO} is equivalent to the following optimization problem. 
\begingroup 
\footnotesize
\allowdisplaybreaks
\begin{subequations} \label{eq:reformulation} 
\begin{align} 
  \min_{\substack{\boldsymbol a_{[T]},\\ {\boldsymbol b}_{[T]} (\boldsymbol \zeta_{[T]}), \\  {\boldsymbol x}_{[T]} (\boldsymbol \zeta_{[T]}),\\ {\boldsymbol y}_{[T]} (\boldsymbol \zeta_{[T]})   }   } &   { \sum_{t=1}^{T}  \sum_{i=1}^L\sum_{j=1}^L (p \mu^{ij}(t) +\tau^{ij} ) a_{t}^{ij}} + \frac{1}{N}    \sum_{n=1}^N  \sup_{ {\boldsymbol \zeta}_{[T]}  \in \mathcal{U}_{N}^{n } }  \bigg(  \sum_{t=1}^{T}\sum_{i=1}^{L} \sum_{j=1}^L  \bigg({  (\theta_{t} p \mu^{ij}(t) + \tau^{ij} ) x_{t}^{ij} (\boldsymbol \zeta_{[t]})   }  + \notag\\ 
 &  (\eta -1) { (p \mu^{ij}(t) +\tau^{ij} ) (x_{t}^{ij} (\boldsymbol \zeta_{[t]})    - b_{t}^{ij} (\boldsymbol \zeta_{[t]})    +  a_{t}^{ij}  )    }  \bigg ) +  \sum_{t=1}^{T}\sum_{i=1}^{L}   {  s_{t}^i y_t^i  (\boldsymbol \zeta_{[t]})  }   \bigg)  \label{eq:ldr_singlepolicy1}   \\  
\mathrm{s.t.} \quad & \sum_{j=1}^L\sum_{k=(t-\omega^{ij} +1)\vee 1 }^{t  } a_{k}^{ij}  \le K^i {,\quad}{  1 \le i \le L, 1 \le t  \le T,  }  \label{eq:ldr_singlepolicy2}  \\ 
& \sum_{j=1}^L\sum_{k=({t}-\omega^{ij}+1)\vee 1 }^{t}  b_{k}^{ij}  (\boldsymbol \zeta_{[k]})     \le K^i {,\quad}{  1 \le i \le L,   }  \label{eq:ldr_singlepolicy3}  \\
 & \zeta_{t}^i - \bigg( K^i - \sum_{j=1}^L\sum_{k=({t}-\omega^{ij}+1) \vee  1 }^{t } b_{k}^{ij}  (\boldsymbol \zeta_{[k]})   +  \sum_{j=1}^L\sum_{k=({t}-\omega^{ji} +1) \vee 1  }^{t} b_{{k}}^{ji}   (\boldsymbol \zeta_{[k]})      \bigg)  \le y_{t}^{i} (\boldsymbol \zeta_{[t]}) {,} 
  { 1 \le t \le T, 1\le i \le L,}   \label{eq:ldr_singlepolicy4} \\  
   & b_{t}^{ij} (\boldsymbol \zeta_{[t]}) - a_{t}^{ij} \le x_{t}^{ij}  (\boldsymbol \zeta_{[t]})  {,\quad}{ 1 \le t \le T, 1\le i ,j\le L, }  \label{eq:ldr_singlepolicy5}  \\ 
 & a_{t}^{ij}, b_{t}^{ij}(\boldsymbol \zeta_{[t]}) ,x_{t}^{ij} (\boldsymbol \zeta_{[t]}) \ge 0 {,\quad}{ 1\le i ,j\le L, 1 \le t \le T,}  \label{eq:ldr_singlepolicy6} \\
 & y_{t}^{i} (\boldsymbol \zeta_{[t]})  \ge 0 {,\quad}{  1\le i \le L, 1 \le t \le T,}   \label{eq:ldr_singlepolicy7} \\ 
 & \forall {\boldsymbol \zeta}_{[T]}   \in \mathcal{U}_N^{n}, 1 \le n \le N.  \label{eq:ldr_singlepolicy10} 
 \end{align}
 \end{subequations}  
 \endgroup 
\end{proposition}
{The proof of Proposition~\ref{prop:reformulation} is in Appendix~\ref{proof:prop1}.}

\subsection{Linear Decision Rule Reformulation } 
\label{sec:Linear_decision_rules}

The adaptive decisions $b_t^{ij} (\boldsymbol \zeta_{[t]})  $, $x_t^{ij} (\boldsymbol \zeta_{[t]})  $ and $y_t^{i} (\boldsymbol \zeta_{[t]})  $ are all general functions of the past demand until day $t$, namely, ${\boldsymbol \zeta}_{[t]}$. 
%More specifically, $b_t^{ij}$, $x_t^{ij}$ and $y_t^{i}$ can be regarded as functions $b_t^{ij}({\boldsymbol \zeta}_{[t]} )$, $x_t^{ij}({\boldsymbol \zeta}_{[t]})$ and $y_t^{i}({\boldsymbol \zeta}_{[t]})$. 
%This characteristic makes the sample robust optimization problem challenging to solve. A common approach to addressing this issue is to constrain the general space of decision rules to the space of linear decision rules (LDR). 
To achieve tractable, adaptive optimization, we employ the linear decision rule (LDR) approach; that is, we restrict the adaptive decisions to the class of affine functions of the demand realizations. 
In our problem, this means that the decisions $b_t^{ij}({\boldsymbol \zeta}_{[t]} )$, $x_t^{ij}({\boldsymbol \zeta}_{[t]})$, and $y_t^{i}({\boldsymbol \zeta}_{[t]})$ are linear functions of ${\boldsymbol \zeta}_{[t]}$.
%The sample robust optimization model with linear decision rule can be reformulated as linear optimization
%This restriction enables the problem to be efficiently solved using commercial solvers such as Gurobi. 
%
%There are two possible choices of linear decision rules, namely, single policy and multi-policy approximations, proposed by \cite{bertsimas2022two}.  
%Under the single policy approximation, a single linear decision rule is employed for all uncertainty sets for each decision. 
%We show the details of  linear decision rules under single policy as follows. 

The LDR of the deployment decision ${\boldsymbol b}_{[T]}$ is given by the following class of functions: 
\begin{equation}
\begin{aligned}
  {\mathscr{L}_b }^{ L^2 \times T }  = \left\{   {b}_{t}^{ij}  ({ {\boldsymbol \zeta}_{[t]} } )  \in \mathscr{R}^{L \times t,1} ,   t\in [T], i,j\in[L] \begin{array}{c|l}
             & \exists { b}_{t}^{0,ij},  { b}_{t,ml}^{1,ij}  \in \mathbb {R}, m\in[t] ,  l\in[L]: \\
            &  b_{t}^{ij}  ({ {\boldsymbol \zeta}_{[t]} } ) = { b}_{t}^{0,ij} +  \sum_{m=1}^{t}  \sum_{l=1}^L  {b_{t,ml}^{1,ij} } \zeta_{m}^{l} 
             \end{array}   
              \right\}.
\end{aligned}
\end{equation}  
In period $t$, from location $i$ to $j$, $\mathscr{R}^{L \times t,1} $ denotes the space of all measurable functions from $\mathbb{R}^{ L \times t } $ to $\mathbb{R}^{ 1 } $,   $ { b}_{t}^{0,ij} $ is the intercept, and $ {b_{t,ml}^{1,ij} } $ is the linear coefficient of the demand realization at location $l$ in period $m$.
%,  and  $ \sum_{m=1}^{t}  \sum_{l=1}^L  {b_{t,ml}^{1,ij} } \zeta_{m}^{l}  $  is the sum over all locations up to day~$t$.  
Similarly, we have  the LDR of  ${\boldsymbol x}_{[T]}$ and ${\boldsymbol y}_{[T]}$ as follows.
\begin{equation}
%\begin{aligned}
     {\mathscr{L}_x }^{ L^2 \times T } = \left\{   {x}_{t}^{ij}  ({ {\boldsymbol \zeta}_{[t]} } )  \in  \mathscr{R}^{L \times t,1} ,   t\in [T], i,j\in[L]
     \begin{array}{c|l}
             & \exists { x}_{t}^{0,ij},  { x}_{t,ml}^{1,ij}  \in \mathbb {R}, m\in[t],  l\in[L]:  \\
             &x_{t}^{ij}  ({ {\boldsymbol \xi}_{[t]} } ) = { x}_{t}^{0,ij} +  \sum_{m=1}^{t}  \sum_{l=1}^L  {x_{t,ml}^{1,ij} } \zeta_{m}^{l} 
    \end{array}   
    \right\},
%\end{aligned}
\end{equation}  
\begin{equation}
%\begin{aligned}
     {\mathscr{L}_y }^{ L \times T } = \left\{   {y}_{t}^{i}  ({ {\boldsymbol \zeta}_{[t]} } )  \in  \mathscr{R}^{L \times t,1} ,   t\in [T], i\in[L]
     \begin{array}{c|l}
             & \exists { y}_{t}^{0,i},  { y}_{t,ml}^{1,i}  \in \mathbb {R}, m\in[t],  l\in[L]:\\
            &  y_{t}^{i}  ({ {\boldsymbol \zeta}_{[t]} } ) = { y}_{t}^{0,i} +  \sum_{m=1}^{t}  \sum_{l=1}^L  {y_{t,ml}^{1,i} } \zeta_{m}^{l} 
    \end{array}   
    \right\}. 
%\end{align}
%\end{aligned}
\end{equation}

Using LDR, Problem \eqref{eq:reformulation} can be reformulated as a linear optimization problem, as outlined in the following proposition. 
\begin{proposition}\label{prop:ldr_singlepolicy_reformulation}
Applying the linear decision rules, Problem  \eqref{eq:reformulation} is equivalent to the following linear optimization problem. 
\begingroup 
\allowdisplaybreaks
\begin{equation*} \label{eq:ldr_singlepolicy_reformulation}
\begin{aligned} 
  \min \quad  &  { \sum_{t=1}^{T}  \sum_{i=1}^L\sum_{j=1}^L (p \mu^{ij}(t) +\tau^{ij} ) a_{t}^{ij}} + \frac{1}{N}    \sum_{n=1}^N \bigg(  \sum_{t=1}^{T}  \bigg(  \sum_{l=1}^{L}   \nu_{tl}^{epi}  \bar{ \zeta}_{t}^{ln}  -  \sum_{l=1}^{L}  \psi_{tl}^{epi}  \underline{ \zeta}_{t}^{ln}     \bigg ) + \\ 
 &   \sum_{t=1}^{T}\sum_{i=1}^{L}  \bigg(  \sum_{j=1}^L  \bigg({  \big(\theta_{t} p \mu^{ij}(t) + \tau^{ij} \big){ x}_{t}^{0,ij}  } +   (\eta -1)  (p \mu^{ij}(t) +\tau^{ij} )  ({ x}_{t}^{0,ij}    - { b}_{t}^{0,ij}  +  a_{t}^{ij}  )      \bigg ) -   {  s_{t}^i { y}_{t}^{0,i}     } \bigg)     \bigg)\\  
\mathrm{s.t.}  & \quad   \nu_{tl}^{epi}   - \psi_{tl}^{epi}  =  \sum_{k=t}^{T}  \sum_{i=1}^{L} \bigg(  \sum_{j=1}^L  \bigg({  \big(\theta_{k} p \mu^{ij}(k) + \tau^{ij} \big)   {x_{k,tl}^{1,ij} }   }+  (\eta -1)  (p \mu^{ij}(k) +\tau^{ij} )    ({x_{k,tl}^{1,ij} }    -   {b_{k,tl}^{1,ij} } )      \bigg )   \\
     & \quad \quad\quad\quad\quad \quad\quad  \quad \quad\quad  +   {  s_{k}^i   {y_{k,tl}^{1,i} }  }  \bigg) ,\quad 1 \le l \le L, 1 \le t \le T, \\  
 &\quad\textnormal{capacity constraints \eqref{eq:ldr_singlepolicy_cap},}   \\  
  &\quad\textnormal{shortage constraints \eqref{eq:ldr_singlepolicy_shortage},}   \\  
    &\quad\textnormal{emgency transfer constraints \eqref{eq:ldr_singlepolicy_emgency},}   \\  
 & \quad\textnormal{nonnegativity constraints \eqref{eq:ldr_singlepolicy_noneg},}  \\ 
 \end{aligned} 
\end{equation*}   
\endgroup  
%\label{eq:ldr_singlepolicy_epigraph}  
where the capacity constraints are given by
{\small 
\begin{equation}   \label{eq:ldr_singlepolicy_cap}
\begin{aligned}  
  & \sum_{j=1}^L \sum_{k=(t-\omega^{ij} +1)\vee 1 }^{t  } a_{k}^{ij}  \le K^i {,\quad}{  1 \le i \le L, 1 \le t  \le T,  }  \\ 
& \sum_{k=1}^{T}  \bigg(  \sum_{l=1}^{L}   \nu^{cap}_{klti} \bar{ \zeta}_{k}^{ln}  -  \sum_{l=1}^{L}  \psi^{cap}_{klti}  \underline{ \zeta}_{k}^{ln}     \bigg )     + \sum_{j=1}^L\sum_{k=({t}-\omega^{ij}+1)\vee 1 }^{t} { b}_{k}^{0,ij}  \le K^i  {,\quad}{1 \le i \le L,  1 \le t \le T,  1 \le n \le N, } \\ 
& \nu^{cap}_{klti}   -  \psi^{cap}_{klti}       =   \sum_{m=k }^{T}  \sum_{j=1}^L   \mathbb{1}_{\{ ({t}-\omega^{ij}+1)\vee 1) \le m \le  t \}}    {b_{m,kl}^{1,ij} } \zeta_{k}^{l}    {,\quad}{1 \le i \le L,  1 \le t \le T, 1 \le l \le L,  1 \le k \le T ; } \\  
 \end{aligned}  
\end{equation}  } 
the shortage constraints are given by
{\small
\begin{equation}   \label{eq:ldr_singlepolicy_shortage}
\begin{aligned}  
 &  \sum_{k=1}^{T}  \bigg(  \sum_{l=1}^{L}   \nu^{sho}_{klti} \bar{ \zeta}_{k}^{ln}  -  \sum_{l=1}^{L}  \psi^{sho}_{klti}  \underline{ \zeta}_{k}^{ln}     \bigg )  + \zeta_{t}^i  - \bigg( K^i  - \sum_{j=1}^L\sum_{k=({t}-\omega^{ij}+1) \vee  1 }^{t } { b}_{k}^{0,ij}   +   \sum_{j=1}^L\sum_{k=({t}-\omega^{ji} +1) \vee 1  }^{t}  { b}_{k}^{0,ji}  \bigg)  \le { y}_{t}^{0,i} , \\ 
 & \quad \quad\quad\quad\quad \quad \quad\quad\quad\quad\quad\quad\quad\quad\quad  \quad  { 1\le i \le L,1 \le t \le T, 1 \le n \le N,}   \\  
   &     \nu^{sho}_{klti}   -  \psi^{sho}_{klti}    =  \sum_{m=k }^{T}  \bigg( \sum_{j=1}^L  \mathbb{1}_{\{ ({t}-\omega^{ij}+1)\vee 1) \le m \le  t \}}  \sum_{l=1}^L  {b_{m,kl}^{1,ij} } \zeta_{k}^{l}    -   \sum_{j=1}^L \mathbb{1}_{\{ ({t}-\omega^{ji}+1)\vee 1) \le m \le  t \}}   \sum_{l=1}^L  {b_{m,kl}^{1,ji} } \zeta_{k}^{l}      -   \\ 
   &  \quad \quad\quad\quad\quad  \quad \quad      \mathbb{1}_{\{ {m} =t  \}} \sum_{l=1}^L    {y_{m,kl}^{1,i} } \zeta_{k}^{l}  \bigg)    {,\quad} 
  { 1\le i,l \le L,1 \le t \le T,  1 \le k \le T ; }   \\  
 \end{aligned}  
\end{equation} }  
the emergency transfer constraints are given by
{\small
\begin{equation}     \label{eq:ldr_singlepolicy_emgency}
\begin{aligned}  
   & \sum_{k=1}^{T}  \bigg(  \sum_{l=1}^{L}   \nu^{eme}_{kltij} \bar{ \zeta}_{k}^{ln}  -  \sum_{l=1}^{L}  \psi^{eme}_{kltij}  \underline{ \zeta}_{k}^{ln}     \bigg ) + { b}_{t}^{0,ij} - a_{t}^{ij}    \le   { x}_{t}^{0,ij}    {,\quad}{  1\le i ,j\le L,1 \le t \le T, 1 \le n \le N,}  \\ 
      &  \nu^{eme}_{kltij}   -  \psi^{eme}_{kltij}      =   \sum_{m=k }^{T}   \mathbb{1}_{\{ {m} =t  \}}   ({b_{m,kl}^{1,ij} } \zeta_{m}^{l}   -   {x_{m,kl}^{1,ij} } \zeta_{m}^{l} )      {,\quad}{  1\le i ,j,l\le L,1 \le t \le T,   1 \le k \le T ;}  \\  
 \end{aligned}  
\end{equation} }  
and the nonnegativity constraints are given by
{\small
\begin{equation}    \label{eq:ldr_singlepolicy_noneg}
\begin{aligned}  
& a_{t}^{ij}  \ge 0 {,\quad}{ 1\le i ,j\le L, 1 \le t \le T,}  \\
 &   \sum_{k=1}^{T}  \bigg(  \sum_{l=1}^{L}   \nu^{nnb}_{kltij} \bar{ \zeta}_{k}^{ln}  -  \sum_{l=1}^{L}  \psi^{nnb}_{kltij}  \underline{ \zeta}_{k}^{ln}     \bigg )   \le  { b}_{t}^{0,ij}   {,\quad}{ 1\le i ,j\le L, 1 \le t \le T, 1 \le n \le N,}  \\
    &   \nu^{nnb}_{kltij}   -  \psi^{nnb}_{kltij}      =   - \sum_{m=k }^{T}   \mathbb{1}_{\{ {m} =t  \}}   {b_{m,kl}^{1,ij} } \zeta_{m}^{l}     {,\quad}{ 1\le i ,j,l\le L, 1 \le t \le T,   1 \le k \le T,  }  \\  
         &   \sum_{k=1}^{T}  \bigg(  \sum_{l=1}^{L}   \nu^{nnx}_{kltij} \bar{ \zeta}_{k}^{ln}  -  \sum_{l=1}^{L}  \psi^{nnx}_{kltij}  \underline{ \zeta}_{k}^{ln}     \bigg )   \le  { x}_{t}^{0,ij}   {,\quad}{ 1\le i ,j\le L, 1 \le t \le T, 1 \le n \le N,}  \\
    &   \nu^{nnx}_{kltij}   -  \psi^{nnx}_{kltij}      =   - \sum_{m=k }^{T}   \mathbb{1}_{\{ {m} =t  \}}   {x_{m,kl}^{1,ij} } \zeta_{m}^{l}     {,\quad}{ 1\le i ,j,l\le L, 1 \le t \le T,   1 \le k \le T,  }  \\   
 &   \sum_{k=1}^{T}    \bigg(  \sum_{l=1}^{L}   \nu^{nny}_{klti} \bar{ \zeta}_{k}^{ln}  -  \sum_{l=1}^{L}  \psi^{nny}_{klti}  \underline{ \zeta}_{k}^{ln}     \bigg )     \le  { y}_{t}^{0,i}   {,\quad}{ 1\le i \le L, 1 \le t \le T,   1 \le n \le N, }   \\    
  &   \nu^{nny}_{klti}   -  \psi^{nny}_{klti}        =  - \sum_{m=k }^{T}   \mathbb{1}_{\{ {m} =t  \}} \sum_{l=1}^L    {y_{m,kl}^{1,i} } \zeta_{k}^{l}     {,\quad}{ 1\le i,l \le L, 1 \le t \le T,   1 \le k \le T ,}   \\ 
    & \nu^{epi} , \psi^{epi}  \in \mathbb{R}_+^{T\times L}, \nu^{cap}, \psi^{cap}, \nu^{sho},\psi^{sho},\nu^{nny}, \psi^{nny},  \in \mathbb{R}_+^{T^2\times L^2}  , \nu^{eme}, \psi^{eme} ,\nu^{nnb}, \psi^{nnb},\in \mathbb{R}_+^{T^2\times L^3}  . 
 \end{aligned}  
\end{equation}  }  
\end{proposition} 
{The proof of Proposition~\ref{prop:ldr_singlepolicy_reformulation} is in Appendix~\ref{proof:prop2}. The main idea of the proof is based on \cite{bertsimas2023dynamic}.}

It should be noted that solutions generated using LDR are not necessarily feasible. 
Moreover, LDRs are typically not used directly as implementable policies in multi-stage optimization problems (\citealt{bertsimas2019}). However, the LDR reformulation offers significant computational advantages. This enables one to re-optimize dynamically as new information becomes available, rather than committing to the LDR policy itself. In the next section, we introduce a rolling-horizon approach to implement this dynamic decision-making framework.

%Furthermore, people usually do not implement LDR as a policy for multi-stage optimization problems (\citealt{bertsimas2019}). Nonetheless, the LDR reformulation can be solved efficiently, which provides the basis for us to resolve it when information is dynamically updated (so that we do not need to rely on the LDR but rather treat it as a way to approximate the future actions. Next, we present a rolling-horizon approach for dynamic implementation. 

\subsection{Rolling-Horizon Implementation}\label{sec:Rolling_horizon}
%\red{This section is to connect the SAA and SRO problem back to the original setting, making them suitable for the practical setting where we cannot foresee all the future demand.} 
%we connect the SRO approach back to the original setting of the multi-period nurse transfer problem, where we need to optimize the problem after observing the nurse demand on each day. 
%This makes the approach suitable for the practical setting where demand from nurses is observed on a daily basis.
%In the multi-period nurse transshipment problem, demand from nurses is observed on a daily basis. 

% since it is not a good choice to use the linear decision rule as the policy for implementation  (\citealt{bertsimas2019}). 
The rolling-horizon approach is a commonly used technique in practice that leverages the observed demand at the beginning of each day to re-optimize decisions over a fixed planning horizon, leading to improved outcomes. This method has been discussed in prior literature, including \cite{chand2002forecast} and \cite{glomb2022rolling}. In this section, we describe the rolling-horizon approach for our specific context.

\begin{figure}[h]
\centering
\includegraphics[width=.8\textwidth]{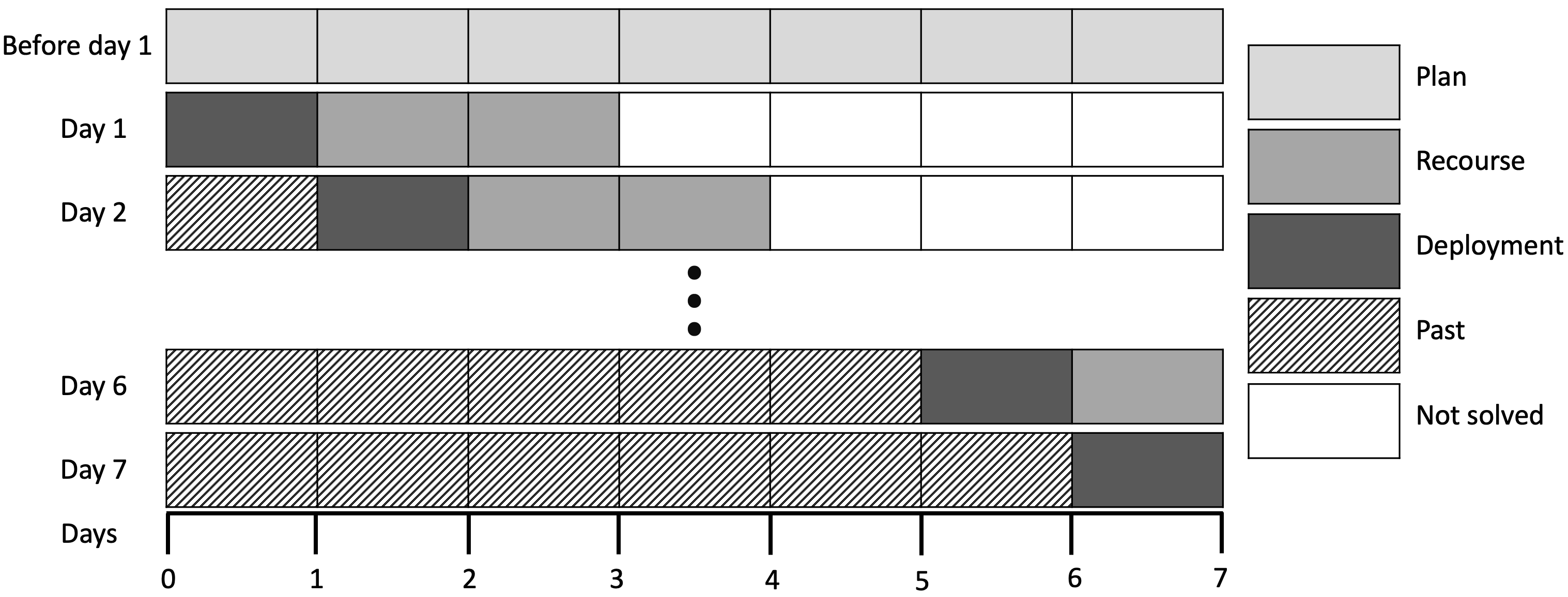}
\caption{ Schematic description of the rolling-horizon approach ($T=7$ and $\omega=3$)} 
\label{timeline:rolling-horizon}
\end{figure}

Figure \ref{timeline:rolling-horizon} illustrates the decision-making timeline.
At the beginning of the planning horizon, we solve the reformulated problem in Proposition~\ref{prop:ldr_singlepolicy_reformulation}  to derive the planned decision.
We denote the resulting optimal planned decision over the entire planning horizon by ${\boldsymbol a}_{[T]}^*$. These specify the planned transfer decisions for each day $t \in [T]$ and correspond to the light gray (Plan) blocks in Figure~\ref{timeline:rolling-horizon}. 
Then, on each day $t$, 
we observe the realized demand ${\boldsymbol \xi}_t$ and re-solve the problem for the remaining horizon $[t, t+S(t)-1]$ using updated demand information. This produces new recourse decisions $[\boldsymbol{b}^*_t,\cdots,\boldsymbol{b}^*_{t+S(t)-1}]$. 
Only the first decision $\boldsymbol{b}^*_t$ is implemented, while the rest are discarded. The dark gray (Deployment) and medium gray (Recourse) blocks in Figure~\ref{timeline:rolling-horizon} illustrate this process: each implemented deployment decision corresponds to the first block in the re-optimized sub-horizon, and the remaining blocks are not used.
This rolling-horizon procedure continues until day $T$, at which point the problem becomes a deterministic optimization problem and the final-day recourse decision is directly implemented.

Feasibility is always guaranteed in this approach because only the first decision $\boldsymbol{b}^*_t$, optimized based on observed demand, is implemented. While future-stage decisions (e.g., $\boldsymbol{b}^*_{t+1}$, $\boldsymbol{b}^*_{t+2}$) are linear functions of yet-to-be-realized demand and may not remain feasible, they are not executed. This iterative process ensures that each deployed decision is always feasible given the current state.

We conclude this section by two remarks. First, we need to choose $S(t)$, the decision horizon, when implementing the approach. 
%Since the secondment $\mu^{ij}(t)= \omega^{ij} \wedge ( T- t+1 )$, 
To balance between the increased computational cost due to a longer horizon and the myopia caused by a shorter decision horizon, we set $S(t) =\omega \wedge ( T- t+1 )$, the secondment length. 
Second, it is necessary to acquire an integer solution for the nurse transfer decisions in practice. One strategy involves enforcing integer constraints on the variables $a_t^{ij}$ and $b_t^{ij}$. However, the computational complexity increases significantly, especially when dealing with a large number of locations or an extended planning horizon. 
%Existing commercial solvers, such as Gurobi, face inefficiencies in solving such problems. 
Thus, we employ the randomized rounding approach \citep{raghavan1987randomized} to round the fractional solutions ${\boldsymbol a}_{[T]}^*$, ${\boldsymbol b}_{1}^*$, $\cdots$, ${\boldsymbol b}_{T}^*$.
This approach keeps the computation tractable while maintaining a close approximation to the mixed-integer optimization problem (\citealt{hao2020robust}).
%This approach bypasses the computational challenges associated with the mixed-integer nature of the problem, making the solution more tractable while maintaining a close approximation to the problem with the integer constraint included (\citealt{hao2020robust}).

\section{Case Study}   
\label{sec:CaseStudy}

In this section, we present a case study to evaluate the effectiveness of our analytical framework.
The study setting is inspired by the configuration of IU Health’s service regions, shown in Figure~\ref{fig_hospital_region}.
To maintain generalizability and avoid attributing specific views, we refer to the system as an ``anonymous hospital network'' in the case study, with four hospitals distributed across four regions: West, East, South, and Central.
% To maintain the generalizability of our case study, we refer to the hospital network in this section as the ``anonymous hospital network,'' with four hospitals distributed across four regions: West, East, South, and Central. 
%This setup corresponds to the four major regions shown in Figure~\ref{fig_hospital_region}. 
Section~\ref{sec:platform} outlines the simulation platform used for performance evaluation, along with the parameter calibration process. Sections \ref{sec:network_design} and \ref{sec:Secondment_effect} focus on analyzing the impact of network design and secondment, respectively. In Section~\ref{sec:ComparisonofPolicies}, we demonstrate the value of robust optimization by comparing solutions from the sample robust optimization with those from the sample average approximation. 
%through which we highlight scenarios where the former approach provides superior results.   

 \subsection{Simulation Evaluation Platform} \label{sec:platform}

Figure~\ref{fig:Data_pipeline_for_simulation_evaluation} illustrates the data pipeline of the simulation evaluation platform, which consists of three key modules: a simulator that generates data, an optimization module that makes nurse transfer decisions and interacts with the simulator to execute those decisions, and an evaluation module that assesses the performance of the deployment strategies. Section~\ref{sec:DataDescription} provides a detailed explanation of each module. Section~\ref{sec:ParameterSetting} outlines the calibration of input parameters.
   
% In this subsection, we first show the details on data generation through a simulation platform and  present the performance evaluation in Section~\ref{sec:DataDescription}. 
%Then we provide a comprehensive overview of the parameter settings employed in our analysis in Section~\ref{sec:ParameterSetting}.  

\subsubsection{Platform Modules} \label{sec:DataDescription}

\quad 

\noindent\textbf{Simulator. }
The first module of our platform simulates daily nurse demand for each hospital for given input parameters (such as daily arrival rate, discharge rate, etc.). The simulator primarily simulates patient movements within the network of hospitals, from which we derive the patient census in each unit at a given time. The nurse demand is then calculated from the patient census using the required nurse-to-patient ratio. 
Specifically, we consider three types of units in each hospital:  medical-surgical unit (MS), progressive care unit (PCU), and intensive care unit (ICU). %, which are defined based on the critical care level of patients. 
To estimate the patient census for each unit, we explicitly model patient transitions between MS, PCU, and ICU, patient discharges, and patient arrivals at each unit; detailed information on these processes is in Appendix \ref{sec:TestingSamplePathGeneration}. Based on the estimated patient census, we apply the following nurse-to-patient ratio to calculate the nurse demand: 5:1 for MS, 3:1 for PCU, and 2:1 for ICU (\citealt{WoltersKluwer2016the}).  
%\red{Wei to confirm these numbers and also see if we can find a reference to support these ratio.}
%These demands form the testing sample paths used for model evaluation.

%This synthetic data generation approach allows for the creation of diverse nurse demand sample paths, enabling a comprehensive assessment of our method. with details provided in Appendix \ref{sec:TestingSamplePathGeneration}. 

% Specifically, we first generate patient arrivals to each hospital by considering both COVID-19 transmission and regular patient arrivals. 
% To further estimate the patient census at each unit (e.g., medical-surgical unit) of the hospital, we present the transition behavior of patients between different units, and detail the generation of patient arrivals, departures, and transitions between these units within each hospital. 
% These data elements are used to estimate nurse demand at each hospital, which serves as the testing sample paths. 

\begin{figure} 
         \centering
         \includegraphics[width=0.8\textwidth]{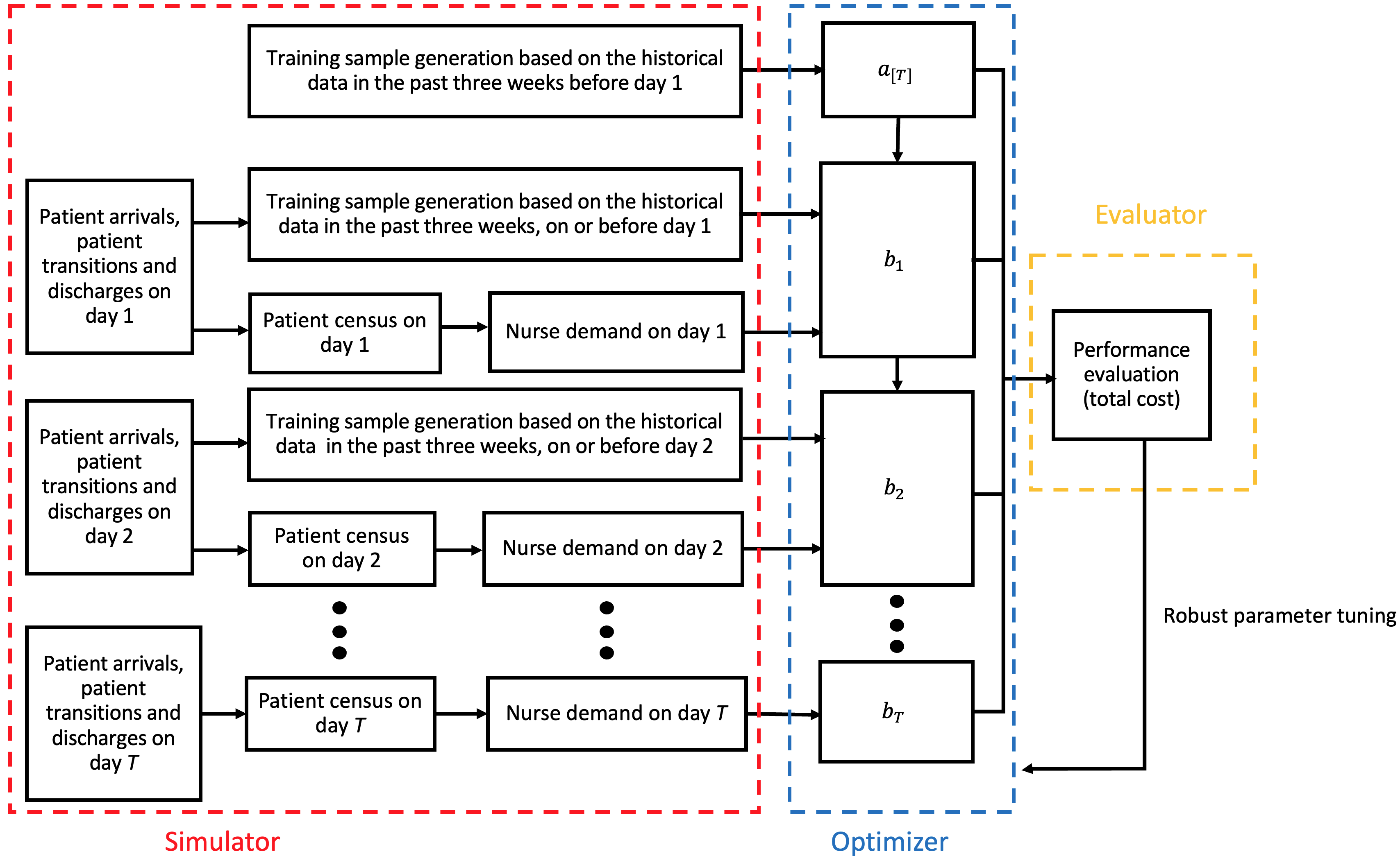}
         \caption{Data pipeline for simulation evaluation (For simplicity, we only show the optimization over one planning horizon for T days. If there are multiple planning horizons, the same approach applies). 
}    
         \label{fig:Data_pipeline_for_simulation_evaluation} 
\end{figure} 

The simulator is used to generate two types of demand data: \textit{testing sample paths}, which are generated using the actual demand distribution (ground truth), and \textit{training sample paths}, which are demand predictions generated based on testing sample paths. Note that the testing and training sample paths have different underlying demand distributions. This distinction allows us to incorporate potential deviations from historical samples or predictions and separate decision-making and decision-evaluation.

\noindent\textbf{Optimizer. }
The optimizer makes planned decisions at the start of each week and deployment decisions on a daily basis. We consider a total horizon of $W=27$ weeks (189 days). Before running the optimizer, we pre-generate $H=30$ testing sample paths, representing sequences of historical demand realized in each of the $L$ hospitals over $W$ weeks.
%We refer to each sequence of the demand data as a \emph{testing sample path}.
The testing data is generated by the simulator calibrated with real hospital data; see details of parameter calibration in  Section~\ref{sec:ParameterSetting}.

The optimizer is executed on each of the $H$ testing sample paths, treating each path as an independent evaluation instance. In each iteration, for every week $w \in[1, W]$ and each day $t$ within that week, we estimate the arrival rates, unit transition probabilities, and discharge probabilities for each hospital using the demand from the previous three weeks. These estimates, which can differ from the true parameters, are then used to generate the training sample paths used in the sample-robust optimization model.

%Note that, depending on the input parameters, the simulator generates demand data from different sets of underlying demand distributions. This distinction allows us to separate training and testing sample paths, even though the same simulator is used for both. We will specify how the input parameters are estimated in generating the training and testing data. 

On each day $t$, we generate $\hat{H}=25$ training sample paths. If day $t$ is the starting day of a week, the weekly planned decisions $\boldsymbol{a}$ are first optimized using the training sample paths. Then, for each day in that week, the deployment decisions $\boldsymbol{b}_t$ are optimized via a rolling-horizon approach, as described in Section~\ref{sec:DRO}. Note that the training sample paths are also generated in a rolling-horizon manner to reflect evolving estimates and maintain consistency with the ``rolling'' decision; see Figure~\ref{fig:Data_pipeline_for_simulation_evaluation} for an illustration. Additional implementation details are provided in Appendix~\ref{sec:TrainingSamplePathGeneration}.

%Due to the potential significant disparities between future and historical sample paths, we predict the possible nurse demand at all hospitals in week $w$ ($1\le w\le W$). %Once we have the simulator, we can use it to generate \emph{training} sample paths of nurse demand, based on which we use
%based on the historical nurse demand on the testing sample path in consideration. 

%Then we can get the planned and deployment decisions for each week under each training set and each testing sample path. 
%This historical data has patient arrivals, departures, and transfers between different units on a daily basis at each hospital, generated in the process of testing sample generation. 
%To estimate the arrival rate, we calculate the sample average of the number of arrivals to each unit of each hospital over a three-week period while accounting for the day-of-the-week effect. 
%Similarly, we estimate the transition probability by computing the sample average of the proportion of customers who moved between different units, and the discharge probability while considering the day-of-the-week effect. 

\noindent\textbf{Evaluator. }
The performance evaluation module is designed to assess the performance of the planned and deployment decisions on the testing data for a given optimization method $\mathcal{A}$. Method $\mathcal{A}$ can refer to either the SRO or the SAA approach (where SAA is a special case of SRO with the robustness parameter set to zero). Evaluation is performed on $H = 30$ testing instances. Under each of the instance, we generate $\hat{H}=25$ training sample paths, each containing nurse demand across $L$ locations and $T = 7$ days. In the SRO method, we group every 5 training sample paths into one training set, resulting in $M = 5$ training sets per instance. For consistency, we apply the same grouping structure to the SAA method.

Let ${\boldsymbol a}_{m,h, w}^{ \mathcal{A} *}$ (${\boldsymbol b}_{m,h,w }^{ \mathcal{A} *} $) denote the optimal planned (deployment) decision in week $w$ using method $\mathcal{A}$ under training data $m$ and testing sample path $h$, where $m\in[M]$, $h\in[H]$, $w\in \{ 1,2,\cdots, W\}$. Then we compute the out-of-sample cost  in week $w$ using method $\mathcal{A}$ as follows:  
 {\small 
\begin{equation}
\begin{aligned}
V_{w} ^{\mathcal{A}} &=  \frac{ 1}{MH } \sum_{m=1}^{M}  \sum_{h=1}^{H }  \bigg(  \sum_{t=1}^{T}   \sum_{i=1}^L\sum_{j=1}^L (p \mu^{ij}(t) +\tau^{ij} ) a_{m,h,w, t}^{\mathcal{A} *, ij}   +  \sum_{t=1}^{T}\sum_{i=1}^{L} \sum_{j=1}^L  \bigg( {  (\theta_{t} p \mu^{ij}(t) + \tau^{ij} ) (b_{m,h,w,t}^{\mathcal{A} *, ij} - a_{m,h,w,t}^{\mathcal{A} *, ij})^+ }    \\ 
& \quad + (\eta -1) { (p \mu^{ij}(t) +\tau^{ij} )  (a_{m,h,w,t}^{\mathcal{A} *, ij} - b_{m,h,w,t}^{\mathcal{A} *, ij}   )^+}  \bigg ) + \sum_{t=1}^{T}\sum_{i=1}^{L}   {  s_{t}^i(\delta_{m, h,w,t}^{\mathcal{A} *,i }  )^+  }  \bigg)  , 
\end{aligned}
\end{equation}
}  
where $\delta_{m,h,w, t}^{\mathcal{A} *,i }$ represents the demand-staffing imbalance on day $t$ in week $w$ using method $\mathcal{A}$, training set $m$, and testing sample path $h$. The average cost over the total horizon with $W$ weeks by using method  $\mathcal{A}$ is given by  
\begin{equation} 
\begin{aligned}
V^{\mathcal{A}} &= \frac{1}{W} \sum_{w=1}^{W}  V_{w} ^{\mathcal{A}}   . 
\end{aligned}
\end{equation}

\subsubsection{Parameter Calibration}  \label{sec:ParameterSetting}    
In this section, we provide the details on the parameters utilized in our analysis. For ease of exposition, we index the four hospitals (West, East, South, and Central) with numeric values -- hospital 1, 2, 3, and 4, respectively.

% it's important to note that each hospital can dynamically adjust its nurse capacity through temporary contracts with external nurses. 

\noindent\textbf{Patient Flow Parameters.}
The key input to the simulator are the number of patient arrivals, discharges, and transitions between units within the hospital. To generate the ground-truth (testing) demand, we begin by estimating baseline arrival rates using historical data from an anonymous hospital system. Temporal and spatial correlations are then introduced through an autoregressive model with noise terms designed to capture these dependencies. Note that the baseline rates and correlation structures are estimated from data, while the arrival model is parameterized to allow tuning of hyperparameters. This enables us to simulate demand surges—such as peak census levels observed during the COVID-19 pandemic—using publicly available benchmarks from the \cite{AmericanHospitalDirectory}. 
%We estimate the probability of discharges and transitions using the hospital data. See the estimation details in Appendix~\ref{sec:TestingSamplePathGeneration}.  
See details in Appendix~\ref{sec:TestingSamplePathGeneration}.

Figure~\ref{nurse_demand_capacity_over_time} illustrates the average nurse demand under the baseline testing scenarios over the $W=27$ weeks at four hospitals, alongside their respective nurse capacities. The demand trajectory in this baseline follows three phases: an increasing demand phase, a decreasing demand phase, and a stable demand phase. Nurse capacity at each hospital is estimated using the provided budgeted staffing ratio and is adjusted weekly to approximate the corresponding trends in demand. Notably, West and Central Hospitals consistently face staffing shortages, while East and South Hospitals maintain surplus nurse capacity. This structural imbalance creates an opportunity for nurse transfers across hospitals, using the redeployment strategies to mitigate system-wide shortages. 
%It is important to note that West and Central Hospitals consistently face a shortage of nurses, whereas East and South Hospitals experience a surplus. This disparity presents an opportunity for nurse transfers between these facilities, potentially alleviating the overall shortage of nursing staff.  

\begin{figure}
     \centering
     \begin{subfigure}[b]{0.23\textwidth}
         \centering
         \includegraphics[width=\textwidth]{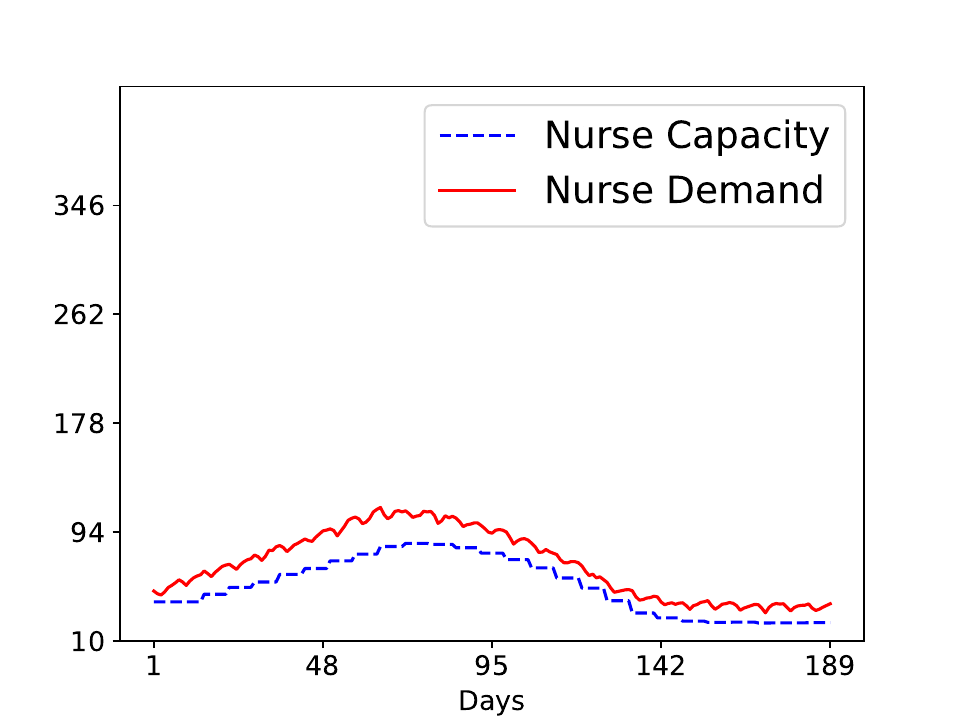}
         \caption{West Hospital} 
         \label{census_Arnett} 
     \end{subfigure}
     \hfill 
     \begin{subfigure}[b]{0.23\textwidth}
         \centering
         \includegraphics[width=\textwidth]{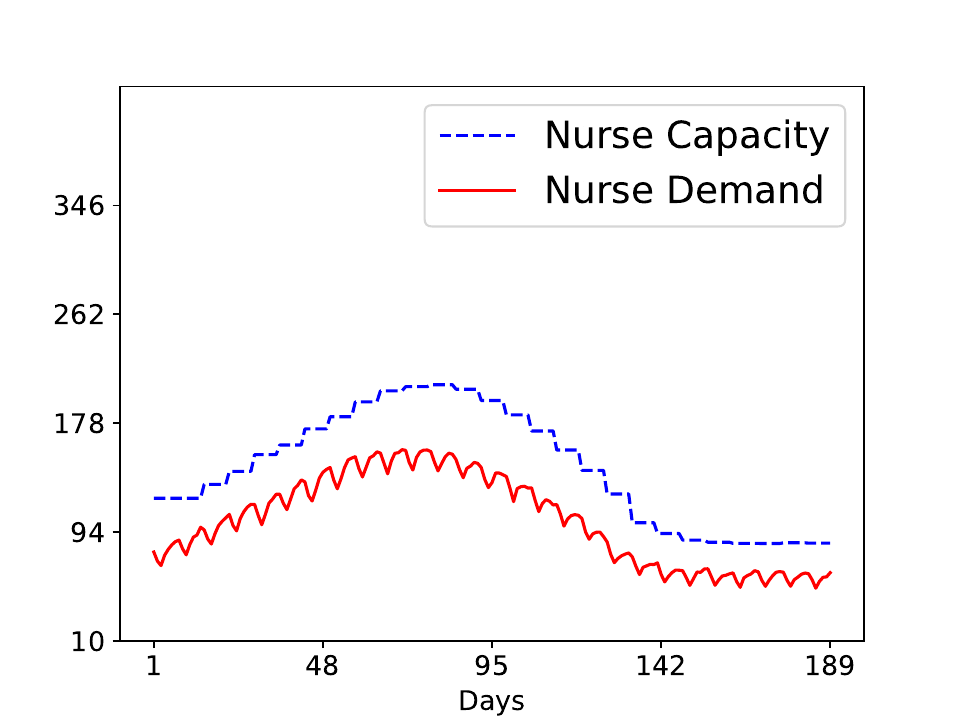}   
         \caption{East Hospital} 
         \label{census_Ball}     
     \end{subfigure} 
          \begin{subfigure}[b]{0.23\textwidth}
         \centering
         \includegraphics[width=\textwidth]{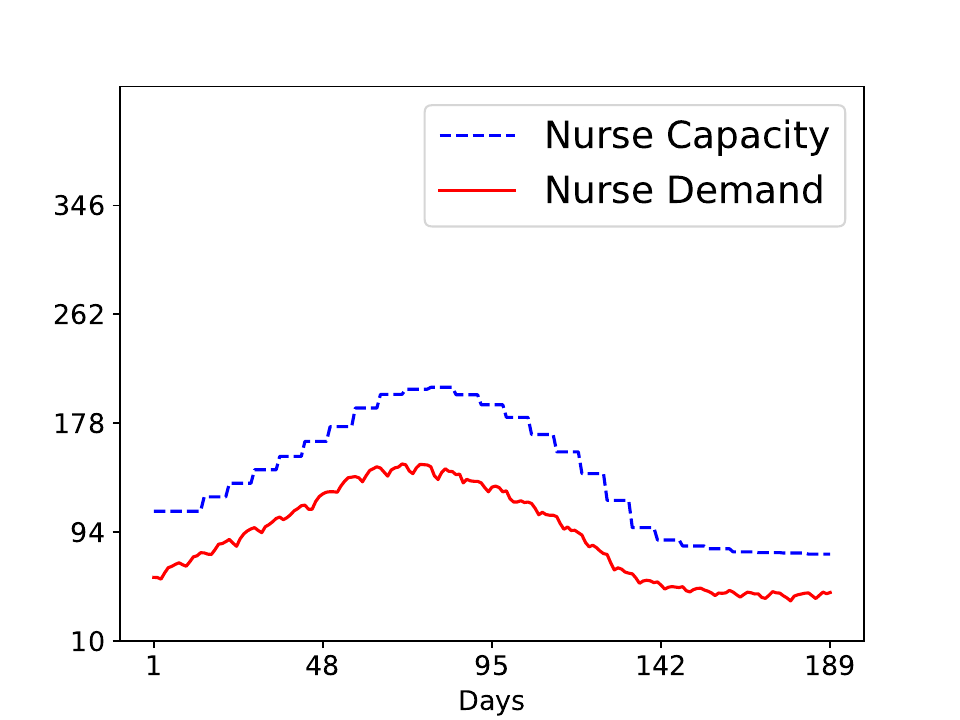}
         \caption{South Hospital} 
         \label{census_Bloomington} 
     \end{subfigure}
     \hfill 
     \begin{subfigure}[b]{0.23\textwidth}
         \centering
         \includegraphics[width=\textwidth]{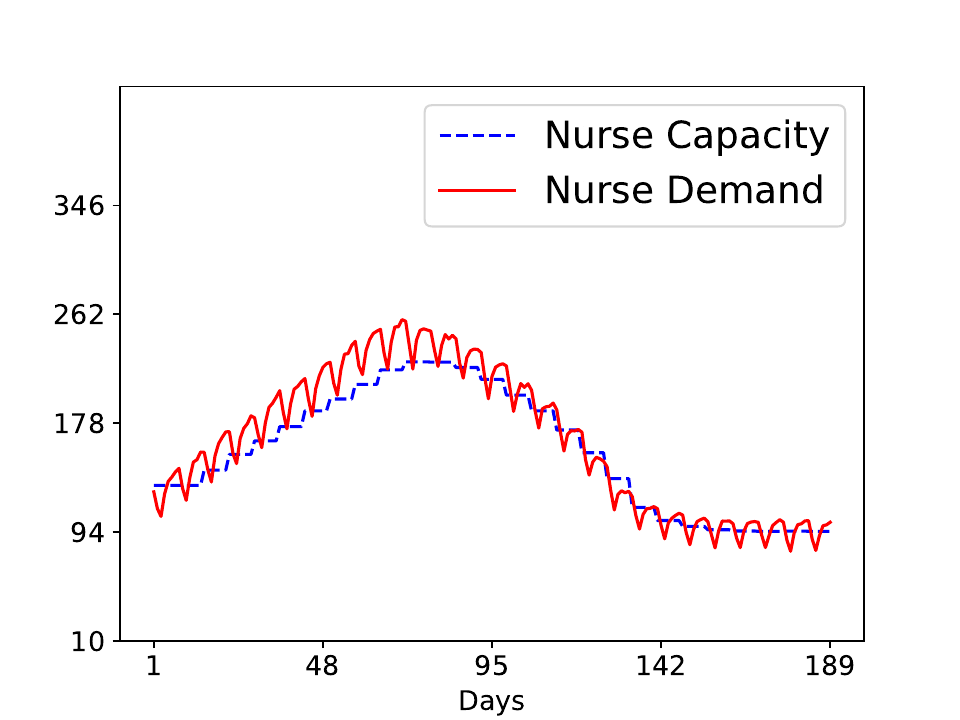}   
         \caption{Central Hospital} 
         \label{census_Methodist}     
     \end{subfigure} 
         \caption{Daily nurse demand and capacity at each hospital}  
        \label{nurse_demand_capacity_over_time} 
\end{figure}

\noindent\textbf{Cost Parameters. } 
We use a daily payment premium $p=1$ as the baseline and normalize other costs accordingly. We set the cancellation fee percentage at 5\%, and the daily salary premium multiplier for emergency transfer $\theta_t$ at $1.6$. We set the non-salary transfer cost (reimbursement for nurses who travel to a remote hospital), $\tau_{ij}$, based on the geographical distance between different hospitals, $i,j\in [L]$. Our secondment setting is also based on the distance between locations, with longer distances resulting in longer secondments. Specifically, secondments between nearby locations (e.g., West and Central, East and Central) have a duration of 1 day, while those between more distant locations (e.g., West and East, West and South) have a duration of 2 days, due to logistical challenges. This is referred to as the base secondment scenario. Table \ref{tab:parameters} summarizes the parameters.  
%see details in Appendix~\ref{sec:SyntheticDataGeneration}. 

% Table generated by Excel2LaTeX from sheet 'Sheet6'      
\begin{table}[htbp]
  \centering 
\scalebox{0.75}{
    \begin{tabular}{|c|l|rrrc|}
    \hline
    Hospital From/To & Parameters & \multicolumn{1}{c}{West} & \multicolumn{1}{c}{East} & \multicolumn{1}{c}{South} & Central \bigstrut\\
    \hline
    \multirow{3}[2]{*}{West} & \multicolumn{1}{c|}{Cost} & \multicolumn{1}{c}{-} & 1.46  & 1.68  & \multicolumn{1}{r|}{1.20} \bigstrut[t]\\
          & Distance & \multicolumn{1}{c}{-} & 88    & 110   & \multicolumn{1}{r|}{62} \\
          & Secondment & \multicolumn{1}{c}{-} & 2     & 2     & \multicolumn{1}{r|}{1} \bigstrut[b]\\
    \hline
    \multirow{3}[2]{*}{East} & \multicolumn{1}{c|}{Cost} & 1.46  & \multicolumn{1}{c}{-} & 1.70  & \multicolumn{1}{r|}{1.14} \bigstrut[t]\\
          & Distance & 88    & \multicolumn{1}{c}{-} & 112   & \multicolumn{1}{r|}{56} \\
          & Secondment & 2     & \multicolumn{1}{c}{-} & 2     & \multicolumn{1}{r|}{1} \bigstrut[b]\\
    \hline
    \multirow{3}[2]{*}{South} & \multicolumn{1}{c|}{Cost} & 1.68  & 1.70  & \multicolumn{1}{c}{-} & \multicolumn{1}{r|}{1.10} \bigstrut[t]\\
          & Distance & 110   & 112   & \multicolumn{1}{c}{-} & \multicolumn{1}{r|}{52} \\
          & Secondment & 2     & 2     & \multicolumn{1}{c}{-} & \multicolumn{1}{r|}{1} \bigstrut[b]\\
    \hline
    \multirow{3}[2]{*}{Central} & \multicolumn{1}{c|}{Cost} & 1.20  & 1.14  & 1.10  & - \bigstrut[t]\\
          & Distance & 62    & 56    & 52    & - \\
          & Secondment & 1     & 1     & 1     & - \bigstrut[b]\\
    \hline
    \end{tabular}%
} 
  \caption{Transfer cost, transfer distance  and secondment setting}  
  \label{tab:parameters}%
\end{table}%

To incentivize hospitals to opt for emergency transfers instead of enduring nurse shortages (particularly when nurses are available at other healthcare facilities), we need to ensure that the unit shortage cost $s_t^i$ exceeds the total emergency cost incurred by the emergency transfer. We set $s_t^i=15$ that is greater than the maximum emergency transfer cost.  

Within our study, we determine the robust parameter for week \( w \) as follows. For the first week, we set the robust parameter \( \epsilon_{N,1} \) to zero. For week \( w \) (\( w \geq 2 \)), we select the robust parameter \( \epsilon_{N,w} \) by choosing the best-performing parameter from the set \( \{ (\epsilon_{N,w-1} - 5\upsilon)^+, (\epsilon_{N,w-1} - 5\upsilon)^+ +5,(\epsilon_{N,w-1} - 5\upsilon)^+ +10 , \ldots, \epsilon_{N,w-1} + 5\upsilon \} \) based on the results from week \( w-1 \). Here, \( \upsilon \) can be adjusted to reflect changes in nurse demand. 
{
For example, we set \( \upsilon = 2 \) for the demand pattern described in Section~\ref{sec:DataDescription}, and increase it to \( \upsilon = 5 \) under scenarios with more pronounced demand changes described in Appendix~\ref{sec:Impact_higher_demand}.}

%when using the previous three weeks for parameter estimation. 
%if we use the past three (or six) weeks' data for prediction. 
%and \( \upsilon = 5 \) for the nurse demand pattern with perturbation. 

%%%%%%%%%%%%%%%%%%%%%%%%%%

\subsection{Network Structure}
%Hub-and-spoke Network and Fully Connected Network    
\label{sec:network_design}  
 
\begin{figure}
     \centering
     \begin{subfigure}[b]{0.3\textwidth}
         \centering
         \includegraphics[width=\textwidth]{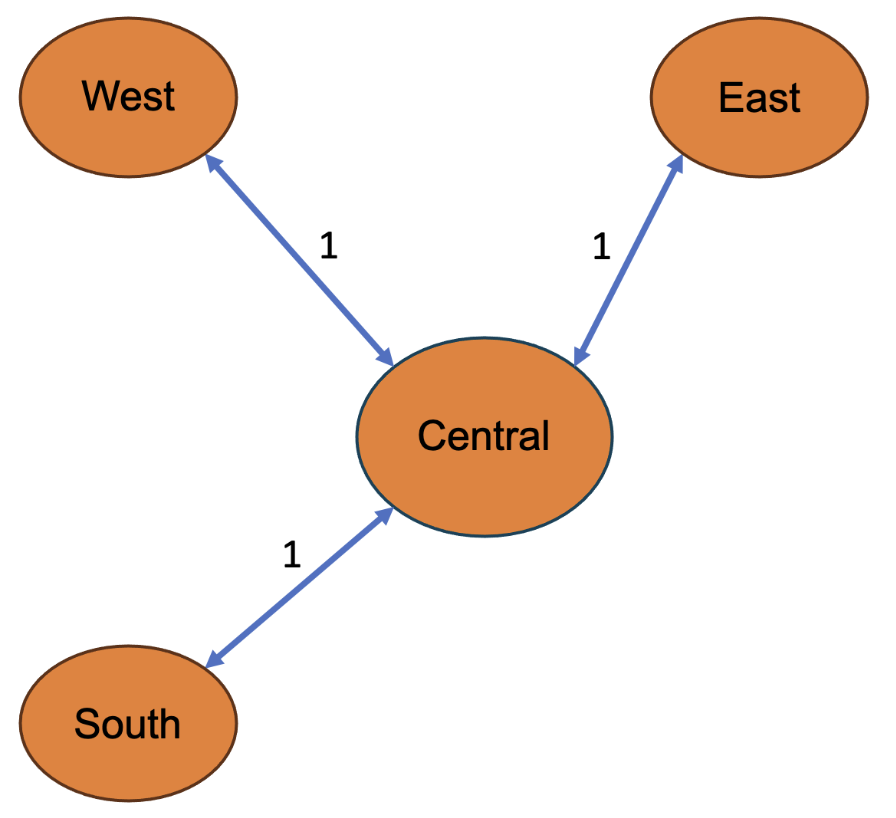} 
         \caption{Hub-and-spoke network }
         \label{fig_HS}    
     \end{subfigure}
     \quad\quad
     \begin{subfigure}[b]{0.3\textwidth}
         \centering
         \includegraphics[width=\textwidth]{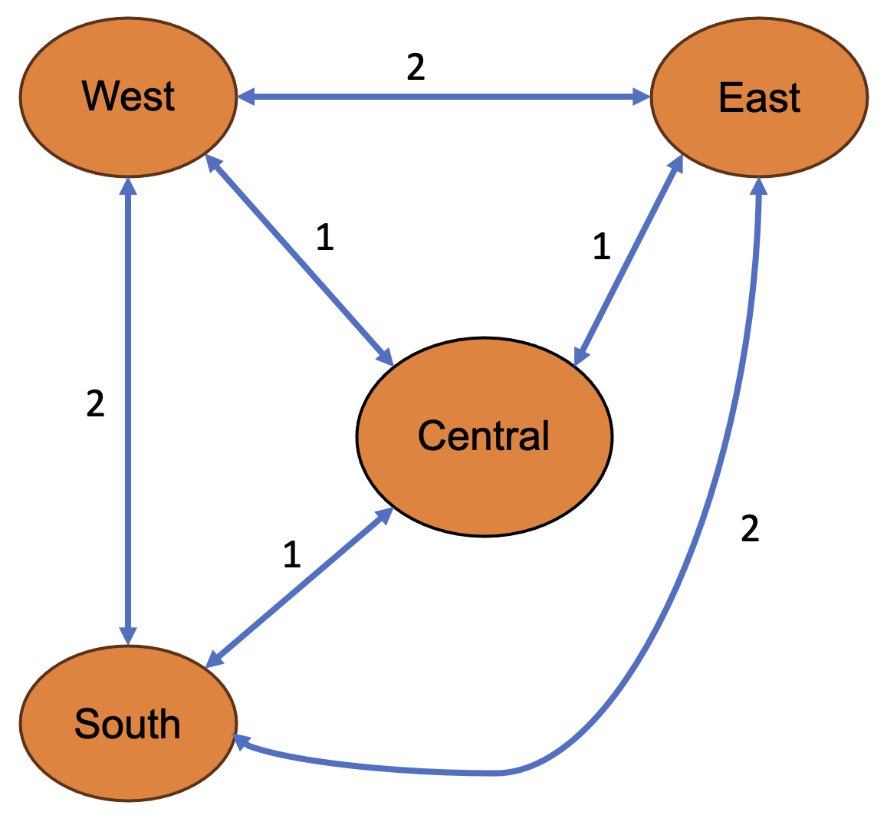}   
         \caption{Fully connected network} 
         \label{fig_FC}    
     \end{subfigure} 
         \caption{Baseline secondment setting under different network designs.}   
         %  (\color{blue} { Draw the figure in square; x-axis: ``Weeks'' (from week 1);  y-axis:  ``Total cost'';  label:   ``One-day secondment scenario '', ``Base secondment scenario'' } ) 
        \label{fig:network_design}
\end{figure}

In this section, we evaluate the impact of network designs on the performance of the nurse redeployment program, comparing two structures: the hub-and-spoke (HS) network and the fully connected (FC) network. In the HS network, nurse transfers are restricted to occur between rural hospitals (West, East, and South -- the spokes) and the Central Hospital (the hub), whereas the FC network allows transfers between any two locations. We use the \textit{baseline secondment} scenario considering the long distances between rural areas, i.e., the secondment duration for transfers between spokes is set to two days, whereas the secondment for transfers between hub and spoke is one day; see Figure \ref{fig:network_design} for the network designs and the secondment settings. 
In addition, we consider two alternative secondment scenarios. 
One is the \textit{one-day secondment} scenario, where all secondments are set to a duration of one day; the other one is the \textit{three-day secondment} scenario, where the secondment between nearby locations is still one day, but the secondment between the locations with long distance is set as three days.
% \red{Did we also do secondment days constant at 2 or 3 days, without the 1-day exception? If yes, we need better names for these secondment scenarios.}
% {\color{blue}We only consider these three secondment scenarios.} 

Table~\ref{tab:Average_metrics} compares the average (weekly) costs under the two designs, over the entire horizon, with different secondment scenarios and the two optimization methods. The results show a significantly lower total cost for the FC network {(reduced by 23\% to 35\%)}, regardless of the method used. We also examine two additional performance metrics: the average number of nurse transfers and the average miles traveled (based on the distance between locations) per week across all locations. Table~\ref{tab:Average_metrics} shows that both metrics decrease significantly under the FC network {(reduced by 13\% to 52\%)}, consistent with the lower costs compared to the HS network. 

% (\red{is this weekly cost or the total cost over all 27 weeks?}{\color{blue}: This is the average weekly cost.})
%This trend is also evident in the weekly cost comparisons in Figures~\ref{total_cost_diff_sec_dro_saa}.

% (Figures~\ref{Deployed_nurse_transfers_Arnett_Methodist_SAA} and \ref{Deployed_nurse_transfers_Arnett_Methodist_DRO}) 

% (Figures~\ref{Deployed_nurse_transfers_Arnett_Methodist_SAA_sec_2} and \ref{Deployed_nurse_transfers_Arnett_Methodist_DRO_sec_2} for the base secondment scenario, and Figures~\ref{Deployed_nurse_transfers_Arnett_Methodist_SAA_sec_1} and \ref{Deployed_nurse_transfers_Arnett_Methodist_DRO_sec_1} for the one-day secondment scenario)

% , and Figures ~\ref{Deployed_nurse_transfers_Arnett_Methodist_SAA} to~\ref{Deployed_nurse_transfers_Arnett_Methodist_DRO_sec_1} are included in Appendix \ref{sec:figures}. 

% Table generated by Excel2LaTeX from sheet 'Sheet5'
\begin{table}[htbp]
  \centering
  \caption{Average cost, deployed transfers, and transferred miles per week for SAA and SRO across various network designs and secondment scenarios.} 
\scalebox{0.8}{
    \begin{tabular}{|c|c|rrcccc|}
    \hline
    \multirow{2}[4]{*}{Metrics} & \multicolumn{1}{c|}{\multirow{2}[4]{*}{Network design}} & \multicolumn{2}{c|}{Baseline secondment} & \multicolumn{2}{c|}{One-day secondment } & \multicolumn{2}{c|}{Three-day secondment } \bigstrut\\
\cline{3-8}          &       & \multicolumn{1}{l|}{SAA} & \multicolumn{1}{l|}{SRO} & \multicolumn{1}{l|}{SAA} & \multicolumn{1}{l|}{SRO} & \multicolumn{1}{l|}{SAA} & \multicolumn{1}{l|}{SRO} \bigstrut\\
    \hline
    \multicolumn{1}{|c|}{\multirow{2}[4]{*}{Cost}} & FC    & 641.90 & 633.05 & \multicolumn{1}{r}{713.30} & \multicolumn{1}{r}{702.20} & \multicolumn{1}{r}{607.52} & \multicolumn{1}{r|}{606.98} \bigstrut\\
\cline{2-2}          & HS    & 933.35 & 914.71 & -     & -     & -     & - \bigstrut\\
    \hline
    \multicolumn{1}{|c|}{\multirow{2}[4]{*}{Deployed transfers}} & FC    & 149.05 & 148.86 & \multicolumn{1}{r}{195.63} & \multicolumn{1}{r}{195.62} & \multicolumn{1}{r}{138.64} & \multicolumn{1}{r|}{139.75} \bigstrut\\
\cline{2-2}          & HS    & 291.69 & 292.04 & -     & -     & -     & - \bigstrut\\
    \hline
    \multicolumn{1}{|c|}{\multirow{2}[4]{*}{Transferred miles}} & FC    & 10137.89 & 10326.54 & \multicolumn{1}{r}{14159.54} & \multicolumn{1}{r}{14289.09} & \multicolumn{1}{r}{9309.35} & \multicolumn{1}{r|}{9529.61} \bigstrut\\
\cline{2-2}          & HS    & 16399.95 & 16438.94 & -     & -     & -     & - \bigstrut\\
    \hline
    \end{tabular}%
} 
  \label{tab:Average_metrics}%
\end{table}%

\begin{figure}
     \centering
     \begin{subfigure}[b]{0.35\textwidth}
         \centering
         \includegraphics[width=\textwidth]{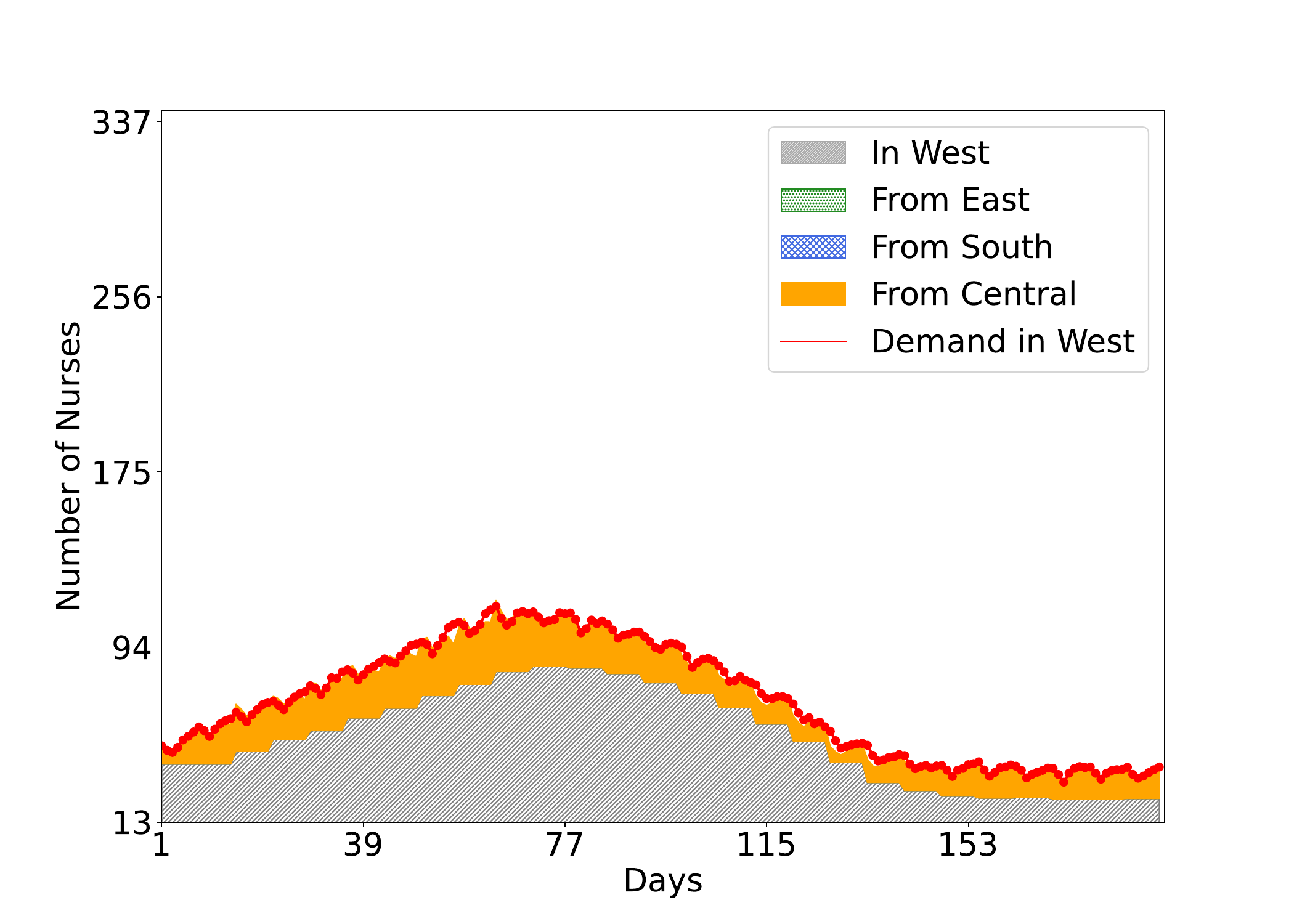}
          \caption{Deployed transfers to West Hospital } 
         % (\color{blue} { x-axis: ``Days'' (from day 1);  y-axis:  ``Number of nurses'';  label: ``Nurse demand at Arnett'', ``Nurses left with home location at Arnett'', ``Nurse transfers from Bloomington'', ``Nurse transfers from Ball'', and ``Nurse transfers from Methodist''; For the labels, use a shaded rectangle for the label; remove the line for nurse transfer and only use the shaded area; for nurse demand, may use red line directly if we a square figure. } )
         \label{arnett_shortage_deployment_dro_sec_1_min_sec_1} 
     \end{subfigure}
     \quad
     \begin{subfigure}[b]{0.35\textwidth} 
         \centering
         \includegraphics[width=\textwidth]{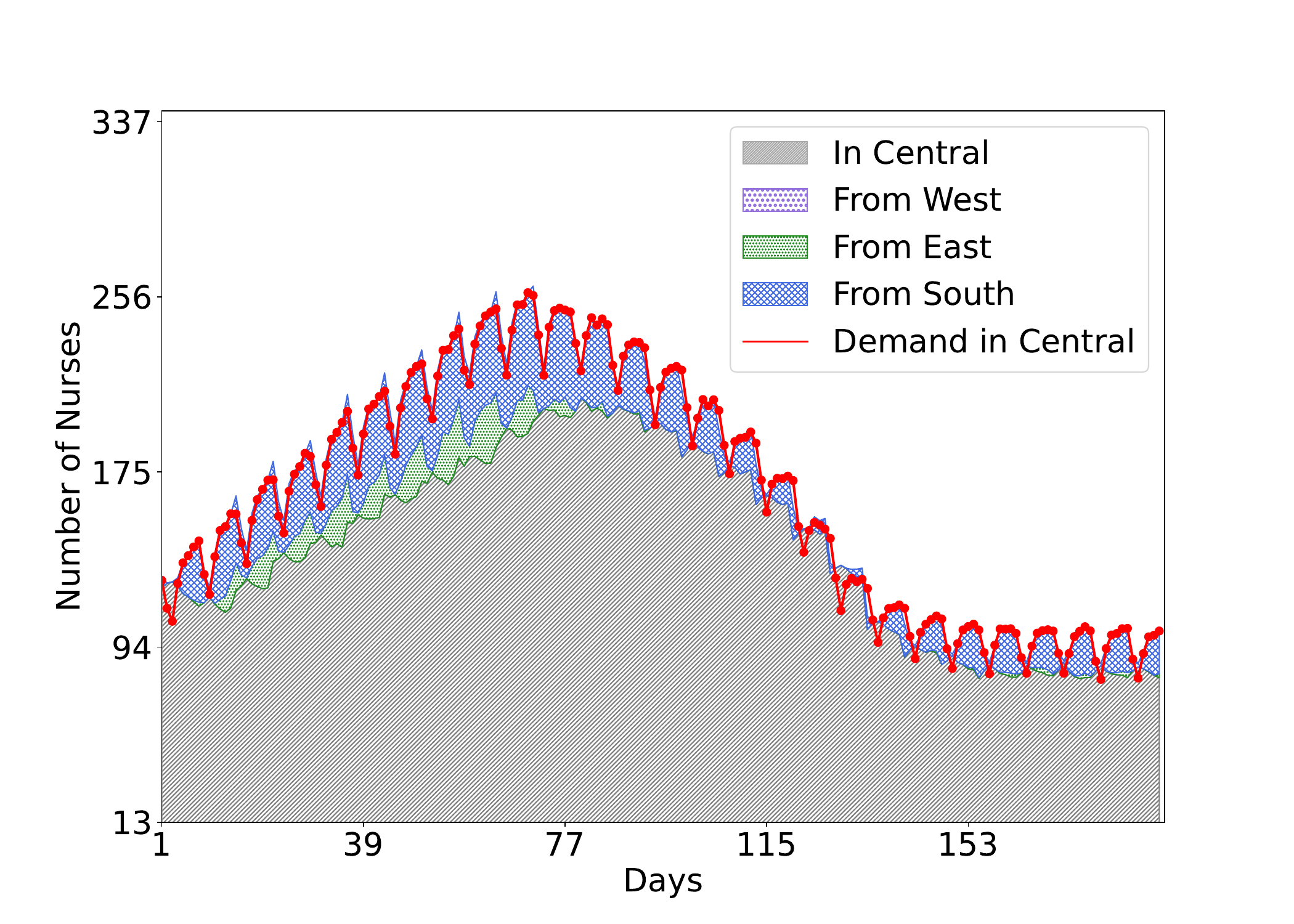}   
          \caption{Deployed transfers to Central Hospital} 
         \label{methodist_shortage_deployment_dro_sec_1_min_sec_1}    
     \end{subfigure}  
         \caption{Daily deployed nurse transfers to West and Central Hospitals by using SRO under the hub-and-spoke (HS) network. }  
        \label{Deployed_nurse_transfers_Arnett_Methodist_DRO} 
\end{figure}

\begin{figure}
     \centering
     \begin{subfigure}[b]{0.35\textwidth}
         \centering
         \includegraphics[width=\textwidth]{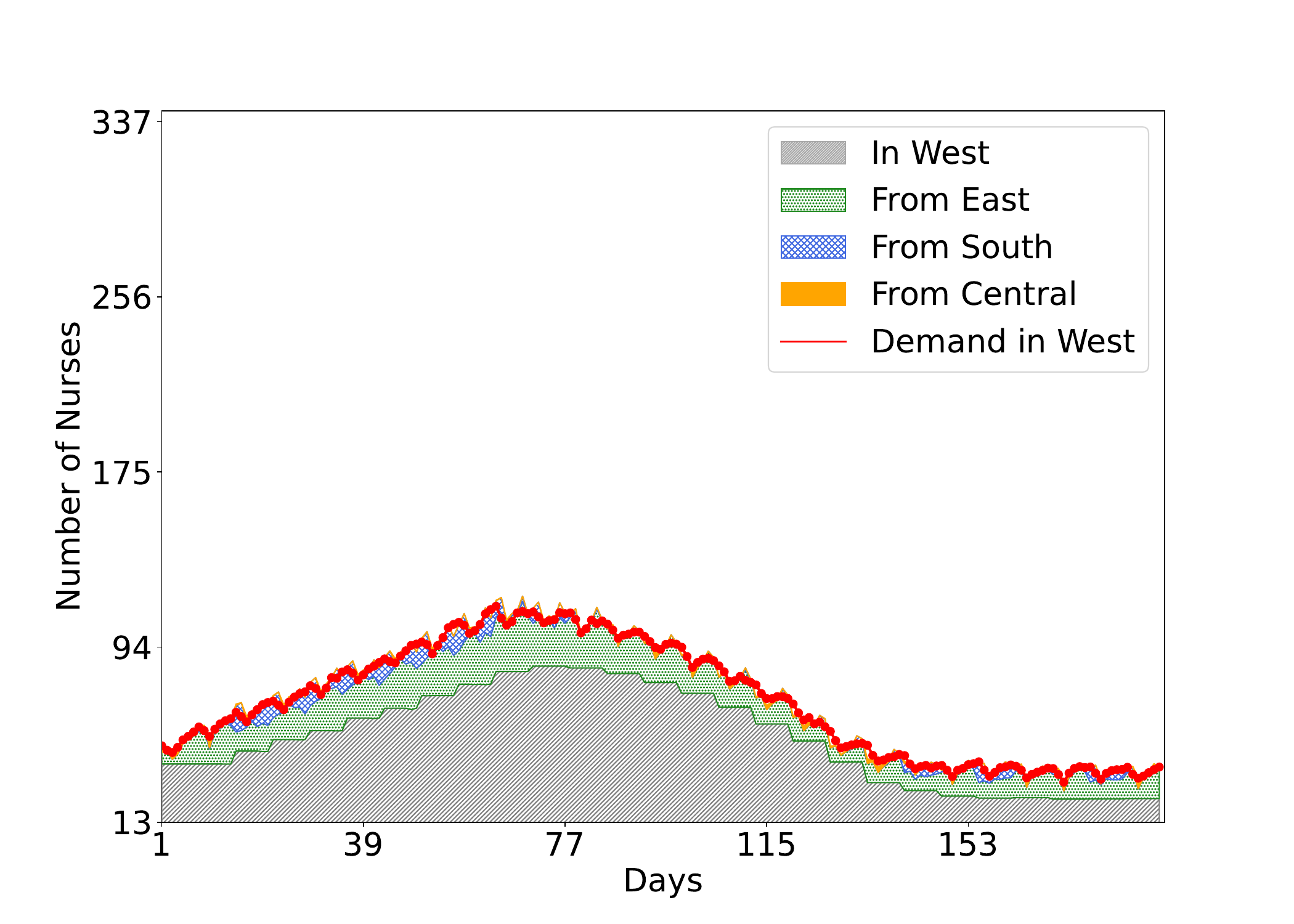}
       \caption{Deployed transfers to West Hospital} 
         % (\color{blue} { x-axis: ``Days'' (from day 1);  y-axis:  ``Number of nurses'';  label: ``Nurse demand at Arnett'', ``Nurses left with home location at Arnett'', ``Nurse transfers from Bloomington'', ``Nurse transfers from Ball'', and ``Nurse transfers from Methodist''; For the labels, use a shaded rectangle for the label; remove the line for nurse transfer and only use the shaded area; for nurse demand, may use red line directly if we a square figure. } )
         \label{arnett_shortage_deployment_dro_sec_2_min_sec_2} 
     \end{subfigure}
     \quad
     \begin{subfigure}[b]{0.35\textwidth}
         \centering
         \includegraphics[width=\textwidth]{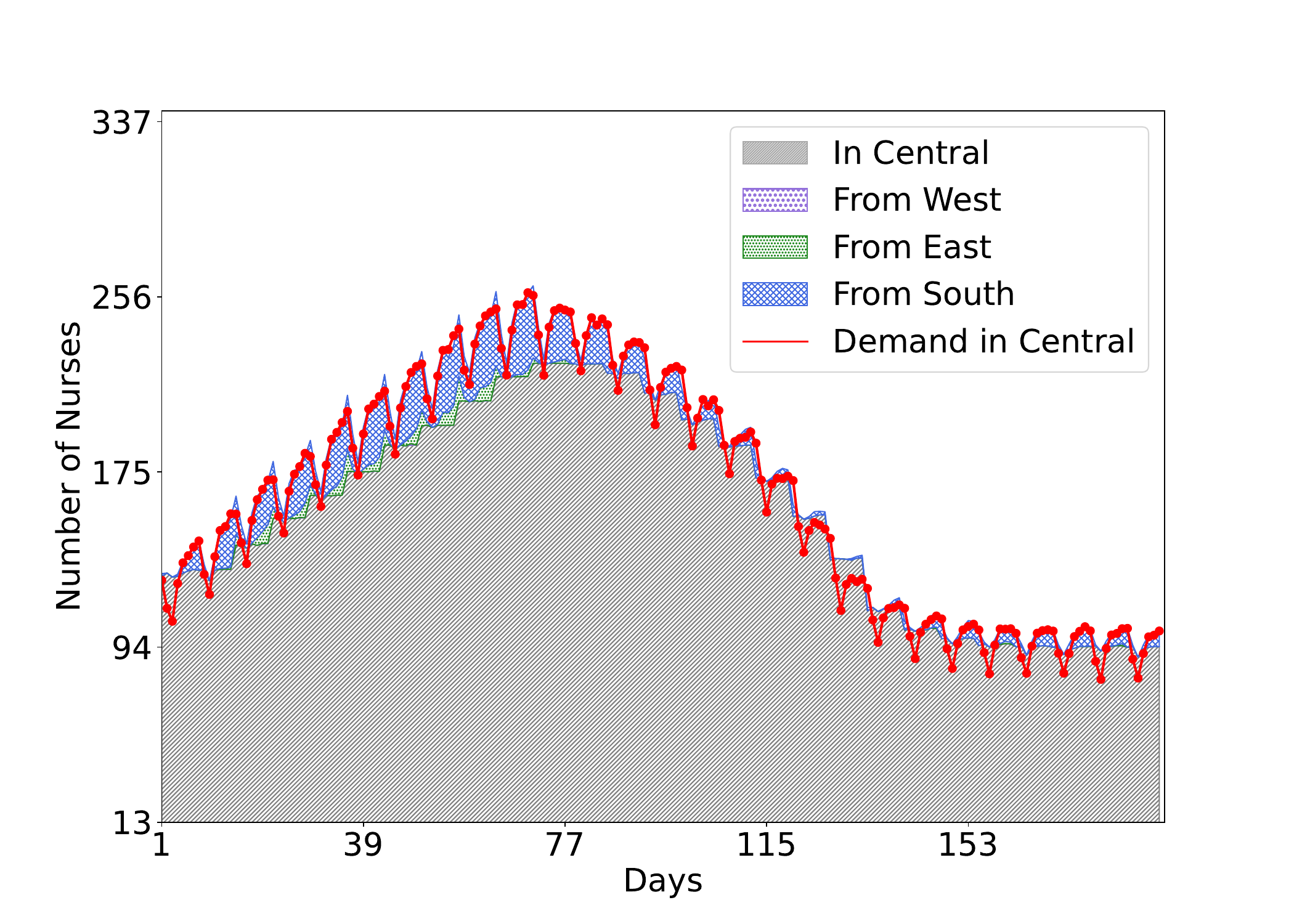}   
          \caption{Deployed transfers to Central Hospital} 
         %(\color{blue} { x-axis: ``Days'' (from day 1);  y-axis:  ``Number of nurses'';  label: ``Nurse demand at Methodist'', ``Nurses left with home location at Methodist'', ``Nurse transfers from Arnett'', ``Nurse transfers from Ball'', and ``Nurse transfers from Bloomington''; For the labels, use a shaded rectangle for the label; remove the line for nurse transfer and only use the shaded area; for nurse demand, may use red line directly if we a square figure. } )
         \label{methodist_shortage_deployment_dro_sec_2_min_sec_2}    
     \end{subfigure} 
         \caption{Daily deployed nurse transfers to West and Central Hospitals  by using SRO under the base secondment scenario and fully connected (FC) network.}  
        \label{Deployed_nurse_transfers_Arnett_Methodist_DRO_sec_2} 
\end{figure}

To explore the reasons behind the advantage demonstrated by the FC design, we analyze (i) the number of transfers from West, South, and East Hospitals to Central, and (ii) the number of transfers from Central, South and East to West, considering the nurse shortages at both West and Central Hospitals. Figure~\ref{Deployed_nurse_transfers_Arnett_Methodist_DRO} shows that in the HS network, transfers occur from Central to West, with Central receiving nurses from South and East. In contrast, Figure~\ref{Deployed_nurse_transfers_Arnett_Methodist_DRO_sec_2} illustrates that in the FC network, nurses can be transferred directly between East and West, and from South to Central. A similar pattern is observed across different secondment scenarios and methods.
%We analyze the comparison under different network designs for SRO and base secondment scenario in detail, and similar phenomenon can be observed for SAA and the other secondment scenario. 
%This reduction is largely due to the more efficient direct transfers from hospitals with surpluses to those with shortages. In other words, the FC network not only reduces costs but also minimizes daily transfers and significantly cuts travel distance and time, highlighting the broader benefits of this design.

To summarize, the FC network outperforms the HS design in this setting for two key reasons. First, it enables direct transfers to hospitals with shortage (e.g., West), which avoids costly intermediary routing through the hub (Central). Second, it leverages longer secondments between distant hospitals (West, South, and East), which reduce average daily transfer costs, as the fixed costs such as relocation bonuses are spread over multiple days.
Notably, the cost gap between the two network designs is more substantial than the differences observed across secondment scenarios or optimization methods, suggesting a dominant role of network structure. However, this does \emph{not} imply that the FC design is universally superior. As we discuss next, its performance depends critically on the alignment between network structure and secondment policies.

%%%%%%%%%%%%%%%%%%%%%%%%%%%%%%%%%
%%%%%%%%%%%%%%%%%%%%%%%%%%%%%%%%%
%%%%%%%%%%%%%%%%%%%%%%%%%%%%%%%%%

\subsection{Effect of Secondment}  \label{sec:Secondment_effect}  

In this section, we study the impact of secondment length. Our main finding is that, if secondment lengths are not properly adapted to the network structure, system performance can deteriorate. Specifically, assigning too short secondments to long-distance transfers may lead to excessive commuting costs or nurse fatigue, which counter-effects the benefit from the improved network connectivity (Section~\ref{sec:ASpecialCase}). Conversely, overly long secondments, while reducing the per-day fixed costs by spreading them out, can reduce responsiveness and limit the ability to reallocate staff when demand shifts (Section~\ref{sec:Secondment_effect1}). Thus, secondment policies must be carefully calibrated to the network structure to balance efficiency, flexibility, and nurse well-being.

\subsubsection{Interplay with Network Structure}  
\label{sec:ASpecialCase}
Under the baseline demand pattern, the FC network consistently outperforms the HS network across all settings. However, this advantage is not universal. In this section, we examine an alternative demand scenario in which only the West Hospital experiences a nurse shortage, while the other three hospitals have surplus capacity (see details in Appendix~\ref{sec:figures}). Moreover, we assume that, under this setting, transfers are restricted to occur only between the West and Central Hospitals in the HS network, and between the West and South Hospitals in the FC network. 
%\red{(move Table~\ref{tab:Average_metrics_special_case} to Appendix~\ref{sec:figures}.)} 
As shown in Table~\ref{tab:Average_metrics_special_case} in Appendix~\ref{sec:figures}, when the secondment length for rural-to-rural transfers is set to one day, the HS network achieves lower total costs—359 (SAA) and 353 (SRO)—compared to the FC network’s 414 (SAA) and 409 (SRO). HS also results in shorter total transfer distances, with total transfer volumes comparable across both network structures. 
However, when the secondment length increases to two or more days, the FC network’s cost drops below 329, once again outperforming the HS network. The FC network also achieves a reduction in total transfer volume, although the total transfer distance may remain slightly higher than that of the HS network at a two-day secondment length.

These findings highlight the importance of allowing longer secondment durations when having rural-to-rural nurse transfers in the FC network. Moreover, the FC design introduces more complexity in coordination, as it requires nurses to travel to and work at more locations (and nurses need to be familiar with multiple workplaces). In Appendix~\ref{sec:coordination_cost}, we further analyze the impact of such coordination complexity—modeled as additional coordination costs—on the performance of both HS and FC networks.

\subsubsection{Impact of Secondment Length} \label{sec:Secondment_effect1}  

In this section, we take a deep dive to investigate the impact of different secondment length, {focusing on the fully connected (FC) network design.} 
Under the baseline setting, we observe a monotonic decrease in total cost as the secondment length increases (Table~\ref{tab:Average_metrics}). This trend is mainly driven by the relatively high base transfer cost (set at 1.1), which makes longer secondments more cost-effective by spreading the cost of relocation over multiple days.
To explore the sensitivity of this trend, we reduce the transfer cost to one-tenth of its original value (0.1), following the calibration in Appendix~\ref{sec:TravelCostSetup}. Under this lower-cost setting, we uncover an interesting \textit{U-shaped} relationship between secondment length and system cost.

\begin{figure}
     \centering
     \begin{subfigure}[b]{0.35\textwidth}
         \centering
         \includegraphics[width=\textwidth]{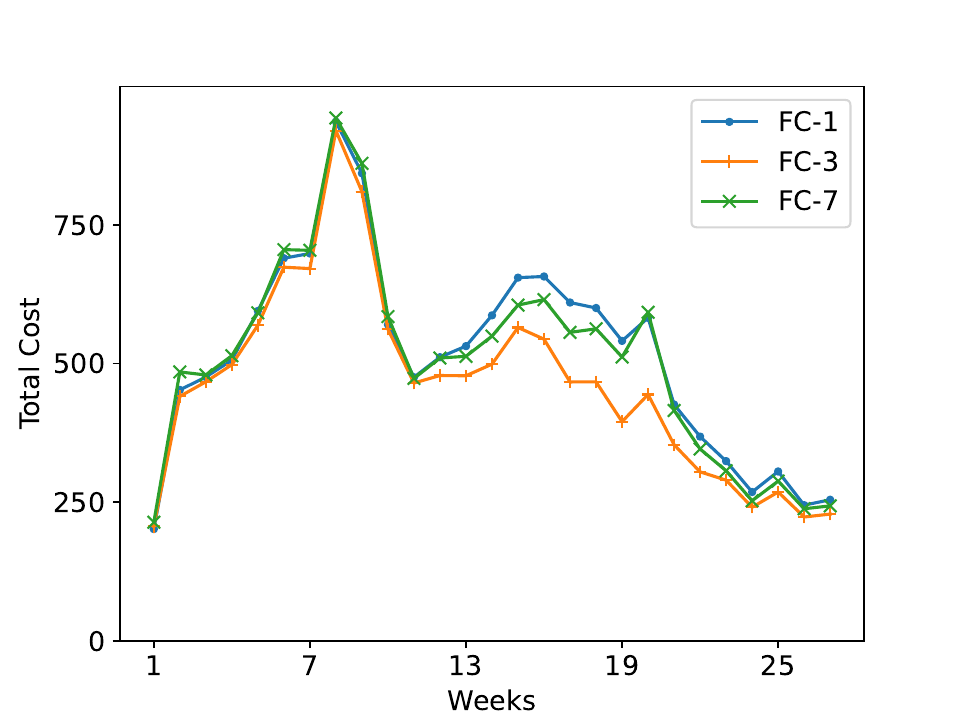}
         \caption{SAA}
         \label{total_cost_diff_sec_robust_0_small_transfer_cost}
     \end{subfigure}
     \small 
     \begin{subfigure}[b]{0.35\textwidth}
         \centering
         \includegraphics[width=\textwidth]{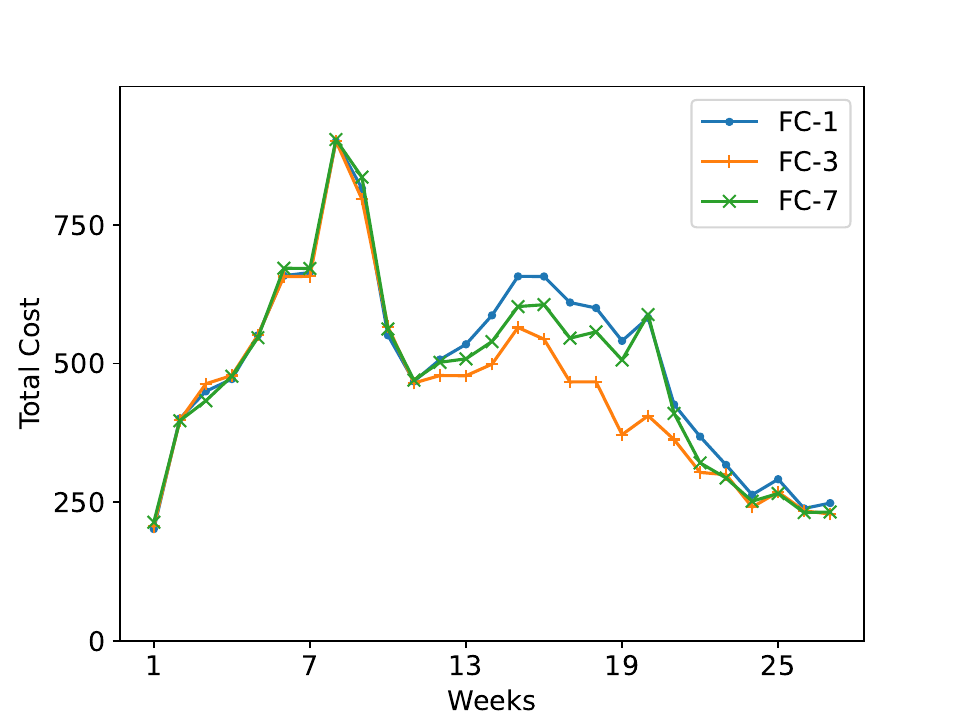}   
         \caption{SRO} 
         \label{total_cost_diff_sec_robust_best_small_transfer_cost} 
     \end{subfigure} 
         \caption{Effect of secondment on the cost as a function of weeks, using base transfer cost 0.1. Here, ``FC-1'' denotes the one-day secondment scenario,  ``FC-3'' represents the three-day secondment scenario, and ``FC-7'' represents the seven-day secondment scenario.} 
         %  (\color{blue} { Draw the figure in square; x-axis: ``Weeks'' (from week 1);  y-axis:  ``Total cost'';  label:   ``One-day secondment scenario '', ``Base secondment scenario'' } ) 
        \label{total_cost_diff_sec_dro_saa_small_transfer_cost} 
\end{figure}

Specifically, Table~\ref{tab:Average_metrics_smaller_transfer_cost} reports the average weekly cost, number of actual nurse transfers, and total transfer mileage under varying secondment scenarios for both SAA and SRO methods (using the lower transfer cost). To isolate the impact of rural-to-rural secondments, we fix the secondment length from central to rural hospitals at one day while varying the rural-to-rural secondment length across 1, 3, and 7 days. We observe from the table that increasing secondment length from 1 to 3 days reduces costs, but extending it further to 7 days leads to an increase. This U-shaped cost pattern is further illustrated in Figure~\ref{total_cost_diff_sec_dro_saa_small_transfer_cost}, which plots weekly cost trends under different secondment lengths for SAA and SRO.

These results indicate the existence of a ``sweet spot'' for optimal secondment duration. The underlying trade-off behind the U-shape is clear: shorter secondments result in high daily transfer costs due to frequent reassignments, while excessively long secondments reduce flexibility and responsiveness in re-deploying nurses to adjust to the changing demand.
%This U-shape trend reflects a trade-off: a shorter secondment duration leads to high daily average transfer costs due to frequent reassignments, whereas excessively long secondments diminish responsiveness and limit flexibility in reallocating staff under demand fluctuations.  
Additionally, we observe that both the average number of transfers and the total mileage decrease as secondment length increases. This is because longer secondments reduce the frequency of travel—each nurse remains at the assigned hospital longer, thereby spreading travel distances across multiple days rather than requiring daily commutes.

% Table generated by Excel2LaTeX from sheet 'Sheet5'
\begin{table}[htbp]
  \centering
  \caption{Average cost, deployed transfers, and transferred miles per week for SAA and SRO across various secondment scenarios, using base transfer cost 0.1.}
    \begin{tabular}{|p{8.915em}|r|r|r|r|r|r|}
    \hline
    \multicolumn{1}{|c|}{\multirow{2}[4]{*}{Metrics}} & \multicolumn{2}{c|}{One-day secondment} & \multicolumn{2}{c|}{Three-day secondment } & \multicolumn{2}{c|}{Seven-day secondment} \bigstrut\\
\cline{2-7}    \multicolumn{1}{|c|}{} & \multicolumn{1}{l|}{SAA} & \multicolumn{1}{l|}{SRO} & \multicolumn{1}{l|}{SAA} & \multicolumn{1}{l|}{SRO} & \multicolumn{1}{l|}{SAA} & \multicolumn{1}{l|}{SRO} \bigstrut\\
    \hline
    Cost  & 515.29 & 501.72 & 464.13 & 456.78 & 506.07 & 486.77 \bigstrut\\
    \hline
    Deployed transfers  & 195.66 & 195.71 & 140.22 & 141.27 & 107.25 & 127.12 \bigstrut\\
    \hline
    Transferred miles & 14131.70 & 14257.93 & 9277.61 & 9553.46 & 7733.71 & 7884.81 \bigstrut\\
    \hline
    \end{tabular}%
  \label{tab:Average_metrics_smaller_transfer_cost}%
\end{table}%

\subsection{Value of Robustness} \label{sec:ComparisonofPolicies}  
From Table~\ref{tab:Average_metrics}, we observe that the SRO method consistently achieves a lower weekly average cost compared to SAA across different network designs and secondment scenarios. To understand when and under what demand patterns this performance gap becomes most significant, we conduct a more detailed phase-by-phase analysis in this section.

\begin{figure} 
     \centering
     \begin{subfigure}[b]{0.3\textwidth} 
         \centering
         \includegraphics[width=\textwidth]{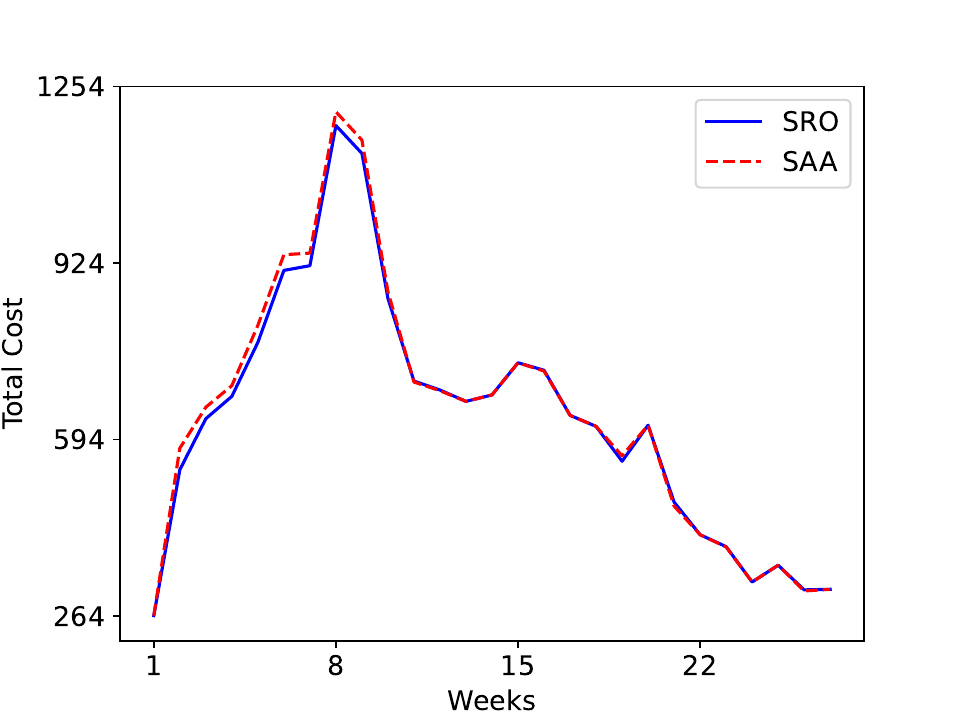} 
         \caption{Comparison on the cost  between SAA and SRO}
         %(\color{blue} { Draw the figure in square; x-axis: ``Weeks'' (from week 1);  y-axis:  ``Total cost'';  label: ``${SAA}$'',   ``${SRO}$'' } ) 
         \label{methodist_shortage_total_cost_sec_3_min_sec_2}
     \end{subfigure}
    \quad 
     \begin{subfigure}[b]{0.3\textwidth} 
         \centering
         \includegraphics[width=\textwidth]{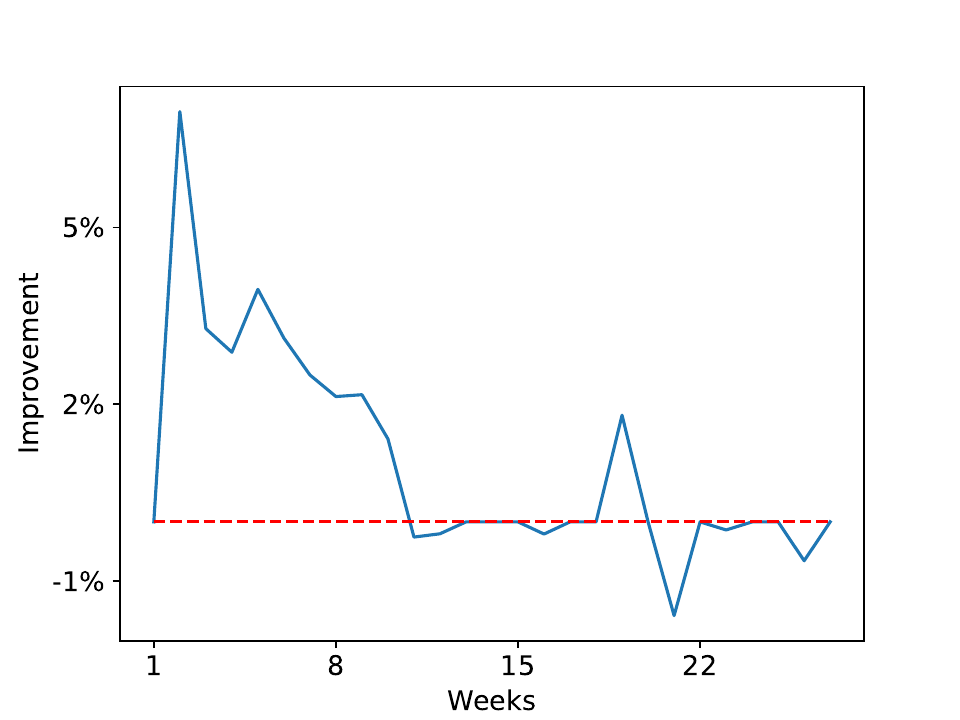}
         \caption{Improvement of SRO over SAA}  
         %  (\color{blue} { Draw the figure in square; x-axis: ``Weeks'' (from week 1);  y-axis:  ``Cost Improvement in \% of DRO Over SAA'';  no label ; draw a dashed line with y-axis$=0$ } )
         \label{methodist_shortage_total_cost_percentage_sec_3_min_sec_2} 
     \end{subfigure}   
    \quad 
     \begin{subfigure}[b]{0.3\textwidth} 
         \centering
         \includegraphics[width=\textwidth]{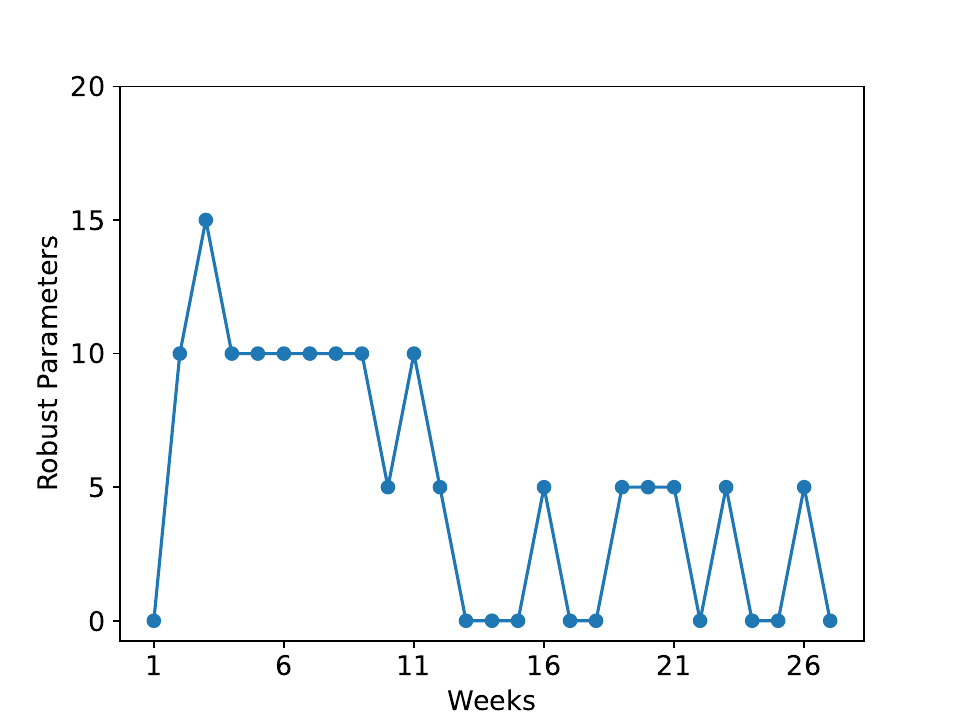}
            \caption{Robust parameter used by SRO}   
         %  (\color{blue} { Draw the figure in square; x-axis: ``Weeks'' (from week 1);  y-axis:  ``Cost Improvement in \% of DRO Over SAA'';  no label ; draw a dashed line with y-axis$=0$ } )
         \label{fully_sec2_weekly_robust_parameter_selection_plot} 
     \end{subfigure}        
        \caption{Weekly comparison between SAA and SRO, and robust parameter used by SRO}   
        \label{Comparison_DRO_SAA}
\end{figure}

We begin by comparing the weekly costs under the FC network and the baseline secondment scenario. Figure~\ref{methodist_shortage_total_cost_sec_3_min_sec_2} shows the weekly costs for both SRO and SAA, and  Figure~\ref{methodist_shortage_total_cost_percentage_sec_3_min_sec_2} illustrates the percentage improvement of SRO over SAA. During the early phase of the pandemic—characterized by increasing demand (before Week 10)—SRO shows a 2–6\% cost reduction. In contrast, during the stable or declining demand phases (Week 10 onward), the performance of SRO and SAA becomes nearly identical (except in weeks 19-20). This latter observation is further supported by the weekly robust parameter selection shown in Figure~\ref{fully_sec2_weekly_robust_parameter_selection_plot}, where SRO selects a robust parameter value of 0 in most of those weeks after Week 12, effectively reducing to SAA.

\begin{figure}
     \centering
     \begin{subfigure}[b]{0.35\textwidth}
         \centering
         \includegraphics[width=\textwidth]{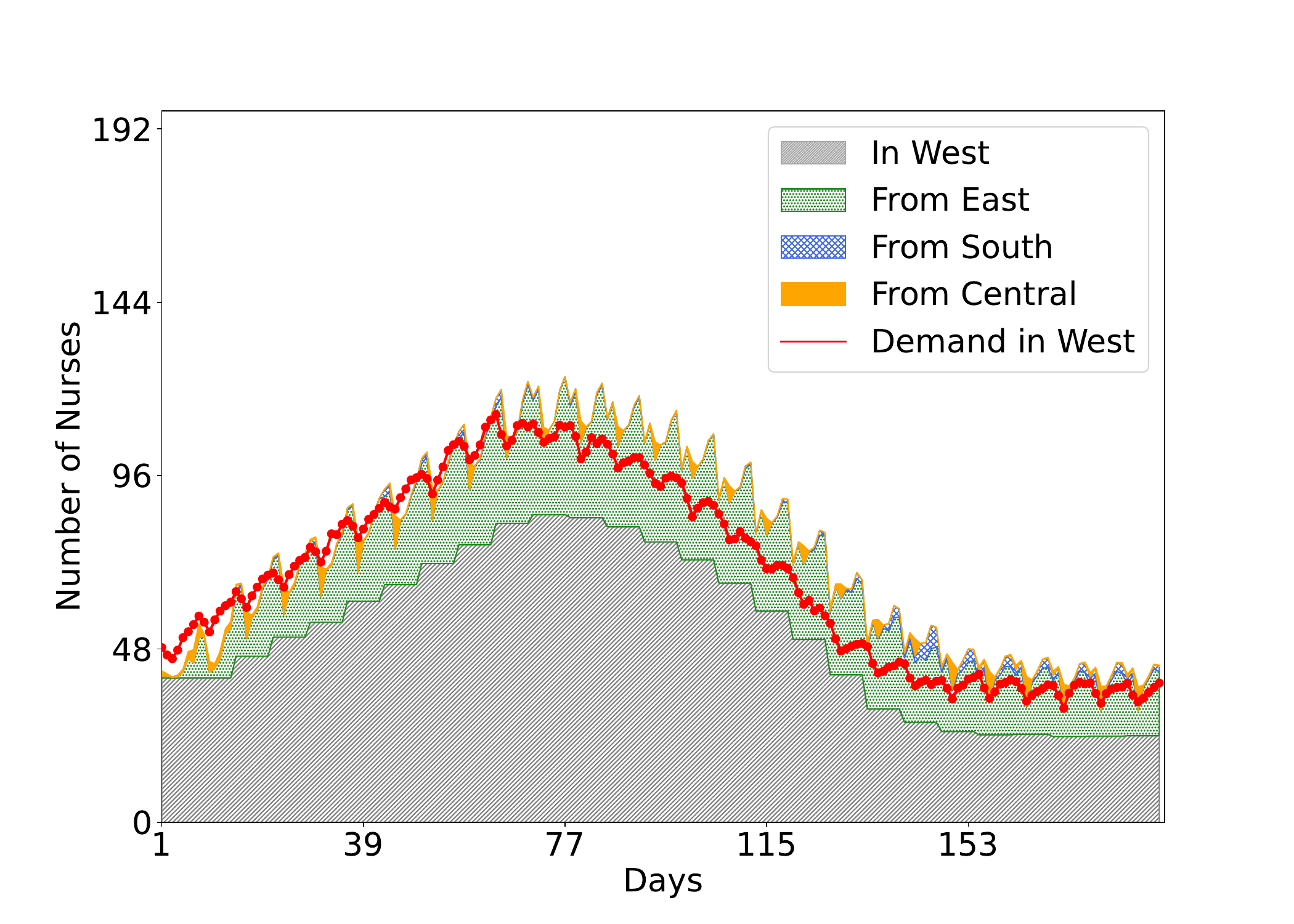}
         \caption{SAA}
         \label{arnett_shortage_staffing_saa_sec_2_min_sec_2} 
     \end{subfigure}
    \quad 
     \begin{subfigure}[b]{0.35\textwidth}
         \centering
         \includegraphics[width=\textwidth]{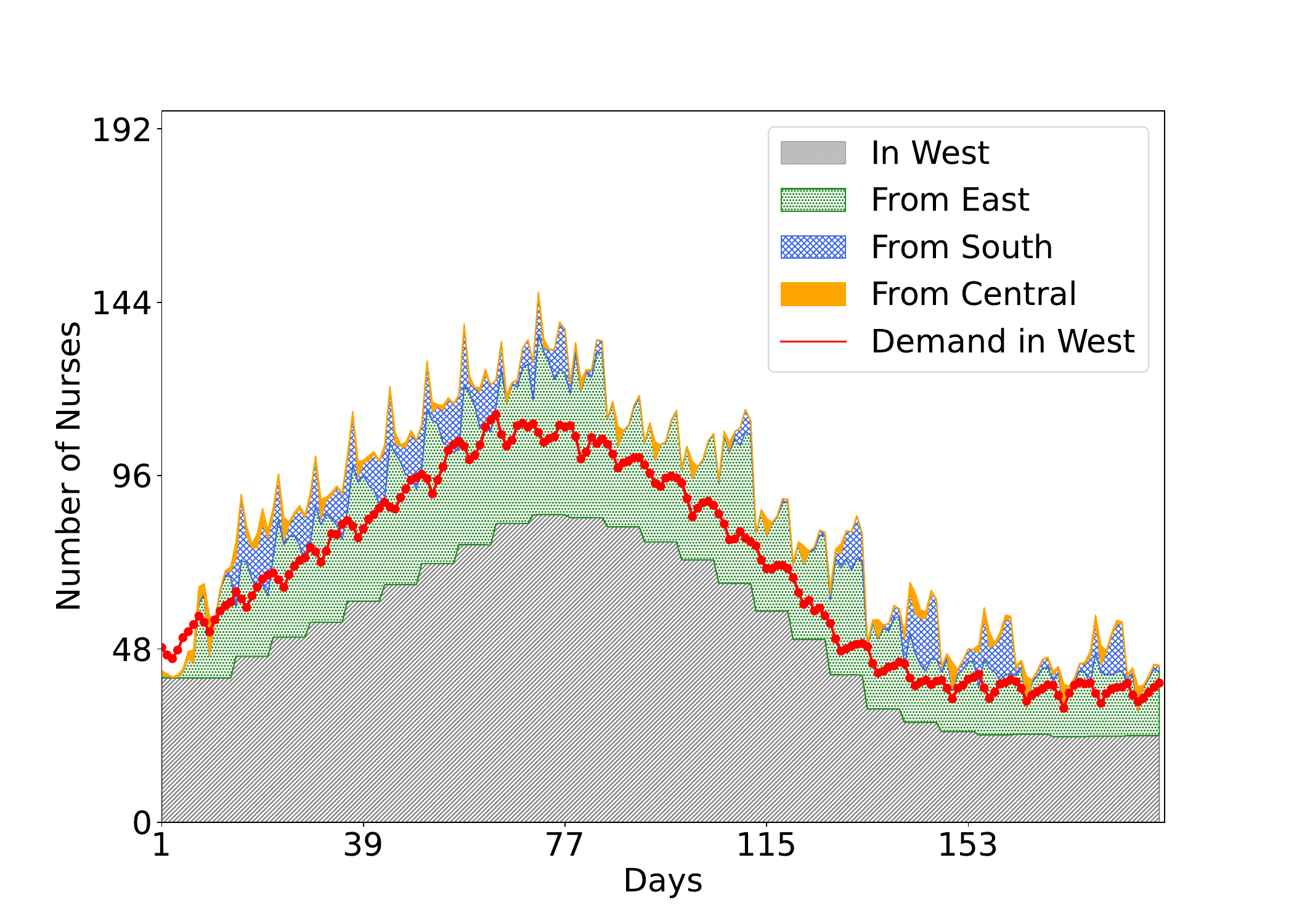}
         \caption{SRO} 
         \label{arnett_shortage_staffing_dro_sec_2_min_sec_2} 
     \end{subfigure} 
        \caption{Daily planned nurse transfers to West Hospital for SAA and SRO}  
        % (\color{blue} { x-axis: ``Days'' (from day 1);  y-axis:  ``Number of nurses'';  label: ``Nurse demand at Methodist'', ``Nurses left with home location at Methodist'', ``Nurse transfers from Arnett'', ``Nurse transfers from Ball'', and ``Nurse transfers from Bloomington''; For the labels, use a shaded rectangle for the label; remove the line for nurse transfer and only use the shaded area; for nurse demand, may use red line directly if we a square figure. } ) 
        \label{Planned_decision_DRO_SAA_arnett} 
\end{figure}

\begin{figure}
     \centering
     \begin{subfigure}[b]{0.35\textwidth}
         \centering
         \includegraphics[width=\textwidth]{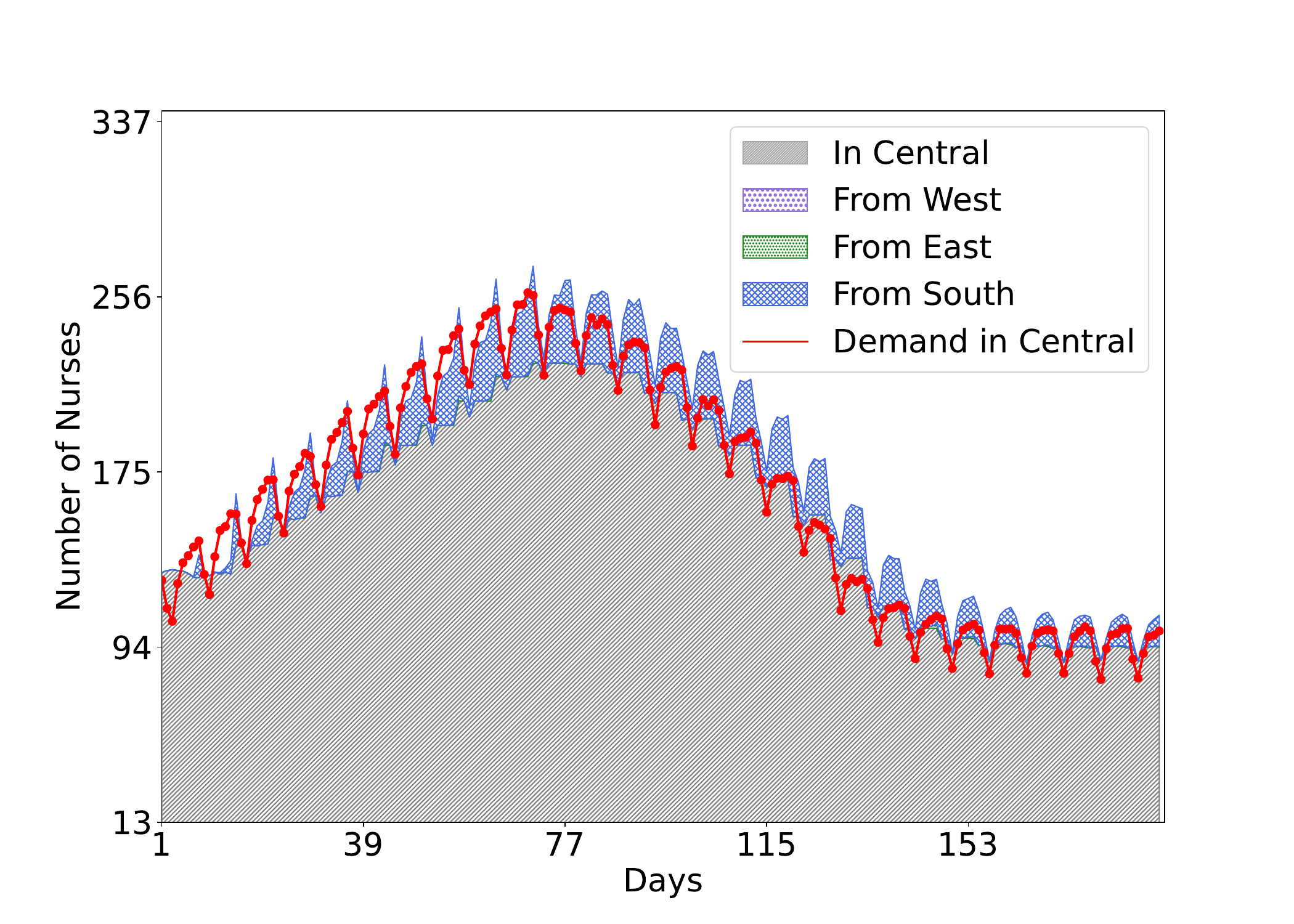}
         \caption{SAA}
         \label{methodist_shortage_staffing_saa_sec_2_min_sec_2} 
     \end{subfigure}
    \quad 
     \begin{subfigure}[b]{0.35\textwidth} 
         \centering
         \includegraphics[width=\textwidth]{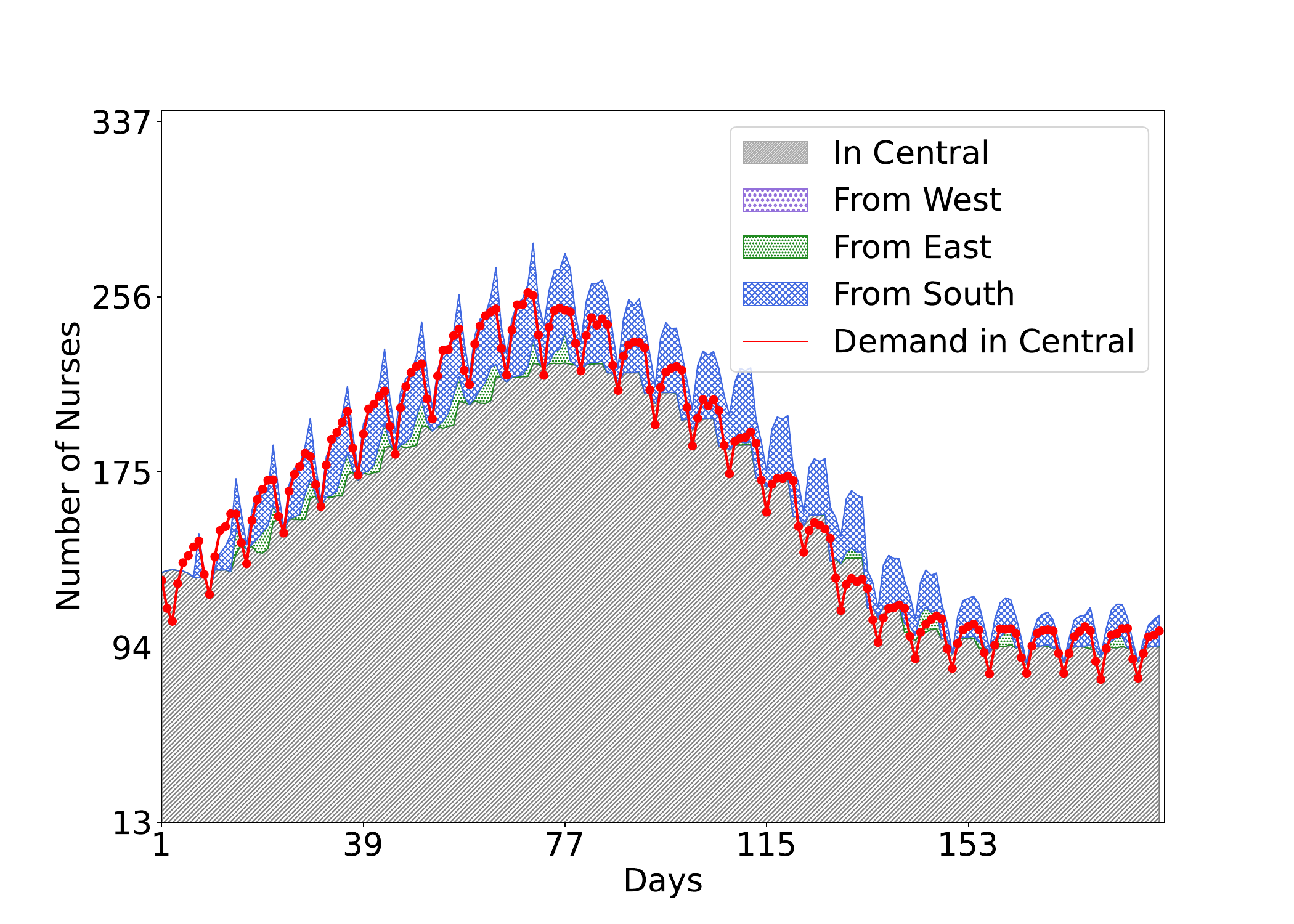} 
         \caption{SRO} 
         \label{methodist_shortage_staffing_dro_sec_2_min_sec_2} 
     \end{subfigure} 
        \caption{Daily planned nurse transfers to Central Hospital for SAA and SRO}   
        %  (\color{blue} { x-axis: ``Days'' (from day 1);  y-axis:  ``Number of nurses'';  label: ``Nurse demand at Methodist'', ``Nurses left with home location at Methodist'', ``Nurse transfers from Arnett'', ``Nurse transfers from Ball'', and ``Nurse transfers from Bloomington''; For the labels, use a shaded rectangle for the label; remove the line for nurse transfer and only use the shaded area; for nurse demand, may use red line directly if we a square figure. } ) 
        \label{Planned_decision_DRO_SAA_methodist}  
\end{figure}

To understand the reasons behind these observations, we examine planned nurse transfers from Central, South, and East Hospitals to West Hospital in Figure~\ref{Planned_decision_DRO_SAA_arnett}, and from West, East, and South to Central Hospital in Figure~\ref{Planned_decision_DRO_SAA_methodist}. We find that, prior to Week 10 (Day 70), SAA consistently under-deploys nurse transfers to high-demand hospitals; whereas SRO tends to allocate more nurses upfront, which reduces the reliance on costly emergency transfers. After Week 10 (Day 71), the planned transfer volumes under SAA become sufficient and align better with demand curve.
%\textcolor{red}{(Wei: To align Figures 11–12 with the weekly timeline, we will clarify the corresponding date ranges in the next revision.)}

\begin{figure}
     \centering
     \begin{subfigure}[b]{0.23\textwidth}
         \centering
         \includegraphics[width=\textwidth]{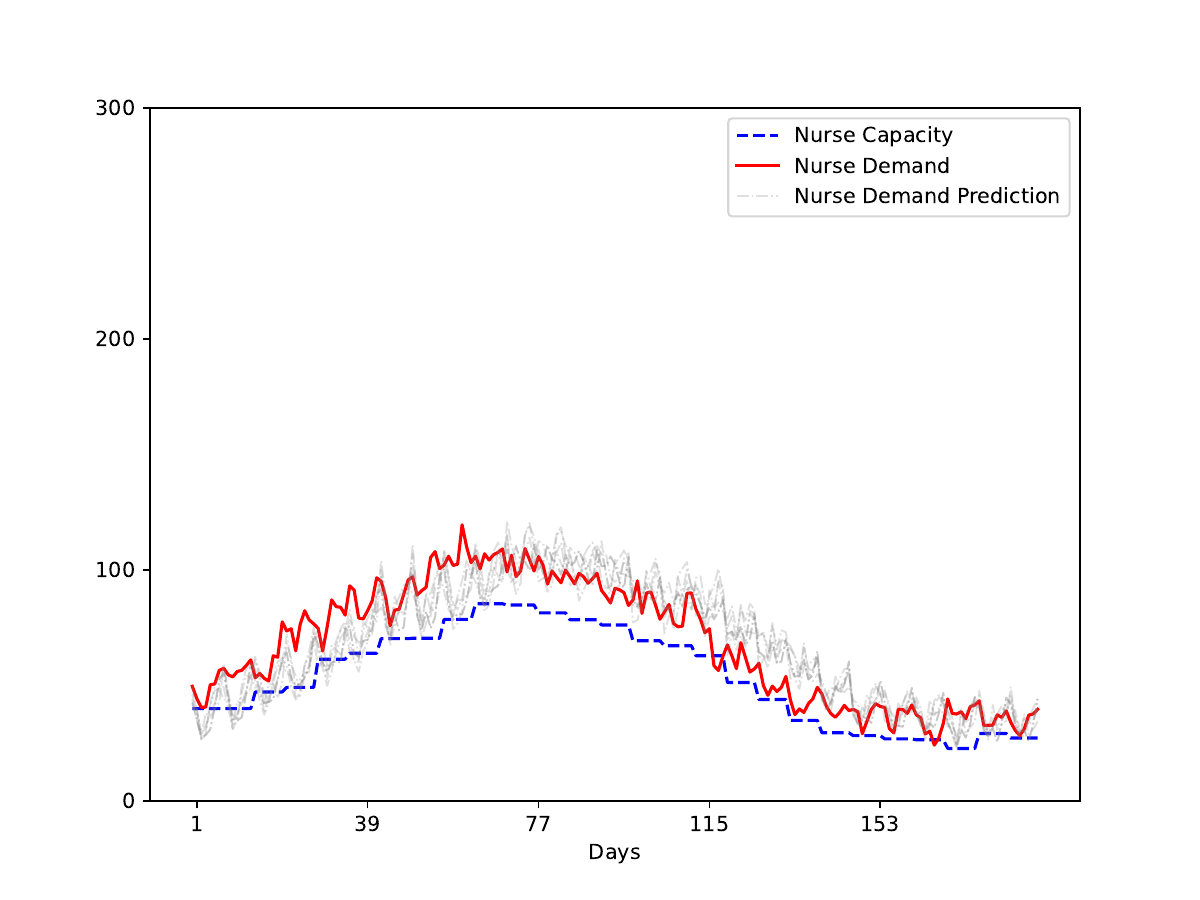}
         \caption{West Hospital} 
         \label{Arnett_demand_vs_capacity_training} 
     \end{subfigure}
     \hfill 
     \begin{subfigure}[b]{0.23\textwidth}
         \centering
         \includegraphics[width=\textwidth]{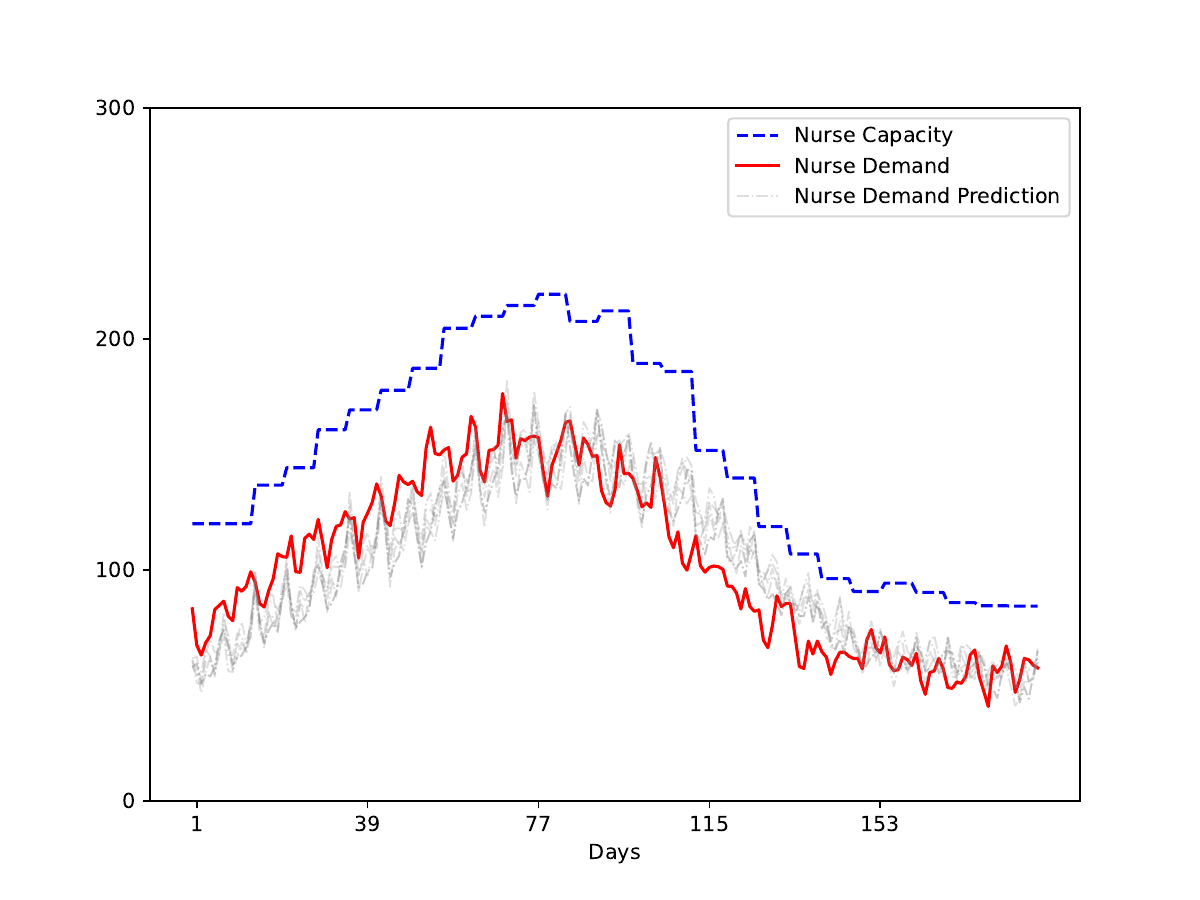}   
         \caption{East Hospital} 
         \label{Ball_Mem_demand_vs_capacity_training}     
     \end{subfigure} 
          \begin{subfigure}[b]{0.23\textwidth}
         \centering
         \includegraphics[width=\textwidth]{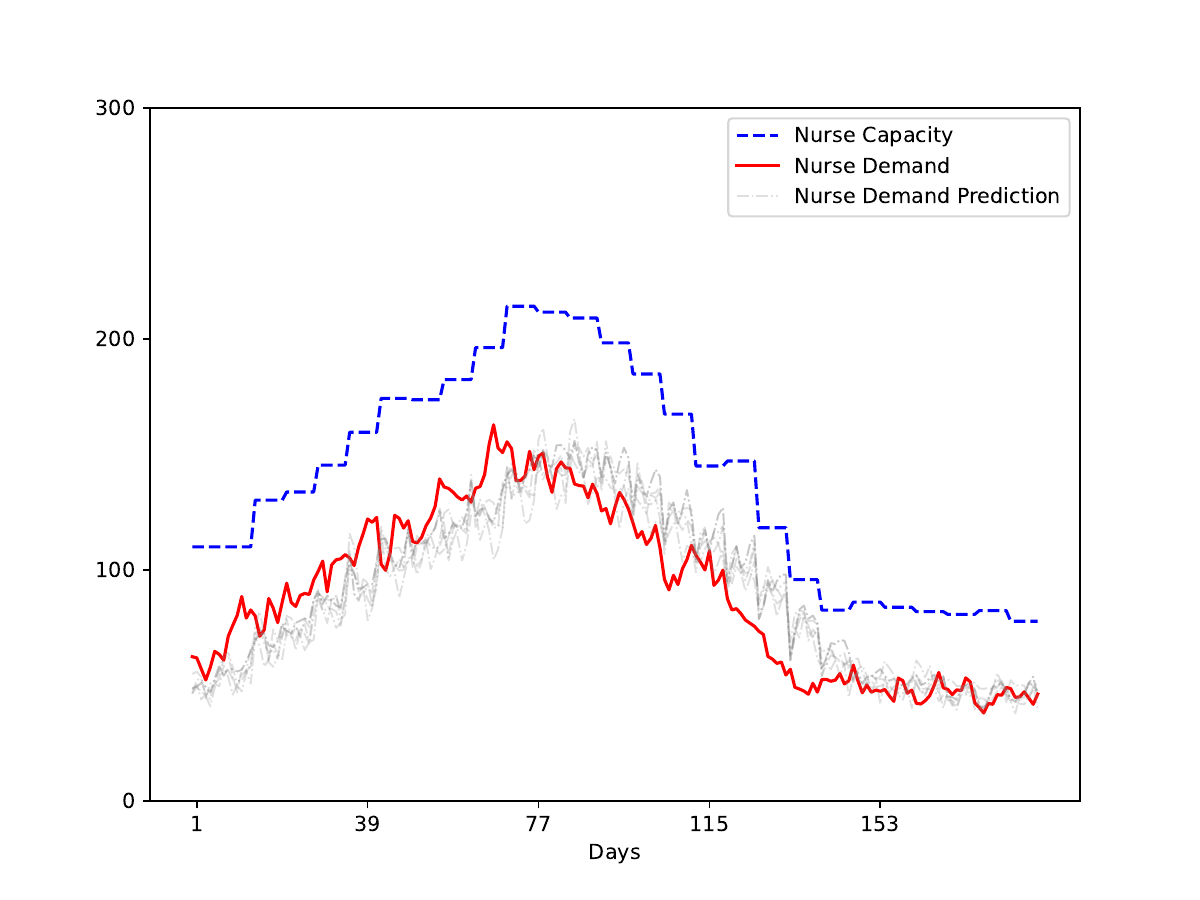}
         \caption{South Hospital} 
         \label{Bloomington_demand_vs_capacity_training} 
     \end{subfigure}
     \hfill 
     \begin{subfigure}[b]{0.23\textwidth}
         \centering
         \includegraphics[width=\textwidth]{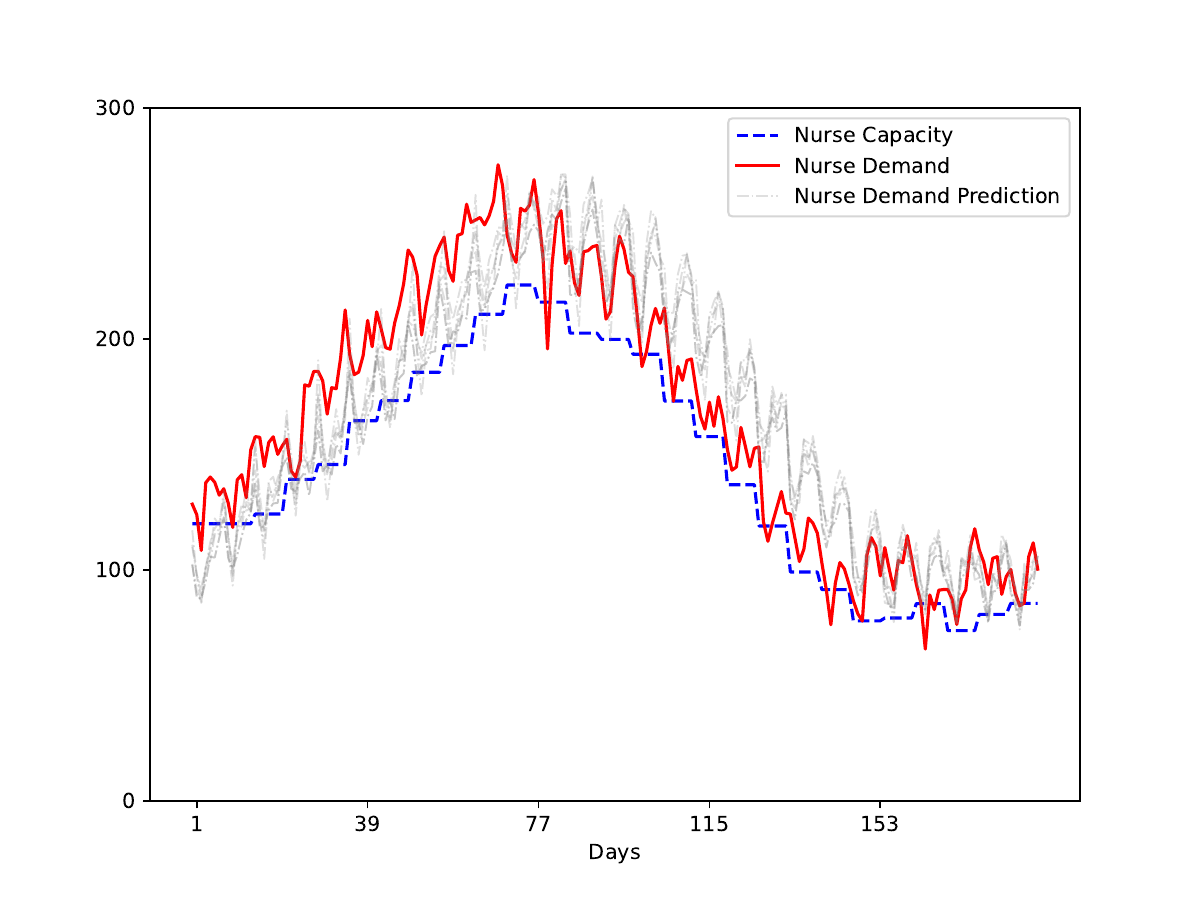}   
         \caption{Central Hospital} 
         \label{Methodist_demand_vs_capacity_training}     
     \end{subfigure} 
         \caption{Daily nurse capacity and demand along with the demand prediction (use the data from the past three weeks) at each hospital for one sample path} 
        \label{nurse_demand_capacity_prediction_over_time} 
\end{figure}

To further explore this phenomenon, Figure~\ref{nurse_demand_capacity_prediction_over_time} plots the actual demand, predicted demand (based on rolling three-week historical data), and nurse capacity for a representative sample path. During the increasing demand phase, the rolling forecast lags behind actual trends, leading to systematic underestimation. Since SAA optimizes over the sample average from these forecasts, it under-allocates staff during this phase. In contrast, SRO hedges against this bias by accounting for demand uncertainty in its framework, which leads to better (upfront) deployment and reduced emergency transfers. In the stable or decreasing demand phase, the forecasts tend to align with or overestimate the actual demand (except in outlier weeks such as Weeks 18–19). Consequently, SAA performs well in these periods, and SRO naturally reduces to SAA by selecting a zero robust parameter. In weeks where forecast errors persist, SRO maintains a performance edge. Additional details on forecast accuracy and demand generation are provided in Appendix~\ref{sec:casestudy_true_transition_prob}.

\textbf{Sensitivity Analysis.}
To assess the robustness of our findings, we conduct two sensitivity analyses—one using an alternative demand forecasts and the other using a demand pattern with higher peaks. First, we evaluate a setting where demand predictions are based on six weeks of historical data (instead of three in the baseline). This longer look-back window introduces greater time-lag, which exacerbates the underestimation during the increasing demand phase and overestimation during the decline phase. As a result, this  amplifies the performance gap between SRO and SAA, for the same reason we explained above. That is, SRO hedges against the forecast bias and proactively transfers more nurses to high-demand hospitals upfront, thereby reducing costly emergency transfers. During stable or decreasing demand phases, both methods behave similarly, and SRO effectively becomes SAA by selecting a robust parameter of zero; see Appendix~\ref{sec:Impact_time_lagging}. 
% \red{Move (original) figure 14 (figure~\ref{nurse_demand_capacity_prediction_over_time_6week}) to appendix; I currently commented it out. and move current figure 14 to appendix. I assume we have an appendix showing results under the 6-week demand forecast setting? refer to that appendix here, so we don't need to show all details here fore sensitivity analysis.}
Second, we explore an alternative demand scenario by increasing the peak of nurse demand without extending the time to reach it—simulating a sharper disease spread. We find that SRO shows even greater improvements over SAA during the increasing demand phase. Specifically, SRO allocates more nurses from surplus regions (East and South) to hospitals with shortages (West and Central) early in the planning stage. SAA, in contrast, fails to respond quickly to the sharp rise, resulting in higher shortage costs. See details in Appendix~\ref{sec:Impact_higher_demand}. 
% \red{Move Figures~\ref{higher_peak_methodist_shortage_total_cost_sec_3_min_sec_2} and \ref{higher_peak_methodist_shortage_total_cost_percentage_sec_3_min_sec_2} to appendix, and put a reference to appendix section here.}

% \blue{  (Do we still include the comparison on work life balance metrics) } 
%  We further conduct a comparative analysis between SRO and SAA by examining the average number of deployed nurse transfer volumes and the average actual transferred miles throughout the week across all locations. As presented in Table \ref{tab:Average_metrics}, there is minimal difference between SAA and SRO regarding these two metrics. A weekly comparison of these metrics also reveals negligible distinctions. This similarity arises because the deployed nurse transfers show almost no variation under both SAA and SRO. We see that compared to SAA, SRO can help reduce costs in the increasing demand pattern, almost without influencing the deployed nurse transfers and actual transferred miles. 

\subsection{Practical Insights}
\label{sec:insights}

Our case study reveals several key insights for designing effective nurse redeployment strategies under uncertainty. First, network design emerges as one of the most influential factors: implementing a fully connected network can significantly reduce total transfer costs, volumes, and distances by allowing direct nurse deployments to hospitals facing shortage while avoiding intermediary routing. This effect is more pronounced than variations in secondment length or optimization method, when the secondment length is appropriately calibrated. However, the fully connected design is not universally better—its benefits diminish or even reverse if secondment policies are poorly aligned with network structure, leading to our next insight.

Second, appropriately calibrated secondments are critical. When secondment lengths are appropriately set (particularly longer secondments for rural-to-rural transfers), they can help further reduce total costs and improve nurse well-being by avoiding excessive daily commuting. However, secondment length must be aligned with network structure. If secondments are too short for long-distance deployments (rural-to-rural under fully connected design), high travel costs are incurred, as well as nurse fatigue. Conversely, overly long secondments can limit flexibility and responsiveness to shifting demand, which also degrade system performance. Thus, secondment policy must be carefully aligned with network design to realize the full benefits of connectivity.

Finally, the SRO approach demonstrates better adaptability to demand uncertainty, especially during periods of increasing or underestimated demand. By proactively deploy nurses upfront, SRO reduces costly emergency transfers and improves system responsiveness. In stable or overestimated demand scenarios, SRO converges to SAA, meaning that it maintains efficiency without unnecessary conservatism. Together, these findings highlight the importance of jointly optimizing network design, secondment policy, and robust planning approach to build a cost-effective healthcare workforce system.

\section{Conclusions}
\label{sec:Conclusions}

This paper addresses the challenge of balancing nurse supply and demand across a network of hospitals by developing a strategy for dynamically deploying nurses from overstaffed to understaffed facilities. 
%We consider a fully connected network design that enables nurse transfers not only between central and rural hospitals, but also directly between rural hospitals. 
Given the substantial distances between rural facilities, transferred nurses are provided with advanced notice and are assigned multi-day secondments at the destination hospital to reduce commute burden, which necessitates a new, multi-stage decision framework that accounts for the temporal correlations caused by the secondments.
Moreover, to address uncertainty in demand forecasts, we adopt an SRO approach. We overcome the associated computational challenges of the SRO approach by leveraging linear decision rules and a rolling-horizon framework. Our extensive case study highlights how network design, secondment policies, and optimization methods (SAA vs. SRO) interact to shape system performance. The analysis yields several practical insights: the fully connected network design provides significant benefits under well-calibrated secondment lengths; secondments help reduce both cost and travel burden; and the SRO method shows more value under rising or uncertain demand.

While our framework is designed for hospital-based nurse redeployment, the underlying principles can be applied more broadly to employee transfer programs in other sectors, such as consulting or emergency response. Adapting the model to new domains may require accounting for industry-specific regulations or organizational constraints, but the core methodology remains relevant. 
An important direction for future work is to integrate individualized nurse staffing into the model, accounting for personal preferences, contractual constraints, and fairness considerations. This extension would enable more tailored and implementable workforce plans, which can further bridge the gap between operational models and real-world staffing needs.

% Acknowledgments here
%\ACKNOWLEDGMENT{The authors gratefully acknowledge the existence of
%the Journal of Irreproducible Results and the support of the Society
%for the Preservation of Inane Research.}

% References here (outcomment the appropriate case)

% CASE 1: BiBTeX used to constantly update the references
%   (while the paper is being written).
\bibliographystyle{informs2014} % outcomment this and next line in Case 1
\bibliography{reference} % if more than one, comma separated

\begin{thebibliography}{47}
\providecommand{\natexlab}[1]{#1}
\providecommand{\url}[1]{\texttt{#1}}
\providecommand{\urlprefix}{URL }

\bibitem[{{American Hospital Directory}(2024)}]{AmericanHospitalDirectory}
{American Hospital Directory} (2024) American hospital directory.
  \urlprefix\url{https://www.ahd.com/}, accessed: 2024-10-02.

\bibitem[{B{\'e}langer et~al.(2019)B{\'e}langer, Ruiz, \protect\BIBand{}
  Soriano}]{belanger2019recent}
B{\'e}langer V, Ruiz A, Soriano P (2019) Recent optimization models and trends
  in location, relocation, and dispatching of emergency medical vehicles.
  \emph{European Journal of Operational Research} 272(1):1--23.

\bibitem[{Benjaafar et~al.(2022)Benjaafar, Jiang, Li, \protect\BIBand{}
  Li}]{benjaafar2022dynamic}
Benjaafar S, Jiang D, Li X, Li X (2022) Dynamic inventory repositioning in
  on-demand rental networks. \emph{Management Science} 68(11):7861--7878.

\bibitem[{Bertsimas et~al.(2023)Bertsimas, McCord, \protect\BIBand{}
  Sturt}]{bertsimas2023dynamic}
Bertsimas D, McCord C, Sturt B (2023) Dynamic optimization with side
  information. \emph{European Journal of Operational Research} 304(2):634--651.

\bibitem[{Bertsimas et~al.(2022)Bertsimas, Shtern, \protect\BIBand{}
  Sturt}]{bertsimas2022two}
Bertsimas D, Shtern S, Sturt B (2022) Two-stage sample robust optimization.
  \emph{Operations Research} 70(1):624--640.

\bibitem[{Bertsimas et~al.(2019)Bertsimas, Sim, \protect\BIBand{}
  Zhang}]{bertsimas2019}
Bertsimas D, Sim M, Zhang M (2019) Adaptive distributionally robust
  optimization. \emph{Management Science} 65(2):604--618.

\bibitem[{Brenner(2025)}]{brenner2025combating}
Brenner A (2025) Combating shortage and burnout with new float pool models.
  \url{https://www.inhouse.health/blog-post/combating-shortage-and-burnout-with-new-float-pool-models},
  blog post.

\bibitem[{Carayon \protect\BIBand{} Gurses(2018)}]{carayon2018nursing}
Carayon P, Gurses A (2018) Nursing workload and patient safety—a human
  factors engineering perspective. \emph{Patient Safety and Quality: An
  Evidence-Based Handbook for Nurses}, chapter~30 (Rockville, MD: Agency for
  Healthcare Research and Quality (US)),
  \urlprefix\url{https://www.ncbi.nlm.nih.gov/books/NBK2657/}.

\bibitem[{Chand et~al.(2002)Chand, Hsu, \protect\BIBand{}
  Sethi}]{chand2002forecast}
Chand S, Hsu VN, Sethi S (2002) Forecast, solution, and rolling horizons in
  operations management problems: A classified bibliography.
  \emph{Manufacturing \& Service Operations Management} 4(1):25--43.

\bibitem[{Cho et~al.(2019)Cho, Bretthauer, Cattani, \protect\BIBand{}
  Mills}]{cho2019behavior}
Cho DD, Bretthauer KM, Cattani KD, Mills AF (2019) Behavior aware service
  staffing. \emph{Production and Operations Management} 28(5):1285--1304.

\bibitem[{Chou et~al.(2023)Chou, Wong, Zhang, Aguiari, Im, Lam, Tse, Tang,
  \protect\BIBand{} Pau}]{chou2023taxi}
Chou KS, Wong KL, Zhang B, Aguiari D, Im SK, Lam CT, Tse R, Tang SK, Pau G
  (2023) Taxi demand and fare prediction with hybrid models: Enhancing
  efficiency and user experience in city transportation. \emph{Applied
  Sciences} 13(18):10192.

\bibitem[{Cook \protect\BIBand{} Goodwin(2008)}]{cook2008airline}
Cook GN, Goodwin J (2008) Airline networks: A comparison of hub-and-spoke and
  point-to-point systems. \emph{Journal of Aviation/Aerospace Education \&
  Research} 17(2):1.

\bibitem[{Ding et~al.(2020)Ding, Nagarajan, \protect\BIBand{}
  Zhang}]{ding2020parallel}
Ding Y, Nagarajan M, Zhang G (2020) Parallel queues with discrete-choice
  arrival pattern: Empirical evidence and asymptotic characterization.
  \emph{Available at SSRN 3584880} .

\bibitem[{Glomb et~al.(2022)Glomb, Liers, \protect\BIBand{}
  R{\"o}sel}]{glomb2022rolling}
Glomb L, Liers F, R{\"o}sel F (2022) A rolling-horizon approach for
  multi-period optimization. \emph{European Journal of Operational Research}
  300(1):189--206.

\bibitem[{Hao et~al.(2020)Hao, He, Hu, \protect\BIBand{} Jiang}]{hao2020robust}
Hao Z, He L, Hu Z, Jiang J (2020) Robust vehicle pre-allocation with uncertain
  covariates. \emph{Production and Operations Management} 29(4):955--972.

\bibitem[{He et~al.(2020)He, Hu, \protect\BIBand{} Zhang}]{he2020robust}
He L, Hu Z, Zhang M (2020) Robust repositioning for vehicle sharing.
  \emph{Manufacturing \& Service Operations Management} 22(2):241--256.

\bibitem[{Helm et~al.(2024)Helm, Shi, Drewes, \protect\BIBand{}
  Cecil}]{helm2024delta}
Helm JE, Shi P, Drewes M, Cecil J (2024) Delta coverage: The analytics journey
  to implement a novel nurse deployment program. \emph{INFORMS Journal on
  Applied Analytics} 54(5):431--454.

\bibitem[{Hoover et~al.(2024)Hoover, Lucy, \protect\BIBand{}
  Mahoney}]{Hoover2024data}
Hoover M, Lucy I, Mahoney K (2024) Data deep dive: A national nursing crisis.
  \urlprefix\url{https://www.uschamber.com/workforce/nursing-workforce-data-center-a-national-nursing-crisis},
  accessed: 2024-05-26.

\bibitem[{Hu et~al.(2025)Hu, Chan, \protect\BIBand{} Dong}]{hu2025prediction}
Hu Y, Chan CW, Dong J (2025) Prediction-driven surge planning with application
  to emergency department nurse staffing. \emph{Management Science}
  71(3):2079--2126.

\bibitem[{Huh et~al.(2013)Huh, Liu, \protect\BIBand{}
  Truong}]{huh2013multiresource}
Huh WT, Liu N, Truong VA (2013) Multiresource allocation scheduling in dynamic
  environments. \emph{Manufacturing \& Service Operations Management}
  15(2):280--291.

\bibitem[{{Indiana State Government}(2023)}]{Indiana_Government_2023}
{Indiana State Government} (2023) Indiana covid-19 home dashboard. Accessed:
  2023-04-22.

\bibitem[{{IU Health}(2023)}]{IU_Health_fiveregions2023}
{IU Health} (2023) Five service regions. Accessed: 2023-04-15.

\bibitem[{Jiang \protect\BIBand{} Guan(2018)}]{jiang2018risk}
Jiang R, Guan Y (2018) Risk-averse two-stage stochastic program with
  distributional ambiguity. \emph{Operations Research} 66(5):1390--1405.

\bibitem[{Lu et~al.(2018)Lu, Chen, \protect\BIBand{} Shen}]{lu2018optimizing}
Lu M, Chen Z, Shen S (2018) Optimizing the profitability and quality of service
  in carshare systems under demand uncertainty. \emph{Manufacturing \& Service
  Operations Management} 20(2):162--180.

\bibitem[{Lyons(2023)}]{Lyons2023nursing}
Lyons J (2023) Nursing shortage by state: A state-by-state breakdown of the
  nursing shortage.
  \urlprefix\url{https://nurse.org/education/nursing-shortage-by-state-analysis/},
  accessed: 2024-05-22.

\bibitem[{Meredith et~al.(2024)Meredith, Turner, Saville, \protect\BIBand{}
  Griffiths}]{meredith2024nurse}
Meredith P, Turner L, Saville C, Griffiths P (2024) Nurse understaffing
  associated with adverse outcomes for surgical admissions. \emph{British
  Journal of Surgery} 111(9):znae215.

\bibitem[{Murphy(2024)}]{PRSG2024navigating}
Murphy K (2024) Navigating the u.s. nurse staffing shortage.
  \urlprefix\url{https://prsglobal.com/blog/navigating-the-us-nurse-staffing-shortage},
  accessed: 2025-05-26.

\bibitem[{Nair \protect\BIBand{} Miller-Hooks(2011)}]{nair2011fleet}
Nair R, Miller-Hooks E (2011) Fleet management for vehicle sharing operations.
  \emph{Transportation Science} 45(4):524--540.

\bibitem[{Nelson(2023)}]{Nelson2023}
Nelson J (2023) What’s the difference between pcu and icu? Accessed:
  2023-04-15.

\bibitem[{Paterson et~al.(2011)Paterson, Kiesm{\"u}ller, Teunter,
  \protect\BIBand{} Glazebrook}]{paterson2011inventory}
Paterson C, Kiesm{\"u}ller G, Teunter R, Glazebrook K (2011) Inventory models
  with lateral transshipments: A review. \emph{European Journal of Operational
  Research} 210(2):125--136.

\bibitem[{Prabhu et~al.(2020)Prabhu, Taaffe, Caglayan, Isik, Song,
  \protect\BIBand{} Hand}]{prabhu2020team}
Prabhu VG, Taaffe K, Caglayan C, Isik T, Song Y, Hand W (2020) Team based, risk
  adjusted staffing during a pandemic: An agent based approach. \emph{2020
  Winter Simulation Conference (WSC)}, 747--758 (IEEE).

\bibitem[{Raghavan \protect\BIBand{} Tompson(1987)}]{raghavan1987randomized}
Raghavan P, Tompson CD (1987) Randomized rounding: a technique for provably
  good algorithms and algorithmic proofs. \emph{Combinatorica} 7(4):365--374.

\bibitem[{Rath et~al.(2023)Rath, Rajaram, Hudson, \protect\BIBand{}
  Mahajan}]{rath2023multilocation}
Rath S, Rajaram K, Hudson M, Mahajan A (2023) Multilocation, dynamic staff
  planning for a healthcare system: Methodology and application. \emph{Dynamic
  Staff Planning for a Healthcare System: Methodology and Application (August
  12, 2023)} .

\bibitem[{Razak et~al.(2020)Razak, Shin, Pogacar, Jung, Pus, Moser,
  Lapointe-Shaw, Tang, Kwan, Weinerman et~al.}]{razak2020modelling}
Razak F, Shin S, Pogacar F, Jung HY, Pus L, Moser A, Lapointe-Shaw L, Tang T,
  Kwan JL, Weinerman A, et~al. (2020) Modelling resource requirements and
  physician staffing to provide virtual urgent medical care for residents of
  long-term care homes: a cross-sectional study. \emph{Canadian Medical
  Association Open Access Journal} 8(3):E514--E521.

\bibitem[{Ryu \protect\BIBand{} Jiang(2025)}]{ryu2025nurse}
Ryu M, Jiang R (2025) Nurse staffing under absenteeism: A distributionally
  robust optimization approach. \emph{Manufacturing \& Service Operations
  Management} .

\bibitem[{Shi et~al.(2022)Shi, Helm, Chen, Lim, Parker, Tinsley,
  \protect\BIBand{} Cecil}]{shi2022operations}
Shi P, Helm JE, Chen C, Lim J, Parker RP, Tinsley T, Cecil J (2022) Operations
  (management) warp speed: Rapid deployment of hospital-focused
  predictive/prescriptive analytics for the covid-19 pandemic. \emph{Production
  and Operations Management} .

\bibitem[{Shu et~al.(2013)Shu, Chou, Liu, Teo, \protect\BIBand{}
  Wang}]{shu2013models}
Shu J, Chou MC, Liu Q, Teo CP, Wang IL (2013) Models for effective deployment
  and redistribution of bicycles within public bicycle-sharing systems.
  \emph{Operations Research} 61(6):1346--1359.

\bibitem[{Song \protect\BIBand{} Dong(2014)}]{song2014empty}
Song DP, Dong JX (2014) Empty container repositioning. \emph{Handbook of Ocean
  Container Transport Logistics: Making Global Supply Chains Effective},
  163--208 (Cham: Springer International Publishing).

\bibitem[{{\v{S}}tev{\'a}rov{\'a} \protect\BIBand{}
  Bad{\'a}nik(2018)}]{vstevarova2018performance}
{\v{S}}tev{\'a}rov{\'a} L, Bad{\'a}nik B (2018) Performance of hub and spoke
  networks of selected airlines. \emph{Transportation research procedia}
  35:240--249.

\bibitem[{{University of St. Augustine for Health
  Sciences}(2024)}]{Augustine_Health_Sciences2021}
{University of St Augustine for Health Sciences} (2024) Nursing shortage: A
  2024 data study reveals key insights.
  \urlprefix\url{https://www.usa.edu/blog/nursing-shortage/}, accessed:
  2025-04-22.

\bibitem[{Wang et~al.(2024)Wang, Nguyen, \protect\BIBand{}
  Hanasusanto}]{wang2024wasserstein}
Wang Y, Nguyen VA, Hanasusanto GA (2024) Wasserstein robust classification with
  fairness constraints. \emph{Manufacturing \& Service Operations Management}
  26(4):1567--1585.

\bibitem[{{Wolters Kluwer}(2016)}]{WoltersKluwer2016the}
{Wolters Kluwer} (2016) The importance of the optimal nurse-to-patient ratio.
  Accessed: 2023-04-16.

\bibitem[{Xie et~al.(2021)Xie, Zhuang, Ang, Chou, Luo, \protect\BIBand{}
  Yao}]{xie2021analytics}
Xie J, Zhuang W, Ang M, Chou MC, Luo L, Yao DD (2021) Analytics for hospital
  resource planning—two case studies. \emph{Production and Operations
  Management} 30(6):1863--1885.

\bibitem[{Yankovic \protect\BIBand{} Green(2011)}]{yankovic2011identifying}
Yankovic N, Green LV (2011) Identifying good nursing levels: A queuing
  approach. \emph{Operations research} 59(4):942--955.

\bibitem[{Yao et~al.(2024)Yao, Shehadeh, \protect\BIBand{}
  Padman}]{yao2024multi}
Yao X, Shehadeh KS, Padman R (2024) Multi-resource allocation and care sequence
  assignment in patient management: a stochastic programming approach.
  \emph{Health Care Management Science} 27(3):352--369.

\bibitem[{Yom-Tov \protect\BIBand{} Mandelbaum(2014)}]{yom2014erlang}
Yom-Tov GB, Mandelbaum A (2014) Erlang-r: A time-varying queue with reentrant
  customers, in support of healthcare staffing. \emph{Manufacturing \& Service
  Operations Management} 16(2):283--299.

\bibitem[{Yuan(2025)}]{yuan2025managing}
Yuan Y (2025) Managing flexible capacity in service systems with worker
  shortages. \emph{Manufacturing \& Service Operations Management}
  27(3):808--824.

\end{thebibliography}

\newpage
\setcounter{page}{1}

\begin{APPENDICES} 
\section{Description of Model Notations} \label{app:table}   
% Table generated by Excel2LaTeX from sheet 'Sheet7'
\begin{table}[htbp]
  \centering
    \caption{Description of model notations}
    \begin{tabular}{lp{31.335em}}
    \hline
    Notation & \multicolumn{1}{l}{Description} \bigstrut\\
    \hline
    $L$   & Number of hospitals \bigstrut[t]\\
    $T$   & Planning horizon length \\
    $\xi_{t}^i$ & Nurse demand on day $t$ at location $i$ \\
        ${\boldsymbol \xi}_t= [\xi_t^{ij}, 1\le i,j \le L]$ & Demand vector at all locations on day $t$\\ 
    ${\boldsymbol \xi}_{[t]} = ({\boldsymbol \xi}_{1} ,\cdots,{\boldsymbol \xi}_{t})$& Historical nurse demands up to day $t$ \\     
    $a_t^{ij}$ & Number of nurses tentatively planned to transfer from location $i$ to location $j$ on day $t$, where   $ i,j \in [L]$  and $t\in [T]$ \\
    ${\boldsymbol a}_t= [a_t^{ij}, 1\le i,j \le L]$ & Vector of on-call decisions for day $t$\\ 
    ${\boldsymbol a}_{[t]} = ({\boldsymbol a}_{1} ,\cdots,{\boldsymbol a}_{t})$&  Sequence of on-call plans from day 1 to day $t$ \\
    $b_t^{ij}$ & Number of nurses deployed from location $i$ to location $j$ on day $t$ \\
    ${\boldsymbol b}_t= [b_t^{ij}, 1\le i,j \le L]$ & Vector of deployment decisions for day $t$\\ 
    ${\boldsymbol b}_{[t]} = ({\boldsymbol b}_{1} ,\cdots,{\boldsymbol b}_{t})$&  Sequence of deployment decisions up to day $t$ \\    
    $\omega^{ij}$ & Minimum required stay from location $i$ to location $j$ when sufficient time remains in the horizon \\
    % $\omega= \max_{i,j}  \omega^{ij} $ & \multicolumn{1}{l}{Maximum secondment} \\
    $\mu^{ij}(t)= \omega^{ij} \wedge ( T- t+1 )$ & Secondment length for a nurse transferred from location $i$ to location $j$ on day $t$ \\
    $p$   & Daily payment premium for working at a non-home location\\
    $\tau^{ij}$ & Additional compensation for traveling  from location $i$ to location $j$ \\ 
    $\eta$ & Percentage fee for canceling a planned transfer \\    
    $\delta_{t}^i$ & Imbalance between demand and available nurses at location $i$ on day $t$ \\
    % $\boldsymbol\delta = [\delta_{t}^i, 1\le i \le  L, 1\le t\le  T] $ &Nurse shortages on day $t$ \\
    $\theta_{t}$ & Premium multiplier applied to the daily wage $p$ for emergency transfers on day t\\
    $K^i$ & Initial nurse capacity at location $i$ \\
    ${ z}_{t}^{ij}  (k) $ & Number of nurses transferred from location $i$ to location $j$, having the number of remaining secondment days as $k$ at the beginning of day $t$ (before the actual transfer), $ 1\le k\le \omega -1  $ \bigstrut[b]\\
    \hline
    \end{tabular}%
  \label{tab:notation}%
\end{table}%

\section{Proof}\label{sec:proof}
\subsection{Proof of Proposition 1}\label{proof:prop1}
%\proof{Proof} 
% From the uncertainty set in Equation \eqref{uncertaintyset}, we have 
%  $$\max_{1\le i\le L,1\le t\le T}  |  \zeta_{t}^i   - \xi_{t}^{i,n}   | \le \epsilon_N  .$$ 
% That is, 
%  $$   \underline{\zeta}_{t}^{i,n}   = \xi_{t}^{i,n} - \epsilon_N  \le  \zeta_{t}^i   \le \xi_{t}^{i,n}  + \epsilon_N =  \bar \zeta_{t}^{i,n} $$ 
% for $1 \le t \le T$, $1\le i\le L$ and $1\le n\le N$.   
% For ease of notation, let $  \underline{\boldsymbol \zeta}_t^n =[\underline{ \zeta}_t^{i,n},  1\le i \le L ]  $, and $  \bar{\boldsymbol \zeta}_t^n =[\bar{ \zeta}_t^{i,n},  1\le i \le L ]  $. 
% Furthermore, we let $   \underline{\boldsymbol \zeta}^n = [  \underline{\boldsymbol \zeta}_1^n,\cdots, \underline{\boldsymbol \zeta}_T^n ] \in \mathbb{R}^{L\times T } $ and $   \bar{\boldsymbol \zeta}^n = [ \bar{\boldsymbol \zeta}_1^n,\cdots, \bar{\boldsymbol \zeta}_T^n ] \in \mathbb{R}^{L\times T } $. 
% These matrices represent the lower and upper bounds of the uncertainty set $ \mathcal{U}_N^{n} $. 

For any sample path ${\boldsymbol \xi}_{[T]}^n$, $1\le n\le N$, the perturbed nurse demand satisfy the following 
 $$  \bar \zeta_{t}^{i,n}  = \xi_{t}^{i,n} - \epsilon_N  \le  \zeta_{t}^i   \le \xi_{t}^{i,n}  + \epsilon_N =  \underline{\zeta}_{t}^{i,n} $$ 
for $1 \le t \le T$ and $1\le i\le L$.   
%
%
%$$ ||  {\boldsymbol \zeta}_{[T]}   - {\boldsymbol \xi}_{[T]}^n ||_\infty \le \epsilon_N  ,$$
%where $n\in [N] $. 
%Then we have 
% $$\max_{1\le i\le L,1\le t\le T}  |  \zeta_{t}^i   - \xi_{t}^{i,n}   | \le \epsilon_N  .$$  
% Then we apply Theorem 12.2.4 in \cite{ben2009robust} for proof. 
% Firstly, it is easy to see that our objective function satisfies Assumption $\mathcal{L}$ in Theorem 12.2.4. 
% 
%
We see that the demand perturbations at each location on each day are independent. 

Let $x_{t}^{ij}   (\boldsymbol \zeta_{[t]})  = (b_{t}^{ij} (\boldsymbol \zeta_{[t]})  - a_{t}^{ij} (\boldsymbol \zeta_{[t]})  )^+ $ for $1 \le t \le T$ and $1\le i\le L$. Then we have 
\begin{equation}
\begin{aligned}
(a_{t}^{ij} - b_{t}^{ij} (\boldsymbol \zeta_{[t]})   )^+   &=x_{t}^{ij} (\boldsymbol \zeta_{[t]})      - b_{t}^{ij}  (\boldsymbol \zeta_{[t]})    +  a_{t}^{ij}  . 
 \end{aligned}
\end{equation}  
Let $y_{t}^{i} (\boldsymbol \zeta_{[t]})  = (\delta_t^i  (\boldsymbol \zeta_{[t]})   )^+ $   for $1 \le t \le T$ and $1\le i\le L$. 
Then we reformulate Problem \eqref{eq:DRO} as Problem \eqref{eq:reformulation}. 

The two formulations are equivalent since all $y_{t}^{i}  (\boldsymbol \zeta_{[t]})  $ can reach their worst case values simultaneously when  $  \bar \zeta_{t}^{i,n}  \le   \zeta_{t}^i   \le  \underline{\zeta}_{t}^{i,n} $ for $1 \le t \le T$ and $1\le i\le L$. 
Consequently, the summation  $ \sum_{t=1}^{T}\sum_{i=1}^{L}   {  s_{t}^i y_t^i (\boldsymbol \zeta_{[t]})   }    $  achieves its worst-case scenario for the same realization of   ${\boldsymbol \zeta}_{[T]}$. 
%\begin{equation}
%\begin{aligned} 
%&    \zeta_{t}^i - \bigg( K^i - \sum_{j=1}^L\sum_{k=({t}-\omega^{ij}+1) \vee  1 }^{t } b_{k}^{ij}   +  \sum_{j=1}^L\sum_{k=({t}-\omega^{ji} +1) \vee 1  }^{t} b_{{k}}^{ji}        \bigg)  \le y_{t}^{i}  {,\quad} 
%  { 1 \le t \le T, 1\le i \le L,}   \\  
% & b_{t}^{ij} - a_{t}^{ij} \le x_{t}^{ij}   {,\quad}{ 1 \le t \le T, 1\le i ,j\le L, }  \\
%& x_{t}^{ij}   \ge 0 {,\quad}{ 1 \le t \le T, 1\le i ,j\le L,} \\
%& y_{t}^{i}   \ge 0 {,\quad}{ 1 \le t \le T, 1\le i \le L.} 
%\end{aligned}
%\end{equation} 
 The result follows. 
\Halmos
%\endproof 

% To conserve space, we present the full proof of Proposition 2 in the online companion.
% (\url{???}). 
%\textcolor{red}{CITE} 

%************************************
\subsection{Proof of Proposition 2} \label{proof:prop2}
%Firstly, we introduce the epigraph variables $v_1$, $v_2$, $\cdots$, $v_N\in \mathbb{R} $. 
%Then we need to consider the following constraints in Problem \eqref{eq:reformulation}: 
%\newpage 
We transform the worst-case objective of Problem \eqref{eq:reformulation} into the following constraints via introducing epigraph variables $v_1$, $v_2$, $\cdots$, $v_N\in \mathbb{R}$, following~\cite{bertsimas2023dynamic}.
{\small
\begin{equation}
\begin{aligned}
& \sum_{t=1}^{T}\sum_{i=1}^{L} \sum_{j=1}^L  \bigg({  \big(\theta_{t} p \mu^{ij}(t) + \tau^{ij} \big) \bigg({ x}_{t}^{0,ij} +  \sum_{m=1}^{t}  \sum_{l=1}^L  {x_{t,ml}^{1,ij} } \zeta_{m}^{l}   }\bigg) +  (\eta -1)  (p \mu^{ij}(t) +\tau^{ij} )  \bigg( \bigg( { x}_{t}^{0,ij} +  \sum_{m=1}^{t}  \sum_{l=1}^L  {x_{t,ml}^{1,ij} } \zeta_{m}^{l}  \bigg)   - \\  
 &    \bigg({ b}_{t}^{0,ij} +  \sum_{m=1}^{t}  \sum_{l=1}^L  {b_{t,ml}^{1,ij} } \zeta_{m}^{l} \bigg)  +  a_{t}^{ij}  \bigg)      \bigg ) +  \sum_{t=1}^{T}\sum_{i=1}^{L}   {  s_{t}^i \bigg( { y}_{t}^{0,i} +  \sum_{m=1}^{t}  \sum_{l=1}^L  {y_{t,ml}^{1,i} } \zeta_{m}^{l} \bigg)    }   \le v_n , 
 \quad \forall {\boldsymbol \zeta}_{[T]}  \in \mathcal{U}_N^{n}, 1 \le n \le N. \\  
 \end{aligned} 
\end{equation}  } 
The constraint can be reformulated as 
{\small
\begin{equation} \label{eq:ldr_singlepolicy_epigraph}  
\begin{aligned} 
& \sum_{t=1}^{T}   \sum_{k=t}^{T}  \sum_{i=1}^{L} \bigg(  \sum_{j=1}^L  \bigg({  \big(\theta_{k} p \mu^{ij}(k) + \tau^{ij} \big)  \sum_{l=1}^L  {x_{k,tl}^{1,ij} } \zeta_{t}^{l}   }+  (\eta -1)  (p \mu^{ij}(k) +\tau^{ij} )   \sum_{l=1}^L  ({x_{k,tl}^{1,ij} }    -   {b_{k,tl}^{1,ij} } ) \zeta_{t}^{l}     \bigg )+   {  s_{k}^i  \sum_{l=1}^L  {y_{k,tl}^{1,i} } \zeta_{t}^{l}     }  \bigg)  \\ 
 &      \le v_n  - \sum_{t=1}^{T}\sum_{i=1}^{L} \bigg( \sum_{j=1}^L  \bigg({  \big(\theta_{t} p \mu^{ij}(t) + \tau^{ij} \big){ x}_{t}^{0,ij}  } +  (\eta -1)  (p \mu^{ij}(t) +\tau^{ij} )  ({ x}_{t}^{0,ij}    - { b}_{t}^{0,ij}  +  a_{t}^{ij}  )      \bigg )  -  s_{t}^i { y}_{t}^{0,i}\bigg) ,  \\ 
 &     \quad\quad  \forall {\boldsymbol \zeta}_{[T]}  \in \mathcal{U}_N^{n}, 1 \le n \le N.  
 \end{aligned} 
\end{equation}  }
Similarly, constraint \eqref{eq:ldr_singlepolicy3} can be reformulated as 
{\small
\begin{equation} \label{eq:ldr_singlepolicy_cap_b}
\begin{aligned} 
& \sum_{k=1}^{T}   \sum_{m=k }^{T}   \mathbb{1}_{\{ ({t}-\omega^{ij}+1)\vee 1) \le m \le  t \}}   \sum_{j=1}^L \sum_{l=1}^L  {b_{m,kl}^{1,ij} } \zeta_{k}^{l}       \le K^i - \sum_{j=1}^L\sum_{k=({t}-\omega^{ij}+1)\vee 1 }^{t} { b}_{k}^{0,ij}, \\ 
&  {\quad \quad \quad\quad\quad\quad\quad\quad\quad\quad\quad\quad\quad\quad\quad\quad\quad\quad\quad}{1 \le i \le L,  1 \le t \le T,  \forall {\boldsymbol \zeta}_{[T]}  \in \mathcal{U}_N^{n}, 1 \le n \le N, } 
 \end{aligned} 
\end{equation}  } 
and constraint \eqref{eq:ldr_singlepolicy4} can be reformulated as 
{\small
\begin{equation} \label{eq:ldr_singlepolicy_demand}
\begin{aligned} 
 &  \sum_{k=1}^{T}   \sum_{m=k }^{T}  \bigg(  \mathbb{1}_{\{ ({t}-\omega^{ij}+1)\vee 1) \le m \le  t \}}   \sum_{j=1}^L \sum_{l=1}^L  {b_{m,kl}^{1,ij} } \zeta_{k}^{l}    -   \mathbb{1}_{\{ ({t}-\omega^{ji}+1)\vee 1) \le m \le  t \}}   \sum_{j=1}^L \sum_{l=1}^L  {b_{m,kl}^{1,ji} } \zeta_{k}^{l}      -    \mathbb{1}_{\{ {m} =t  \}} \sum_{l=1}^L    {y_{m,kl}^{1,i} } \zeta_{k}^{l}  \bigg)   \\  
&  \le { y}_{t}^{0,i} +  K^i - \sum_{j=1}^L\sum_{k=({t}-\omega^{ij}+1) \vee  1 }^{t } { b}_{k}^{0,ij}    +  \sum_{j=1}^L\sum_{k=({t}-\omega^{ji} +1) \vee 1  }^{t}  { b}_{k}^{0,ji}  - \zeta_{t}^i    {,\quad} 
  { 1\le i \le L,1 \le t \le T,  \forall {\boldsymbol \zeta}_{[T]}  \in \mathcal{U}_N^{n}, 1 \le n \le N,}   \\  
 \end{aligned} 
\end{equation}  }

constraint \eqref{eq:ldr_singlepolicy5} can be reformulated as 
\begin{equation}  \label{eq:ldr_singlepolicy_emergency}
\begin{aligned}  
   & \sum_{k=1}^{T}   \sum_{m=k }^{T}   \mathbb{1}_{\{ {m} =t  \}}   \sum_{l=1}^L  ({b_{m,kl}^{1,ij} } \zeta_{m}^{l}   -   {x_{m,kl}^{1,ij} } \zeta_{m}^{l} )    \le a_{t}^{ij}  +  { x}_{t}^{0,ij} - { b}_{t}^{0,ij}     {,\quad}{  1\le i ,j\le L,1 \le t \le T, \forall {\boldsymbol \zeta}_{[T]}  \in \mathcal{U}_N^{n}, 1 \le n \le N,}  \\ 
 \end{aligned} 
\end{equation}  

 constraint \eqref{eq:ldr_singlepolicy6}  can be reformulated as 
\begin{equation} \label{eq:ldr_singlepolicy_nonb}
\begin{aligned} 
  &  -  \sum_{k=1}^{T}   \sum_{m=k }^{T}   \mathbb{1}_{\{ {m} =t  \}}  \sum_{l=1}^L  {b_{m,kl}^{1,ij} } \zeta_{m}^{l}   \le  { b}_{t}^{0,ij}   {,\quad}{ 1\le i ,j\le L, 1 \le t \le T, \forall {\boldsymbol \zeta}_{[T]}  \in \mathcal{U}_N^{n}, 1 \le n \le N,}  \\
  &  -  \sum_{k=1}^{T}   \sum_{m=k }^{T}   \mathbb{1}_{\{ {m} =t  \}}  \sum_{l=1}^L  {x_{m,kl}^{1,ij} } \zeta_{m}^{l}   \le  { x}_{t}^{0,ij}   {,\quad}{ 1\le i ,j\le L, 1 \le t \le T, \forall {\boldsymbol \zeta}_{[T]}  \in \mathcal{U}_N^{n}, 1 \le n \le N,}  
 \end{aligned} 
\end{equation}  
and constraint \eqref{eq:ldr_singlepolicy7}  can be reformulated as 
\begin{equation}  \label{eq:ldr_singlepolicy_nony} 
\begin{aligned} 
 &- \sum_{k=1}^{T}   \sum_{m=k }^{T}   \mathbb{1}_{\{ {m} =t  \}} \sum_{l=1}^L    {y_{m,kl}^{1,i} } \zeta_{k}^{l}     \le  { y}_{t}^{0,i}   {,\quad}{ 1\le i \le L, 1 \le t \le T,  \forall {\boldsymbol \zeta}_{[T]}  \in \mathcal{U}_N^{n}, 1 \le n \le N.}   \\ 
% & b_{t}^{ij}   \in  {\mathscr{L}_b }^{ L^2 \times T }  , x_{t}^{ij}   \in  {\mathscr{L}_x }^{ L^2 \times T }   {,\quad}{  1\le i ,j\le L, 1 \le t \le T, } \\ 
%& y_{t}^{i}    \in  {\mathscr{L}_y }^{ L \times T }    {,\quad}{  1\le i \le L, 1 \le t \le T,}   \\
%& \forall {\boldsymbol \zeta}_{[T]}  \in \mathcal{U}_N^{n}, 1 \le n \le N. 
 \end{aligned} 
\end{equation}

We see that each semi-infinite constraint can be reformulated as 
$$ \max_{{\boldsymbol \zeta}_{[T]} \in  \mathcal{U}_N^{n} }  \sum_{t=1}^T  {\boldsymbol \beta}_t   {\boldsymbol\zeta}_t  \le  \gamma  $$
for some vectors ${\boldsymbol \beta} = ( {\boldsymbol\beta}_1, \cdots,  {\boldsymbol\beta}_T) \in \mathbb{R}^{L\times T } $ and $\gamma \in \mathbb{R}$. 
It follows from strong duality for linear optimization that
\begin{equation} 
\begin{aligned}
  \max_{{\boldsymbol \zeta}_{[T]} \in \mathcal{U}_N^{n} }  \sum_{t=1}^T  {\boldsymbol \beta}_t   {\boldsymbol\zeta}_t = 
\left \{ \begin{array}{lll} 
&\min_{ {\boldsymbol\lambda}_t, \boldsymbol\alpha_t  \in \mathbb{R}_+^{L} }   \sum_{t=1}^T \big(  \bar{\boldsymbol \zeta}_t^n {\boldsymbol\lambda}_t - \underline{\boldsymbol \zeta}_t^n {\boldsymbol\alpha}_t) \\
 & \\
&\mathrm{s.t.},  \quad {\boldsymbol\lambda}_t  - {\boldsymbol\alpha}_t  =  {\boldsymbol \beta}_t   \quad \quad \forall t \in [T] .\\ 
\end{array} \right.  
\end{aligned}
\end{equation}  
The optimal solutions are given by ${\boldsymbol\lambda}_t = [{\boldsymbol \beta}_t ]_+$ and ${\boldsymbol\alpha}_t = [ - { \boldsymbol \beta}_t ]_+ $. 
Importantly, these solutions are independent of the index $n$ of the uncertainty set $\mathcal{U}_N^{n}$.  
 
Consequently, constraint \eqref{eq:ldr_singlepolicy_epigraph} is satisfied if and only if there exist $\nu^{epi} $, $\psi^{epi} \in \mathbb{R}_+^{T\times L}  $  such that  
{\small   
\begin{equation*}   
\begin{aligned}  
& \sum_{t=1}^{T}  \bigg(  \sum_{l=1}^{L}   \nu_{tl}^{epi} \bar{ \zeta}_{t}^{ln}  -  \sum_{l=1}^{L}  \psi_{tl}^{epi}  \underline{ \zeta}_{t}^{ln}     \bigg )    +  \sum_{t=1}^{T}\sum_{i=1}^{L}  \bigg( \sum_{j=1}^L  \bigg({  \big(\theta_{t} p \mu^{ij}(t) + \tau^{ij} \big){ x}_{t}^{0,ij}  } +  (\eta -1)  (p \mu^{ij}(t) +\tau^{ij} )  ({ x}_{t}^{0,ij}    - { b}_{t}^{0,ij}  +  a_{t}^{ij}  )      \bigg ) \\ 
 &\quad\quad\quad \quad\quad\quad\quad\quad\quad\quad\quad\quad \quad\quad\quad \quad\quad\quad \quad  -   {  s_{t}^i { y}_{t}^{0,i}     }    \bigg )   \le v_n      ,   \quad  1 \le n \le N, \\ 
 &   \nu_{tl}^{epi}   - \psi_{tl}^{epi}  =  \sum_{k=t}^{T}  \sum_{i=1}^{L} \bigg(  \sum_{j=1}^L  \bigg({  \big(\theta_{k} p \mu^{ij}(k) + \tau^{ij} \big)   {x_{k,tl}^{1,ij} }   }+  (\eta -1)  (p \mu^{ij}(k) +\tau^{ij} )    ({x_{k,tl}^{1,ij} }    -   {b_{k,tl}^{1,ij} } )      \bigg )  
     +   {  s_{k}^i   {y_{k,tl}^{1,i} }  }  \bigg)   ,  \\
     & \quad \quad\quad\quad\quad 1 \le l \le L, 1 \le t \le T. \\ 
 \end{aligned}  
\end{equation*}  }

Constraint \eqref{eq:ldr_singlepolicy_cap_b} is satisfied if and only if there exist $\nu^{cap} $, $\psi^{cap} \in \mathbb{R}_+^{T^2\times L^2}  $  such that    
{\small 
\begin{equation*}   
\begin{aligned}  
& \sum_{k=1}^{T}  \bigg(  \sum_{l=1}^{L}   \nu^{cap}_{klti} \bar{ \zeta}_{k}^{ln}  -  \sum_{l=1}^{L}  \psi^{cap}_{klti}  \underline{ \zeta}_{k}^{ln}     \bigg )     + \sum_{j=1}^L\sum_{k=({t}-\omega^{ij}+1)\vee 1 }^{t} { b}_{k}^{0,ij}  \le K^i  {,\quad}{1 \le i \le L,  1 \le t \le T,  1 \le n \le N, } \\ 
& \nu^{cap}_{klti}   -  \psi^{cap}_{klti}       =   \sum_{m=k }^{T}  \sum_{j=1}^L   \mathbb{1}_{\{ ({t}-\omega^{ij}+1)\vee 1) \le m \le  t \}}    {b_{m,kl}^{1,ij} } \zeta_{k}^{l}    {,\quad}{1 \le i \le L,  1 \le t \le T, 1 \le l \le L,  1 \le k \le T.  } \\  
 \end{aligned}  
\end{equation*}  } 

Constraint \eqref{eq:ldr_singlepolicy_demand} is satisfied if and only if there exist $\nu^{sho} $, $\psi^{sho} \in \mathbb{R}_+^{T^2\times L^2}  $    such that 
{\small
\begin{equation*}   
\begin{aligned}  
 &  \sum_{k=1}^{T}  \bigg(  \sum_{l=1}^{L}   \nu^{sho}_{klti} \bar{ \zeta}_{k}^{ln}  -  \sum_{l=1}^{L}  \psi^{sho}_{klti}  \underline{ \zeta}_{k}^{ln}     \bigg )  + \zeta_{t}^i  - \bigg( K^i  - \sum_{j=1}^L\sum_{k=({t}-\omega^{ij}+1) \vee  1 }^{t } { b}_{k}^{0,ij}   +   \sum_{j=1}^L\sum_{k=({t}-\omega^{ji} +1) \vee 1  }^{t}  { b}_{k}^{0,ji}  \bigg)  \le { y}_{t}^{0,i} , \\ 
 & \quad \quad\quad\quad\quad \quad \quad\quad\quad\quad\quad\quad\quad\quad\quad    { 1\le i \le L,1 \le t \le T, 1 \le n \le N,}   \\  
   &     \nu^{sho}_{klti}   -  \psi^{sho}_{klti}    =  \sum_{m=k }^{T}  \bigg( \sum_{j=1}^L  \mathbb{1}_{\{ ({t}-\omega^{ij}+1)\vee 1) \le m \le  t \}}  \sum_{l=1}^L  {b_{m,kl}^{1,ij} } \zeta_{k}^{l}    -   \sum_{j=1}^L \mathbb{1}_{\{ ({t}-\omega^{ji}+1)\vee 1) \le m \le  t \}}   \sum_{l=1}^L  {b_{m,kl}^{1,ji} } \zeta_{k}^{l}      -   \\ 
   &  \quad \quad\quad\quad\quad  \quad \quad      \mathbb{1}_{\{ {m} =t  \}} \sum_{l=1}^L    {y_{m,kl}^{1,i} } \zeta_{k}^{l}  \bigg)    {,\quad} 
  { 1\le i, l  \le L,1 \le t \le T,  1 \le k \le T. }   \\  
 \end{aligned}  
\end{equation*} }  

Constraint \eqref{eq:ldr_singlepolicy_emergency} is satisfied if and only if there exist $\nu^{eme} $, $\psi^{eme} \in \mathbb{R}_+^{T^2\times L^3}  $   such that 
{\small
\begin{equation*}   
\begin{aligned}  
   & \sum_{k=1}^{T}  \bigg(  \sum_{l=1}^{L}   \nu^{eme}_{kltij} \bar{ \zeta}_{k}^{ln}  -  \sum_{l=1}^{L}  \psi^{eme}_{kltij}  \underline{ \zeta}_{k}^{ln}     \bigg ) + { b}_{t}^{0,ij} - a_{t}^{ij}    \le   { x}_{t}^{0,ij}    {,\quad}{  1\le i ,j\le L,1 \le t \le T, 1 \le n \le N,}  \\ 
      &  \nu^{eme}_{kltij}   -  \psi^{eme}_{kltij}      =   \sum_{m=k }^{T}   \mathbb{1}_{\{ {m} =t  \}}   ({b_{m,kl}^{1,ij} } \zeta_{m}^{l}   -   {x_{m,kl}^{1,ij} } \zeta_{m}^{l} )      {,\quad}{  1\le i ,j,l\le L,1 \le t \le T,   1 \le k \le T. }  \\  
% & a_{t}^{ij}  \ge 0 {,\quad}{ 1\le i ,j\le L, 1 \le t \le T.}  
 \end{aligned}  
\end{equation*} }  
Constraints \eqref{eq:ldr_singlepolicy_nonb} and \eqref{eq:ldr_singlepolicy_nony}  are satisfied if and only if there exist $\nu^{nnb} $, $\psi^{nnb} \in \mathbb{R}_+^{T^2\times L^3}  $, $ \nu^{nny} $, $\psi^{nny} \in \mathbb{R}_+^{T^2\times L^2}  $   such that  
{\small
\begin{equation*}   
\begin{aligned}  
 &   \sum_{k=1}^{T}  \bigg(  \sum_{l=1}^{L}   \nu^{nnb}_{kltij} \bar{ \zeta}_{k}^{ln}  -  \sum_{l=1}^{L}  \psi^{nnb}_{kltij}  \underline{ \zeta}_{k}^{ln}     \bigg )   \le  { b}_{t}^{0,ij}   {,\quad}{ 1\le i ,j\le L, 1 \le t \le T, 1 \le n \le N,}  \\
    &   \nu^{nnb}_{kltij}   -  \psi^{nnb}_{kltij}      =   - \sum_{m=k }^{T}   \mathbb{1}_{\{ {m} =t  \}}   {b_{m,kl}^{1,ij} } \zeta_{m}^{l}     {,\quad}{ 1\le i ,j,l\le L, 1 \le t \le T,   1 \le k \le T,  }  \\  
     &   \sum_{k=1}^{T}  \bigg(  \sum_{l=1}^{L}   \nu^{nnx}_{kltij} \bar{ \zeta}_{k}^{ln}  -  \sum_{l=1}^{L}  \psi^{nnx}_{kltij}  \underline{ \zeta}_{k}^{ln}     \bigg )   \le  { x}_{t}^{0,ij}   {,\quad}{ 1\le i ,j\le L, 1 \le t \le T, 1 \le n \le N,}  \\
    &   \nu^{nnx}_{kltij}   -  \psi^{nnx}_{kltij}      =   - \sum_{m=k }^{T}   \mathbb{1}_{\{ {m} =t  \}}   {x_{m,kl}^{1,ij} } \zeta_{m}^{l}     {,\quad}{ 1\le i ,j,l\le L, 1 \le t \le T,   1 \le k \le T,  }  \\  
 &   \sum_{k=1}^{T}    \bigg(  \sum_{l=1}^{L}   \nu^{nny}_{klti} \bar{ \zeta}_{k}^{ln}  -  \sum_{l=1}^{L}  \psi^{nny}_{klti}  \underline{ \zeta}_{k}^{ln}     \bigg )     \le  { y}_{t}^{0,i}   {,\quad}{ 1\le i \le L, 1 \le t \le T,   1 \le n \le N, }   \\    
  &   \nu^{nny}_{klti}   -  \psi^{nny}_{klti}        =  - \sum_{m=k }^{T}   \mathbb{1}_{\{ {m} =t  \}} \sum_{l=1}^L    {y_{m,kl}^{1,i} } \zeta_{k}^{l}     {,\quad}{ 1\le i,l \le L, 1 \le t \le T,   1 \le k \le T .}   
 \end{aligned}  
\end{equation*}  }  
Thus, we can remove the epigraph decision variables and the result follows. 
\Halmos
%\endproof  

%************************************

\section{Additional Details for Simulation Platform}
\label{sec:SyntheticDataGeneration}  

In this section, we describe the simulation platform used in our case study.  
Section~\ref{sec:TestingSamplePathGeneration} details the process of calibrating the simulator with COVID-19 data and generating the testing sample paths. 
Section~\ref{sec:TrainingSamplePathGeneration}  presents the procedure for generating training sample paths using rolling demand predictions. 
Section~\ref{sec:TravelCostSetup} explains how we calibrate the transfer cost based on geographical distance between different hospitals.  Section~\ref{sec:NurseCapacityGeneration} specifies the setup of nurse capacity. 

% we first describe the process of calibrating the simulator with COVID-19 data. The calibrated simulator is used to generate the testing sample paths. We then describe the generation of the training sample paths, which uses the same simulator but with different input data. 

 \subsection{Simulation Calibration and Testing Path Generation} \label{sec:TestingSamplePathGeneration} 
 
 In this section, we describe the process of calibrating the simulator with COVID-19 data, which is then used to generate the testing sample paths. 
We first introduce the method for generating patient arrivals at each hospital in Appendix~\ref{sec:Patient_Arrivals}. 
To estimate the patient census at each hospital unit, we introduce the patient flow model within one hospital in Appendix~\ref{sec:Patient_Transitions}, which captures the transition of patients between different units in the hospital. 
Subsequently, we present the census estimation in Appendix~\ref{sec:PatientsCensusGeneration}. 
Finally, in Appendix \ref{sec:NurseDemandEstimation} we describe how to use the patient census to generate estimates for nurse demand at each hospital, which are used as the testing sample paths.

\subsubsection{Patient Arrivals to Each Hospital}  \label{sec:Patient_Arrivals}

We assume that the number of patient arrivals on day $t$ ($t \in [\mathcal{T}]$) at location $i$ ($i \in [L]$) follows a Poisson distribution with an arrival rate of $\lambda_{i,t}$ per day, where $\lambda_{i,t}$'s are temporally correlated across days and spatially correlated across locations. We show the details below; Algorithm \ref{alg:patient_arr_generation} provides the corresponding pseudocode. 

\begin{algorithm}
    \caption{Patient Arrival Generation}
    \label{alg:patient_arr_generation}
    \begin{itemize}
        \item \textbf{Input}: $p$, $Y$, $t_{\text{lag}}$, $r$, $\mathcal{T}$, $L$, $\boldsymbol{\alpha}$, $\boldsymbol{\beta}$, $\boldsymbol{\zeta}$, $\boldsymbol{\gamma}$, $\boldsymbol{\phi}$, $\boldsymbol{\theta}$, $\boldsymbol{\kappa}$, $\boldsymbol{\epsilon}$, $[\lambda_{i,t} , 1\le i \le L, -p+1 \le t \le 0]$
        \item \textbf{Main Iteration}: For each $t \in [\mathcal{T}]$ and $i \in [L]$, perform:
        \begin{enumerate}
            \item Generate the arrival rate $\lambda_{i,t}$ on day $t$ at location $i$ from
            \begin{align*}
                \lambda_{i,t} &= \lambda_{i,t}^{1} + \lambda_{i,t}^{2}, \\
                \text{where} \\
                \lambda_{i,t}^{1} &= \beta_{i,t} \sum_{l=1}^{p} \phi_{i,l} \lambda_{i,t-l}^{1} + \zeta_{i,t}\kappa_{i, (t-1)\% Y +1 } +  \epsilon_t,\; \epsilon_t \sim N(0, \gamma_{i,t}\kappa_{i, (t-1)\% Y +1 })\\
                \lambda_{i,t}^{2} &= \sum_{\substack{j=1 \\ j \neq i}}^{L} \alpha_{ji,t} \theta_{ji} \sum_{l=1+t_{\text{lag}}}^{t_{\text{lag}}+r} \lambda_{j,t-l}^{2}. 
               % \lambda_{i,t}^{3} &\sim \text{Poisson}(\zeta_{i,t} \kappa_{i, (t-1)\%Y+1}).
            \end{align*}
            \item Generate the patient arrivals $A_{i,t}$ on day $t$ at location $i$ from
            \begin{align*}
                A_{i,t} \sim \text{Poisson}(\lambda_{i,t}).
            \end{align*}
        \end{enumerate}
        \item \textbf{Output}: $A_{i,t}$ 
    \end{itemize}
\end{algorithm}

\begin{itemize}
    \item \textbf{\textit{Temporal Correlation:}} Inspired by the autoregressive model, 
    we recursively generate the first component of the arrival rate, $\lambda_{i,t}^{1}$,  on day $t \in \{1,2,\cdots,\mathcal{T} \}$ at location $i \in \{1,2,\cdots,L\}$ as  
    \begin{align*}
        \lambda_{i,t}^{1} = \beta_{i,t}\sum_{l=1}^p \phi_{i,l}\lambda^{1}_{i,t-l} + \zeta_{i,t}\kappa_{i, (t-1)\% Y +1 }+ \epsilon_t,
    \end{align*} 
    where $p$ is the number of days from the past that correlates with patient arrivals on day $t$, $ \phi_{i,l}$ is the temporal correlation factor,  $ \beta_{i,t}$ is is a scaling parameter reflecting the strength of disease spread, $\kappa_{i, (t-1)\%Y+1}$ is the day-of-the-week effect, $\zeta_{i,t} \kappa_{i, (t-1)\%Y+1} $ is the patient arrival by considering the day-of-the-week effect,   and $\epsilon_t \sim N(0, \gamma_{i,t}\kappa_{i, (t-1)\% Y +1 })$ is a white noise term. 
   % Here, $\gamma_{i,t}$  is the effect of patient arrival strength on the variance of the white noise. 
   We let     ${\boldsymbol\phi}= [\phi_{i,t}, 1\le i \le L, 1\le t\le p] $,  ${\boldsymbol\beta}= [\beta_{i,t},1\le i \le L, 1\le t\le p]$, $\boldsymbol{\kappa} = [\kappa_{i,t},1\le i \le L, 1\le t\le Y ]$ and $\boldsymbol{\zeta} = [\zeta_{i,t}, 1\le i \le L, 1\le t\le p]$, and $\boldsymbol{\gamma}= [\gamma_{i,t}, 1\le i \le L, 1\le t\le p]$.  
% {\color{blue} We estimate the temporal correlation factor ${\boldsymbol\phi} $ by fitting the real dataset from an anonymous hospital through the autoregressive model.}  
% \red{What is the ``real data'' we are using?
% and does this sentence mean we use the ARIMA model to fit the $\phi$ part?} 
% and ${\boldsymbol\beta}$, ${\boldsymbol\kappa}$,  $\boldsymbol{\zeta}$ and $\boldsymbol{\gamma}$  can be adjusted based on the spread of the pandemic.  

    \item \textbf{\textit{Spatial Correlation:}} To incorporate the spatial factor of the disease spread, we recursively generate the second component of the arrival rate, $\lambda^{2}_{i,t}$, on day $t \in \{1,2,\cdots,\mathcal{T}\}$ at location $i \in \{1,2,\cdots,L\}$ as
    \begin{align*}
       \lambda^{2}_{i,t} &= \sum_{{j=1,j \neq i}}^{L} \alpha_{ji,t} \theta_{ji} \sum_{l=1+t_{lag}}^{t_{lag}+r} \lambda^{2}_{j,t-l}  ,
    \end{align*}
    where $r$ is the length of the period during which disease spread from other locations would affect the current location, $t_{lag}$ denotes the delay effects in disease spread, $\theta_{ij}$ is the spatial correlation factor between location $i$ and location $j$,  $\alpha_{ij,t}$ is the infectious strength from location $i$ to location $j$. We let $\boldsymbol{\theta}  = [ \theta_{ij}, 1\le i,j\le L]$  and $\boldsymbol{\alpha} =  [ \alpha_{ij,t}, 1\le i,j\le L, 1\le t\le T]$.  
%    
%    
    % \item \textbf{\textit{Other Factors:}}  Except for the arrivals due to the temporal and spatial correlation, the arrival rate due to other factors is generated based on  
    % \begin{align*}
    %     \lambda_{i,t}^{3} \sim Piosson(\zeta_{i,t}\kappa_{i, (t-1)\% Y +1 })  , 
    % \end{align*} 
    % where $\kappa_{i, (t-1)\%Y+1}$ is the arrival rate by considering the day-of-week effect, and the $\zeta_{i,t}$ is the patient arrival strength because of other factors. We let $\boldsymbol{\kappa} = [\kappa_{i,t},1\le i \le L, 1\le t\le Y ]$ and $\boldsymbol{\zeta}_i = [\zeta_{i,t}, 1\le i \le L, 1\le t\le p]$. 
    \end{itemize} 

%   { \color{red}  Change 7 to parameter. } 

With the two components, the arrival rate at location $i$ on day $t$ is calculated as \begin{align*}
    \lambda_{i,t} = \lambda^{1}_{i,t} + \lambda^{2}_{i,t} . 
\end{align*}

% \textbf{Arrival Generation. } Based on the generated arrival rate, we generate arrivals to each hospital $a_{i,t} \sim Poisson(\lambda_{i,t})$. 

For our case study, we set \(\mathcal{T} =189\),  \(Y=7\), \(p=7\),  \(L=4\).
We estimate the parameters of the autoregressive model by fitting it to historical arrival data from an anonymous hospital. The resulting coefficients are
\([\phi_{i,1}, \phi_{i,2}, \ldots, \phi_{i,p}]=[0.061, -0.165, -0.042, -0.072, -0.148, 0.035, 0.588]\) for each \(i \in [L]\). 
%Furthermore, \(\lambda_{i,t}\) is the historical arrival rate at location $i$ on day $t$ (\(1 \leq i \leq L\)) for \(-p+1 \leq t \leq 0\). 
We set ${\alpha}_{ij,t}$, ${\beta}_{i,t}$, $\zeta_{i,t}$, and $\gamma_{i,t}$ as 
\[
\begin{cases}
    {\alpha}_{ij,t} = f_t,\\
    {\beta}_{i,t} = f_t,\\
    \zeta_{i,t} = c_{i} f_t,\\
    \gamma_{i,t} = c_{i}f_t , 
\end{cases}
\]
where 
\[
f_t
=
\begin{cases} 
(c_{peak} - 1) \cdot \sqrt{\frac{t - t_{start}}{t_{peak} - t_{start} }} + 1, & \text{if } t_{start} \leq t < t_{peak}, \\
(c_{peak} - 1) \cdot \sqrt{\frac{t_{end} - t}{t_{end} - t_{peak}}} + 1, & \text{if } t_{peak} \leq t \leq t_{end},\\ 
1, & \text{if }  t > t_{end}.   
\end{cases}
\] 
Here, we set $t_{start}=1$, $t_{peak} = 49$, $t_{end}=119$. We choose $[c_1, c_2, c_3,c_4] =[0.3,0.4,0.5,1]$ to reflect the relative magnitude of patient arrivals at the four hospitals (i.e., West, East, South, and Central) under consideration. 
We set $c_{peak}=1.5$ for the experiment discussed in Section~\ref{sec:DataDescription} and $c_{peak}=1.7$ for the experiment discussed in Appendix~\ref{sec:Impact_higher_demand}.

Finally, for $\theta_{ij}$, we have  
    \begin{align*}
        \theta_{ij} = \exp{\left(-z \frac{d_{ij}}{\max_{i,j \in \{1,...,L\}, i\neq j} d_{ij}}\right)} , 
    \end{align*}
where $z$ is a positive coefficient. 
Here, we set $z = 6.5$, along with the transition and discharge probabilities to be introduced in Section~\ref{sec:Patient_Transitions}, to ensure that the peak patient census across the four hospitals approximates the total staffed bed capacity of the anonymous hospital network, which operated at full capacity during the height of the pandemic \citep{AmericanHospitalDirectory}. 
%\red{How did you choose this $exponetial$ form with $z$ form? Any explanation? How did you choose the value of $-6.5$?}
%four IU Health hospitals (Arnett, Ball Memorial, Bloomington, and Methodist) 

Figure~\ref{fig:arrival_rate} displays the historical and simulated arrival rates over 30 weeks for each hospital. The red line denotes the expected arrival rate, while the gray lines show 30 randomly generated realizations. Days -20 to 0 correspond to the three-week historical input that uses in the rolling horizon approach. 

\begin{figure}[htbp]
    \centering
    \begin{subfigure}[t]{0.23\textwidth}
        \includegraphics[width=\textwidth]{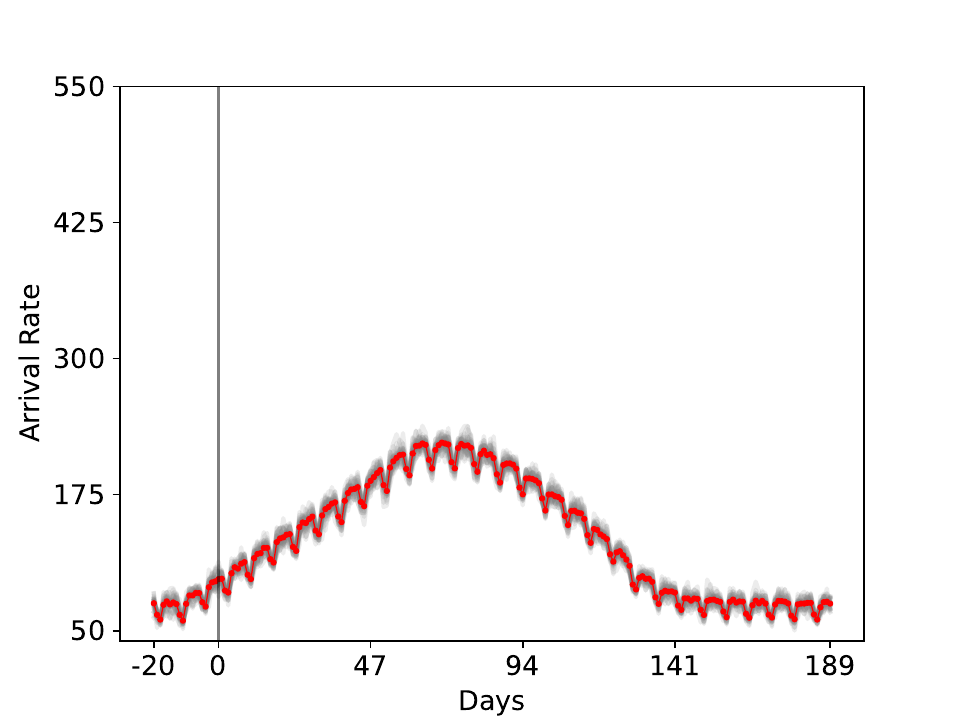}
        \caption{West Hospital}
    \end{subfigure}
        \begin{subfigure}[t]{0.23\textwidth}
        \includegraphics[width=\textwidth]{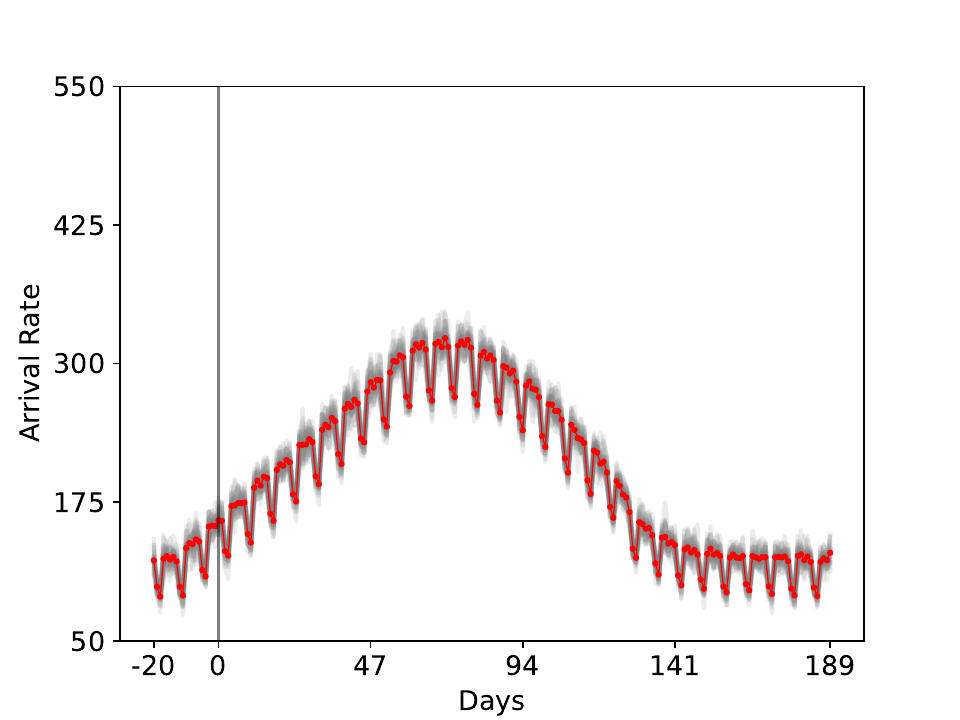}
        \caption{East Hospital}
    \end{subfigure}
        \begin{subfigure}[t]{0.23\textwidth}
        \includegraphics[width=\textwidth]{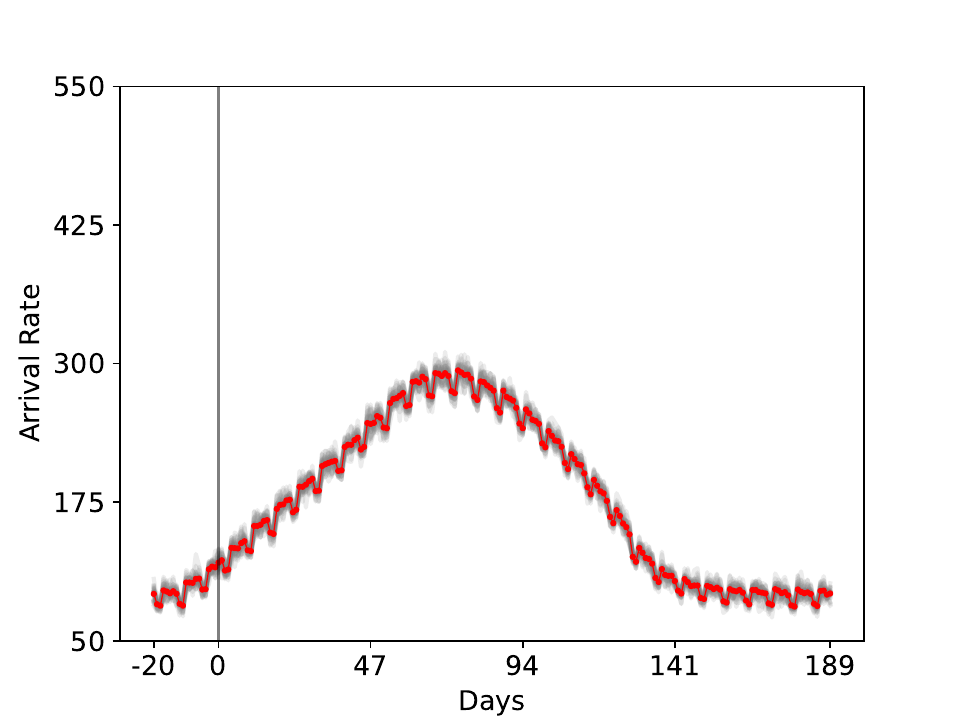}
        \caption{South Hospital}
    \end{subfigure}
        \begin{subfigure}[t]{0.23\textwidth}
        \includegraphics[width=\textwidth]{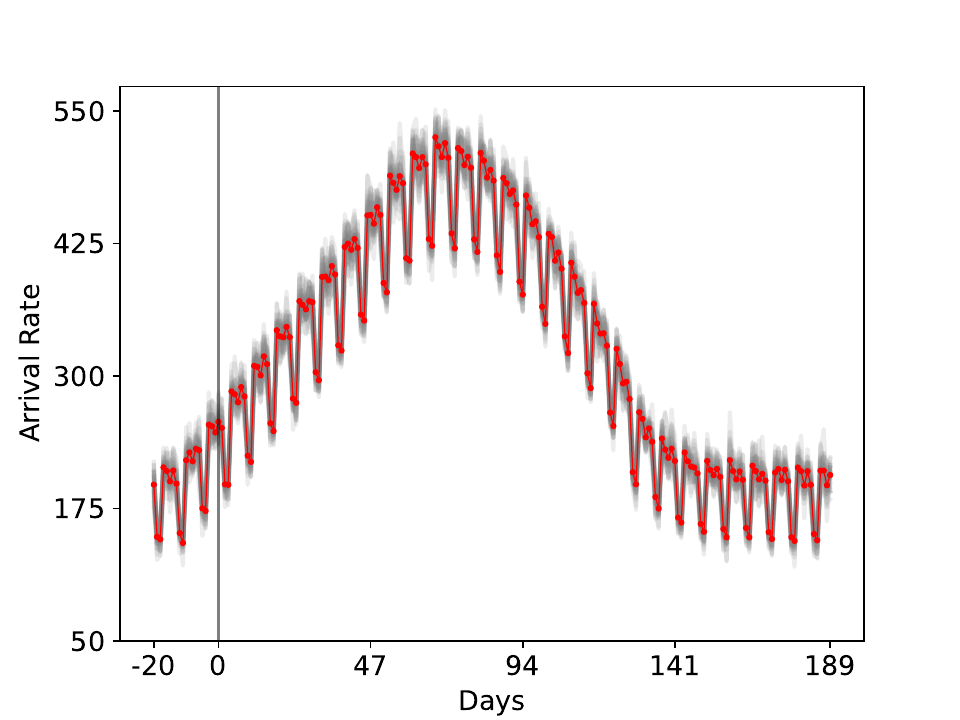}
        \caption{Central Hospital}
    \end{subfigure}
    \caption{Daily arrival rates for each hospital} 
    \label{fig:arrival_rate} 
\end{figure}

\subsubsection{Patient Transitions Within the Hospital}    
\label{sec:Patient_Transitions}  

In this section, we describe the patient transition dynamics across hospital units, which are critical for estimating patient census and, consequently, nurse demand. Each hospital under study contains three unit types: medical-surgical (MS), progressive care (PCU), and intensive care (ICU), differentiated by the level of care required. ICU patients require the most intensive treatment, PCU patients need moderate monitoring, and MS patients typically receive standard inpatient care.

Figure~\ref{fig_transition_diagram} illustrates the transition behavior of patients among these units. New patients are admitted daily to MS, PCU, or ICU at rates $\lambda_M$, $\lambda_P$, and $\lambda_I$, respectively. We use indices $M$, $P$, and $I$ to represent the units. Once admitted to unit $u \in {M, P, I}$, patients may remain, transfer to another unit, or be discharged (we use index $D$ for discharge), depending on their condition and progression of care needs; see \citealt{Nelson2023} for clinical context. The probability of transitioning from unit $u$ to $v$ is denoted by $p_{u,v}$ (if $u=v$, $p_{u,u}$ corresponds to staying in the current unit), and the discharge probability from unit $u$ is $p_{u,D}$ for $u \in {M, P, I}$. These probabilities satisfy the normalization condition $\sum_{v\in\{M,P,I\}} p_{u,v} + p_{u,D} = 1$.

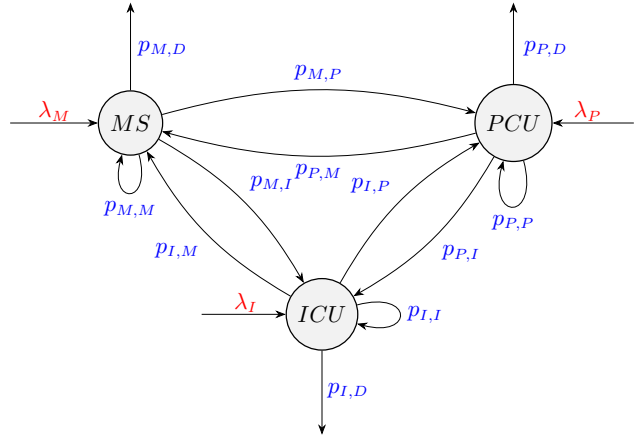
\begin{figure}[htbp]
\centering  
\scalebox{0.8}{
    \begin{tikzpicture} [ on grid,
every edge/.style = {draw, -{Stealth[scale=1]}, bend left=15}]
%every edge quotes/.append style = {auto, inner sep=2pt, font=\footnotesize}
%                        ]
\node (n1)  [state] {$MS$}; 
\node (n3)  [state,below right=of n1]   {$ICU$};  
\node (n2)  [state,above right=of n3]   {$PCU$};

\path   (n1)    edge ["$ \textcolor{blue}{p_{M,P}} $"] (n2)
                edge ["$ \textcolor{blue}{p_{M,I}}$"] (n3) 
              %  edge [loop below] (n1)  
        (n2)    edge ["$ \textcolor{blue}{p_{P,M}}$"] (n1)
                edge ["$  \textcolor{blue}{p_{P,I} }$"] (n3)
        (n3)    edge ["$ \textcolor{blue}{p_{I,M} }$"] (n1)
                edge["$ \textcolor{blue}{p_{I,P} }$"] (n2); 
                
 \draw[-{Stealth}] (n1) to[loop below] node{$  \textcolor{blue}{p_{M,M}}  $} ();                 
 \draw[-{Stealth}] (n2) to[loop below] node{$  \textcolor{blue}{p_{P,P}}  $} ();     
  \draw[-{Stealth}] (n3) to[loop right] node{$  \textcolor{blue}{p_{I,I}}  $} ();                     
                
         \draw [-{Stealth}] ($(n1)+(-2,0)$)  -- (n1);
         \draw [-{Stealth}]  (n1) --  ($(n1)+(0,2)$)  ; 
         
          \draw [-{Stealth}] (n2) --  ($(n2)+(0,2)$)   ; 
         \draw [-{Stealth}]  ($(n2)+(2,0)$)    -- (n2)  ;   
         
          \draw [-{Stealth}] (n3) --  ($(n3)+(0,-2)$)   ; 
         \draw [-{Stealth}]  ($(n3)+(-2,0)$)   -- (n3)  ;

         \node[anchor=north] at (-1.2,.5) { \textcolor{red}{$\lambda_M$ }  };       
          \node[anchor=north] at (0.6,1.5) { \textcolor{blue}{$p_{M,D}$  }  };       
          
           \node[anchor=north] at ($(n2)+(1.3,0.5)$)  { \textcolor{red}{$\lambda_P$ }  };       
          \node[anchor=north] at ($(n2)+(0.6,1.5)$)   { \textcolor{blue}{$p_{P,D}$  }  };     
          
           \node[anchor=north] at ($(n3)+(-1.2,0.5)$)  { \textcolor{red}{$\lambda_I$ }  };       
          \node[anchor=north] at ($(n3)+(0.5,-1)$)   { \textcolor{blue}{$p_{I,D}$  }  };
%https://tex.stackexchange.com/questions/18389/tikz-node-at-same-x-coordinate-as-another-node-but-specified-y-coordinate 
    \end{tikzpicture}
    }
    \caption{ Transition diagram in MS, PCU and ICU}
  \label{fig_transition_diagram} 
 \end{figure}

Table~\ref{tab:transition_prob} presents the set of transition and discharge probabilities used to generate the testing sample paths for each hospital. As noted earlier, these probabilities are calibrated to ensure that the peak patient census aligns with the total staffed bed capacity of the anonymous hospital network. These values are used throughout the numerical study in Section~\ref{sec:CaseStudy}. It is important to note that these probabilities differ from those estimated directly from the data. For reference, Table~\ref{tab:transition_prob_original_monday} shows the empirically estimated transition probabilities for a representative weekday (Monday). Despite these differences, the managerial insights derived in Section~\ref{sec:CaseStudy} remain robust. For completeness, we provide results using the empirical probabilities in Section~\ref{sec:casestudy_true_transition_prob}.

% The table reveals some interesting insights. Firstly, we observe that patients in ICU, PCU, and MS have an increasing probability of being discharged, which is consistent with the fact that patients in these units have decreasing levels of severity. Moreover, the patients in MS have a much smaller probability of being transferred to ICU compared to being transferred to PCU since their situation may deteriorate gradually. Conversely, the patients in ICU have the same probability of continuing to stay in ICU or being transferred to PCU, which is much higher than the probability of being transferred to MS directly. Similarly, the patients in PCU have a much larger probability of being transferred to MS compared to the probability of staying in PCU or being transferred to ICU. 
% These observations provide valuable insights into the transition behavior of patients, which will inform our nurse demand estimation model. 

 \begin{table}
\parbox{.45\linewidth}{
\centering
    \begin{tabular}{|l|r|r|r|r|} 
    \hline
    Unit  & \multicolumn{1}{l|}{MS} & \multicolumn{1}{l|}{PCU} & \multicolumn{1}{l|}{ICU} & \multicolumn{1}{l|}{Discharge} \bigstrut\\
    \hline
    MS    & 0.05  & 0.2   & 0.1  & 0.65 \bigstrut\\
    \hline
    PCU   & 0.25   & 0.05  & 0.1  & 0.6 \bigstrut\\
    \hline
    ICU   & 0.5   & 0.4  & 0.05  & 0.05 \bigstrut\\
    \hline
    \end{tabular}% 
  \caption{Adjusted probabilities} 
      \label{tab:transition_prob} 
}
\hfill
\parbox{.45\linewidth}{
\centering
    \begin{tabular}{|l|r|r|r|r|} 
    \hline
    Unit  & \multicolumn{1}{l|}{MS} & \multicolumn{1}{l|}{PCU} & \multicolumn{1}{l|}{ICU} & \multicolumn{1}{l|}{Discharge} \bigstrut\\
    \hline
    MS    & 0.7675 & 0.0125 & 0.0124 & 0.2076 \bigstrut\\
    \hline
    PCU   & 0.0852 & 0.7463 & 0.0258 & 0.1427 \bigstrut \\
    \hline
    ICU   & 0.1044 & 0.0596 & 0.7849 & 0.0511 \bigstrut \\
    \hline
    \end{tabular}% 
  \caption{Estimated probabilities} 
      \label{tab:transition_prob_original_monday} 
}
\end{table}

\subsubsection{Patient Census Estimation}  \label{sec:PatientsCensusGeneration} 
% Combined with the pseudocode 

Let $N_{i, u}^t$ be the patient census at unit $u$ in hospital $i$ at the beginning of day $t$. The census on day $t+1$ can be recursively determined based on the census at day $t$, accounting for patient inflows (arrivals and transfers-in) and outflows (discharges and transfers-out).

Specifically, we assume that the number of patients staying in the same unit ($n_{i,uu}^t$), are discharged ($n_{i,uD}^t$), or are transferred to other units ($n_{i,uv}^t$ and $n_{i,uw}^t$) in hospital $i$ on day $t$ follows a multinomial distribution, with parameters $N_{i,u}^t$ and ($ p_{i,uu}^{(t-1)\%Y+1}$, $ p_{i,uD}^{(t-1)\%Y+1}$, $ p_{i,uv}^{(t-1)\%Y+1}$, $ p_{i,uw}^{(t-1)\%Y+1}$). 
Here, $u,v,w\in\{M,P,I\}$ (with $v\neq u$ or $w\neq u$ representing the transfer to a different unit) and the superscript reflects the day-of-week effect. The transition and discharge probabilities are specified in Appendix~\ref{sec:Patient_Transitions}.

For patient arrivals, let $\Lambda_{i,u}^t$ be the number of new patients arriving at unit $u$ in hospital $i$ on day $t$. We first generate the total number of arrivals, $A_{i}^t$, as described in Section~\ref{sec:Patient_Arrivals}, and then allocate them across units using a multinomial distribution with parameters $A_i^t$ and probabilities ($q_{i,M}$, $q_{i,P}$, $q_{i,I}$), estimated as (0.7659, 0.153, 0.0811) from the anonymous hospital data.

The patient census at each unit at the end of day t is then updated as:
$$N_{i,u}^{t+1} = N_{i,u}^{t } + {\Lambda}_{i,u}^t +  \sum_{v \in \{ M,P,I\} } n_{i,vu}^t  - \sum_{v \neq u,v \in \{ M,P,I\}  } n_{i,uv}^t  -n_{i,uD}^t.$$  
The corresponding pseudocode is provided in Algorithm~\ref{alg:Nurse_Demand_Estimation}. Estimated patient censuses for each unit in our case study are shown in Figure~\ref{fig:patient_census}, with the shaded areas indicating 30 randomly generated trajectories. Days -20 to 0 correspond to the three-week historical input used in the rolling horizon approach.

\vspace{-0.1in}

\begin{algorithm}
    \caption{Patient Census Estimation}
    \label{alg:Nurse_Demand_Estimation}
    \begin{itemize}
        \item \textbf{Input}:  $\mathcal{T} $, $L$, $[A_{i,t} , 1\le i \le L, 1 \le t \le \mathcal{T} ]$, $[q_{i,u}, 1\le i \le L, u \in \{M, P, I\} ]$, $[p_{i,uv}^t, 1\le i \le L, 1\le t \le Y,  u,v \in \{M, P, I, D\} ]$, $[N_{i,u}^{1} ,1\le i \le L,  u \in \{M, P, I\} ]$.  
        \item \textbf{Main Iteration}: For each $t \in [\mathcal{T}-1]$ and $u \in \{M, P, I\}$, perform:
        \begin{enumerate}
            \item Estimate the patient census at unit $u$ in hospital $i$ on day $t+1$ from  
            \begin{align*}
                        N_{i,u}^{t+1} &= N_{i,u}^{t } + {\Lambda}_{i,u}^t +  \sum_{v \in \{ M,P,I\} } n_{i,vu}^t  - \sum_{v \neq u,v \in \{ M,P,I\}  } n_{i,uv}^t  -n_{i,uD}^t, \\
                \text{where} \\
                          & (\Lambda^t_{i,M}, \Lambda^t_{i,P}, \Lambda^t_{i,I}) \sim Multinomial(A_{i}^t, q_{i,M}, q_{i,P}, q_{i,I}),  \\
&(n_{i,uu}^t, n_{i,uD}^t, n_{i,uv}^t, n_{i,uw}^t) \sim Multinomial(N_{i,u}^t, p_{i,uu}^{(t-1)\%Y+1}, p_{i,uD}^{(t-1)\%Y+1}, p_{i,uv}^{(t-1)\%Y+1}, p_{i,uw}^{(t-1)\%Y+1}). 
            \end{align*}
        \end{enumerate}
        \item \textbf{Output}: $N_{i,u}^{t+1}$ 
    \end{itemize}
\end{algorithm}

       \begin{figure}
    \centering
    \begin{subfigure}[t]{0.23\textwidth}
        \includegraphics[width=\textwidth]{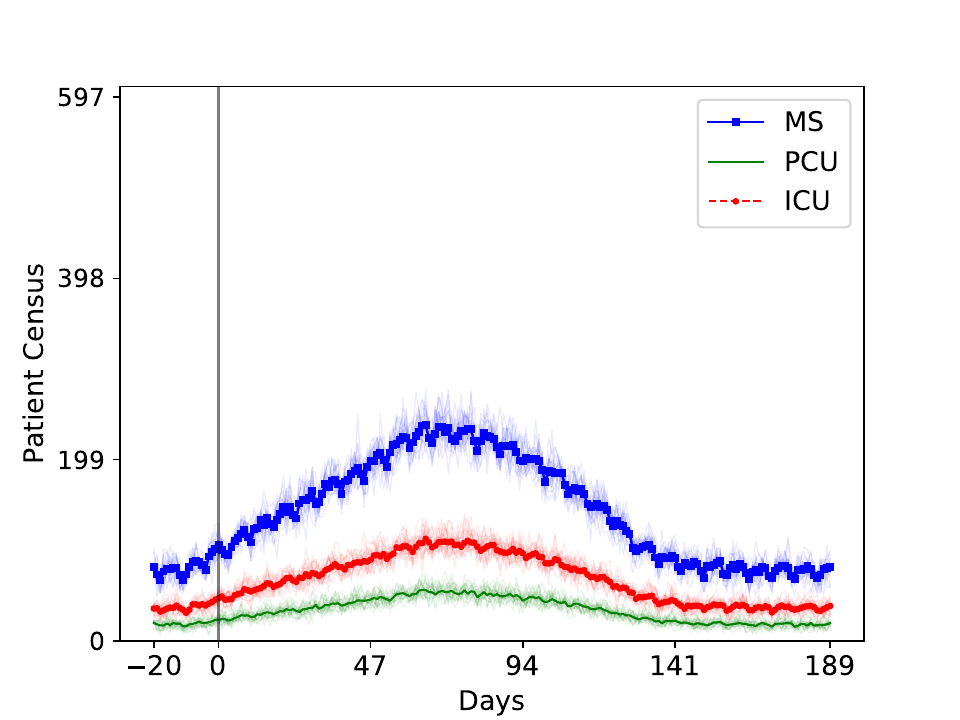}
        \caption{West Hospital}
    \end{subfigure}
        \begin{subfigure}[t]{0.23\textwidth}
        \includegraphics[width=\textwidth]{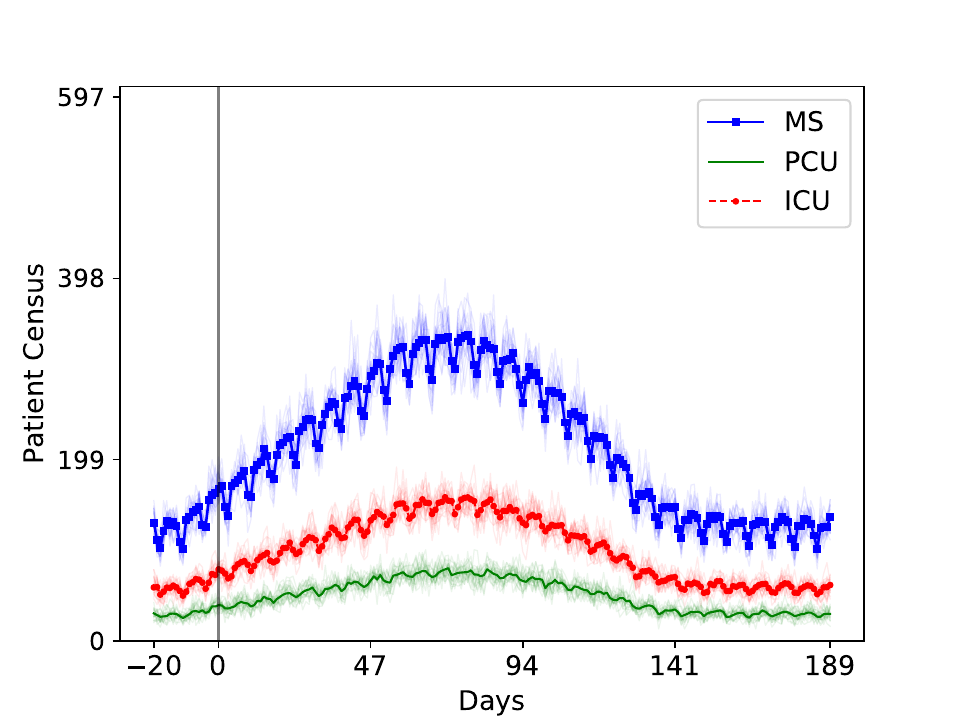}
        \caption{East Hospital}
    \end{subfigure}
        \begin{subfigure}[t]{0.23\textwidth}
        \includegraphics[width=\textwidth]{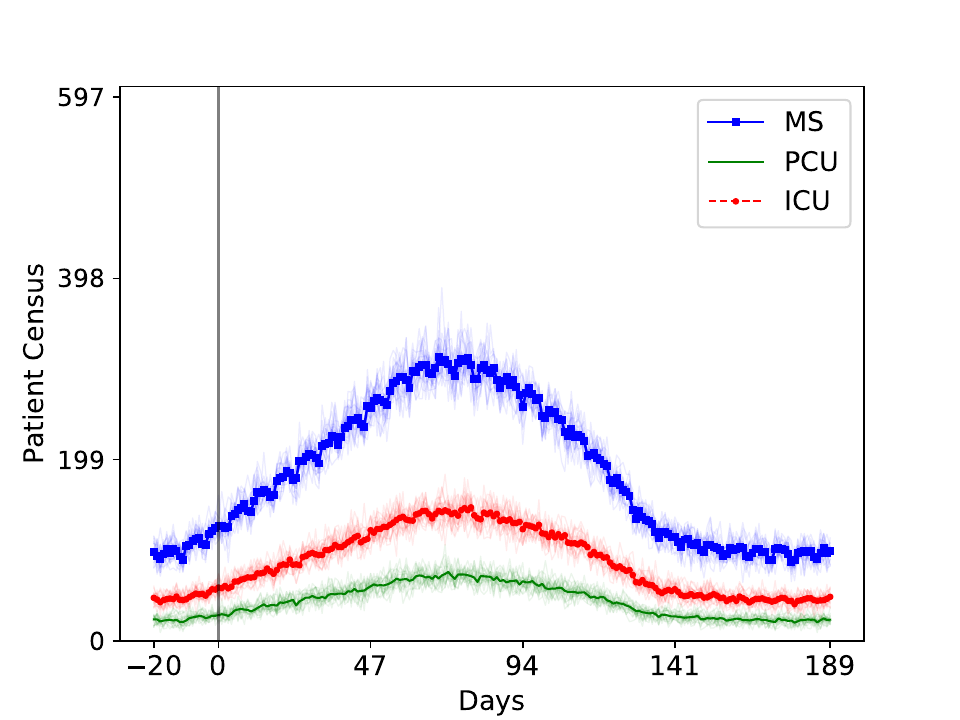}
        \caption{South Hospital}
    \end{subfigure}
        \begin{subfigure}[t]{0.23\textwidth}
        \includegraphics[width=\textwidth]{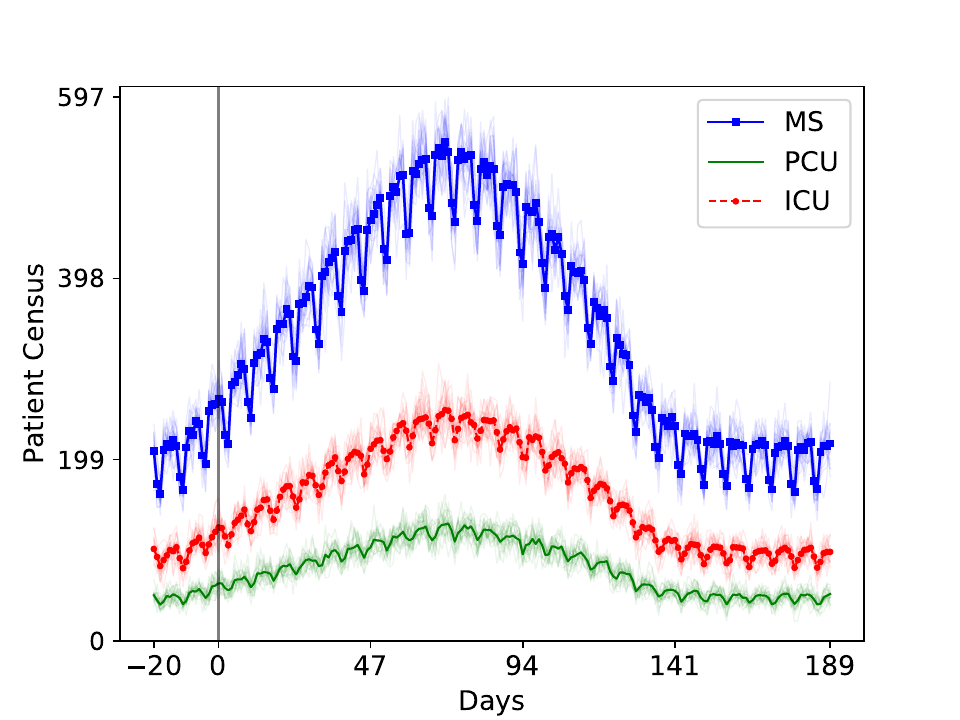}
        \caption{Central Hospital}
    \end{subfigure}
    \caption{Daily patient census at each unit for each hospital}
    \label{fig:patient_census}
\end{figure}

\subsubsection{Nurse Demand Estimation}  \label{sec:NurseDemandEstimation}

We adopt the nurse-to-patient staffing ratios used in the anonymous hospital system, which align with standard practices across U.S. hospitals (\citealt{yankovic2011identifying}): 1:5 for MS units, 1:3 for PCU, and 1:2 for ICU. Based on these ratios, the nurse demand at hospital $i$ on day $t$ is estimated by aggregating the demand across all units as follows: 
$$  \xi_i^t = \frac{N_{i,M}^{t } } {5}  +\frac{N_{i,P}^{t } } {3}  + \frac{N_{i,I}^{t } } {2}.$$   
Figure~\ref{fig:nurse_demand} shows the estimated nurse demand for each hospital in our case study. The red line indicates the expected demand, while the gray lines represent 30 randomly generated nurse demand trajectories. Days -20 to 0 correspond to the three-week historical input used in the rolling-horizon forecasting process.

\begin{figure}[htbp]
    \centering
    \begin{subfigure}[t]{0.23\textwidth}
        \includegraphics[width=\textwidth]{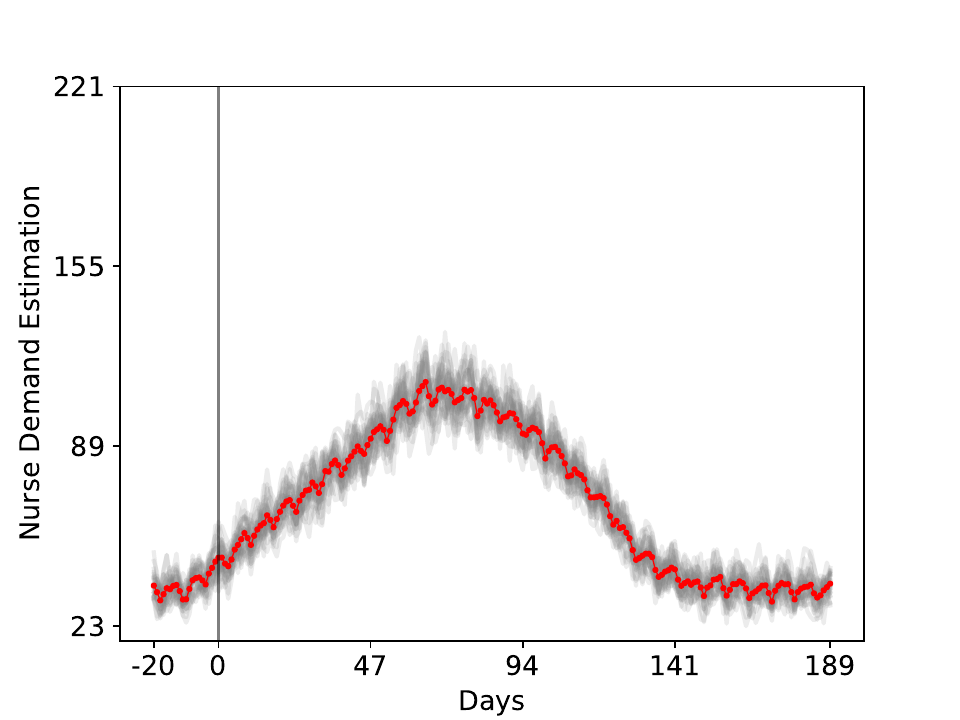}
        \caption{West Hospital}
    \end{subfigure}
        \begin{subfigure}[t]{0.23\textwidth}
        \includegraphics[width=\textwidth]{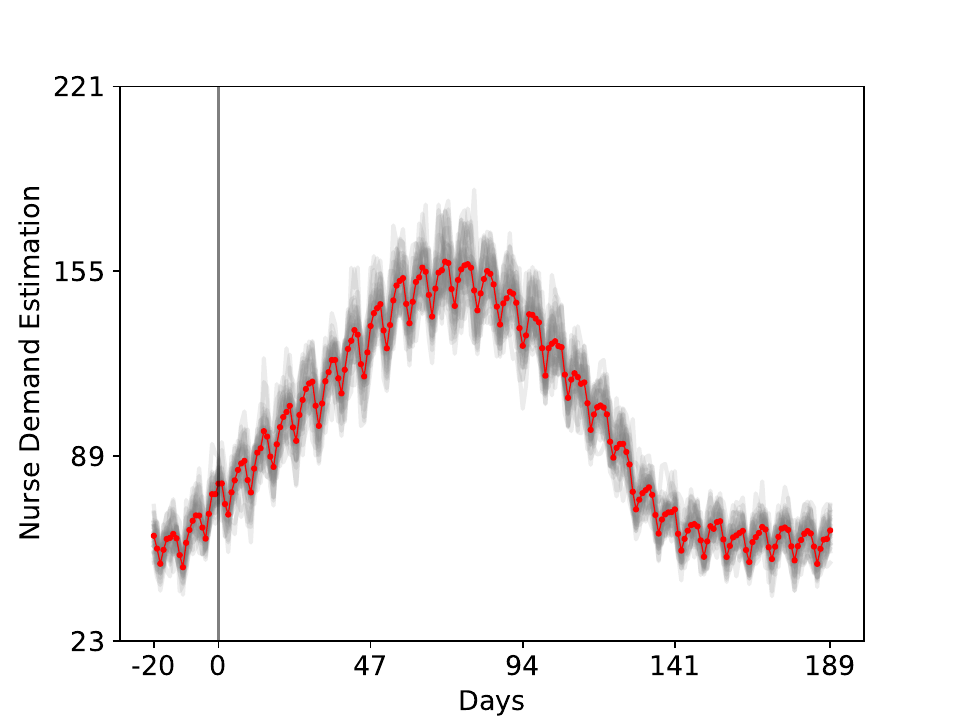}
        \caption{East Hospital}
    \end{subfigure}
        \begin{subfigure}[t]{0.23\textwidth}
        \includegraphics[width=\textwidth]{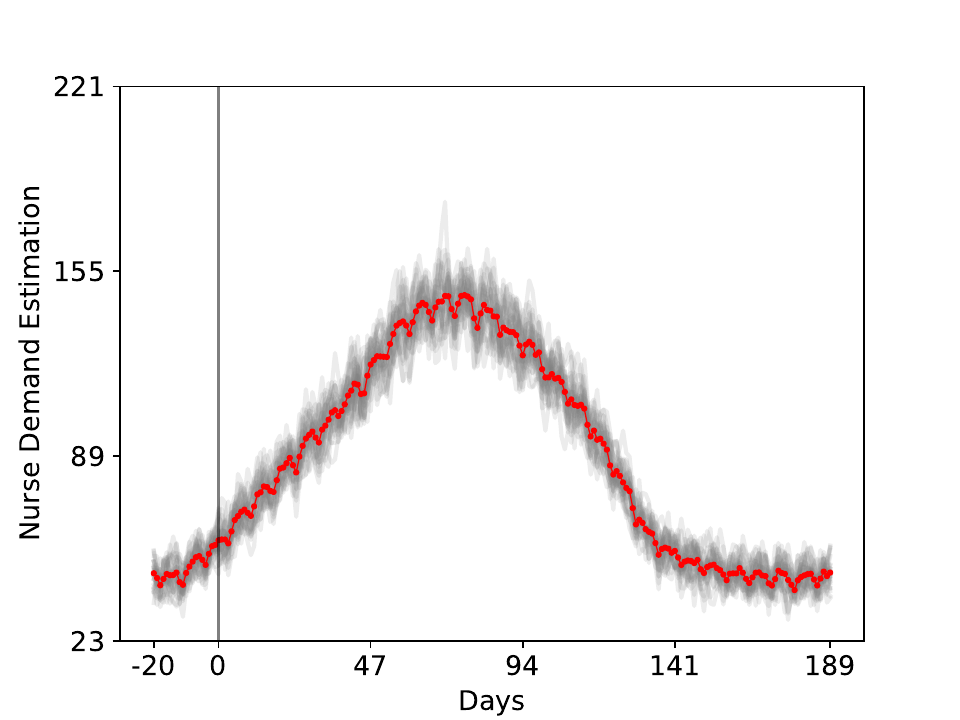}
        \caption{South Hospital}
    \end{subfigure}
        \begin{subfigure}[t]{0.23\textwidth}
        \includegraphics[width=\textwidth]{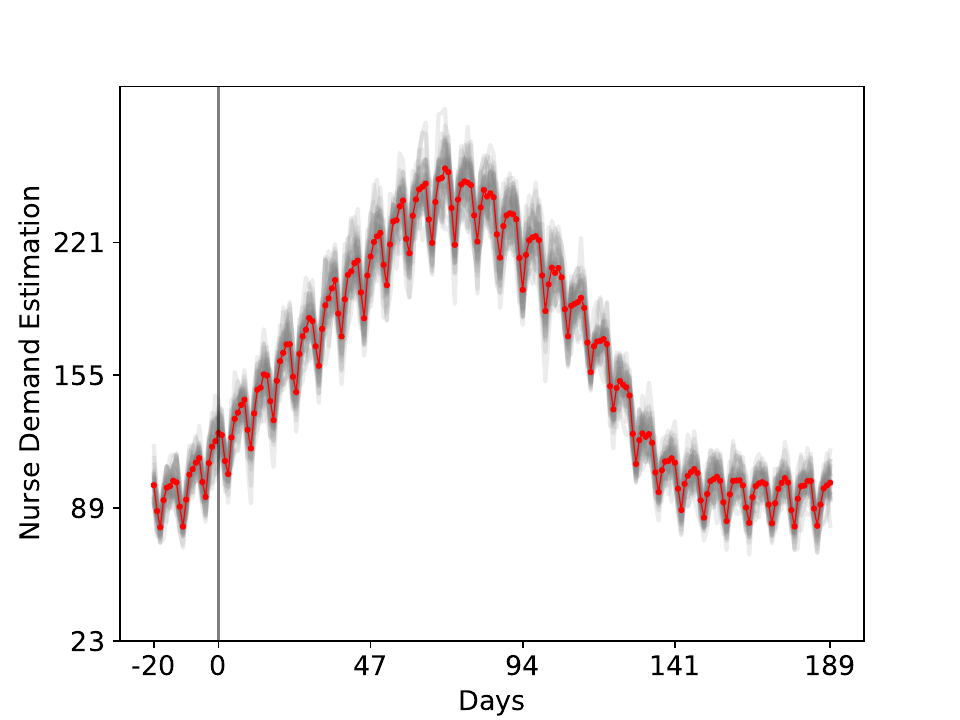}
        \caption{Central Hospital}
    \end{subfigure}
    \caption{Daily nurse demand estimation for each hospital}
    \label{fig:nurse_demand}
\end{figure}

%{\bf Remark:}  This comprehensive approach enables the generation of nurse demand at each hospital over an extended period. These demands constitute the testing sample path for nurse demand. It is essential to note that we estimate nurse demand using the patient census rather than historical nurse demand data. This choice is influenced not only by the lack of access to historical nurse demand data but also by the recognition that nurse demand is primarily driven by patient census. The explicit modeling of the patient transition process between MS, PCU, and ICU accurately captures patient transitions and discharges in the hospital, providing a highly credible estimation of the patient census at each unit. Using this estimated patient census, we then apply the nurse-to-patient ratio to calculate nurse demand.

%{\bf Give the figure for nurse demand} 

\subsection{Training Sample Path Generation}    \label{sec:TrainingSamplePathGeneration}      

In this section, we describe the procedure for generating the training sample paths, which differs from the one used for testing. This distinction is necessary due to the nonstationary nature of patient arrivals and nurse demand, which can lead to significant discrepancies between historical and future patterns. To address this, we adopt a rolling prediction approach: using data from the past three weeks to estimate key patient flow parameters and generate demand forecasts for the upcoming weeks.

Specifically, we estimate the arrival rates, transition probabilities, and discharge probabilities using daily data on patient arrivals, discharges, and inter-unit transfers from the prior three weeks at each hospital. Arrival rates are estimated by taking the day-of-week-specific sample average over three observations (e.g., the past three Mondays), denoted by $\hat {\lambda}_{i,u}^t$ for unit $u$ on day $t$ in hospital $i$. Transition probabilities—$\hat p_{i,uu}^t$, $\hat p_{i,uv}^i$, and $\hat p_{i,uw}^t$—are calculated as the average proportion of patients moving between units, also conditioned on the day of the week. Similarly, discharge probabilities $\hat p_{i,uD}^t$ are estimated using the same logic. With these estimated parameters, we simulate patient flows and compute nurse demand following the procedures in Sections~\ref{sec:PatientsCensusGeneration} and~\ref{sec:NurseDemandEstimation}; the output forms the training sample paths.

% {\bf Remark:}  It is important to recognize that our model, which predicts patient census and nurse demand utilizing historical transfer data, is one of numerous potential prediction models. Our modeling framework is designed for adaptability and is capable of incorporating a variety of prediction models. Furthermore, in the absence of historical transfer data, where only past census data are available, a direct predictive model can be employed to forecast patient census, subsequently facilitating the estimation of nurse demand. 

\subsection{Travel Cost Setting}\label{sec:TravelCostSetup}

We set the non-salary transfer cost (i.e., the reimbursement provided to nurses who travel to a remote hospital), denoted by $\tau_{ij}$, based on the geographical distance $d_{ij}$ between hospitals $i$ and $j$, for all $i, j \in [L]$. The transfer bonuses, distances, and secondment settings across locations were provided in Table~\ref{tab:parameters} in the main paper, and we explain here on how we estimate the values.

To build the connection between the distance and the non-salary bonus, we first establish a baseline by identifying the hospital pair ($\hat{i}, \hat{j}$) with the shortest distance:
$(\hat i, \hat j) = \argmin_{i,j}  d_{ij}$. 
We then set the minimum non-salary transfer cost for this pair as $\tau_{\hat i, \hat j}  = \tau_{\hat j, \hat i }  = 1.1$.
For all other hospital pairs, we define the non-salary transfer cost as linearly increasing with distance:
$$\tau_{ i,  j}  = \tau_{\hat i, \hat j} +   k (d_{i,j} - {d_{\hat i, \hat j}  }  ), $$ 
where $k = 0.01$ governs the cost increment per unit distance. This linear formulation is consistent with assumptions commonly made in the literature (e.g., \citealt{chou2023taxi}).

\subsection{Nurse Capacity Setting}
\label{sec:NurseCapacityGeneration}     

We adjust the nurse capacity as outlined in Algorithm~\ref{alg:semi-capacity-alg}, using prior nurse capacity levels and historical nurse demand values ($\xi$'s) to reflect the practice that hospitals adjust staffing levels adaptively based on demand. In our setup, we set $W = 27, L = 4, Y = 7, m = 2$, and $n = 0.8$. The initial nurse capacity for the first two weeks is given by $[C_{1,1}, C_{2,1}, C_{3,1}, C_{4,1}] = [C_{1,2}, C_{2,2}, C_{3,2}, C_{4,2}] = [40, 120, 110, 130]$, and the adjustment scale for each hospital is determined by the hyperparameter vector $[b_1, b_2, b_3, b_4] = [0.11, 0.17, 0.17, 0.115]$.
The resulting nurse capacity trajectories for each hospital are illustrated in Figure~\ref{fig:capacity}, where the red line indicates the expected capacity and the gray lines represent 30 simulated trajectories.

%The historical and estimated nurse demand for each hospital are displayed in Figure \ref{fig:nurse_demand}, where the red line represents the expected demand and the grey lines represent 30 random nurse demands. The data from day -20 to day 0 represents historical input, while data from day 1 to day 189 represent the estimated nurse demand.

\begin{figure}[htbp]
    \centering
    \begin{subfigure}[t]{0.23\textwidth}
        \includegraphics[width=\textwidth]{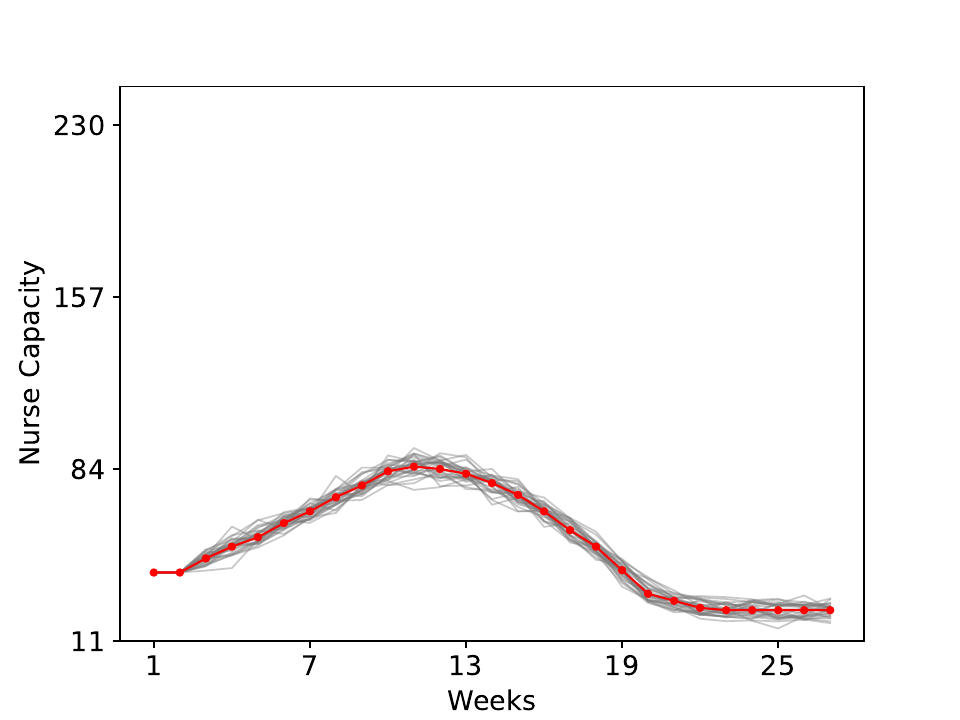}
        \caption{West Hospital}
    \end{subfigure}
        \begin{subfigure}[t]{0.23\textwidth}
        \includegraphics[width=\textwidth]{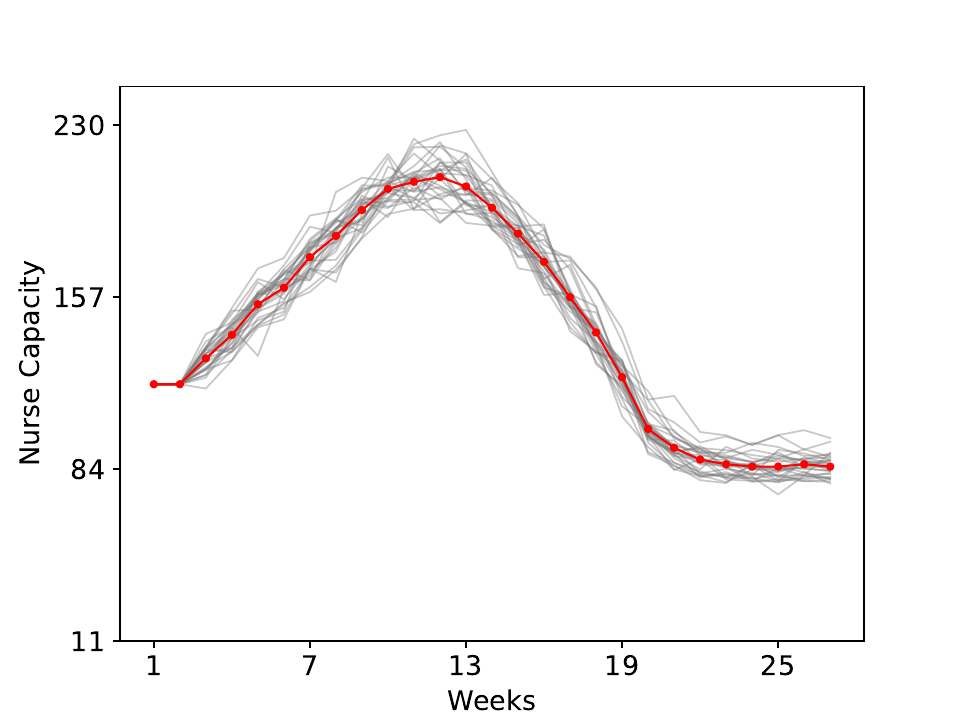}
        \caption{East Hospital}
    \end{subfigure}
        \begin{subfigure}[t]{0.23\textwidth}
        \includegraphics[width=\textwidth]{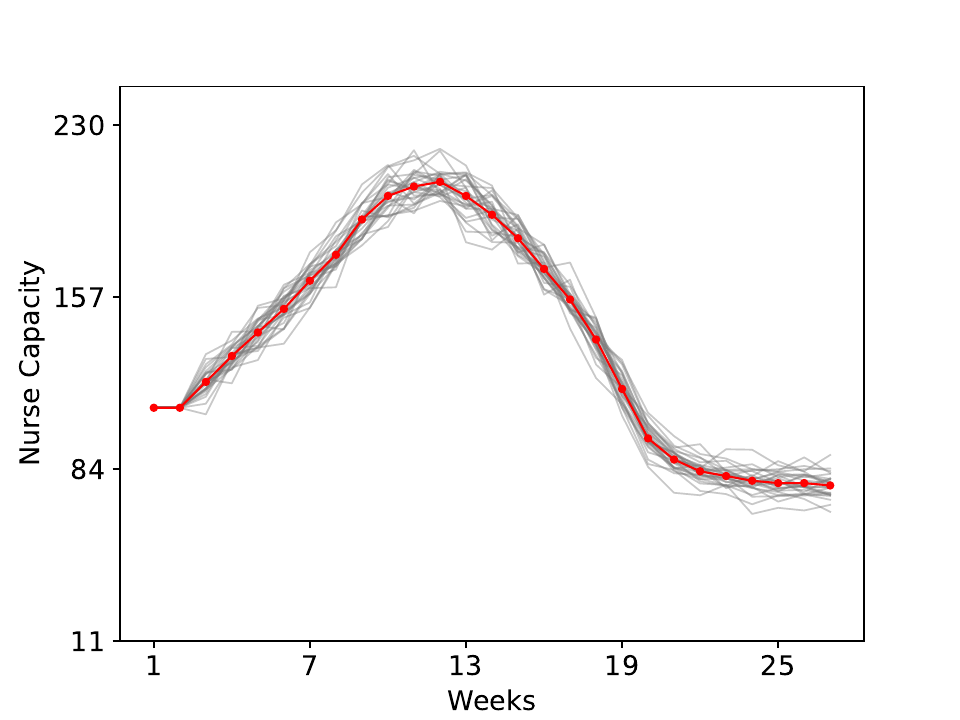}
        \caption{South Hospital}
    \end{subfigure}
        \begin{subfigure}[t]{0.23\textwidth}
        \includegraphics[width=\textwidth]{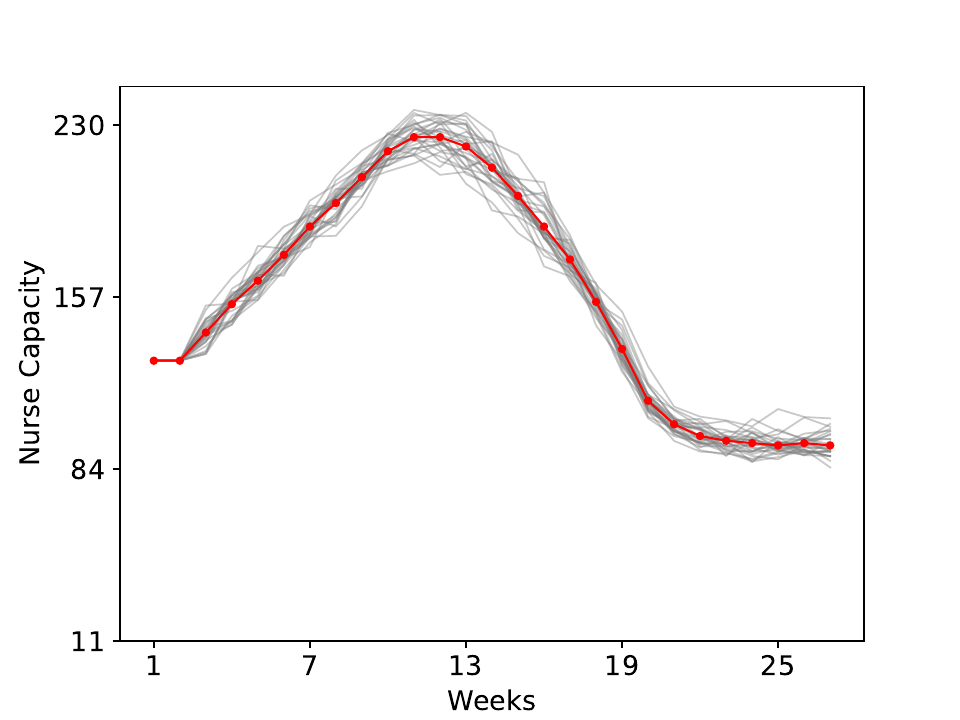}
        \caption{Central Hospital}
    \end{subfigure}
    \caption{Weekly nurse capacity for each hospital}
    \label{fig:capacity}
\end{figure}

% {\bf Show the figure for capacity} 

\begin{algorithm}[btbp]
    \caption{Nurse Capacity Estimation} 
    \label{alg:semi-capacity-alg}  
    \begin{itemize}
        \item \textbf{Input}: $W$, $L$, $Y$,  $m$, $n$, $[C_{i,w}, 1\le i \le L] $, $[b_1, \cdots, b_L]$, $[{\xi}_{i,t}, { 1\le i \le L, Y(w-3)+1\le t\le  WY } ] $ 
        \item \textbf{Main Iteration}: For $w\in \{3,4,\cdots, W\}$ and $i\in [L]$, perform:      
     \begin{enumerate}
        \item Calculate 
            \begin{align*}
                D_{i,w} &=  \frac{\sum_{t=Y(w-2)+1}^{Y(w-1)}\xi_{i,t}}{Y} - \frac{\sum_{t=Y(w-3)+1}^{Y(w-2)}\xi_{i,t}}{Y},\\ 
                C_{i,w} &= C_{i,w-1} + \phi_w(D_{i,w}) b_i D_{i,w}, 
            \end{align*}
        \quad    where 
            \begin{align*}
    \phi_w(D_{i,w}) = 
    \begin{cases}
        m, \quad \text{ if } D_{i,w} \geq 0,\\
        n, \quad \text{\; if } D_{i,w} < 0. 
    \end{cases}
\end{align*}
\end{enumerate}
        \item \textbf{Output}: $C_{i,w}$ 
    \end{itemize}
\end{algorithm}

\section{Supplementary Analysis for Case Study}
\label{sec:Additional_results_case_study}     
\subsection{Alternative Demand Pattern for Section~\ref{sec:Secondment_effect}}
\label{sec:figures}   

To investigate the effect of secondment beyond the baseline setting, we consider an alternative demand pattern, where only the West Hospital experiences a nurse shortage (the other three hospitals have surplus capacity). The nurse demand and capacity for each hospital are shown in Figure~\ref{nurse_demand_capacity_over_time_special_demand}. 
The corresponding results are summarized in Table~\ref{tab:Average_metrics_special_case}.

\begin{figure}
     \centering
     \begin{subfigure}[b]{0.23\textwidth}
         \centering
         \includegraphics[width=\textwidth]{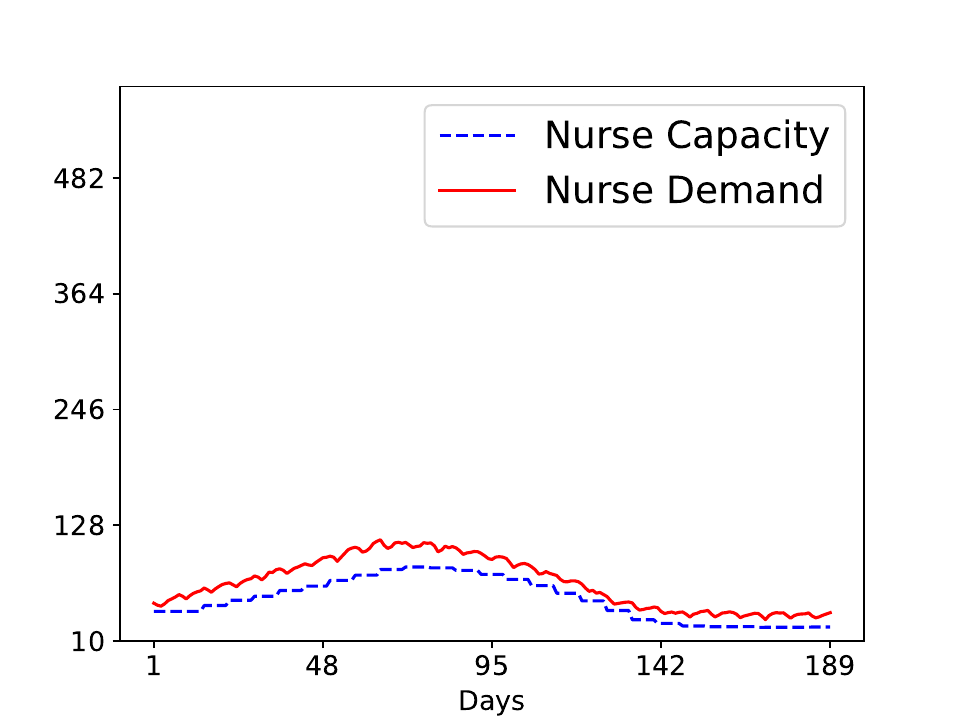}
         \caption{West Hospital} 
         \label{census_Arnett_special_demand} 
     \end{subfigure}
     \hfill 
     \begin{subfigure}[b]{0.23\textwidth}
         \centering
         \includegraphics[width=\textwidth]{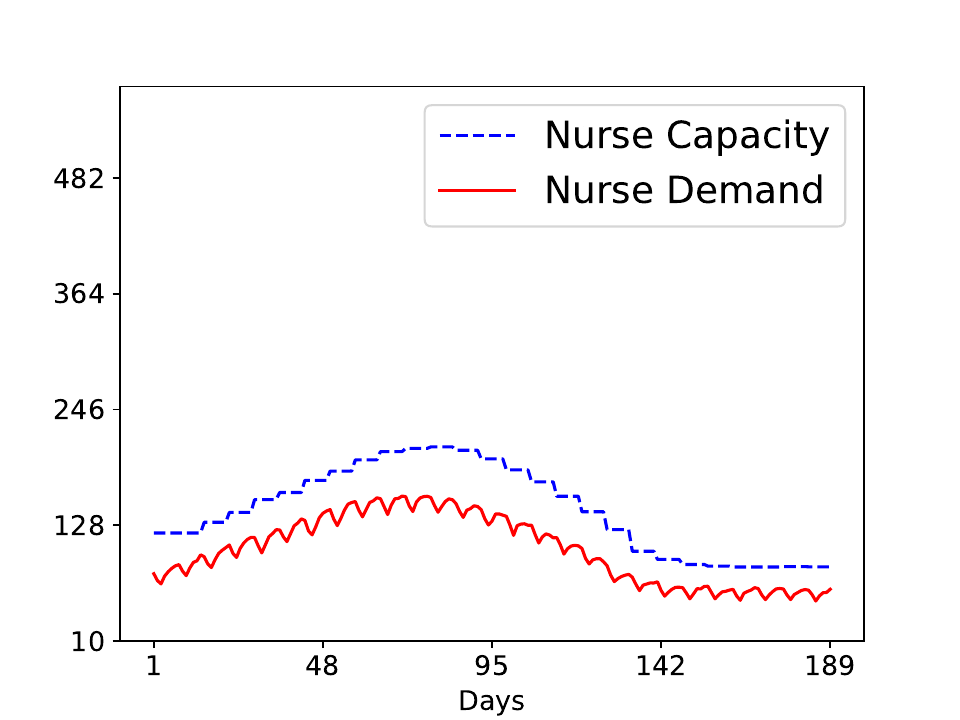}   
         \caption{East Hospital} 
         \label{census_Ball_special_demand}     
     \end{subfigure} 
          \begin{subfigure}[b]{0.23\textwidth}
         \centering
         \includegraphics[width=\textwidth]{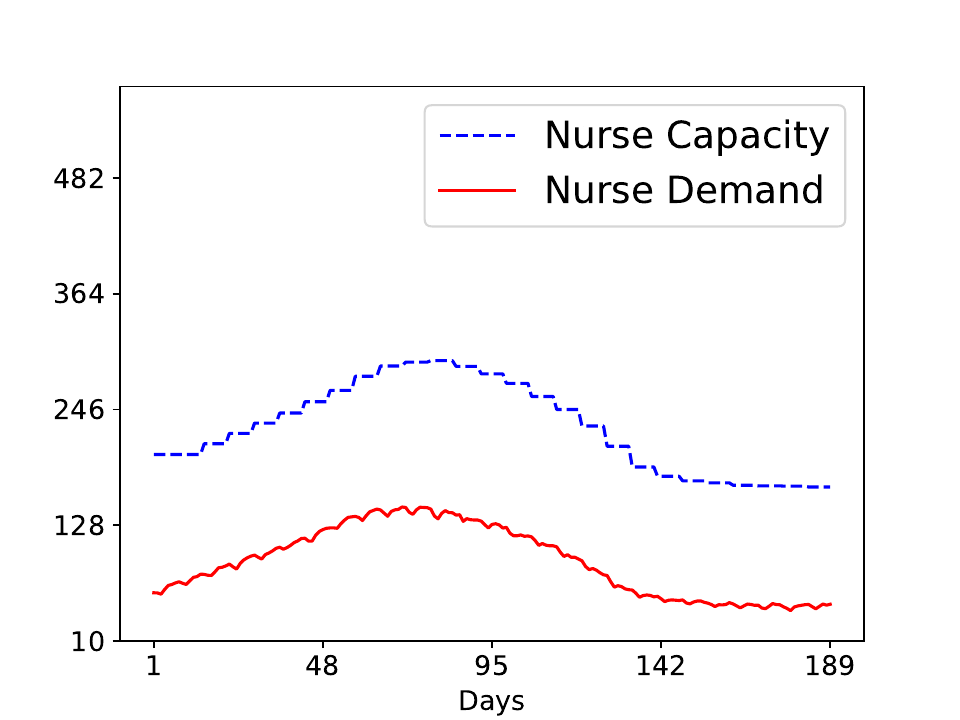}
         \caption{South Hospital} 
         \label{census_Bloomington_special_demand} 
     \end{subfigure}
     \hfill 
     \begin{subfigure}[b]{0.23\textwidth}
         \centering
         \includegraphics[width=\textwidth]{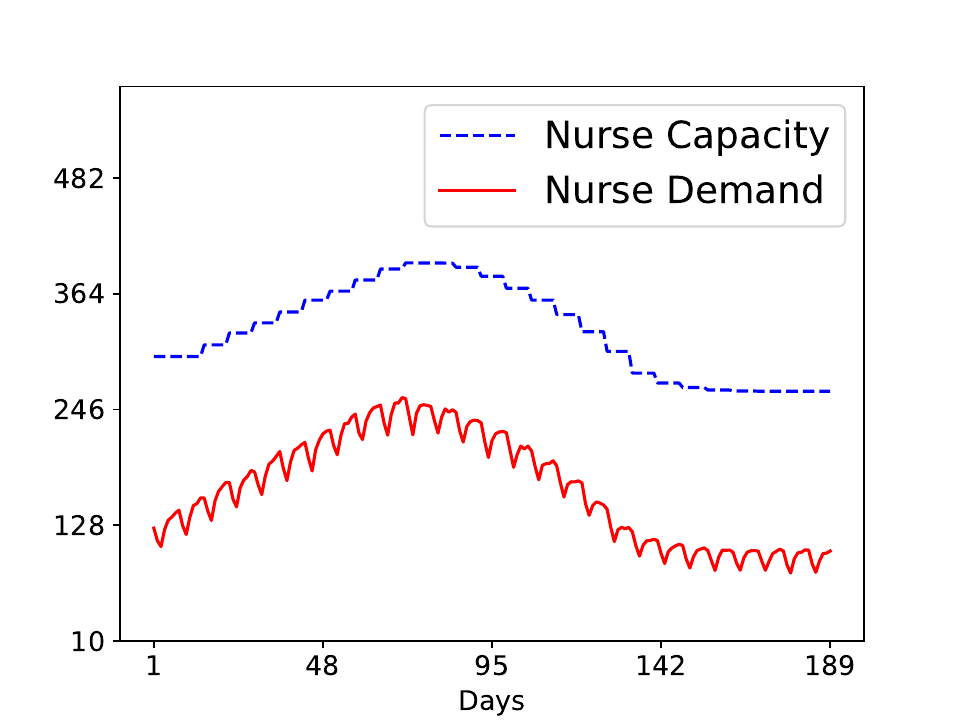}   
         \caption{Central Hospital} 
         \label{census_Methodist_special_demand}     
     \end{subfigure} 
         \caption{Daily nurse demand and capacity at each hospital under the special demand pattern}  
        \label{nurse_demand_capacity_over_time_special_demand}
\end{figure} 

% Table generated by Excel2LaTeX from sheet 'Sheet5'
\begin{table}[htbp]
  \centering
  \caption{Average cost, deployed transfers, and transferred miles per week for SAA and SRO across various network designs and secondment scenarios under the alternative demand pattern 
 % \red{after a second thought, let's move this table to Appendix~\ref{sec:figures} to save space for the main paper. I added the numbers in the main text.}
  }
  \scalebox{0.8}{
    \begin{tabular}{|c|c|rrcccc|}
    \hline
    \multirow{2}[4]{*}{Metrics} & \multicolumn{1}{c|}{\multirow{2}[4]{*}{Network design}} & \multicolumn{2}{c|}{Baseline secondment} & \multicolumn{2}{c|}{One-day secondment } & \multicolumn{2}{c|}{Three-day secondment } \bigstrut\\
\cline{3-8}          &       & \multicolumn{1}{l|}{SAA} & \multicolumn{1}{l|}{SRO} & \multicolumn{1}{l|}{SAA} & \multicolumn{1}{l|}{SRO} & \multicolumn{1}{l|}{SAA} & \multicolumn{1}{l|}{SRO} \bigstrut\\
    \hline
    \multicolumn{1}{|c|}{\multirow{2}[4]{*}{Cost}} & FC    & 329.46 & 325.78 & \multicolumn{1}{r}{413.50} & \multicolumn{1}{r}{408.52} & \multicolumn{1}{r}{302.29} & \multicolumn{1}{r|}{302.29} \bigstrut\\
\cline{2-2}          & HS    & 358.95 & 353.25 & -     & -     & -     & - \bigstrut\\
    \hline
    \multicolumn{1}{|c|}{\multirow{2}[4]{*}{Deployed transfers}} & FC    & 67.81 & 67.83 & \multicolumn{1}{r}{115.95} & \multicolumn{1}{r}{116.02} & \multicolumn{1}{r}{58.43} & \multicolumn{1}{r|}{58.43} \bigstrut\\
\cline{2-2}          & HS    & 116.08 & 116.08 & -     & -     & -     & - \bigstrut\\
    \hline
    \multicolumn{1}{|c|}{\multirow{2}[4]{*}{Transferred miles}} & FC    & 7459.06 & 7461.62 & \multicolumn{1}{r}{12754.37} & \multicolumn{1}{r}{12761.94} & \multicolumn{1}{r}{6427.52} & \multicolumn{1}{r|}{6427.52} \bigstrut\\
\cline{2-2}          & HS    & 7196.70 & 7196.70 & -     & -     & -     & - \bigstrut\\
    \hline
    \end{tabular}%
    }
  \label{tab:Average_metrics_special_case}%
\end{table}%

\subsection{Effect of Coordination Costs on Network Performance}\label{sec:coordination_cost}    
It is important to note that under the FC network, nurses from rural hospitals must be prepared to work in all other rural hospitals, in addition to the Central Hospital. 
This broader assignment scope may incur a coordination cost -- capturing the nurse dissatisfaction stemming from unfamiliarity with multiple workplaces. If we incorporate this cost, the performance of the FC and HS networks then vary depending on the magnitude of this coordination cost. 
For instance, when the coordination cost is set to 250, the HS network yields a lower overall cost 933.35 than the FC network {(with a cost of 963.30)} when the secondment length between rural regions is one, under the SAA approach. However, the FC network outperforms HS when the secondment length is two or more. For example, when the secondment length is two days, the cost is 891.90. 
A similar phenomenon can be observed under the SRO approach. 
%\red{Wei to add the number of the cost here for secondment length of two days. }
%This finding underscores the importance of incorporating longer secondment durations when enabling nurse transfers between rural hospitals, particularly in the presence of coordination-related frictions.

\subsection{Impact of Time-Lagging in Demand Prediction}\label{sec:Impact_time_lagging}   
% \noindent\textbf{Impact of time-lagging in demand prediction. }
In this section, we analyze the case where nurse demand predictions are based on data from the past six weeks (instead of three weeks in the baseline). 
%This can help us further evaluate the impact of underestimation. 
We show the nurse capacity, and nurse demand along with the nurse demand prediction for a representative sample path in Figure~\ref{nurse_demand_capacity_prediction_over_time_6week}. 
When the demand forecast is based on the past 6-week data, it shows an even greater underestimation of nurse demand during periods of increasing demand, compared to using the three-week data in the baseline. It also shows an even greater overestimation during periods of decreasing demand, due to the time-lagging effect. As a result, we see a even more significant improvements with SRO compared to SAA in Figure~\ref{Comparison_DRO_SAA_6week}. 
In particular, SRO adapts to the potential misestimation of the demand by assigning additional nurses, thereby reducing emergency transfers and achieving superior performance. 

\begin{figure}
     \centering
     \begin{subfigure}[b]{0.23\textwidth}
         \centering
         \includegraphics[width=\textwidth]{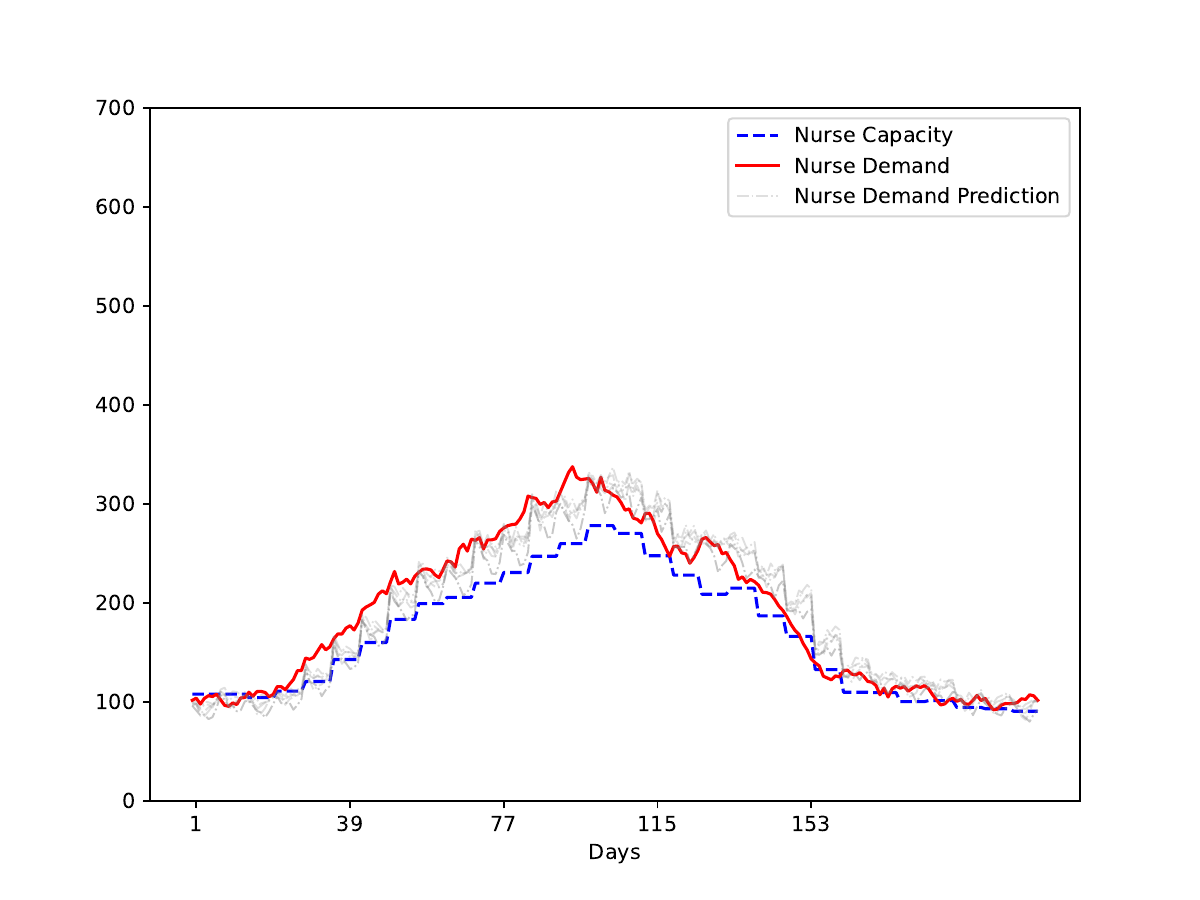}
         \caption{West Hospital} 
         \label{Arnett_demand_vs_capacity_training_6week} 
     \end{subfigure}
     \hfill 
     \begin{subfigure}[b]{0.23\textwidth}
         \centering
         \includegraphics[width=\textwidth]{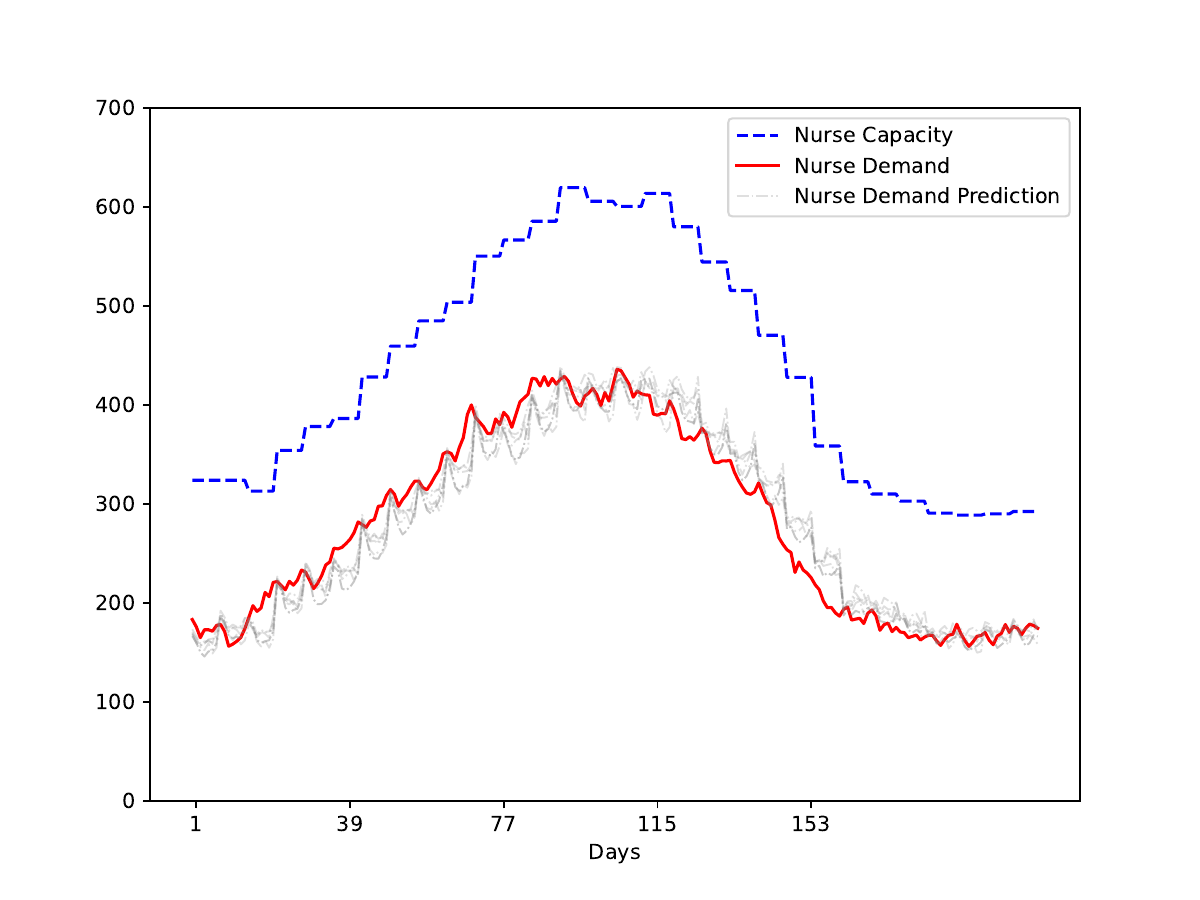}   
         \caption{East Hospital} 
         \label{Ball_Mem_demand_vs_capacity_training_6week}     
     \end{subfigure} 
          \begin{subfigure}[b]{0.23\textwidth}
         \centering
         \includegraphics[width=\textwidth]{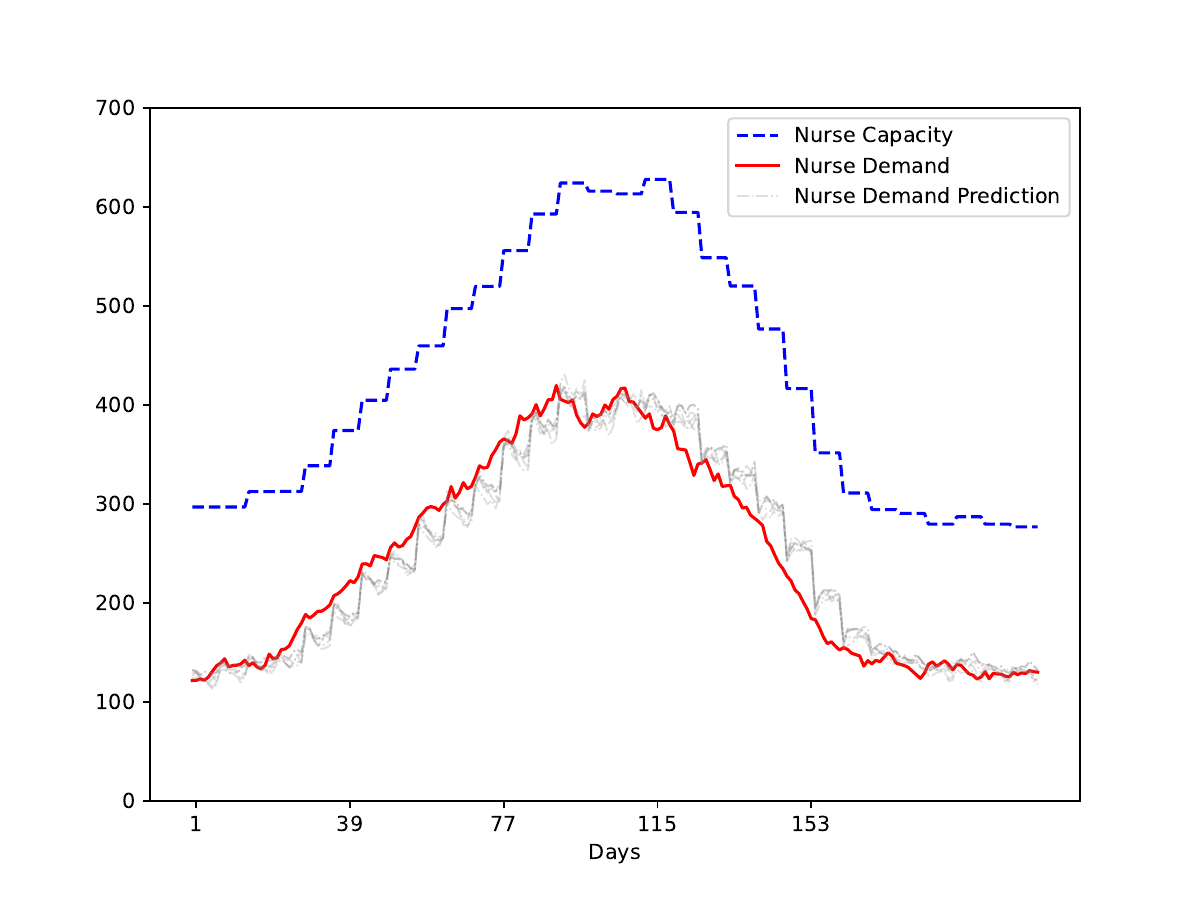}
         \caption{South Hospital} 
         \label{Bloomington_demand_vs_capacity_training_6week} 
     \end{subfigure}
     \hfill 
     \begin{subfigure}[b]{0.23\textwidth}
         \centering
         \includegraphics[width=\textwidth]{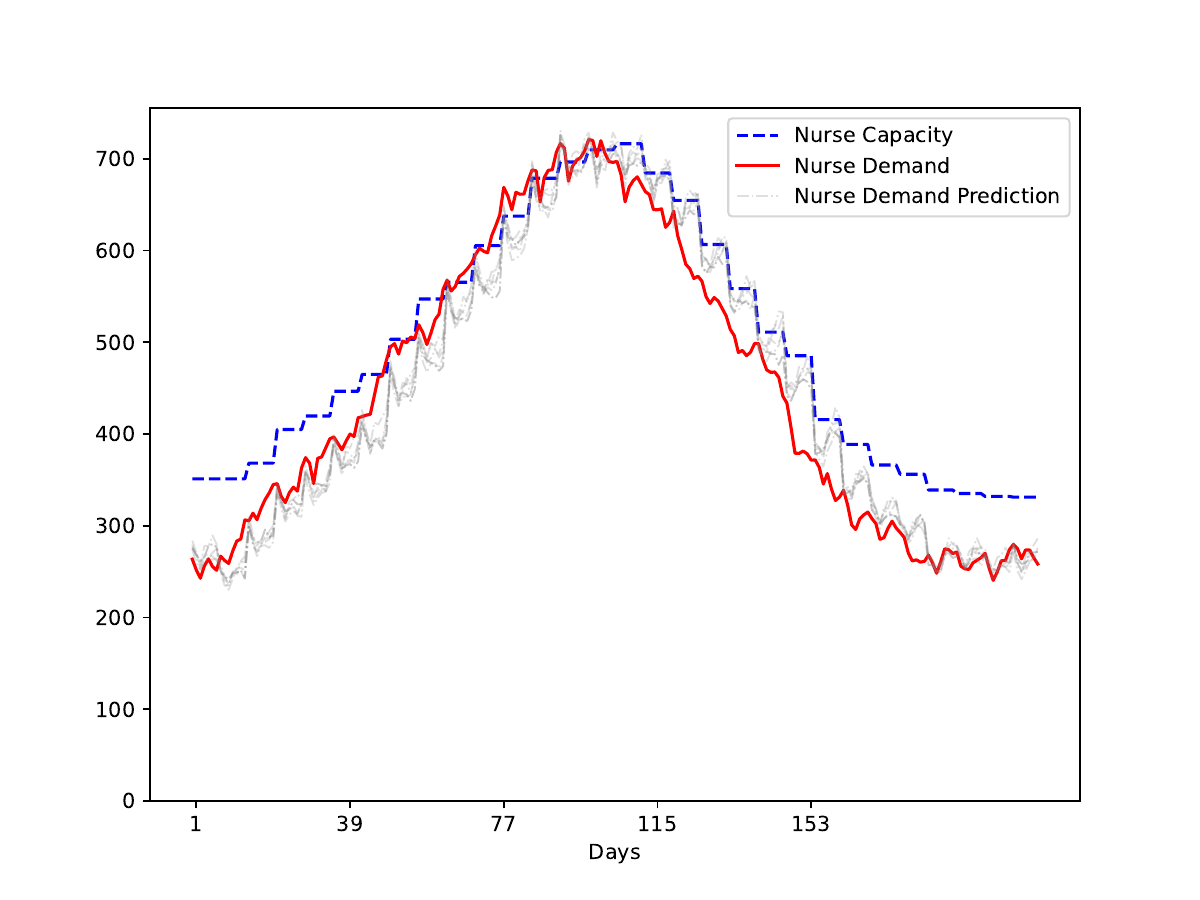}   
         \caption{Central Hospital} 
         \label{Methodist_demand_vs_capacity_training_6week}     
     \end{subfigure} 
         \caption{Daily nurse capacity and demand along with the demand prediction (use the data from the past six weeks) for one sample path} 
        \label{nurse_demand_capacity_prediction_over_time_6week} 
\end{figure}

\begin{figure} 
     \centering
     \begin{subfigure}[b]{0.3\textwidth} 
         \centering
         \includegraphics[width=\textwidth]{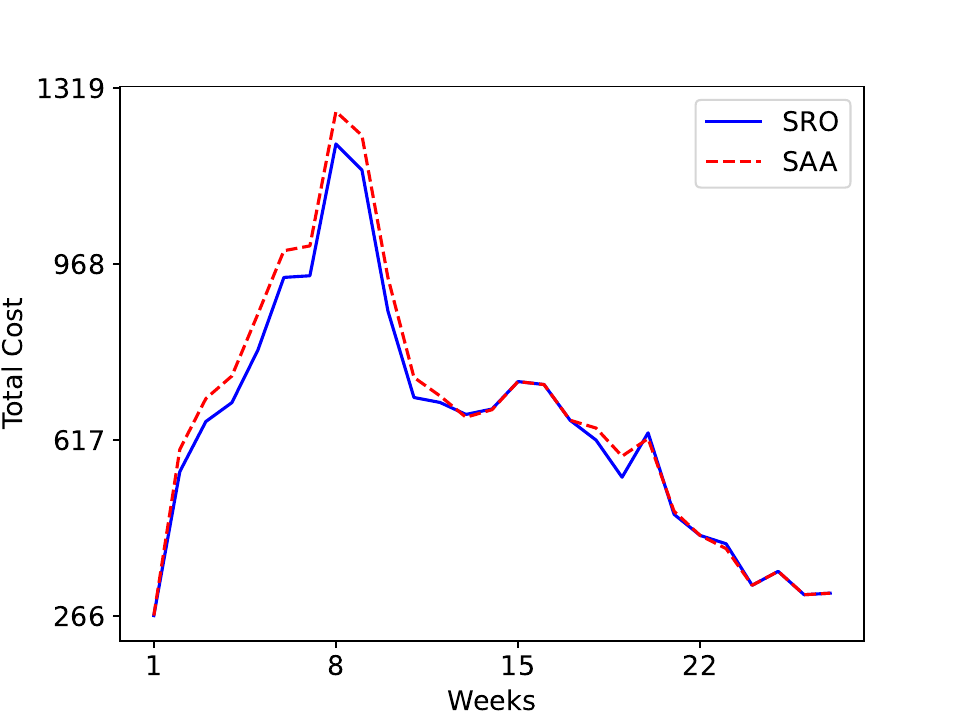} 
         \caption{Comparison on cost}
         %(\color{blue} { Draw the figure in square; x-axis: ``Weeks'' (from week 1);  y-axis:  ``Total cost'';  label: ``${SAA}$'',   ``${SRO}$'' } ) 
         \label{cost_constrast_6week}
     \end{subfigure}
     \hfill
     \begin{subfigure}[b]{0.3\textwidth} 
         \centering
         \includegraphics[width=\textwidth]{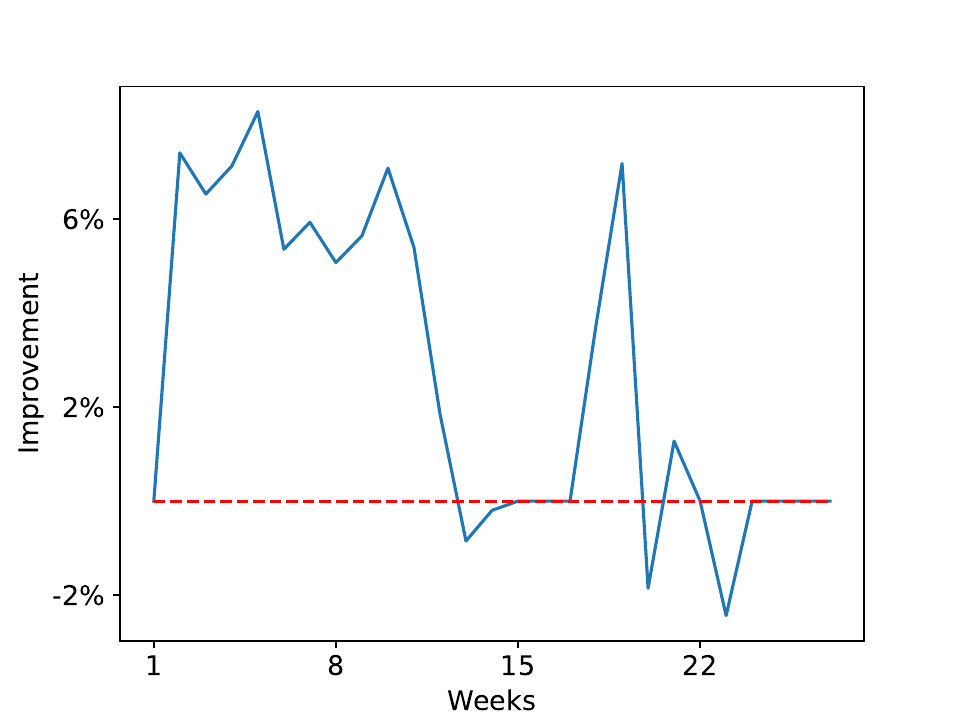}
         \caption{Improvement of SRO}  
         %  (\color{blue} { Draw the figure in square; x-axis: ``Weeks'' (from week 1);  y-axis:  ``Cost Improvement in \% of DRO Over SAA'';  no label ; draw a dashed line with y-axis$=0$ } )
         \label{cost_constrast_percentage_6week} 
     \end{subfigure}   
     \hfill
     \begin{subfigure}[b]{0.3\textwidth} 
         \centering
         \includegraphics[width=\textwidth]{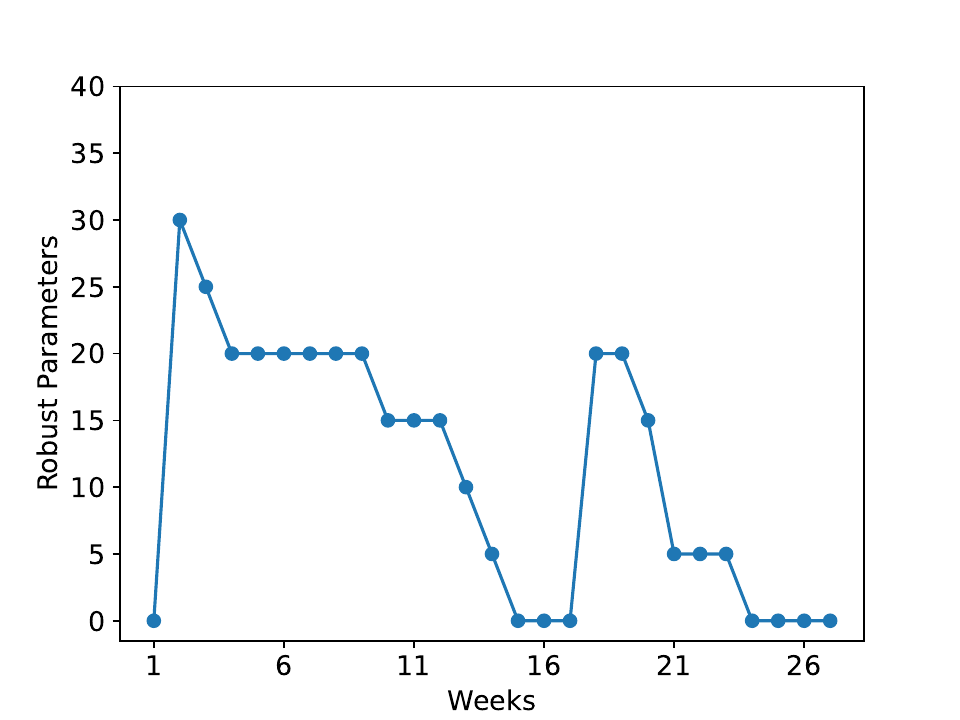}
            \caption{Robust parameter of SRO}   
         %  (\color{blue} { Draw the figure in square; x-axis: ``Weeks'' (from week 1);  y-axis:  ``Cost Improvement in \% of DRO Over SAA'';  no label ; draw a dashed line with y-axis$=0$ } )
         \label{fully_sec2_weekly_robust_parameter_selection_plot_6week} 
     \end{subfigure}        
        \caption{Weekly comparison between SAA and SRO, and robust parameter used by SRO (use the data from the past six weeks)
       % \red{Move this figure to appendix.}
        }   
        \label{Comparison_DRO_SAA_6week}
\end{figure}

\subsection{Impact of Higher Demand Peak}\label{sec:Impact_higher_demand}  
% \noindent\textbf{Impact of higher demand peak. } 
In this section, we explore the impact of variations in disease spread on the efficacy of our sample robust optimization approach. 
Specifically, we modify the disease spread pattern to increase the peak nurse demand while maintaining the same duration to reach this peak, resulting in a sharper increase in demand. 
The nurse demand and capacity are displayed in Figure~\ref{higher_peak_nurse_demand_capacity_over_time}, and the detailed setting is shown in Appendix \ref{sec:TestingSamplePathGeneration}.  
This adjustment enables us to assess how changes in the intensity of demand affect our optimization strategy. 
%Through this analysis, our goal is to thoroughly understand how the dynamics of disease spread influence the performance of the sample robust optimization approach we propose.

We compare SRO with SAA under the HC network and baseline secondment scenario. 
Specifically, we analyze the costs associated with SAA and SRO in Figure \ref{higher_peak_methodist_shortage_total_cost_sec_3_min_sec_2}. We show the improvement of SRO over SAA in Figure \ref{higher_peak_methodist_shortage_total_cost_percentage_sec_3_min_sec_2},
 and the chosen robust parameters for SRO in each week are shown in Figure \ref{higher_peak_fully_sec2_weekly_robust_parameter_selection_plot}. 
We see that the improvement of SRO over SAA is greater in the increasing demand phase.  
This is because West and Central Hospitals experience an even sharper nurse demand and shortage (compared to the baseline), prompting SRO to allocate more nurses from East and South Hospitals to these hospitals during the planning stage, as observed in Figures \ref{higher_peak_Planned_decision_DRO_SAA_arnett} and \ref{higher_peak_Planned_decision_DRO_SAA_methodist}. 
That is, SRO is more adaptable to the demand change.
In contrast, SAA lacks the flexibility to increase nurse transfers quickly in response to the sudden rise in demand. 
Consequently, SRO significantly reduces the overall nurse shortage compared to SAA.

  \begin{figure}
     \centering
     \begin{subfigure}[b]{0.23\textwidth}
         \centering
         \includegraphics[width=\textwidth]{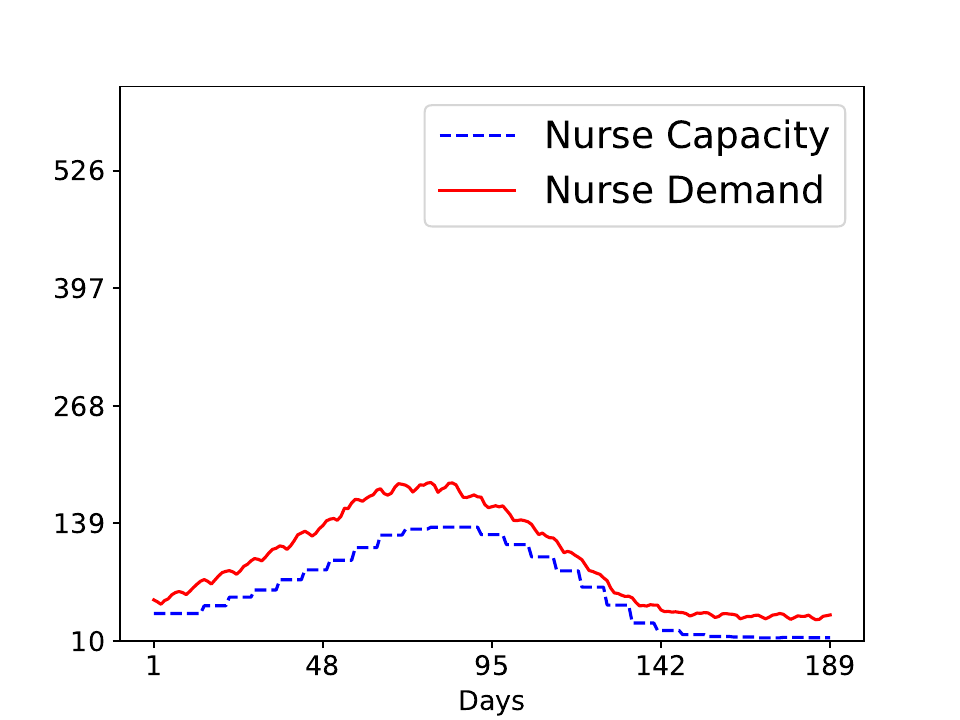}
         \caption{West Hospital} 
         \label{higher_peak_census_Arnett} 
     \end{subfigure}
     \hfill 
     \begin{subfigure}[b]{0.23\textwidth}
         \centering
         \includegraphics[width=\textwidth]{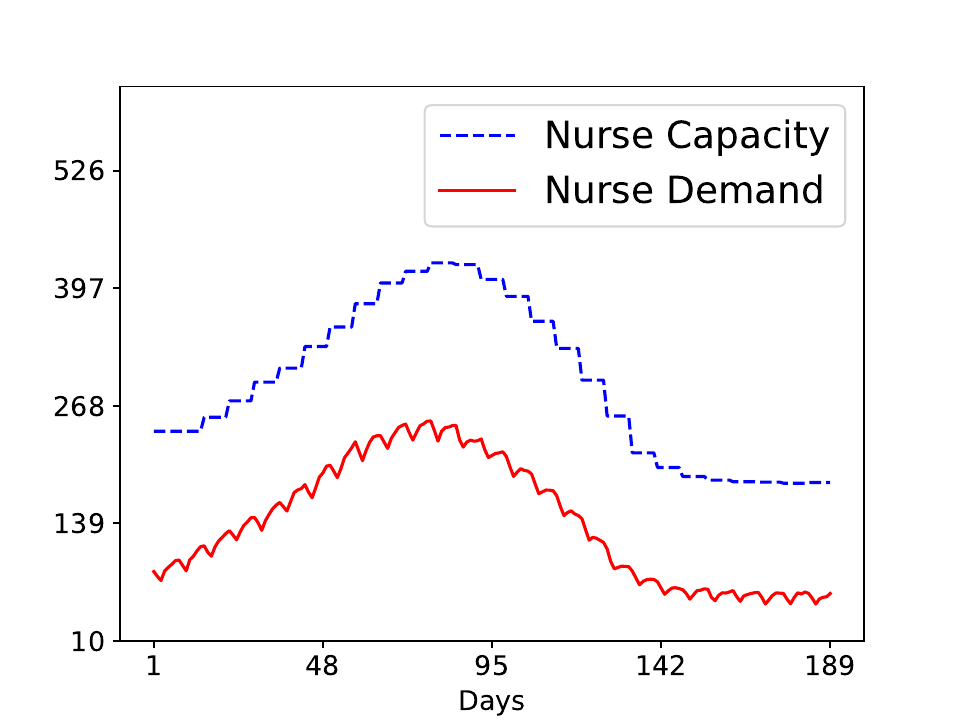}   
         \caption{East Hospital} 
         \label{higher_peak_census_Ball}    
     \end{subfigure} 
          \begin{subfigure}[b]{0.23\textwidth}
         \centering
         \includegraphics[width=\textwidth]{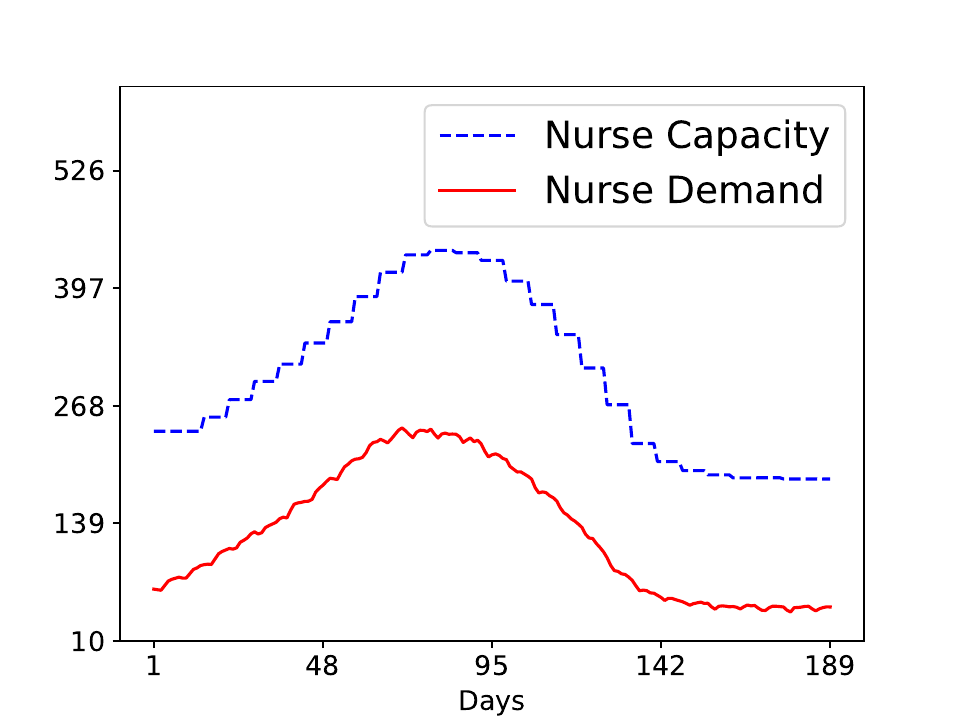}
         \caption{South Hospital} 
         \label{higher_peak_census_Bloomington} 
     \end{subfigure}
     \hfill 
     \begin{subfigure}[b]{0.23\textwidth}
         \centering
         \includegraphics[width=\textwidth]{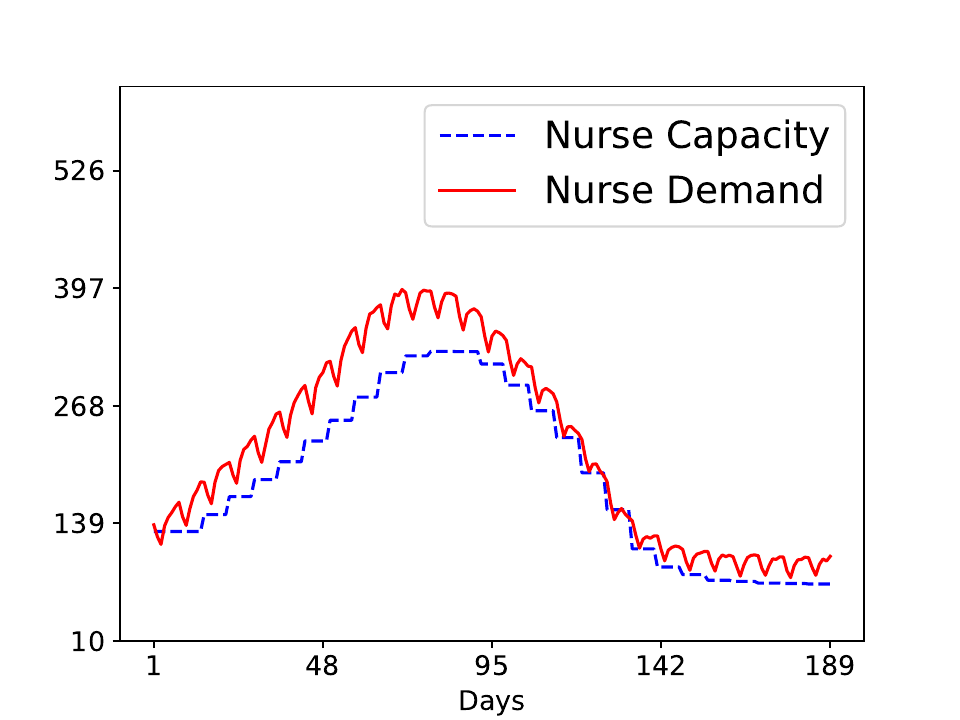}   
         \caption{Central Hospital}   
         \label{higher_peak_census_Methodist}     
     \end{subfigure}  
         \caption{Daily nurse demand and capacity at each hospital (under a higher demand peak)}
         % \red{Move this figure to appendix.}
         % \red{this has no label on the x-axis - should be ``days'', right? also, figure 10 uses week on the x-axis, while most other figures uses days; need some clarification in the caption of figure.}
         % }  
        \label{higher_peak_nurse_demand_capacity_over_time}  
\end{figure} 

\begin{figure} 
     \centering
     \begin{subfigure}[b]{0.3\textwidth} 
         \centering
         \includegraphics[width=\textwidth]{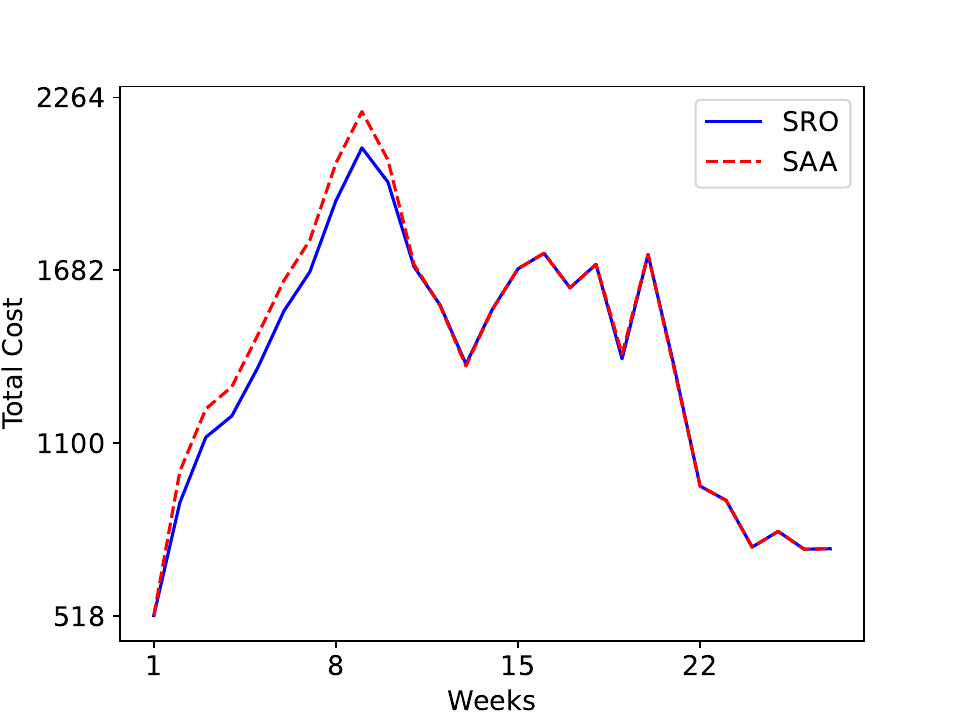} 
         \caption{Comparison on cost}
         %(\color{blue} { Draw the figure in square; x-axis: ``Weeks'' (from week 1);  y-axis:  ``Total cost'';  label: ``${SAA}$'',   ``${SRO}$'' } ) 
         \label{higher_peak_methodist_shortage_total_cost_sec_3_min_sec_2}
     \end{subfigure}
     \hfill
     \begin{subfigure}[b]{0.3\textwidth} 
         \centering
         \includegraphics[width=\textwidth]{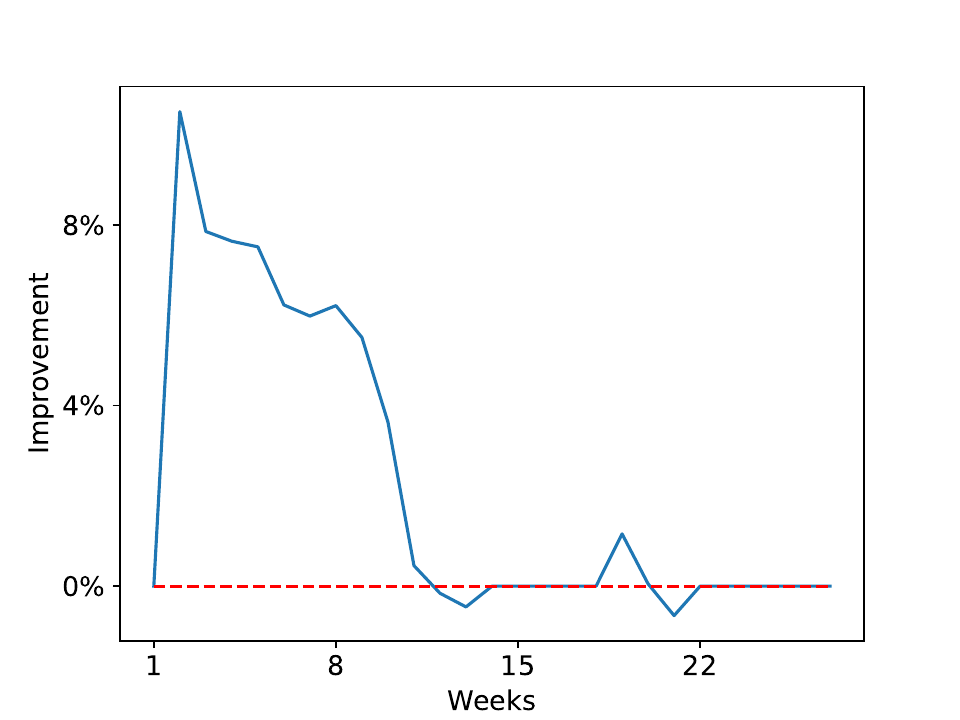}
         \caption{Improvement of SRO}  
         %  (\color{blue} { Draw the figure in square; x-axis: ``Weeks'' (from week 1);  y-axis:  ``Cost Improvement in \% of DRO Over SAA'';  no label ; draw a dashed line with y-axis$=0$ } )
         \label{higher_peak_methodist_shortage_total_cost_percentage_sec_3_min_sec_2} 
     \end{subfigure}   
     \hfill
     \begin{subfigure}[b]{0.3\textwidth} 
         \centering
         \includegraphics[width=\textwidth]{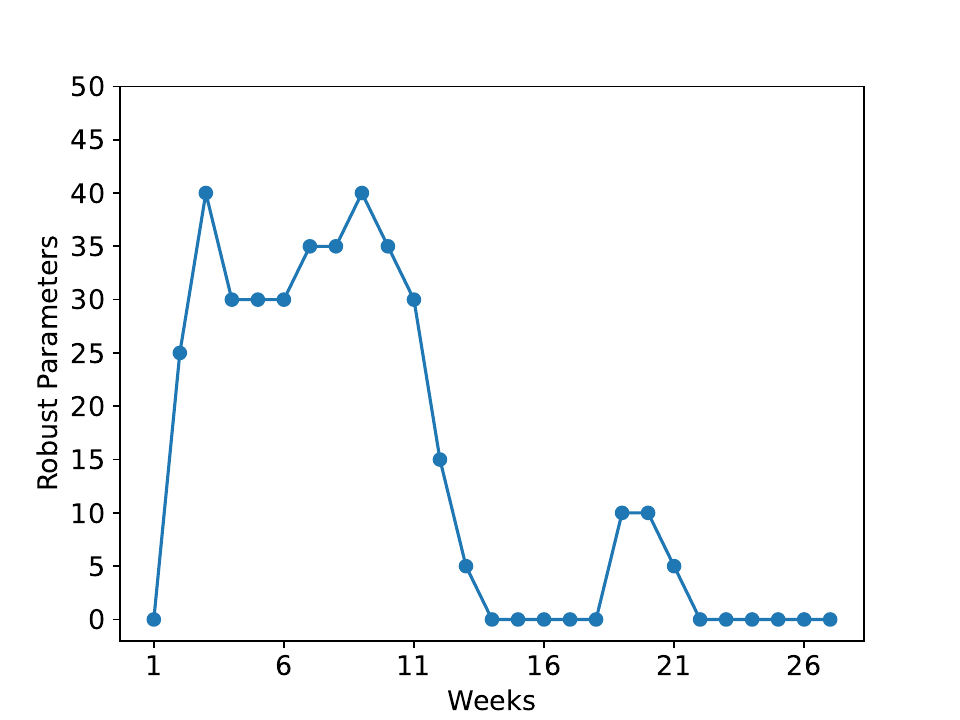}
            \caption{Robust parameter of SRO}   
         %  (\color{blue} { Draw the figure in square; x-axis: ``Weeks'' (from week 1);  y-axis:  ``Cost Improvement in \% of DRO Over SAA'';  no label ; draw a dashed line with y-axis$=0$ } )
         \label{higher_peak_fully_sec2_weekly_robust_parameter_selection_plot} 
     \end{subfigure}        
        \caption{Weekly comparison between SAA and SRO, and robust parameter used by SRO (under a higher demand peak)}   
        \label{higher_peak_Comparison_DRO_SAA}
\end{figure} 

\begin{figure}
     \centering
     \begin{subfigure}[b]{0.35\textwidth}
         \centering
         \includegraphics[width=\textwidth]{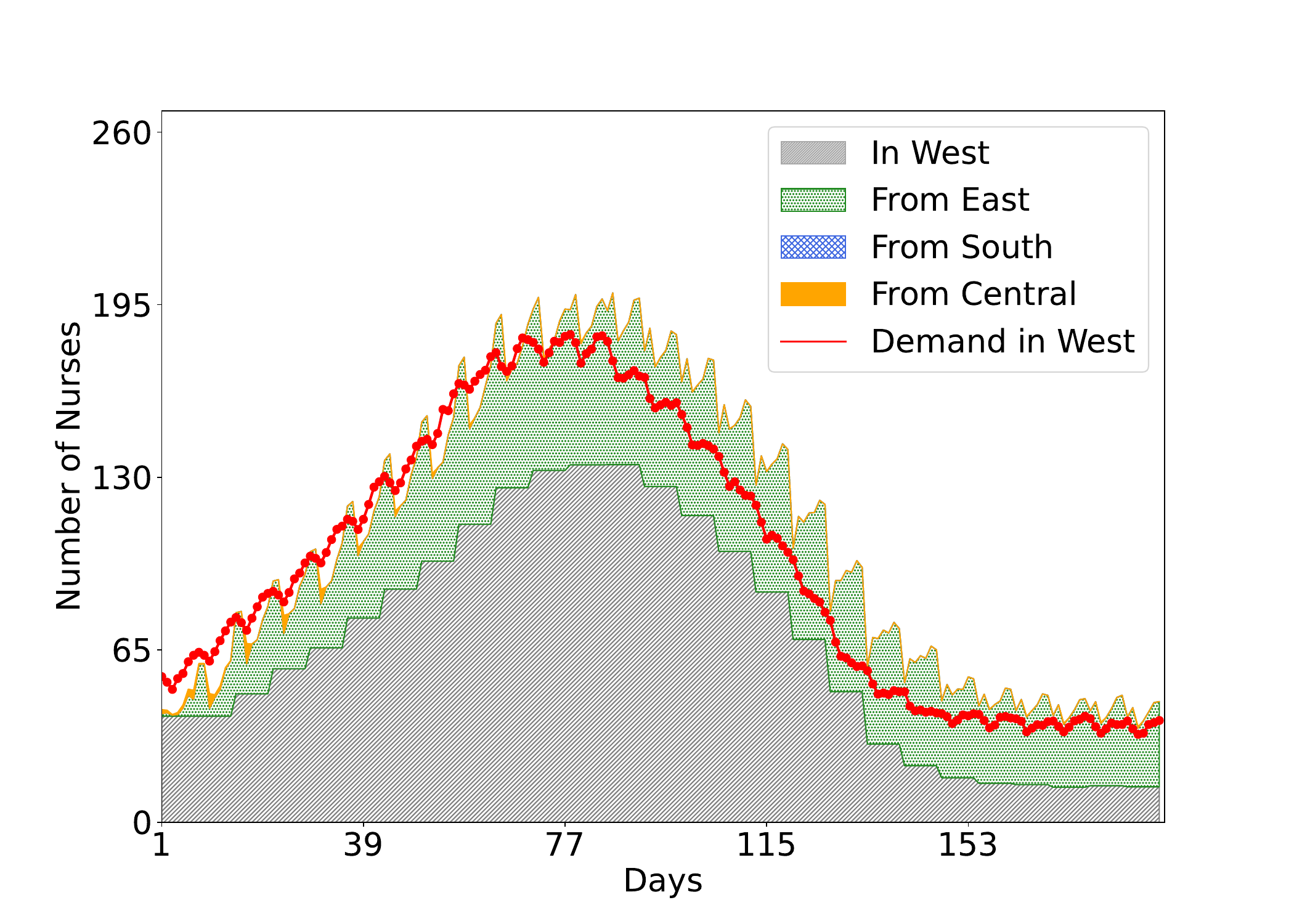}
         \caption{SAA}
         \label{higher_peak_arnett_shortage_staffing_saa_sec_2_min_sec_2}  
     \end{subfigure}
     \quad
     \begin{subfigure}[b]{0.35\textwidth}
         \centering
         \includegraphics[width=\textwidth]{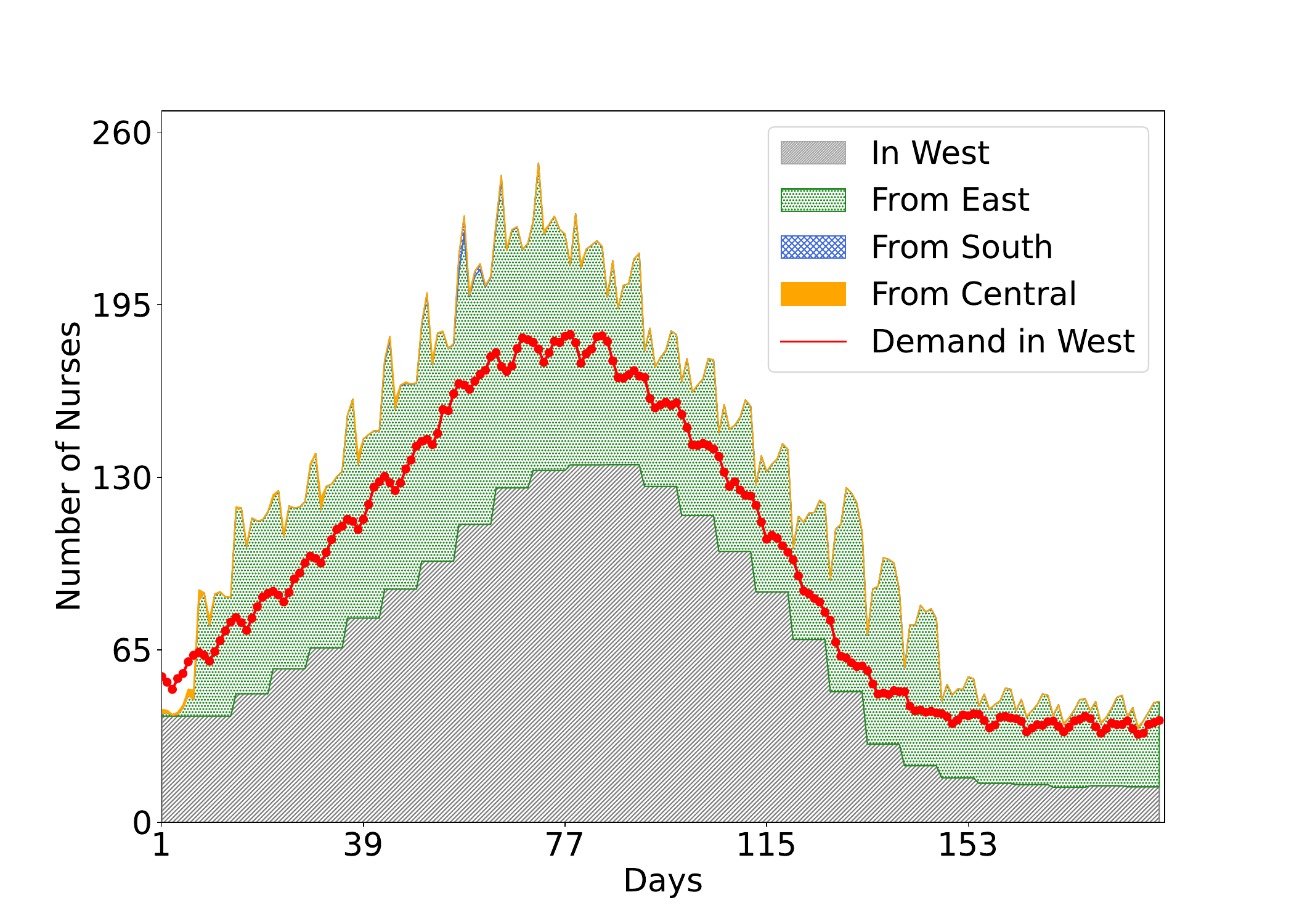}
         \caption{SRO} 
         \label{higher_peak_arnett_shortage_staffing_dro_sec_2_min_sec_2} 
     \end{subfigure} 
        \caption{Daily planned nurse transfers to West Hospital for  SAA and SRO(under a higher demand peak)}  
        % (\color{blue} { x-axis: ``Days'' (from day 1);  y-axis:  ``Number of nurses'';  label: ``Nurse demand at Methodist'', ``Nurses left with home location at Methodist'', ``Nurse transfers from Arnett'', ``Nurse transfers from Ball'', and ``Nurse transfers from Bloomington''; For the labels, use a shaded rectangle for the label; remove the line for nurse transfer and only use the shaded area; for nurse demand, may use red line directly if we a square figure. } ) 
        \label{higher_peak_Planned_decision_DRO_SAA_arnett} 
\end{figure}

\begin{figure}
     \centering
     \begin{subfigure}[b]{0.35\textwidth}
         \centering
         \includegraphics[width=\textwidth]{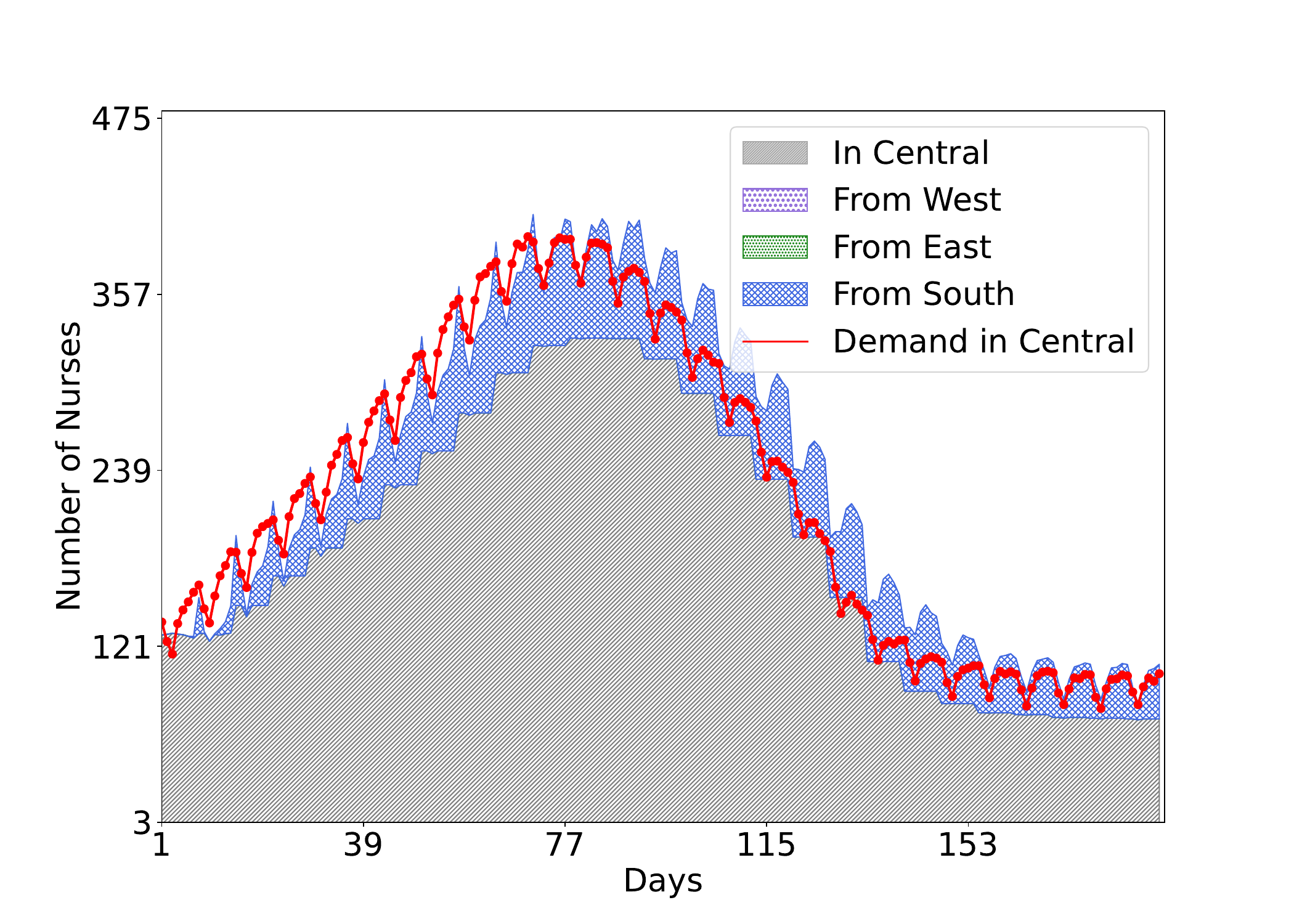}
         \caption{SAA}
         \label{higher_peak_methodist_shortage_staffing_saa_sec_2_min_sec_2} 
     \end{subfigure}
     \quad
     \begin{subfigure}[b]{0.35\textwidth} 
         \centering
         \includegraphics[width=\textwidth]{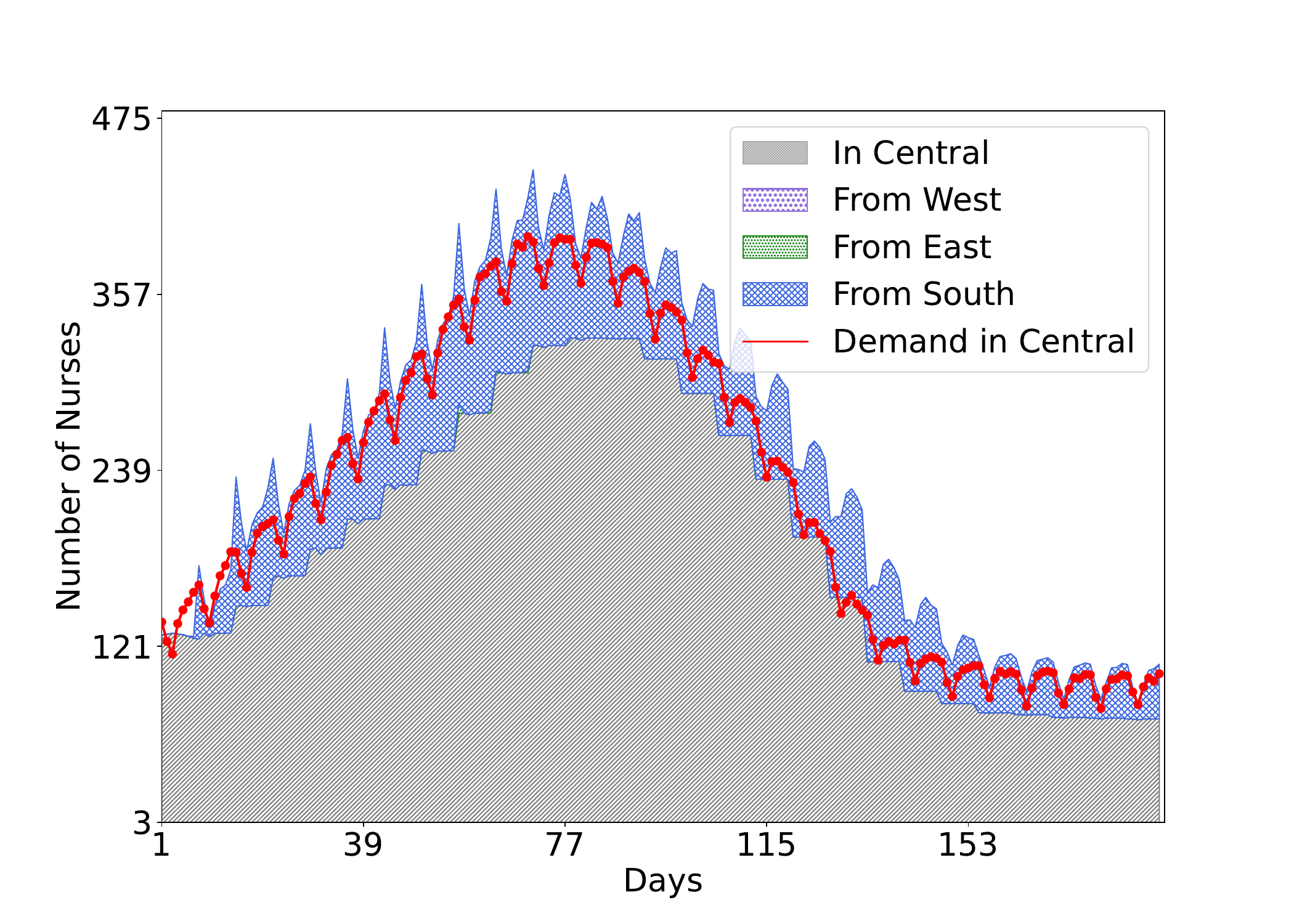} 
         \caption{SRO} 
         \label{higher_peak_methodist_shortage_staffing_dro_sec_2_min_sec_2} 
     \end{subfigure} 
        \caption{Daily planned nurse transfers to Central Hospital for SAA and SRO (under a higher demand peak)}     
        %  (\color{blue} { x-axis: ``Days'' (from day 1);  y-axis:  ``Number of nurses'';  label: ``Nurse demand at Methodist'', ``Nurses left with home location at Methodist'', ``Nurse transfers from Arnett'', ``Nurse transfers from Ball'', and ``Nurse transfers from Bloomington''; For the labels, use a shaded rectangle for the label; remove the line for nurse transfer and only use the shaded area; for nurse demand, may use red line directly if we a square figure. } ) 
        \label{higher_peak_Planned_decision_DRO_SAA_methodist}  
\end{figure}

%\newpage{}
\section{Additional Results Using Estimated Transfer Probabilities}
\label{sec:casestudy_true_transition_prob}  

In all previous case study experiments, we use the transfer probabilities listed in Table~\ref{tab:transition_prob}. 
In this section, we show that the managerial insights of the case study in Section~\ref{sec:CaseStudy} still hold when we change to use the probabilities in Table~\ref{tab:transition_prob_original_monday}. 
%We first show the average nurse demand along with the corresponding nurse capacity at the four hospitals in Section~\ref{sec:nurse_capacity_demand_origin}.
The average nurse demand along with the corresponding nurse capacity under this setting are shown in Figure~\ref{nurse_demand_capacity_over_time_true_transition_prob}, which have 
similar patterns as those in Figure~\ref{nurse_demand_capacity_over_time}.
We show the effect of network design in Section~\ref{sec:app_network}. Then we examine the effect of secondment in Section~\ref{sec:app_secondment}. 
Finally, we show the value of robustness in Section~\ref{sec:app_robustness}. 

%\subsection{Nurse Demand and Capacity}\label{sec:nurse_capacity_demand_origin}

%Specifically, the nurse demand still follows sequential patterns of increasing, decreasing, and stable phases. 
%Additionally, we continue to observe nurse shortages at West and Central Hospitals and surpluses at East and West Hospitals.

\begin{figure}
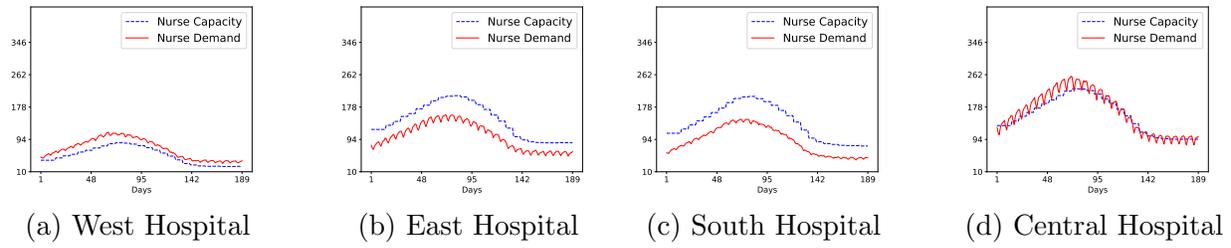

     \centering
     \begin{subfigure}[b]{0.23\textwidth}
         \centering
         \includegraphics[width=\textwidth]{census_Arnett}
         \caption{West Hospital} 
         \label{census_Arnett_true_transition_prob} 
     \end{subfigure}
     \hfill 
     \begin{subfigure}[b]{0.23\textwidth}
         \centering
         \includegraphics[width=\textwidth]{census_Ball}   
         \caption{East Hospital} 
         \label{census_Ball_true_transition_prob}     
     \end{subfigure} 
          \begin{subfigure}[b]{0.23\textwidth}
         \centering
         \includegraphics[width=\textwidth]{census_Bloomington}
         \caption{South Hospital} 
         \label{census_Bloomington_true_transition_prob} 
     \end{subfigure}
     \hfill 
     \begin{subfigure}[b]{0.23\textwidth}
         \centering
         \includegraphics[width=\textwidth]{census_Methodist}   
         \caption{Central Hospital} 
         \label{census_Methodist_true_transition_prob}     
     \end{subfigure} 
         \caption{Daily nurse demand and capacity at each hospital using probabilities in Table~\ref{tab:transition_prob_original_monday}}  
        \label{nurse_demand_capacity_over_time_true_transition_prob}
\end{figure} 

%Include the nurse capacicity and nurse demand at each hospital. 

\subsection{Network Structure}
\label{sec:app_network}

The results under different network design are in Table~\ref{tab:Average_cost_original}. 
We see that the total cost under the FC network is much lower than that under the HS network, and the network design has the most significant impact on the total cost compared to the secondment duration and approaches used. This is consistent with the observations in Section~\ref{sec:network_design}. 

% Table generated by Excel2LaTeX from sheet 'Sheet5'
\begin{table}[htbp]
  \centering
  \caption{Average cost per week for SAA and DRO across various network designs and secondment scenarios using probabilities in Table~\ref{tab:transition_prob_original_monday}} 
    \begin{tabular}{|c|rrcccc|} 
    \hline
    \multicolumn{1}{|c|}{\multirow{2}[4]{*}{Network design}} & \multicolumn{2}{c|}{Baseline secondment} & \multicolumn{2}{c|}{One-day secondment } & \multicolumn{2}{c|}{Three-day secondment } \bigstrut\\
\cline{2-7}          & \multicolumn{1}{l|}{SAA} & \multicolumn{1}{l|}{SRO} & \multicolumn{1}{l|}{SAA} & \multicolumn{1}{l|}{SRO} & \multicolumn{1}{l|}{SAA} & \multicolumn{1}{l|}{SRO} \bigstrut\\
    \hline
    FC    & 1432.8163 & 1403.16 & \multicolumn{1}{r}{1658.14} & \multicolumn{1}{r}{1591.98} & \multicolumn{1}{r}{1293.89} & \multicolumn{1}{r|}{1282.03} \bigstrut\\
\cline{1-1}    HS    & 2137.01 & 2130.38 & -     & -     & -     & - \bigstrut\\
    \hline
    \end{tabular}%
  \label{tab:Average_cost_original}%
\end{table}%

\subsection{Effect of Secondment}
\label{sec:app_secondment}
Next, we examine the effect of secondment length. 
As shown in Table~\ref{tab:Average_cost_original}, the total cost decreases with increasing secondment length, as the minimum transfer cost is set to 1.1 in this setting. 
This makes longer secondments relatively more cost-effective than shorter ones.
When we reduce the minimum transfer cost from 1.1 to 0.1, we can observe similar tradeoff from the secondment in Section~\ref{sec:Secondment_effect}.

% We show the total cost in each week under different secondment scenarios by using the SAA and SRO in Figure~\ref{total_cost_diff_sec_dro_saa_original}. 
%   \begin{figure}
%      \centering
%      \begin{subfigure}[b]{0.45\textwidth}
%          \centering
%          \includegraphics[width=\textwidth]{total_cost_diff_sec_robust_0_original}
%          \caption{SAA}
%          \label{total_cost_diff_sec_saa_true_transition_prob}
%      \end{subfigure}
%      \hfill
%      \begin{subfigure}[b]{0.45\textwidth}
%          \centering
%          \includegraphics[width=\textwidth]{total_cost_diff_sec_robust_best_original}   
%          \caption{SRO} 
%          \label{total_cost_diff_sec_dro_true_transition_prob}    
%      \end{subfigure} 
%          \caption{Effect of secondment on the cost as a function of weeks (use probabilities in Table~\ref{tab:transition_prob_original_monday})}   
%          %  (\color{blue} { Draw the figure in square; x-axis: ``Weeks'' (from week 1);  y-axis:  ``Total cost'';  label:   ``One-day secondment scenario '', ``Base secondment scenario'' } ) 
%         \label{total_cost_diff_sec_dro_saa_original}
% \end{figure} 

%Include the figure on the comparison between different secondment scenarios. 

\subsection{Value of Robustness}
\label{sec:app_robustness}

We present the comparison along with the robust parameter used by SRO in Figure~\ref{Comparison_DRO_SAA_true_transition_prob}.
The results indicate that SRO generally outperforms SAA under increasing, decreasing, and stable demand patterns, with only a few exceptions.
To understand this difference, we analyze the planned nurse transfers to East and Methodist Hospitals, shown in Figures~\ref{Planned_decision_DRO_SAA_arnett_true_transition_prob} and~\ref{Planned_decision_DRO_SAA_methodist_true_transition_prob}. 
We see that SRO assigns more nurses during the planning stage, which helps reduce the need for emergency transfers. 
In contrast, SAA often fails to assign enough nurses to meet demand. 
Specifically, SAA generally does not allocate enough nurses during periods of increasing demand and fails to assign adequate staff at the beginning of each week under decreasing and stable demand patterns.
This results in a higher number of emergency transfers, which aligns with our findings in Section~\ref{sec:ComparisonofPolicies}.

% To further investigate, we examine nurse capacity, actual nurse demand, and predicted nurse demand for a representative sample path, as shown in Figure~\ref{nurse_demand_capacity_prediction_over_time_true_transition_prob}. 
% We observe that predicted nurse demand underestimates actual demand during periods of increasing demand.
% Additionally, underestimation may occur at the beginning of each week during decreasing and stable demand patterns, despite more frequent overestimations at other times. 
% This results in the superior performance of SRO, as SAA makes decisions based on the sample average of the objective function across all predicted sample paths.
% These findings not only highlight that SRO outperforms SAA when there is some level of underestimation in nurse demand but also demonstrate that SRO exhibits greater adaptability to demand changes. This aligns with our findings in Section~\ref{sec:ComparisonofPolicies}.

\begin{figure} 
     \centering
     \begin{subfigure}[b]{0.3\textwidth} 
         \centering
         \includegraphics[width=\textwidth]{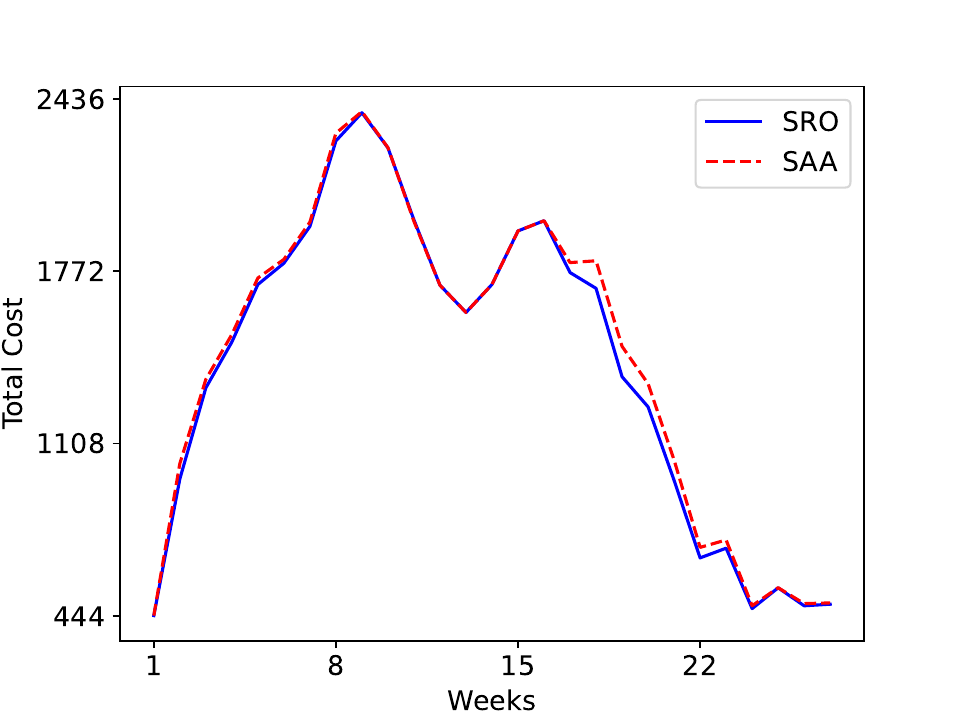} 
         \caption{Comparison on cost}
         %(\color{blue} { Draw the figure in square; x-axis: ``Weeks'' (from week 1);  y-axis:  ``Total cost'';  label: ``${SAA}$'',   ``${SRO}$'' } ) 
         \label{methodist_shortage_total_cost_sec_3_min_sec_2_true_transition_prob}
     \end{subfigure}
     \hfill
     \begin{subfigure}[b]{0.3\textwidth} 
         \centering
         \includegraphics[width=\textwidth]{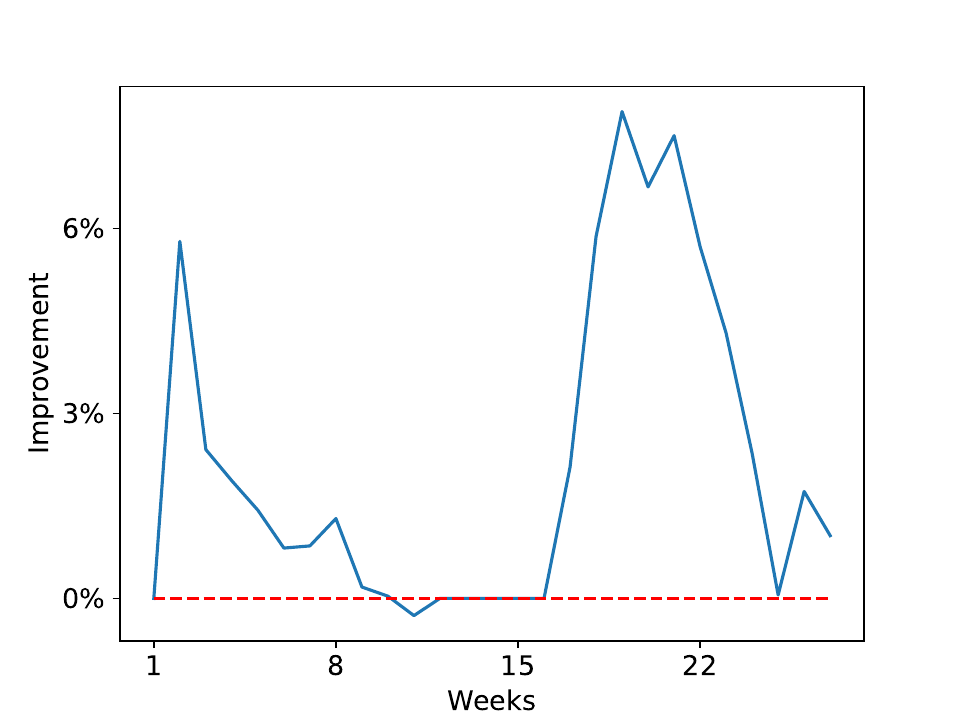}
         \caption{Improvement of SRO}  
         %  (\color{blue} { Draw the figure in square; x-axis: ``Weeks'' (from week 1);  y-axis:  ``Cost Improvement in \% of DRO Over SAA'';  no label ; draw a dashed line with y-axis$=0$ } )
         \label{methodist_shortage_total_cost_percentage_sec_3_min_sec_2_true_transition_prob} 
     \end{subfigure}   
     \hfill
     \begin{subfigure}[b]{0.3\textwidth} 
         \centering
         \includegraphics[width=\textwidth]{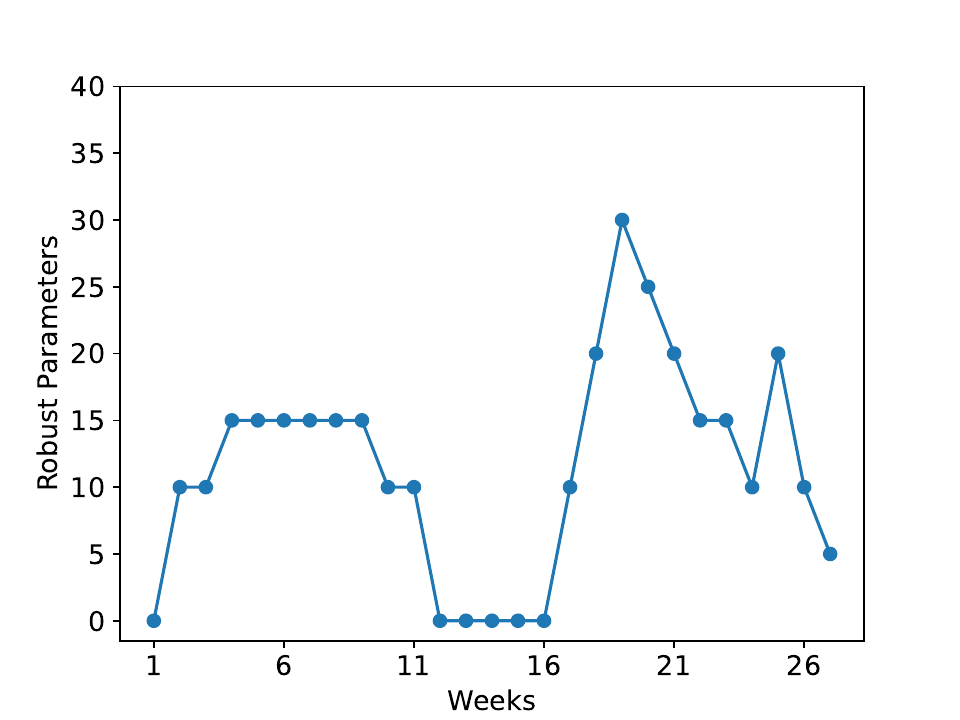}
            \caption{Robust parameter of SRO}   
         %  (\color{blue} { Draw the figure in square; x-axis: ``Weeks'' (from week 1);  y-axis:  ``Cost Improvement in \% of DRO Over SAA'';  no label ; draw a dashed line with y-axis$=0$ } )
         \label{fully_sec2_weekly_robust_parameter_selection_plot_true_transition_prob} 
     \end{subfigure}        
        \caption{Weekly comparison between SAA and SRO, and robust parameter used by SRO, using probabilities in Table~\ref{tab:transition_prob_original_monday}}   
        \label{Comparison_DRO_SAA_true_transition_prob}
\end{figure} 

% \begin{figure} 
%          \centering
%          \includegraphics[width=0.45\textwidth]{fully_sec2_weekly_robust_parameter_selection_plot}
%          \caption{Robust parameter used by SRO under the fully connected network and base secondment scenario}    
%          \label{fully_sec2_weekly_robust_parameter_selection_plot} 
% \end{figure} 

\begin{figure}
     \centering
     \begin{subfigure}[b]{0.35\textwidth}
         \centering
         \includegraphics[width=\textwidth]{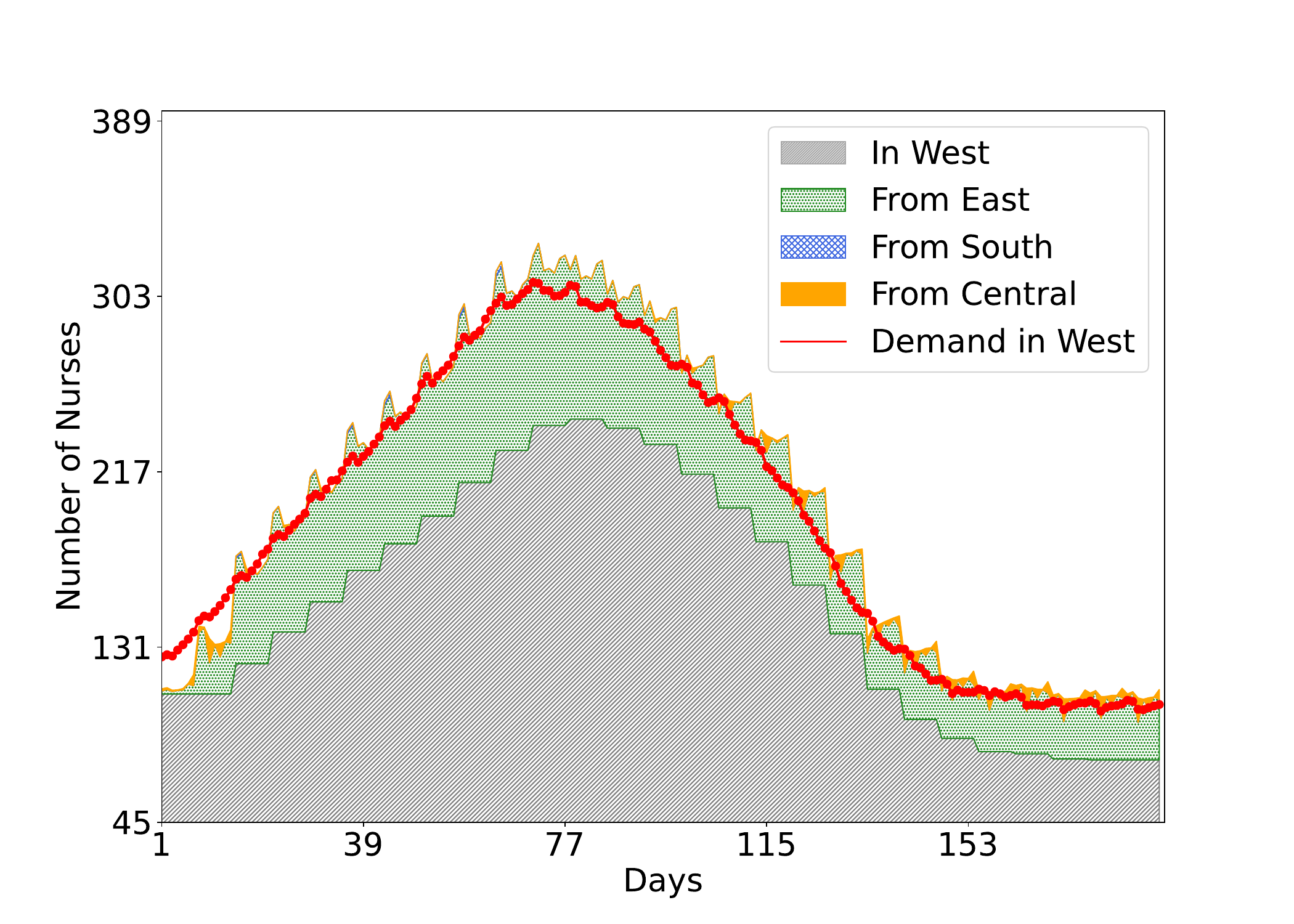}
         \caption{SAA}
         \label{arnett_shortage_staffing_saa_sec_2_min_sec_2_true_transition_prob} 
     \end{subfigure}
     \quad
     \begin{subfigure}[b]{0.35\textwidth}
         \centering
         \includegraphics[width=\textwidth]{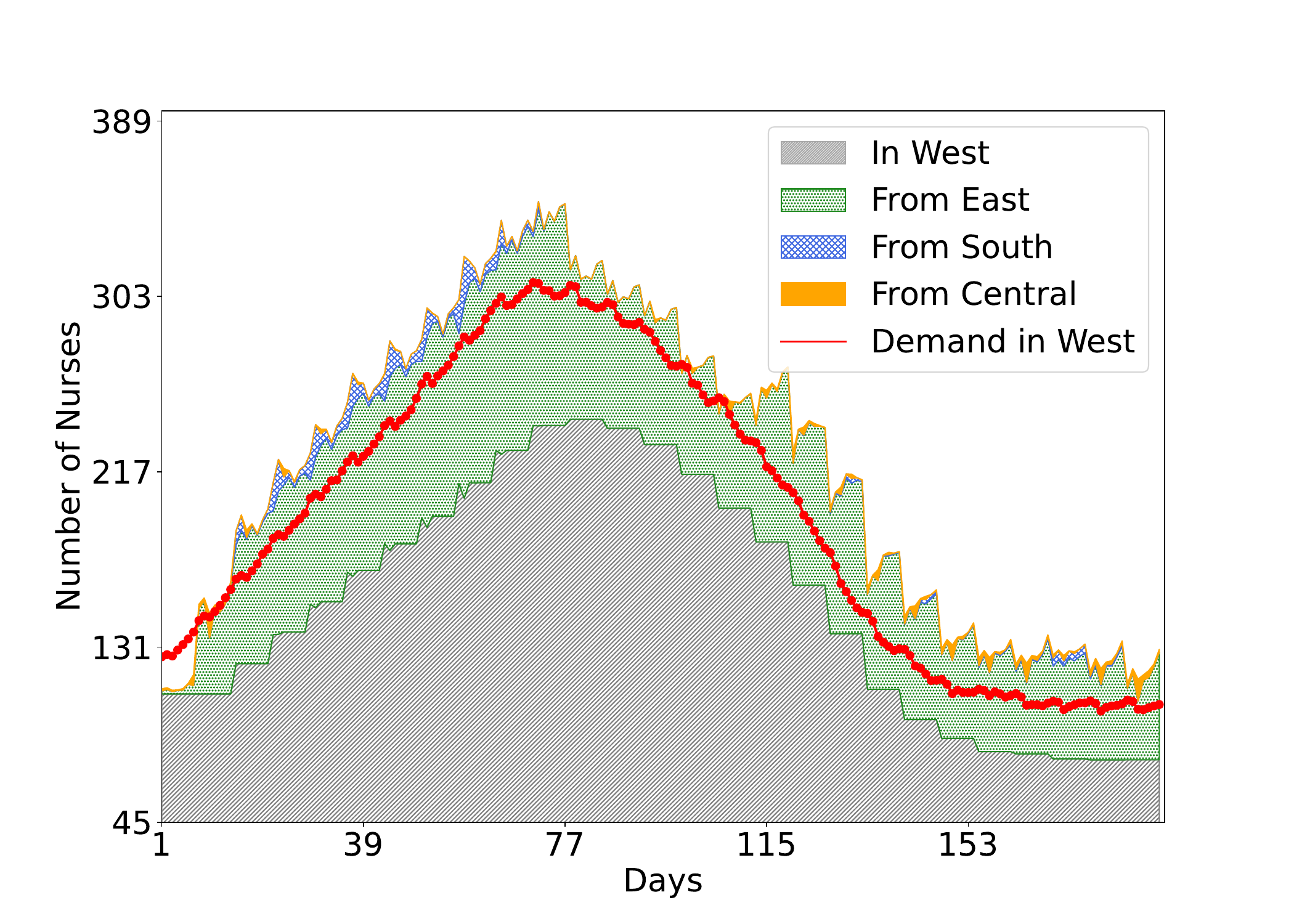}
         \caption{SRO} 
         \label{arnett_shortage_staffing_dro_sec_2_min_sec_2_true_transition_prob} 
     \end{subfigure} 
        \caption{Daily planned nurse transfers to West Hospital for SAA and SRO, using probabilities in Table~\ref{tab:transition_prob_original_monday}}  
        % (\color{blue} { x-axis: ``Days'' (from day 1);  y-axis:  ``Number of nurses'';  label: ``Nurse demand at Methodist'', ``Nurses left with home location at Methodist'', ``Nurse transfers from Arnett'', ``Nurse transfers from Ball'', and ``Nurse transfers from Bloomington''; For the labels, use a shaded rectangle for the label; remove the line for nurse transfer and only use the shaded area; for nurse demand, may use red line directly if we a square figure. } ) 
        \label{Planned_decision_DRO_SAA_arnett_true_transition_prob} 
\end{figure}

\begin{figure}
     \centering
     \begin{subfigure}[b]{0.35\textwidth}
         \centering
         \includegraphics[width=\textwidth]{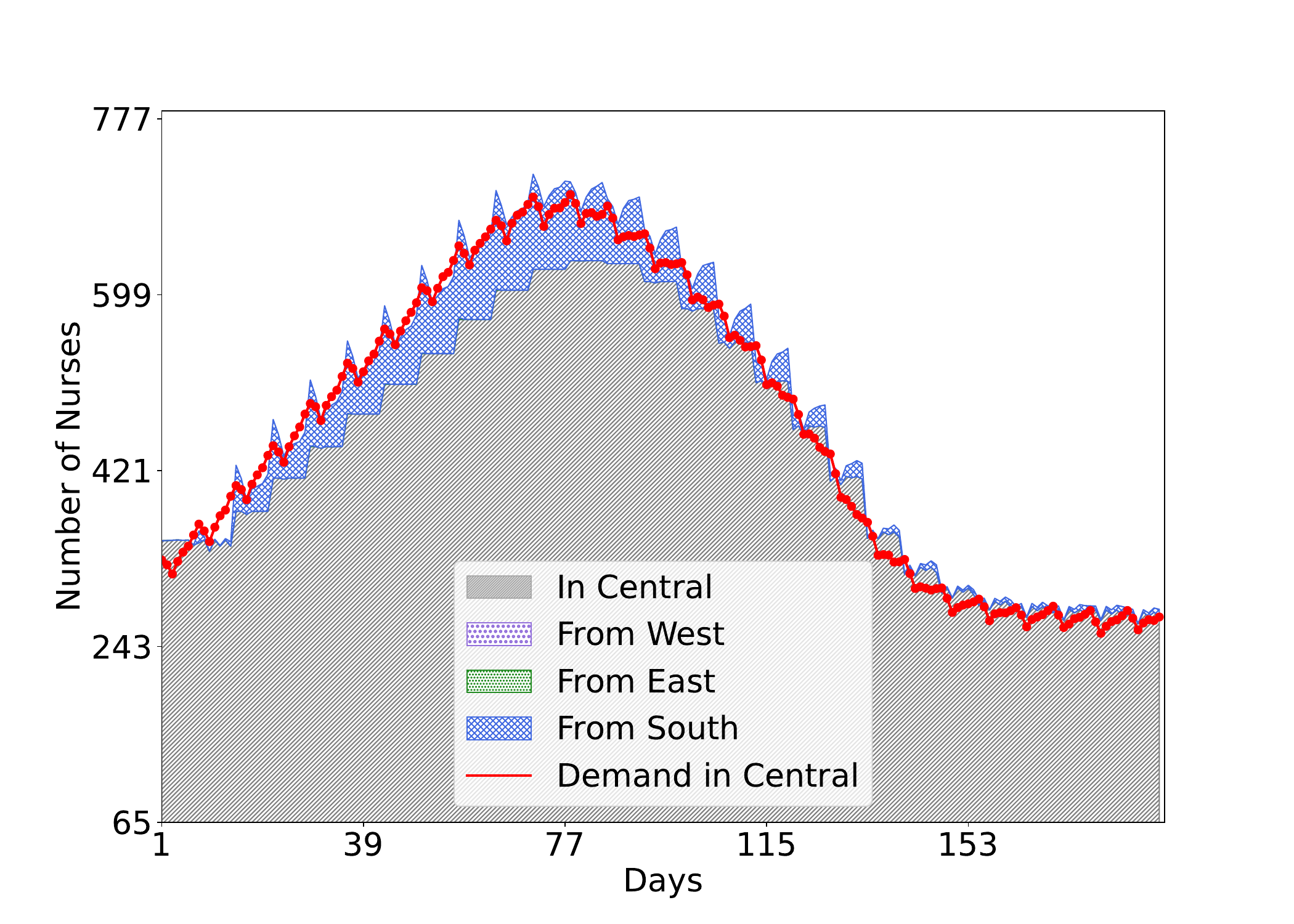}
         \caption{SAA}
         \label{methodist_shortage_staffing_saa_sec_2_min_sec_2_true_transition_prob} 
     \end{subfigure}
     \quad
     \begin{subfigure}[b]{0.35\textwidth} 
         \centering
         \includegraphics[width=\textwidth]{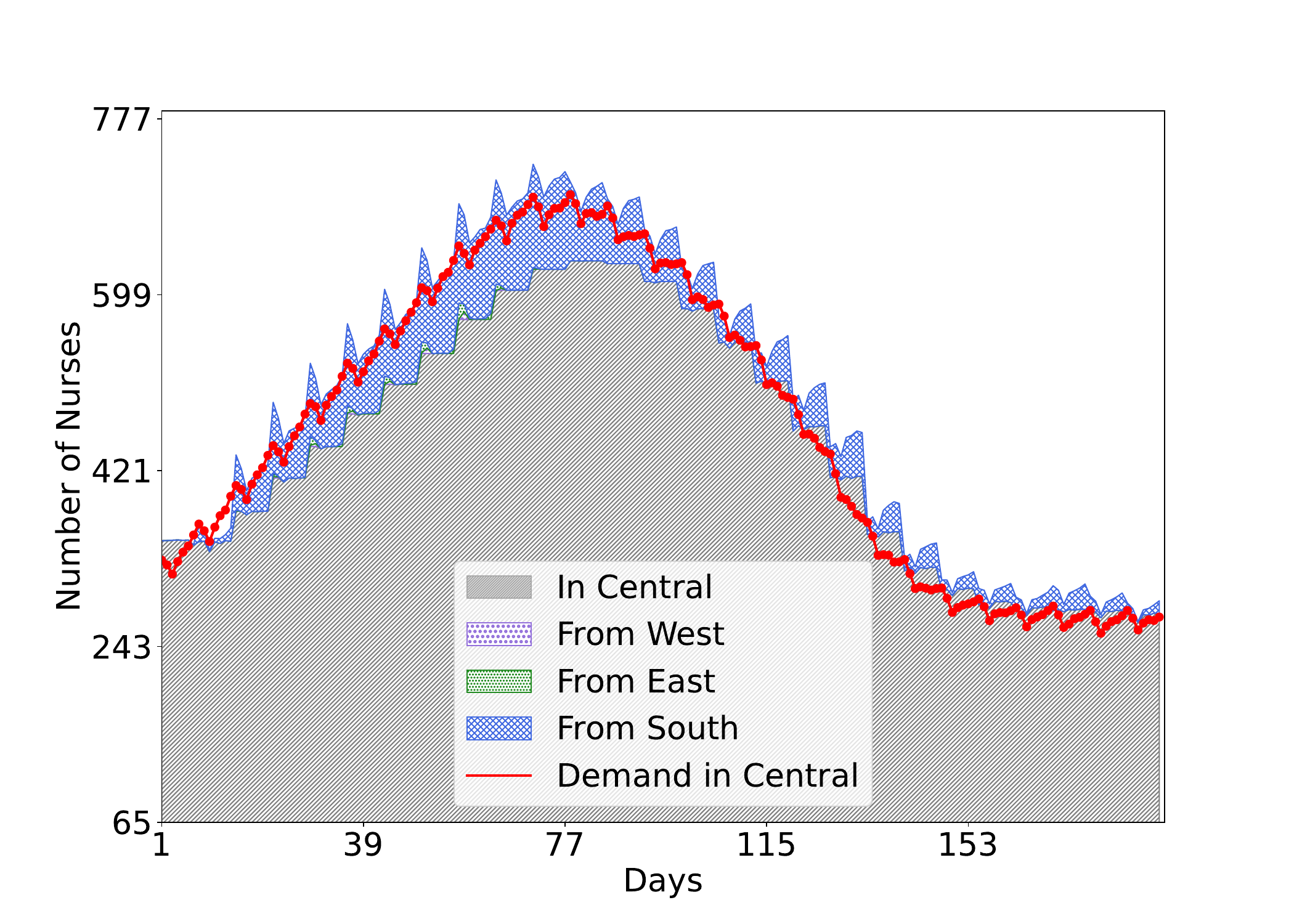} 
         \caption{SRO} 
         \label{methodist_shortage_staffing_dro_sec_2_min_sec_2_true_transition_prob} 
     \end{subfigure} 
        \caption{Daily planned nurse transfers to Central Hospital for SAA and SRO, using probabilities in Table~\ref{tab:transition_prob_original_monday}}   
        %  (\color{blue} { x-axis: ``Days'' (from day 1);  y-axis:  ``Number of nurses'';  label: ``Nurse demand at Methodist'', ``Nurses left with home location at Methodist'', ``Nurse transfers from Arnett'', ``Nurse transfers from Ball'', and ``Nurse transfers from Bloomington''; For the labels, use a shaded rectangle for the label; remove the line for nurse transfer and only use the shaded area; for nurse demand, may use red line directly if we a square figure. } ) 
        \label{Planned_decision_DRO_SAA_methodist_true_transition_prob}  
\end{figure}  

% \begin{figure}
%      \centering
%      \begin{subfigure}[b]{0.23\textwidth}
%          \centering
%          \includegraphics[width=\textwidth]{Arnett_demand_vs_capacity_training_true_transition_prob}
%          \caption{West Hospital} 
%          \label{Arnett_demand_vs_capacity_training_true_transition_prob} 
%      \end{subfigure}
%      \hfill 
%      \begin{subfigure}[b]{0.23\textwidth}
%          \centering
%          \includegraphics[width=\textwidth]{Ball_Mem_demand_vs_capacity_training_true_transition_prob}   
%          \caption{East Hospital} 
%          \label{Ball_Mem_demand_vs_capacity_training}     
%      \end{subfigure} 
%           \begin{subfigure}[b]{0.23\textwidth}
%          \centering
%          \includegraphics[width=\textwidth]{Bloomington_demand_vs_capacity_training_true_transition_prob}
%          \caption{South Hospital} 
%          \label{Bloomington_demand_vs_capacity_training_true_transition_prob} 
%      \end{subfigure}
%      \hfill 
%      \begin{subfigure}[b]{0.23\textwidth}
%          \centering
%          \includegraphics[width=\textwidth]{Methodist_demand_vs_capacity_training_true_transition_prob}   
%          \caption{Central Hospital} 
%          \label{Methodist_demand_vs_capacity_training_true_transition_prob}     
%      \end{subfigure} 
%          \caption{Daily nurse capacity and demand along with the demand prediction   for one sample path (underTable~\ref{tab:transition_prob_original_monday})} 
%         \label{nurse_demand_capacity_prediction_over_time_true_transition_prob} 
% \end{figure} 

\end{APPENDICES} 

%%%%%%%%%%%%%%%%%
\end{document}